%% file: TESE.tex
\documentclass[dsc,numbers]{pgim}
\usepackage{amsmath,amsthm,amsfonts,amssymb,amscd,amstext}
\usepackage[ansinew]{inputenc}
\usepackage[titletoc]{appendix}
\usepackage{textcomp}
\usepackage{flowchart}
\usepackage{float}

\usepackage{euscript}

\usepackage{epsfig,psfrag}

\usepackage[active]{srcltx}
\usepackage{xcolor}
\usepackage[colorlinks=true, linkcolor=blue, urlcolor=blue]{hyperref}

\usepackage{hyperref}
\usepackage[english]{babel}

\usepackage{verbatim}

\newtheorem{theorem}{Theorem}[section]

\newtheorem{corollary}[theorem]{Corollary}
\newtheorem{lemma}[theorem]{Lemma}
\newtheorem{algorithm}[theorem]{Algorithm}
\newtheorem{definition}[theorem]{Definition}
\newtheorem{remark}[theorem]{Remark}
\newtheorem{example}[theorem]{Example}

\let\oldmarginpar\marginpar
\renewcommand\marginpar[1]{\-\oldmarginpar[\raggedleft\footnotesize #1]%
{\raggedright\footnotesize #1}}

\makelosymbols
\makeloabbreviations

\begin{document}
  \title{\bf{Tensor decompositions and algorithms, with applications to tensor learning}}
  \foreigntitle{Tensor decompositions and algorithms, with applications to tensor learning}
  \author{Felipe}{Bottega Diniz}
  \advisor{Prof.}{Gregorio}{Malajovich}{Munoz}


  \examiner{Prof.}{gdhjagsd}{D.Sc.(presidente)}
  \examiner{Prof.}{}{D.Sc.}
  \examiner{Prof.}{}{D.Sc.}
  \examiner{Prof.}{}{D.Sc.}
  \examiner{Prof.}{}{D.Sc.}
  \examiner{Prof.}{}{D.Sc.}
  \examiner{Prof.}{}{D.Sc.(suplente)}
  \department{PGMAT}
  \date{11}{2019}

  \keyword{Tensors}
  \keyword{Canonical polyadic decomposition}
  \keyword{Multilinear singular value decomposition}
  \keyword{Nonlinear least squares}
  \keyword{Machine learning}

  \maketitle
  \include{thanks}
  \frontmatter
  \include{resumo}
  \include{abstract}
  \include{dedic} 
  \tableofcontents 
  \listoffigures 
  \listoftables 
  \printlosymbols 
  \printloabbreviations 

  \mainmatter 
  \include{CHAPTER_Intro}

  \include{CHAPTER_1}
  \include{CHAPTER_2}
  \include{CHAPTER_3}
  \include{CHAPTER_4}
  \include{CHAPTER_5}

  \include{Conclusion}

  \begin{appendices}
  \include{appenA}
  \include{appenB}

  \end{appendices}

  \backmatter
  \bibliographystyle{coppe-unsrt}


\end{document}

%% file: thanks.tex
\section*{ }
	\pagestyle{empty}
	\begin{center}
 		{\Large{Tensor decompositions and algorithms, with applications to tensor learning}}
	\bigskip
	\\
	{Felipe Bottega Diniz}\\
	{Orientador: Gregorio Malajovich Mu\~noz}
	\end{center}
	\bigskip\bigskip\bigskip
	Tese de doutorado apresentada ao Programa de P\'os-Gradua\c c\~ao
	em Matem\'atica, Instituto de Matem\'atica da Universidade
	Federal do Rio de Janeiro (UFRJ), como parte dos requisitos
	necess\'arios \`a obten\c c\~ao do t\'itulo de doutor em Matem\'atica.\\\\
	Aprovada por:\\\\
	
	Presidente, Prof. Gregorio Malajovich Mu\~noz\\
	
	Prof. Bernardo Freitas Paulo da Costa\\
	
	Prof. Nick Vannieuwenhoven\\
	
	Prof. Andr\'e Lima Ferrer de Almeida\\
	
	Prof. Amit Bhaya\\

	\bigskip\bigskip\bigskip\bigskip\bigskip

	\begin{center}
	{Rio de Janeiro}\\
	{Novembro de 2019}
	\end{center}

%% file: resumo.tex
\section*{\hspace{6cm} Resumo}
	Neste trabalho \'e apresentada uma nova implementa\c c\~ao da canonical polyadic decomposition (CPD). Ela possui uma menor complexidade computacional e menor uso de mem\'oria do que as implementa\c c\~oes estado da arte dispon\'iveis.

	Come\c camos com alguns exemplos de aplica\c c\~oes da CPD para problemas do mundo real. Um breve resumo das principais contribui\c c\~oes deste trabalho \'e o seguinte. No cap\'itulo 1, revisamos a \'algebra e geometria cl\'assicas de tensores, com foco na CPD. O cap\'itulo 2 \'e focado na compress\~ao tensorial, que \'e considerada (neste trabalho) uma das partes mais importantes do algoritmo CPD. No cap\'itulo 3, falamos sobre o m\'etodo de Gauss-Newton, que \'e um m\'etodo de m\'inimos quadrados n\~ao-lineares usado para minimizar fun\c c\~oes n\~ao-lineares. O cap\'itulo 4 \'e o mais longo deste trabalho. Neste cap\'itulo, apresentamos o personagem principal desta tese: Tensor Fox. Basicamente, \'e um pacote tensorial que inclui um CPD solver. Ap\'os a introdu\c c\~ao do Tensor Fox, realizaremos muitas experi\^encias computacionais comparando esse solver com v\'arios outros. No final deste cap\'itulo, apresentamos a decomposi\c c\~ao Tensor Train e mostramos como us\'a-la para calcular CPDs de ordem superior. Tamb\'em discutimos alguns detalhes importantes, como regulariza\c c\~ao, pr\'e-condicionamento, condicionamento, paralelismo, etc. No cap\'itulo 5, consideramos a interse\c c\~ao entre decomposi\c c\~oes tensoriais e machine learning. \'E introduzido um novo modelo, que funciona como uma vers\~ao tensorial de redes neurais. Finalmente, no cap\'itulo 6, fazemos as conclus\~oes finais e introduzimos nossas expectativas de trabalho futuro.

%% file: abstract.tex
\section*{\hspace{6cm} Abstract}
	A new algorithm of the canonical polyadic decomposition (CPD) presented here. It features lower computational complexity and memory usage than the available state of the art implementations. 

	We begin with some examples of CPD applications to real world problems. A short summary of the main contributions in this work follows. In chapter 1 we review classical tensor algebra and geometry, with focus on the CPD. Chapter 2 focuses on tensor compression, which is considered (in this work) to be one of the most important parts of the CPD algorithm. In chapter 3 we talk about the Gauss-Newton method, which is a nonlinear least squares method used to minimize nonlinear functions. Chapter 4 is the longest one of this thesis. In this chapter we introduce the main character of this thesis: Tensor Fox. Basically it is a tensor package which includes a CPD solver. After introducing Tensor Fox we will conduct lots of computational experiments comparing  this solver with several others. At the end of this chapter we introduce the Tensor Train decomposition and show how to use it to compute higher order CPDs. We also discuss some important details such as regularization, preconditioning, conditioning, parallelism, etc. In chapter 5 we consider the intersection between tensor decompositions and machine learning. A novel model is introduced, which works as a tensor version of neural networks. Finally, in chapter 6 we reach the final conclusions and introduce our expectations for future developments. 

%% file: dedic.tex
\section*{\hspace{4.8cm} Acknowledgments}
	This work is the result of a long journey, which started more than 10 years ago, when I decided to become a mathematician. At first, I was just a graduate student who liked math; however, little by little, its infinite beauty began to reveal itself to me. Today it is a big honor to be a mathematician. I owe many thanks to some wonderful people who helped me in one way or another, and I dedicate this space to them.
	
	First I thank to my mother, F\'atima Val\'eria, who provided me a comfortable and stimulating environment where I could focus on my studies. There is no library comparable to my home and I wouldn't be here if it weren't for her.  
	
	All my family were very supporting and believed in me. In particular I thank for my young brother, Bruno, my grandmother Dalba and grandfather Walder, my stepfather Fernando who is now gone (God bless him), and my father Nirvando. 
	
	During all my doctorate I was very motivated by my girlfriend, Eliane (Lili). She was very important to me and helped me a lot in difficult times. I admire she as a person and as a mathematician. She helped me to overcome this difficult chapter of my life. She taught me the meaning of love.  
	
	My friends suffered a little with my constant absence from social meetings, because I always said that `` I have to study '', but they supported and I was very fortunate to have such good friends. Thank you very much for your patience and for believing in me, Jorge, Thiago, Rodolpho, Karina, J\'ulio, Jose Hugo, Bernardo, Deodato, Yuri and Taynara. You all are the best and know that your friendship means everything to me!
	
	I would like to thank my advisor Gregorio Malajovich. In the past years we developed a healthy relation which turned out to be very productive, I learned a lot, and most importantly, thanks to him I was able to find myself in mathematics. Gregorio was the responsible for showing me the world of tensors. At the beginning I wasn't sure if I wanted to pursue this path, but as I progress it become clear that he knew better than me that this was a good path. Thank you very much for believing in me , Gregorio.
	
	I also want to thank professor Bernardo Freitas, who accompanied my whole journey during the doctorate. He is an amazing person who is always paying attention and challenging you. His valuable observations and advice helped me a lot be where I am now.   
	
	Nilson da C. Bernardes, it was privilege to be present at your lectures of Functional Analysis and even more to have you in my qualification exam, also in Functional Analysis. Your way of teaching is inspiring, I hope to be able to teach like this in the future.
	
	This study was financed in part by the Coordena\c c\~ao de Aperfei\c coamento de Pessoal de N\'ivel Superior - Brasil (CAPES) - Finance Code 001. 
	
	Finally, I like to thank Gregorio Malajovich, Bernardo Freitas, Nick Vannieuwenhoven, Andr\'e de Almeida and Amit Bhaya for being part of my doctoral defense committee.  

%% file: CHAPTER_Intro.tex
\chapter*{Introduction}\label{intro}
\begin{figure}[h] 
		\hspace*{-.5cm}
		\includegraphics[scale=.5]{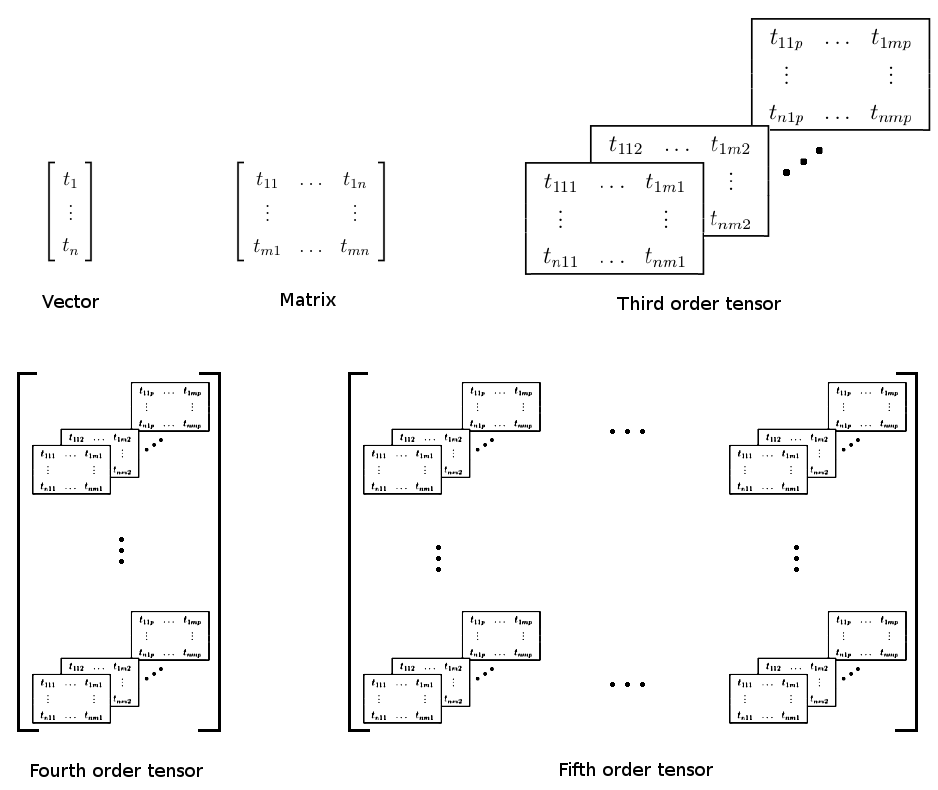}
		\caption{\footnotesize{Shapes of tensors for the first five orders.}}
		\label{tensor-intuition}
	\end{figure}

\addcontentsline{toc}{chapter}{Introduction} 
	A vector can be thought as data arranged in an unidimensional fashion, that is, an ordered sequence of numbers, strings, or any other kind of information. In the same way a matrix can be thought as bidimensional data, which also is a sequence of vectors with the same length. Tensors are a natural generalization of this process. One can dispose several matrices of same shape in an ordered sequence, forming a 3D-block of data, see figure~\ref{tensor-intuition}. This is what is called a \emph{third order tensor}. Analogously, matrices are second order tensors and vectors are first order tensors. It is useful to define scalars as zero order tensors. Recursively, one can define a \emph{$L$-th order tensor} as an ordered sequence of $(L-1)$-th order tensors of same shape.  
	
	Tensors can be defined rigorously as mathematical objects but, for the moment, it will be convenient to think of tensors just as multidimensional arrays of data. Given a field $\mathbb{K}$ and positive integers $I_1, \ldots, I_L$, the set of tensors with shape $I_1 \times \ldots \times I_L$ is defined as
		$$\mathbb{K}^{I_1 \times \ldots \times I_L} = \{ (t_{i_1 \ldots i_L})_{1 \leq i_1 \leq I_1, \ldots, 1 \leq i_L \leq I_L} | \ t_{i_1 \ldots i_L} \in \mathbb{K} \}.$$
This is an informal definition of a tensor space since it depends on coordinates, although it is a useful way to visualize and think of tensors. Later we will give a new definition more aligned with multilinear algebra. In this work we always use an upper index to indicate a sequence of vectors and a lower index to indicate their coordinates. 
	
	Given any vectors $\textbf{v}^{(1)} \in \mathbb{K}^{I_1}, \ldots, \textbf{v}^{(1)} \in \mathbb{K}^{I_L}$, define the tensor $\textbf{v}^{(1)} \otimes \ldots \otimes \textbf{v}^{(L)} \in \mathbb{K}^{I_1} \otimes \ldots \otimes \mathbb{K}^{I_L}$ by
	
	$$(\textbf{v}^{(1)} \otimes \ldots \otimes \textbf{v}^{(L)})_{i_1 \ldots i_L} = v_{i_1}^{(1)} \cdot \ldots \cdot v_{i_L}^{(L)}.$$	
Any tensor of this form is called a \emph{rank one tensor}. One can also say the tensor \emph{has rank one}. Notice that in the case of second order tensors (matrices), this definition agrees with the definition of rank one matrices. It is not hard to see that any tensor can be written as a sum of rank one tensors. It is of interest to find the minimum number of rank one terms necessary to construct such a sum, and this minimum number is called the \emph{rank} of the tensor. The decomposition of a tensor as a sum of rank one terms is called a \emph{canonical polyadic decomposition} (CPD). For decades, tensor decompositions have been applied to general multidimensional data with success. Today they are excel in several applications, including blind source separation, dimensionality reduction, pattern/image recognition, machine learning and data mining \cite{bro, cichoki2009, cichoki2014, savas, smilde, anandkumar}. One of the reasons for the successfulness of tensor decompositions comes from its uniqueness, which occurs for all higher order tensors. This is a desired property which is not found by matrices. 
	
	We start giving some motivational examples which highlight the applicability of tensor decompositions, the main topic of this work. In particular, some attention will be given to machine learning applications, a theme to be explored in more details only in the last chapter, when we will have developed the necessary machinery for such.
	
	\section{First example: Gaussian mixtures}	
		Consider a mixture of $K$ Gaussian distributions with identical covariance matrices. We have lots of data with unknown averages and unknown covariance matrices. The problem at hand is to design an algorithm to \emph{learn} these parameters from the data given. We use $\mathbb{P}$ to denote probability and $\mathbb{E}$ to denote expectation (which we may also call the \emph{mean} or \emph{average}).
		
\begin{figure}[h] 
			\centering
			\includegraphics[scale=.42]{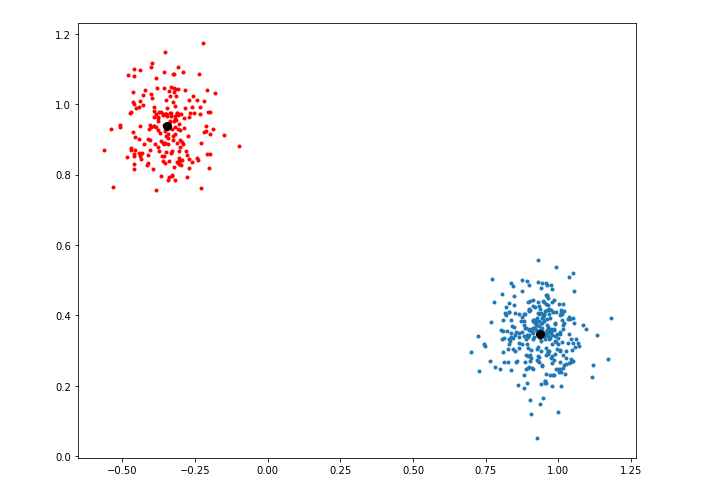}
			\caption{\footnotesize{Gaussian mixture in the plane with 2 clusters. The first cluster has mean $\textbf{u}^{(1)} = [-0.34, \  0.93]^T$ and the second has mean $\textbf{u}^{(2)} = [0.93, \ 0.34]^T$. The variance is $\sigma^2 = 0.0059$.}}
			\label{gaussian_mix}
		\end{figure}

		Let $\textbf{x}^{(1)}, \ldots, \textbf{x}^{(N)} \in \mathbb{R}^d$ be a set of collected data sample. Let $h$ be a discrete random variable with values in $\{1, 2, \ldots, K\}$ such that $\mathbb{P}[h = i]$ is the probability that a sample $\textbf{x}$ is a member of the $i$-th cluster. We denote $w^{(i)} = \mathbb{P}[h = i]$ and $\textbf{w} = [ w^{(1)}, \ldots, w^{(K)}]^T$, the vector of probabilities. Let $\textbf{u}^{(i)} \in \mathbb{R}^d$ be the mean of the $i$-th distribution and assume that all distributions have the same covariance matrix $\sigma^2 \textbf{I}_d$ for $\sigma > 0$. See figure~\ref{gaussian_mix} for an illustration of a Gaussian mixture in the case where $d = 2$ and $K = 2$.
		
		Given a sample point $\textbf{x}$, note that we can write
		$$\textbf{x} = \textbf{u}_h + \textbf{z},$$
where $\textbf{z}$ is a random vector with mean 0 and covariance $\sigma^2 \textbf{I}_d$. We summarise the main results in the next theorem whose proof can be found in \cite{anandkumar}.

		\begin{theorem}[Hsu and Kakade, 2013] \label{moments} 
			Assume $d \geq K$. The variance $\sigma^2$ is the smallest eigenvalue of the covariance matrix $\mathbb{E}[ \textbf{x} \otimes \textbf{x} ] - \mathbb{E}[ \textbf{x} ] \otimes \mathbb{E}[ \textbf{x} ]$. Furthermore, if		
			\begin{flalign*}
				& M_1 = \mathbb{E}[ \textbf{x} ],\\
				& M_2 = \mathbb{E}[ \textbf{x} \otimes \textbf{x} ] - \sigma^2 \textbf{I}_d,\\
				& M_3 = \mathbb{E}[ \textbf{x} \otimes \textbf{x} \otimes \textbf{x} ] - \sigma^2 \sum_{i=1}^d \left( \mathbb{E}[ \textbf{x} ] \otimes \textbf{e}_i \otimes \textbf{e}_i + \textbf{e}_i \otimes \mathbb{E}[ \textbf{x} ] \otimes \textbf{e}_i + \textbf{e}_i \otimes \textbf{e}_i \otimes \mathbb{E}[ \textbf{x} ] \right),
			\end{flalign*}
then
			\begin{flalign*}
				& M_1 = \sum_{i=1}^K w^{(i)}\ \textbf{u}^{(i)},\\
				& M_2 = \sum_{i=1}^K w^{(i)}\ \textbf{u}^{(i)} \otimes \textbf{u}^{(i)},\\
				& M_3 = \sum_{i=1}^K w^{(i)}\ \textbf{u}^{(i)} \otimes \textbf{u}^{(i)} \otimes \textbf{u}^{(i)}.
			\end{flalign*}
		\end{theorem} 
		
		Theorem~\ref{moments} allows us to use the method of moments, which is a classical parameter estimation technique from statistics. This method consists in computing certain statistics of the data (often empirical moments) and use it to find model parameters that give rise to (nearly) the same corresponding population quantities. Now suppose that $N$ is large enough so we have a reasonable number of sample points to make useful statistics. First we compute the empirical mean 
		
		\begin{equation}		
			\hat{\mu} := \frac{1}{N} \sum_{j=1}^N \textbf{x}^{(j)} \approx \mathbb{E}[ \textbf{x} ]. \label{empirical_mean}
		\end{equation}
		
		Now use this result to compute the empirical covariance matrix
		
		\begin{equation}
			\hat{\textbf{S}} := \frac{1}{N} \sum_{j=1}^N ( \textbf{x}^{(j)} \otimes \textbf{x}^{(j)} - \hat{\mu} \otimes \hat{\mu} ) \approx \mathbb{E}[ \textbf{x} \otimes \textbf{x} ] - \mathbb{E}[ \textbf{x} ] \otimes \mathbb{E}[ \textbf{x} ]. \label{empirical_cov}
		\end{equation}
		
		The smallest eigenvalue of $\hat{\textbf{S}}$ is the empirical variance $\hat{\sigma}^2 \approx \sigma^2$. Now we compute the empirical third moment (empirical skewness)
		
		 \begin{equation}
			\hat{\mathcal{S}} := \frac{1}{N} \sum_{j=1}^N \textbf{x}^{(j)} \otimes \textbf{x}^{(j)} \otimes \textbf{x}^{(j)} \approx \mathbb{E}[ \textbf{x} \otimes \textbf{x} \otimes \textbf{x} ] \label{empirical_skew}
		\end{equation}
and use it to get the empirical value of $M_3$,
		\begin{equation}
			\hat{\mathcal{M}}_3 := \hat{\mathcal{S}} - \hat{\sigma}^2 \sum_{i=1}^d \left( \hat{\mu} \otimes \textbf{e}_i \otimes \textbf{e}_i + \textbf{e}_i \otimes \hat{\mu} \otimes \textbf{e}_i + \textbf{e}_i \otimes \textbf{e}_i \otimes \hat{\mu} \right) \approx M_3.\label{empirical_M3}
		\end{equation}
		
		By theorem~\ref{moments}, $M_3 = \displaystyle \sum_{i=1}^K w^{(i)}\ \textbf{u}^{(i)} \otimes \textbf{u}^{(i)} \otimes \textbf{u}^{(i)}$, which is a symmetric tensor containing all parameter information we want to find. The idea is, after computing a symmetric CPD for $\hat{\mathcal{M}}_3$, normalize the factors so each vector has unit norm. By doing this we have a tensor of the form
		$$\sum_{i=1}^K \hat{w}^{(i)}\ \hat{\textbf{u}}^{(i)} \otimes \hat{\textbf{u}}^{(i)} \otimes \hat{\textbf{u}}^{(i)}$$
as a candidate to solution. Note that it is easy to make all $\hat{w}^{(i)}$ positive. If some of them is negative, just multiply it by $-1$ and multiply one of the associated vectors also by $-1$. The final tensor is unchanged but all $\hat{w}^{(i)}$ now are positive. For more on this subject we recommend reading \cite{anandkumar}. 

	\section{Second example: Topic models}		
		Consider a set of documents (texts) together with a set of $K$ possible topics, this structured set of documents is called a \emph{corpus}. Each topic can be represented by a number $1 \leq j \leq K$. Additionally, consider that this corpus has $d$ distinct words in its vocabulary and that each document has $L \geq 3$ words. The words are labelled as numbers between 1 and $d$.  
		
		\begin{figure}[h] 
			\centering
			\includegraphics[scale=.42]{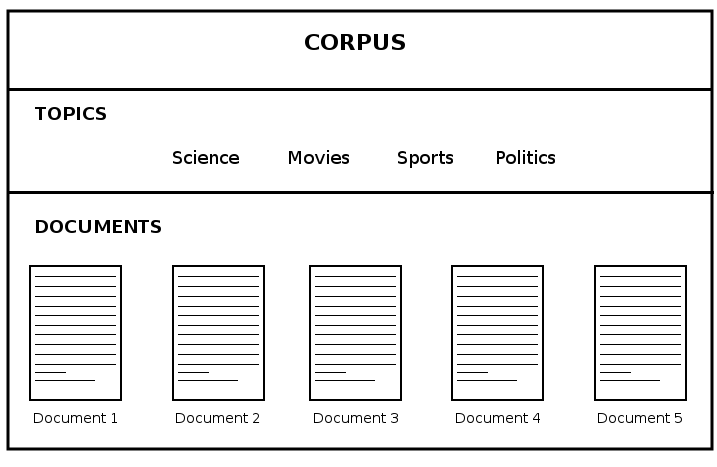}
			\caption{\footnotesize{Example of a corpus structure.}}
			\label{topics}
		\end{figure}	
		
		The \emph{bag-of-words} model \cite{bag_of_words, bag_of_words_2} is a system of representation used in natural language processing (NLP). In this model, any text is represented as the multiset of its words, disregarding grammar and word ordering but keeping multiplicity. This multiset is the ``bag'' containing the words. We consider the bag-of-words model in the problem of document classification, that is, given any text we want a method to classify it into some topic in a prescribed list of topics. 	
		
		Now suppose that the topics of the documents follows a (discrete) probability distribution such that $\mathbb{P}[h = j] = w^{(j)}$ is the probability that a random document belongs to topic $j$ (the random variable $h$ is a latent variable, responsible for the topics assignment). Let $\textbf{w} = [w^{(1)}, \ldots, w^{(K)}]^T$ be the vector of probabilities of the topics. Given the topic $h$ and a random document, the words are assumed to follow a probability distribution such that $\mathbb{P}[\textbf{x} = i|h = j] = u_j^{(i)}$. In other words, given the topic $h = j$, $u_j^{(i)}$ is the probability that a random word $\textbf{x}$ drawn in a document is the word $i$. Denote $\textbf{u}_j = [u_j^{(1)}, \ldots, u_j^{(L)}]^T$ for the vector of probabilities of the words in a document, given the topic $h = j$. Additionally, suppose that randomly drawing sample points from this distribution generates i.i.d. (independent and identically distributed) random variables.
		
		As first observation, we have that $\displaystyle \sum_{j=1}^K w_j = 1$ and $\displaystyle \sum_{i=1}^L u_j^{(i)} = 1$ for any topic $j$. We also will convert the words into vectors, that is, each word $i$ is represented by the canonical basis vector $\textbf{e}_i \in \mathbb{R}^d$. One advantage of this encoding is that the moments of these random vectors correspond to probabilities over words. Consider a document with words $\textbf{x}^{(1)}, \textbf{x}^{(2)}, \ldots, \textbf{x}^{(L)}$, then we have that
		
		$$\mathbb{E}[\textbf{x}^{(i')} \otimes \textbf{x}^{(j')}] = \sum_{i,j=1}^d \mathbb{P}[ \textbf{x}^{(i')} = \textbf{e}^{(i)}, \ \textbf{x}^{(j')} = \textbf{e}^{(j)} ] \ \textbf{e}^{(i)} \otimes \textbf{e}^{(j)} = $$
		
		$$ = \sum_{i,j=1}^d \mathbb{P}[ i'\text{-th word} = i, \ j'\text{-th word} = j ] \ \textbf{e}^{(i)} \otimes \textbf{e}^{(j)} = $$
		
		$$ = 
		\left[ \begin{array}{ccc}
			\mathbb{P}[ i'\text{-th word} = 1, \ j'\text{-th word} = 1 ] & \ldots & \mathbb{P}[i'\text{-th word} = 1, \ j'\text{-th word} = d ] \\
			\vdots & & \vdots \\
			\mathbb{P}[ i'\text{-th word} = d, \ j'\text{-th word} = 1 ] & \ldots & \mathbb{P}[ i'\text{-th word} = d, \ j'\text{-th word} = d ] 
		\end{array} \right].$$
		
		More generally, the entry $(i_1, i_2, \ldots, i_L)$ of the tensor $\mathbb{E}[\textbf{x}^{(1)} \otimes \textbf{x}^{(2)} \otimes \ldots \otimes \textbf{x}^{(L)}]$ is 
		$$\mathbb{P}[ \text{first word} = i_1, \ \text{second word} = i_2, \ldots, \ L\text{-th word} = i_L ].$$
We also remark that the conditional expectation of $\textbf{x}^{(i')}$ given $h$ is simply $\textbf{u}^{(j)}$. More precisely,
		
		$$\mathbb{E}[\textbf{x}^{(i')} | h = j] = \sum_{i=1}^d \mathbb{P}[i' \text{-th word} = i | h=j] \textbf{e}^{(i)} = \sum_{i=1}^d u_i^{(j)} \textbf{e}^{(i)} = \textbf{u}^{(j)}.$$
Since the words are conditionally independent given the topic, we can use this property with conditional moments. More precisely, we have that
		
		$$\mathbb{E}[ \textbf{x}^{(i')} \otimes \textbf{x}^{(j')} | h = j ] = \mathbb{E}[ \textbf{x}^{(i')} | h = j]  \otimes \mathbb{E}[ \textbf{x}^{(j')} | h = j ] = \textbf{u}^{(j)} \otimes \textbf{u}^{(j)}.$$
		
		With similar calculations we obtain the following theorem \cite{anandkumar, anandkumar2}.
				
		\begin{theorem}[Anandkumar et al., 2012]
			Let $1 \leq i', j', k' \leq L$. If		
			\begin{flalign*}
				& M_1 = \mathbb{E}[ \textbf{x}^{(i')} ]\\
				& M_2 = \mathbb{E}[ \textbf{x}^{(i')} \otimes \textbf{x}^{(j')} ]\\
				& M_3 = \mathbb{E}[ \textbf{x}^{(i')} \otimes \textbf{x}^{(j')} \otimes \textbf{x}^{(k')} ],
			\end{flalign*}
then
			\begin{flalign*}
				& M_1 = \sum_{i=1}^K w^{(i)}\ \textbf{u}^{(i)}\\
				& M_2 = \sum_{i=1}^K w^{(i)}\ \textbf{u}^{(i)} \otimes \textbf{u}^{(i)}\\
				& M_3 = \sum_{i=1}^K w^{(i)}\ \textbf{u}^{(i)} \otimes \textbf{u}^{(i)} \otimes \textbf{u}^{(i)}.
			\end{flalign*}
		\end{theorem}
		
		We may estimate these moments using the actual data at hand. By using any method to compute a CPD for $M_3$ we get estimates for the latent variables $w^{(i)}$ and $\textbf{u}^{(i)}$. In many aspects this model resembles the Gaussian mixture model, both use the method of moments to construct a third order tensor which we try to approximate with the CPD.
		
	\section{Third example: Approximation of functions} \label{func-grid}
		Consider a multivariate function $\varphi: \mathbb{R}^L \to \mathbb{R}$ which are difficult to handle analytically, but evaluating $\varphi$ is feasible. Hence, it is possible to take numerical data to study this function. One may try to consider evaluating $\varphi(x^{(1)}, \ldots, x^{(L)})$ in a closed grid of points such that each direction is partitioned in $I_\ell$ parts. This mean we have the points $x_1^{(\ell)}, \ldots, x_{I_\ell}^{(\ell)}$ to consider in the direction of the $\ell$-th coordinate. Overall we will compute $\varphi(x_{i_1}^{(1)}, \ldots, x_{i_L}^{(L)})$ for all $i_1 = 1 \ldots I_1, \ \ldots, \ i_L = 1 \ldots I_L$. 
		
		The results of the computation can be stored in a tensor $\mathcal{T} \in \mathbb{R}^{I_1 \times \ldots \times I_L}$ such that $t_{i_1 \ldots i_L} = \varphi(x_{i_1}^{(1)}, \ldots, x_{i_L}^{(L)})$. Although this is possible, note that as $L$ increases, the \emph{curse of dimensionality}\footnote{\url{https://en.wikipedia.org/wiki/Curse_of_dimensionality}} becomes apparent so storing the results this way requires too much memory. For instance, if $L = 50$ and the dimensions are small, $I_1 = I_2 = \ldots = I_{50} = 2$, storing $\mathcal{T}$ would require 9 petabytes. To overcome this problem we must to store $\mathcal{T}$ in a more economic form, and this is possible with a low rank CPD approximation. Assume that $\varphi$ is \emph{separable}, in the sense that there are functions $\varphi_r^{(\ell)}: \mathbb{R}^{I_\ell} \to \mathbb{R}$, for $\ell=1 \ldots L$ and $r = 1 \ldots R$, such that
		
		$$\varphi(x_{i_1}^{(1)}, \ldots, x_{i_L}^{(L)}) = \sum_{r=1}^R \varphi_r^{(1)}(x_{i_1}^{(1)}) \cdot \varphi_r^{(2)}(x_{i_2}^{(2)}) \cdot \ldots \cdot \varphi_r^{(1)}(x_{i_L}^{(L)}).$$		
Define $\textbf{W}^{(\ell)} = \left[ \textbf{w}_1^{(\ell)}, \ldots, \textbf{w}_R^{(L)} \right] \in \mathbb{R}^{I_\ell \times R}$, where $\textbf{w}_r^{(\ell)} = \left[ \varphi_r^{(\ell)}(x_1^{(\ell)}), \ldots, \varphi_r^{(\ell)}(x_{I_\ell}^{(\ell)}) \right]^T \in \mathbb{R}^{I_\ell}$ is the $r$-th column of $\textbf{W}^{(\ell)}$. Then we have the equality $\mathcal{T} = \displaystyle \sum_{r=1}^R \textbf{w}_r^{(1)} \otimes \ldots \otimes \textbf{w}_r^{(L)}$, which come from a CPD for $\mathcal{T}$. Storing this CPD costs $R \displaystyle \sum_{\ell=1}^L I_\ell$ floats, which is much better than $\displaystyle \prod_{\ell=1}^L I_\ell$ floats necessary to store the tensor in coordinate-wise format. The formulation of this problem and the proposed approach to solve it is based on \cite{grasedyck, hackbusch}. Note that this approach still works if $\varphi$ is not exactly separable but can be well approximated by a separable function. In fact this is what usually happens, one wants to approximate a function of $L$ variables by a finite sum of products of functions of one variable.

	\section{Results}
		In order to obtain a CPD for the tensor it is clear that one would want to know its rank in the first place. Unfortunately this problem is known to be NP-hard \cite{lim2}. Usually one already have some prior knowledge of the problem or, in the worst case, one have to compute several CPDs for different ranks to find the best fit. This second approach seems reasonable, however it is of limited use due to the \emph{border rank} phenomenon, which will be further discussed. 
		
		In this work we always suppose the rank is known in advance, or at least a decent estimate for the rank is known. In this case all we have to worry is with the computation of the CPD. In the past years several algorithms were proposed and implemented \cite{cichoki2013, comon2, tensorlab, tensorlab2, tensortoolbox, cp_opt, tensorly, breiding, rgh} so, today, we have a better understanding about how each approach performs. In particular, Gauss-Newton algorithms are proven to have better convergence properties, also verified experimentally. Algorithms based on this approach usually were much slower \cite{tomasi, tomasi2}, but this is not the reality today. Exploiting the structure of the approximated Hessian matrix lead to algorithms competitive in terms of speed, and with better accuracy \cite{cichoki2013, tensorlab2}. Following this path, in chapter 3 we exploit this structure towards the goal to speeding up conjugate gradient iterations of subproblems faced at each iteration of the Gauss-Newton algorithm. 
		
		While researching about previous implementations for the CPD, I tried to spot the parts overlooked by others. Aspects as damping parameter (when there is regularization), number of conjugate gradient iterations, compression - preprocessing, etc, does not always receives the due attention. With this in mind I designed a new program taking all these little details in account. The result is a tensor package called \emph{Tensor Fox}.\footnote{This package is free for download at \url{https://github.com/felipebottega/Tensor-Fox}. } In chapter 4 we give a detailed description of this package with respect to the CPD computation.	
		
		Let $\mathcal{T}$ be an $L$-order tensor with shape $\underbrace{n \times n \times \ldots \times n}_{L \text{ times}}$. Below we show the cost per iteration (in \emph{flops} - floating point operations) of state of art implementations and Tensor Fox,  to compute a rank-$R$ CPD for $\mathcal{T}$. The constants $c_1, c_2$ are positive integers with little influence on the costs. At first sight it seems that Tensor Toolbox and Tensorly are better, since their cost per iteration is cheaper. In fact it is the opposite, the other solvers indeed make slower iterations, but their iterations have more quality, which leads to faster convergence. Alternating Least squares, for instance, can take thousands of iterations to converge, whereas a Gauss-Newton based algorithm may finish within less than a hundred iterations. The quality of the steps counts. Furthermore, with the exception of Tensor Fox, all the other Gauss-Newton based solvers are costly in the rank. Tensorlab has a factor of $R^3$ and Tensor Box has a factor of $R^6$, whereas Tensor Fox is quadratic on $R$. Computational experiments reinforces these observations.\bigskip
		
		\footnotesize
		\begin{table}
			$\hspace{-1cm}$
		 	\begin{tabular}{|c|c|c|}
				\hline
				\textbf{Package} & \textbf{Algorithm} & \textbf{Computational cost}\\
				\hline
				\footnotesize Tensorlab & \footnotesize Gauss-Newton & \footnotesize$\mathcal{O}\left( 2(L + 1)Rn^L + c_1(\frac{5}{2}L^2R^2 +8LR^2n + \frac{1}{3}LR^3) \right)$\\
				\hline
				\footnotesize Tensor Toolbox & \footnotesize Gradient-based optimization & \footnotesize$\mathcal{O}\left( LRn^L \right)$\\
				\hline
				\footnotesize Tensorly & \footnotesize Alternating Least Squares & \footnotesize$\mathcal{O}\left( LRn^L \right)$\\
				\hline
				\footnotesize Tensor Box & \footnotesize Gauss-Newton & \footnotesize$\mathcal{O}\left( L^3R^2 + LR^3 + LRn^L + L^3n + L^3R^6 \right)$\\
				\hline
				\footnotesize Tensor Fox & \footnotesize Gauss-Newton & \footnotesize$\mathcal{O}\left( (LR+L-1)n^L + LR^2 ( 1 + n ) + 3LRn + c_2 (9 LR + L^2 R^2)n \right)$\\
				\hline
			\end{tabular}
			\caption{\footnotesize{Cost per iteration of several CPD solvers.}} \bigskip
		\end{table} 		
		\normalsize		
		
		The computational complexity to compute CPDs makes it hard to aim at really big tensors, because of the factor $n^L$ present in all algorithms. It is the curse of dimensionality in action. In the era of Big Data it is not enough to just have good tensor models, they also need to be computable within a reasonable time. In chapter 4, section 4.6, we show that it is enough to use the Gauss-Newton approach only for third order tensors. With the ideas of \cite{cpd_tensortrain} we are able to compute higher order CPDs while avoiding the curse of dimensionality. The \emph{tensor train decomposition} (TTD) was recently linked to the CPD, providing a way to compute higher order CPDs much faster than any previous algorithm. These new ideas are implemented in Tensor Fox, which leads to a cost of
		$$\mathcal{O}\left( (3R+2)n^3 + 3R^2(1+n) + 9Rn + 9c_2 (3R + R^2)n \right)$$
flops per iteration for higher order tensors. Of course this is not all. This is the cost per iteration of one third order CPD to be computed, between $L$ of them. Additionally, before the computation of these third order CPDs we have to compute $L-$ SVDs, which adds a cost of 
		$$\mathcal{O}\left( 2 \left( n^{L+1} + R^2 \frac{n^L - n^2}{n-1} + (1 + R^3) n^3 \right) \right) \text{ flops}.$$ 
Therefore we didn't avoid completely the curse of dimensionality. On the other hand, we remark that this cost has a low constant and it is added only once to the overall cost, while the costs of the previous showed algorithms are added at each iteration. The tensor train approach performed much better than any other algorithm in our tests.  

		We remark that only Tensorlab performs tensor compression before the iterations. In chapter 2 we take a closer look at tensor compression and show why this is crucial to alleviate the curse of dimensionality too. At the moment, most solvers consider compression as an optional action to take, but this should be default. For example, if we have a tensor of shape $n \times n \times n$ and want to compute an approximate CPD of rank $R \ll n$, then it is possible to compress it to a tensor of shape $R \times R \times R$ and use this one to find the CPD. There is virtually no loss in precision and the cost of doing that is of $\mathcal{O} \left( \displaystyle \sum_{\ell=1}^3 \min\{R, n\} \prod_{\ell=1}^3 n \right) = \mathcal{O} \left( 3 R n^3 \right)$ flops. Compare this to the previous costs, where we have something of $\mathcal{O} \left( 3 R n^3 \right)$ flops \textbf{at each iteration}. We try to stress this point here with the known results of the area and computational experiments.  
		
		Tensor decompositions are amazing tools to model multidimensional data, and that is why developing new algorithms to compute the CPD is necessary in this area. This work is an attempt to improve the state of art overall performance. The main contributions of this work are the following: 
		\begin{itemize}
			\item The diagonal regularization introduced at~\ref{diag_reg} reduces substantially the condition number compared with the other regularization approaches used in the literature.
			\item The approximated Hessian of the problem has a block structure which is exploited in theorem~\ref{Hv} to accelerate any algorithm based on Krylov methods to solve the normal equations of the Gauss-Newton step.
			\item Tensor Fox performs specially better for higher order tensors. This is possible with the CPD Tensor Train technique developed at \cite{cpd_tensortrain} and implemented in Tensor Fox.
			\item A development of a new algorithm combining the best parts of several state of art algorithms was implemented in tensor Fox. Improvements includes a method to exploit the approximated Hessian.
			\item A new tensor package software, Tensor Fox, which is competitive and freely available to the interested practitioners and researches. I remark that all routines of Tensor Fox were written from scratch, that is, not a single part of other tensor package was copied. Routines are optimized for speed.
			\item Several benchmarks are introduced as an attempt to obtain a fair comparison between all the state of art CPD solvers for a range of distinct problems. Usually the papers makes comparisons between the one they are introducing and just one or two outside solvers. This seems to be the first time such a broad comparison is made. 
			\item New tensor models for machine learning problems are introduced and their potential is experimentally validated.
		\end{itemize}

%% file: CHAPTER_1.tex
\chapter{Basic notions}\label{cap-1} 
	We introduce the necessary notations and preliminary results of multilinear algebra. After that we will talk about tensor products, decompositions and the geometry of tensor spaces. In this chapter we also formalize the main challenge of this work, which is to compute the canonical polyadic decomposition. 

	\section{Notations} 
		 Scalars will be denoted by lower case letters, including greek letters, e.g., $a$ or $\lambda$. Sometimes we can use capital letters for natural numbers, e.g., $L$ or $C$. Vectors are denoted by bold lower case letters, e.g., $\textbf{x}$. Matrices are denoted by bold capital letters, e.g., $\textbf{X}$. Tensors are denoted by calligraphic capital letters, e.g., $\mathcal{T}$. Capital greek letters will be more flexible, appearing sometimes as matrices, sometimes as tensors, and sometimes as sets. The $i$-th entry of a vector $\textbf{x}$ is denoted by $x_i$, the entry $(i,j)$ of a matrix $\textbf{X}$ is denoted by $x_{ij}$, and the entry $(i_1, i_2, \ldots, i_L)$ of a tensor $\mathcal{T}$ with $L$ indexes is denoted by $t_{i_1  i_2 \ldots i_L}$. Sometimes it will be necessary to denote the entry $(i,j)$ of $\textbf{X}$ by $(\textbf{X})_{ij}$, and the same may happen to a tensor. Any kind of sequence will be indicated by superscripts. For example, we write $\textbf{x}^{(1)}, \textbf{x}^{(2)}, \ldots$ for a sequence of vectors. The $n \times n$ identity matrix will be denoted by $\textbf{I}_n$. 
		 
		 In the case we have a function $f$ with $n$ scalar arguments, we denote them by $x_1, x_2, \ldots, x_n$ and write $f(x_1, x_2, \ldots, x_n)$. If there is just three arguments we prefer the classical $f(x, y, z)$, and similar considerations for two or just one argument, where we will use $f(x, y)$, and $f(x)$, respectively. 
		 
		 Vector spaces, groups and subsets in general will be denoted by blackboard bold capital letters or just capital letters, e.g., $\mathbb{V}$ or $S$. An important particular case is of a field, which will be denoted by $\mathbb{K}$. However, this work is limited to the cases $\mathbb{K} = \mathbb{R}$ (real numbers) and $\mathbb{K} = \mathbb{C}$ (complex numbers). In this work it will convenient to define the set of natural numbers as being the set $\mathbb{N} = \{1, 2, 3, \ldots \}$. The symbols $\mathbb{P}$ and $\mathbb{E}$ are reserved for probability and expectation, respectively. Every time we introduce a basis, we will be using the set notation with the implicit understanding it is an ordered set. Given two tuples $(r_1, \ldots, r_L), (R_1, \ldots, R_L)$, we write $(r_1, \ldots, r_L) < (R_1, \ldots, R_L)$ if $r_i < R_i$ for all $i = 1 \ldots L$.
		 
		  We also adopt the Matlab notational style when it is desired to take slices\footnote{Since the author uses much more Python/Numpy than Matlab there is a chance to appear some slight differences.} or to fix a subset of indexes while varying the others. For example, if $\textbf{X}$ is a $m \times n$ matrix, then $X_{i:}$ is the $i$-th row of $\textbf{X}$, while $\textbf{X}_{:j}$ is its $j$-th column. 
		 
		The symbol $^T$ denotes the transpose of a vector or matrix, $^\ast$ denotes the conjugate transpose of a vector or matrix, $\dagger$ denotes the pseudoinverse of a matrix. If we use $^\ast$ for a vector space, then it means the dual space. For example, if $\mathbb{V}$ is a vector space, then $\mathbb{V}^\ast$ is the dual space of $\mathbb{V}$. The symbol $\cong$ will be used to denote isomorphism between vector spaces.	Let $\left\{ \textbf{e}_1, \ldots, \textbf{e}_n \right\}$ be a basis of $\mathbb{V}$ and let $\textbf{x} = x_1 \textbf{e}_1 + \ldots + x_n \textbf{e}_n, \textbf{y} = y_1 \textbf{e}_1 + \ldots y_n \textbf{e}_n$ be two vectors in this space. The \emph{Hermitian inner product} between $\textbf{x}$ and $\textbf{y}$ is defined by $\langle \textbf{x}, \textbf{y} \rangle = x_1 \overline{y_1} + \ldots + x_n \overline{y_n}$. We also adopt the definition $\langle \textbf{x}, \textbf{y} \rangle_\mathbb{R} = x_1 y_1 + \ldots + x_n y_n$ and call it the \emph{Euclidean inner product}. 
		
		Finally, consider the Euclidean vector space $\mathbb{K}^n$ with basis $\{ \textbf{e}_1, \ldots, \textbf{e}_n \}$ and dual basis $\{ \textbf{f}_1^\ast, \ldots, \textbf{f}_n^\ast \}$, where $\textbf{f}_i \in \mathbb{K}^n$ for each $i = 1 \ldots n$. In this context, for any vector $\textbf{x} \in \mathbb{K}^n$ we define $\textbf{f}_i^\ast(\textbf{x}) = \textbf{f}_i^\ast \cdot \textbf{x} = \langle \textbf{x}, \textbf{f}_i \rangle$. If the basis is orthonormal, then $\textbf{f}_i^\ast = \textbf{e}_i^\ast$.

	\section{Multilinear maps}
		Let $\mathbb{V}^{(1)}, \ldots, \mathbb{V}^{(L)}, \mathbb{V}, \mathbb{W}$ be vector spaces over the same field $\mathbb{K}$ such that $\dim (\mathbb{V}^{(\ell)}) = I_\ell$ for each $\ell = 1 \ldots L$, and $\dim (\mathbb{W}) = J$. A map $\mathcal{T}:\mathbb{V}^{(1)} \times \ldots \times \mathbb{V}^{(L)} \to \mathbb{W}$ is said to be \emph{$L$-linear} if $\mathcal{T}$ is linear in each coordinate. More precisely, for all $(\textbf{x}^{(1)}, \ldots, \textbf{x}^{(L)}) \in \mathbb{V}^{(1)} \times \ldots \times \mathbb{V}^{(L)}$ and all $\alpha, \beta \in \mathbb{K}$ we have that		
	$$\mathcal{T}(\textbf{x}^{(1)}, \ldots, \alpha \textbf{x}^{(i)} + \beta \textbf{x}^{(i+1)}, \ldots, \textbf{x}^{(L)}) = $$
	$$ = \alpha \mathcal{T}(\textbf{x}^{(1)}, \ldots, \textbf{x}^{(i)}, \ldots, \textbf{x}^{(L)}) + \beta \mathcal{T}(\textbf{x}^{(1)}, \ldots, \textbf{x}^{(i+1)}, \ldots, \textbf{x}^{(L)}).$$  
	
		Sometimes it is not relevant to mention the value $L$ and one can just say that $\mathcal{T}$ is a \emph{multilinear map}. Notice that for $L = 1$, $\mathcal{T}$ is just a classical linear map. We denote the space of $L$-linear maps $\mathbb{V}^{(1)} \times \ldots \times \mathbb{V}^{(L)} \to \mathbb{W}$ by $\mathcal{L}(\mathbb{V}^{(1)}, \ldots, \mathbb{V}^{(L)}; \mathbb{W})$. Some common abbreviations are $\mathcal{L}_L(\mathbb{V}; \mathbb{W}) = \mathcal{L}( \underbrace{\mathbb{V}, \ldots, \mathbb{V}}_{L \text{ times}}; \mathbb{W} )$ and $\mathcal{L}(\mathbb{V}) = \mathcal{L}(\mathbb{V}; \mathbb{V})$. 
	
		\begin{lemma} \label{isomorphism}
			Let $(i_1, \ldots, i_p, i_1', \ldots, i_{L-p}')$ be any permutation of $(1, \ldots, L)$. Then			
			$$\mathcal{L}(\mathbb{V}^{(1)}, \ldots, \mathbb{V}^{(L)}; \mathbb{W}) \cong \mathcal{L}\left( \mathbb{V}^{(i_1)},\ldots, \mathbb{V}^{(i_p)}; \mathcal{L}(\mathbb{V}^{(j_1)}, \ldots, \mathbb{V}^{(j_{L-p})}; \mathbb{W}) \right).$$
			
		\end{lemma}
	
		\begin{corollary} \label{isomorphism2}
			$\mathcal{L}(\mathbb{V}^{(1)}, \ldots, \mathbb{V}^{(L)}; \mathbb{W}) \cong \mathcal{L}(\mathbb{V}^{(1)}, \ldots, \mathbb{V}^{(L)}, \mathbb{W}^\ast; \mathbb{K})$.
		\end{corollary}
	
		This corollary is a direct consequence of lemma~\ref{isomorphism} because 
	
	 	$$\mathcal{L}(\mathbb{V}^{(1)}, \ldots, \mathbb{V}^{(L)}, \mathbb{W}^\ast; \mathbb{K}) \cong \mathcal{L}\left( \mathbb{V}^{(1)}, \ldots, \mathbb{V}^{(L)}; \mathcal{L}(\mathbb{W}^\ast; \mathbb{K}) \right) \cong \mathcal{L}(\mathbb{V}^{(1)}, \ldots, \mathbb{V}^{(L)}; \mathbb{W}).$$
	
		With corollary~\ref{isomorphism2} we are able to concentrate our attention to multilinear maps of the form $\mathcal{T}:\mathbb{V}^{(1)} \times \ldots \times \mathbb{V}^{(L)} \to \mathbb{K}$. In the case of a linear map $\mathcal{T}:\mathbb{K}^n \to \mathbb{K}$ (linear functional), we know there is a vector $\textbf{a} \in \mathbb{K}^n$ such that 
		
		\begin{equation} \label{order1}
			\mathcal{T}(\textbf{x}) =  \textbf{a}^T \textbf{x}
		\end{equation} 		
for all $\textbf{x} \in \mathbb{K}^n$. In the case of a 2-linear (bilinear) map $\mathcal{T}:\mathbb{K}^m \times \mathbb{K}^n \to \mathbb{K}$, there is a matrix $\textbf{A} \in \mathbb{K}^{n \times m}$ such that\footnote{In the complex case there is the notion of a sequilinear form, which is a map $(\textbf{x}, \textbf{y}) \mapsto \textbf{y}^\ast \textbf{A} \textbf{x}$. Sesquilinear forms and bilinear complex forms are not the same thing.} 

	\begin{equation} \label{order2}	
		\mathcal{T}(\textbf{x}, \textbf{y}) = \textbf{y}^T \textbf{A} \textbf{x} 
	\end{equation}
for all $\textbf{x} \in \mathbb{K}^m, \textbf{y} \in \mathbb{K}^n$. If one want a general formula for 3-linear (trilinear) maps $\mathcal{A}:\mathbb{K}^m \times \mathbb{K}^n \times \mathbb{K}^p \to \mathbb{K}$ or more, the concept of tensors is a must. In order to have a better understanding of what is happening it is convenient to work in coordinates after fixing a basis for each $\mathbb{V}^{(\ell)}$. 
	
		\begin{theorem} \label{theorem_coordinates}
			Let $\{ \textbf{e}_1^{(\ell)}, \ldots, \textbf{e}_{I_\ell}^{(\ell)} \}$ be a basis for each $\mathbb{V}^{(\ell)}$ and let $\mathcal{T} \in \mathcal{L}(\mathbb{V}^{(1)}, \ldots, \mathbb{V}^{(L)}; \mathbb{K})$ be a $L$-linear map. Then there exists scalars $t_{i_1 \ldots i_L} \in \mathbb{K}$ such that
			$$\mathcal{T}(\textbf{e}_{i_1}^{(1)}, \ldots, \textbf{e}_{i_L}^{(L)}) = t_{i_1 \ldots i_L},$$
for $i_1 = 1, \ldots, I_1$, $\ i_2 = 1, \ldots, I_2$, $\ \ldots,\ i_L = 1, \ldots, I_L$.
		\end{theorem}
		
		The values $t_{i_1 \ldots i_L}$ are called the \emph{coordinates of $\mathcal{T}$} with respect to the given bases. 	As it happens for linear maps and matrices, once we have fixed bases it is possible to associated the multilinear map $\mathcal{T}$ with the coordinates $t_{i_1 \ldots i_L}$. This is the same thing we do with matrices, considering them as a static table of numbers or as a linear transformation, depending on the context. So one can identify $\mathcal{T}$ with its coordinate representation. In the case $\mathcal{T}:\mathbb{K}^n \to \mathbb{K}$ we have that 
		$$\mathcal{T} = \left[ 
		\begin{array}{c}
			t_1\\ 
			\vdots\\
			t_n
		\end{array}
		\right],$$
in the case $\mathcal{T}:\mathbb{K}^m \times \mathbb{K}^n \to \mathbb{K}$ we have	
		$$\mathcal{T} = \left[ 
		\begin{array}{ccc}
			t_{11} & \ldots & t_{1 n}\\
			\vdots & & \vdots\\
			t_{m 1} & \ldots & t_{m n}
		\end{array}	
		\right],$$
and in the case $\mathcal{T}:\mathbb{K}^m \times \mathbb{K}^n \times \mathbb{K}^p \to \mathbb{K}$ we have the ``rectangular matrix'' below.

		\begin{figure}[h]
			\centering
			\includegraphics[scale=.55]{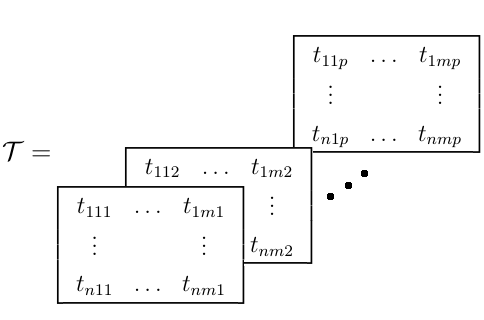}
			\caption{\footnotesize{Trilinear map in coordinates.}}
		\end{figure}
		
		\begin{remark} \label{matrix_tensor_relation}
			In the case $\mathcal{T}: \mathbb{K}^n \to \mathbb{K}$ the vector $\mathcal{T}$ is related with $\textbf{a}$ by the identity $\mathcal{T} = \textbf{a}^T$, that is, $\mathcal{T}(\textbf{x}) = \textbf{a}^T \cdot \textbf{x}$. In the case $\mathcal{T}:\mathbb{K}^m \times \mathbb{K}^n \to \mathbb{K}$ the matrix $\mathcal{T}$ is $m \times n$ while $\textbf{A}$ is $n \times m$. Both matrices are related by the identity $\mathcal{T} = \textbf{A}^T$. Now we have that $\mathcal{T}(\textbf{x}, \textbf{y}) = \textbf{y}^T \textbf{A} \textbf{x}$.
		\end{remark}
	
		Remember that these coordinate representations are tensor of orders 1,2,3, respectively. Now let's see how one can obtain an explicit formula for $\mathcal{T}( \textbf{x}^{(1)}, \ldots, \textbf{x}^{(L)} )$, where $(\textbf{x}^{(1)}, \ldots, \textbf{x}^{(L)}) \in \mathbb{V}^{(1)} \times \ldots \times \mathbb{V}^{(L)}$ is arbitrary and such that $\textbf{x}^{(\ell)} = x_1^{(\ell)} \textbf{e}_1^{(\ell)} + \ldots + x_{I_\ell}^{(\ell)} \textbf{e}_{I_\ell}^{(\ell)}$ for each $\ell =1, \ldots, L$. As consequence of theorem~\ref{theorem_coordinates} we have that
		$$\mathcal{T}(\textbf{x}^{(1)}, \ldots, \textbf{x}^{(L)}) = \mathcal{T}( x_1^{(1)} \textbf{e}_1^{(1)} + \ldots + x_{I_1}^{(1)} \textbf{e}_{I_1}^{(1)}, \ldots, x_1^{(L)} \textbf{e}_1^{(L)} + \ldots + x_{I_L}^{(L)} \textbf{e}_{I_L}^{(L)} ) = $$
		
	$$ = \sum_{i_1=1}^{I_1} \ldots \sum_{i_L=1}^{I_L} x_{i_1}^{(1)} \ldots x_{i_L}^{(L)} \mathcal{T}(\textbf{e}_{i_1}^{(1)}, \ldots, \textbf{e}_{i_L}^{(L)}) = $$
	
		\begin{equation} \label{formula_coordinates}
			= \sum_{i_1=1}^{I_1} \ldots \sum_{i_L=1}^{I_L} x_{i_1}^{(1)} \ldots x_{i_L}^{(L)} t_{i_1 \ldots i_L}.
		\end{equation}  

		The vectors obtained by fixing all dimensions except one are important and have their own name. 
	
		\begin{definition}	
			Let $1 \leq \ell \leq L$. Then, for each choice of indexes $i_1, \ldots, i_{\ell-1}, i_{\ell+1}, \ldots, i_L$, the vector $\mathcal{T}_{i_1 \ldots i_{\ell-1} \ : \ i_{\ell+1} \ldots i_L}$ is called a mode-$\ell$ fiber of $\mathcal{T}$.
		\end{definition}	
	
\begin{figure}[h]
			\centering
			\includegraphics[scale=.15]{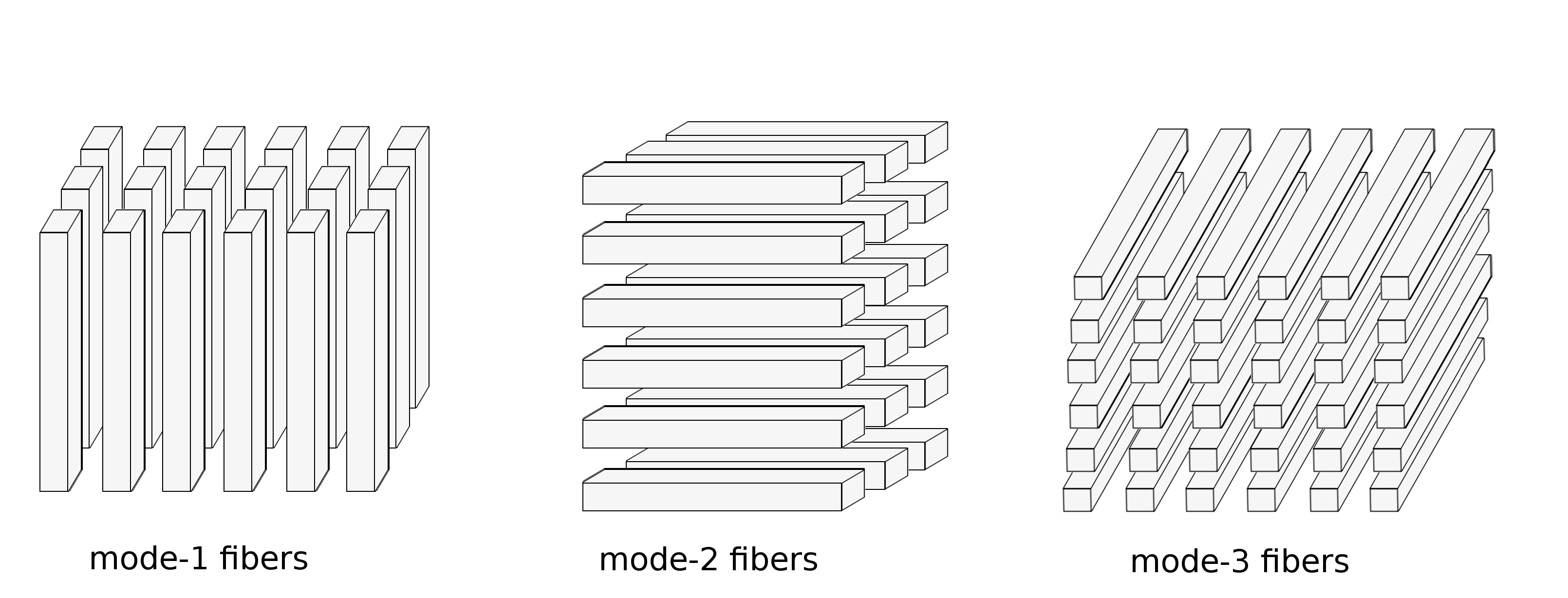}
			\caption{\footnotesize{Fibers of a third order tensor.}}
			\label{fibers}
		\end{figure}

		In the case $L = 1$ the only fiber is $\mathcal{T}$ itself. In the case $L = 2$ the mode-1 fibers are the columns and the mode-2 fibers are the rows of $\mathcal{T}$. The case $L = 3$ is illustrated in figure~\ref{fibers}.
	
		Although we are not going to use this now, it will be important to us the subtensors we obtain by fixing all dimensions except two. This will give rise to bidimensional subtensors, that is, matrices.
		
		\begin{definition}	
			Let $1 \leq \ell < \ell' \leq L$. Then, for each choice of indexes 
			$$i_1, \ldots, i_{\ell-1}, i_{\ell+1}, \ldots, i_{\ell'-1}, i_{\ell'+1}, \ldots, i_L,$$ 
the vector $\mathcal{T}_{i_1 \ldots i_{\ell-1} \ : \ i_{\ell+1} \ldots \ i_{\ell-1} \ : \ i_{\ell+1}, \ldots i_L}$ is called a \emph{slice} of $\mathcal{T}$.
		\end{definition}
		
\begin{figure}[h]
			\centering
			\includegraphics[scale=.15]{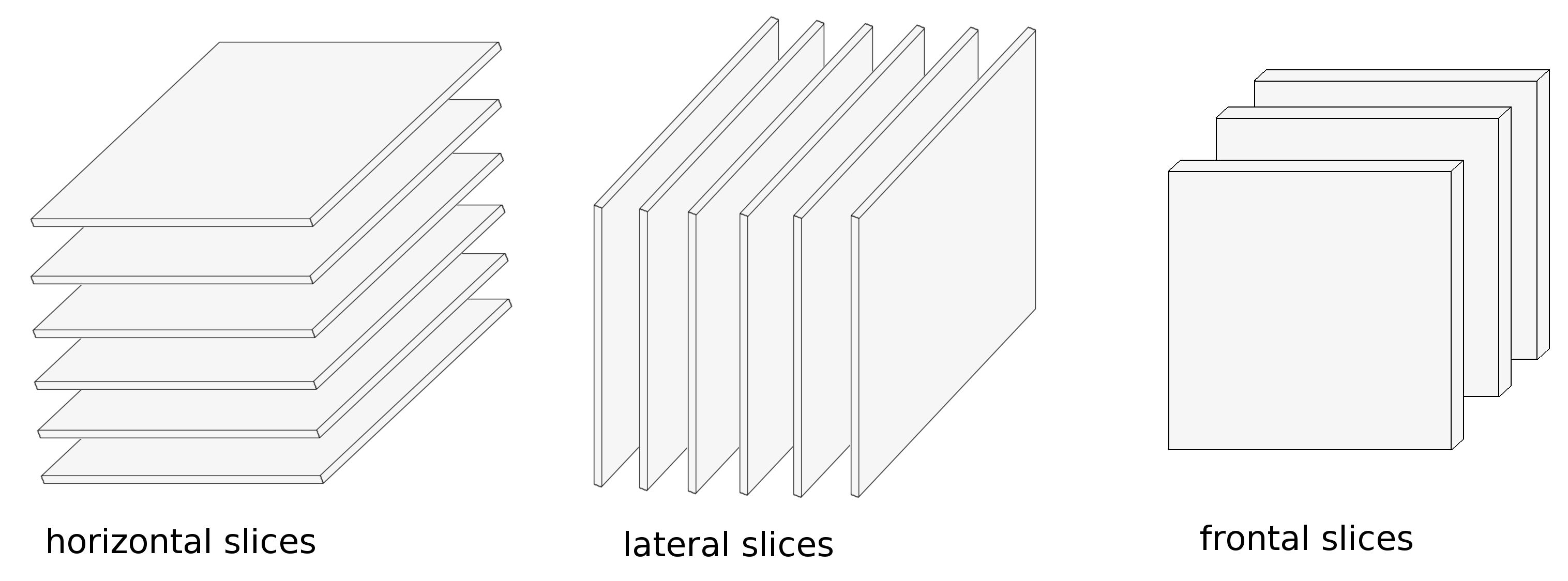}
			\caption{\footnotesize{Slices of a third order tensor.}}
			\label{slices}
		\end{figure}

		In the special case of $\mathcal{T}$ being a third order tensor, we can call the matrices $\mathcal{T}_{i :: }$ the \emph{horizontal slices}, $\mathcal{T}_{: j:}$ the \emph{lateral slices}, and $\mathcal{T}_{::k}$ the \emph{frontal slices}. These types of slices are illustrated in figure~\ref{slices}.
		
	\section{Tensor products}		
		\begin{definition}
			Let $f^{(\ell)} \in (\mathbb{V}^{(\ell)})^\ast$ for each $\ell =1 \ldots L$. The \emph{tensor product} between the functionals $f^{(\ell)}$ is the map $f^{(1)} \otimes \ldots \otimes f^{(L)}: \mathbb{V}^{(1)} \times \ldots \times \mathbb{V}^{(L)} \to \mathbb{K}$ defined as
			
			$$f^{(1)} \otimes \ldots \otimes f^{(L)} (\textbf{x}^{(1)}, \ldots, \textbf{x}^{(L)}) = f^{(1)}(\textbf{x}^{(1)}) \cdot \ldots \cdot f^{(L)}(\textbf{x}^{(L)}).$$
			
		\end{definition} 	
		
		The linear space generated by all tensor products of the form $f^{(1)} \otimes \ldots \otimes f^{(L)}$ is denoted by $(\mathbb{V}^{(1)})^\ast \otimes \ldots \otimes (\mathbb{V}^{(L)})^\ast$ and called the \emph{tensor product} between the spaces $(\mathbb{V}^{(\ell)})^\ast$. An element of $(\mathbb{V}^{(1)})^\ast \otimes \ldots \otimes (\mathbb{V}^{(L)})^\ast$ can be called a \emph{covariant $L$-tensor} or a \emph{covariant $L$-th order tensor}.
		
		\begin{definition}
			Let $\textbf{v}^{(\ell)} \in \mathbb{V}^{(\ell)}$ for each $\ell =1 \ldots L$. The \emph{tensor product} between the vectors $\textbf{v}^{(\ell)}$ is the map $\textbf{v}^{(1)} \otimes \ldots \otimes \textbf{v}^{(L)}: (\mathbb{V}^{(1)})^\ast \times \ldots \times (\mathbb{V}^{(L)})^\ast \to \mathbb{K}$ defined as
			
			$$\textbf{v}^{(1)} \otimes \ldots \otimes \textbf{v}^{(L)} (f^{(1)}, \ldots, f^{(L)}) = f^{(1)}(\textbf{v}^{(1)}) \cdot \ldots \cdot f^{(L)}(\textbf{v}^{(L)}).$$
			
		\end{definition} 	
		
		The linear space generated by all tensor products of the form $\textbf{v}^{(1)} \otimes \ldots \otimes \textbf{v}^{(L)}$ is denoted by $\mathbb{V}^{(1)} \otimes \ldots \otimes \mathbb{V}^{(L)}$ and called the \emph{tensor product} between the spaces $\mathbb{V}^{(\ell)}$. An element of $\mathbb{V}^{(1)} \otimes \ldots \otimes \mathbb{V}^{(L)}$ can be called a \emph{contravariant $L$-tensor} or a \emph{contravariant $L$-th order tensor}.
		
		One may also work with \emph{mixed tensors}, that is, tensors which are the product of vector spaces and dual spaces. Since the ordering is not so important (lemma~\ref{isomorphism}) we can define a mixed tensor to be an element of the space $(\mathbb{V}^{(1)})^\ast \otimes \ldots \otimes (\mathbb{V}^{(L)})^\ast \otimes \mathbb{W}^{(1)} \otimes \ldots \otimes \mathbb{W}^{(M)}$. These tensors are called \emph{tensors of type $(L,M)$}. In particular, a contravariant $L$-th order tensor is a tensor of type $(L,0)$ and a covariant $M$-th order tensor is a tensor of type $(0,M)$. Generally, one may refer to a tensor product of vector spaces just as a \emph{tensor space}. To finish this set of notations and terminology, when we have tensor products between the same space $\mathbb{V}$, it is common to denote $\mathbb{V}^{\otimes L} = \underbrace{\mathbb{V} \otimes \ldots \otimes \mathbb{V}}_{L \text{ times}}$. The next theorem summarizes the main properties of tensor spaces and their relation to multilinear maps. For more details about the algebra of tensor products, consult appendix~\ref{appenB}. 
		
		\begin{theorem}\label{algebric_results}
			The following statements holds true.
			\begin{enumerate}
				\item $(\mathbb{V}^{(1)})^\ast \otimes \ldots \otimes (\mathbb{V}^{(L)})^\ast \cong \mathcal{L}(\mathbb{V}^{(1)}, \ldots, \mathbb{V}^{(L)}; \mathbb{K})$
				\item $\mathbb{V}^{(1)} \otimes \ldots \otimes \mathbb{V}^{(L)} \cong \mathcal{L}\big( (\mathbb{V}^{(1)})^\ast, \ldots, (\mathbb{V}^{(L)})^\ast; \mathbb{K} \big)$
				\item $\dim\big( (\mathbb{V}^{(1)})^\ast \otimes \ldots \otimes (\mathbb{V}^{(L)})^\ast \big) = \dim(\mathbb{V}^{(1)} \otimes \ldots \otimes \mathbb{V}^{(L)}) = \displaystyle \prod_{\ell=1}^L I_\ell$
				\item $\left\{\textbf{e}_{i_1}^{(1)} \otimes \ldots \otimes \textbf{e}_{i_L}^{(L)}: \ i_1=1 \ldots I_1, \ \ldots, \ i_L=1 \ldots I_L \right\}$ is a basis for $\mathbb{V}^{(1)} \otimes \ldots \otimes \mathbb{V}^{(L)}$
				\item Any tensor $\mathcal{T} \in \mathbb{V}^{(1)} \otimes \ldots \otimes \mathbb{V}^{(L)}$ may be written as
				$$\mathcal{T} = \sum_{i_1=1}^{I_1} \ldots \sum_{i_L=1}^{I_L} t_{i_1 \ldots i_L} \ \textbf{e}_{i_1}^{(1)} \otimes \ldots \otimes \textbf{e}_{i_L}^{(L)},$$
where $t_{i_1 \ldots i_L}$ are the coordinates given in theorem~\ref{theorem_coordinates}.
			\end{enumerate}			 
		\end{theorem}
		
		\begin{remark}
			Sometimes it is useful to consider the isomorphism $\mathbb{V}^{(1)} \otimes \ldots \otimes \mathbb{V}^{(L)} \cong \mathcal{L}\big( (\mathbb{V}^{(1)})^\ast, \ldots, (\mathbb{V}^{(L-1)})^\ast;\ \mathbb{V}^{(L)} \big)$ and consider $\textbf{v}^{(1)} \otimes \ldots \otimes \textbf{v}^{(L)} \in \mathbb{V}^{(1)} \otimes \ldots \otimes \mathbb{V}^{(L)}$ as the map given by		
			\begin{equation} \label{product_isomorphism}
			\textbf{v}^{(1)} \otimes \ldots \otimes \textbf{v}^{(L)} (f^{(1)}, \ldots, f^{(L-1)}) = f^{(1)}(\textbf{v}^{(1)}) \cdot \ldots \cdot f^{(L)}(\textbf{v}^{(L-1)}) \cdot \textbf{v}^{(L)}.
			\end{equation}		 
		\end{remark}\bigskip
		
		\begin{example}
			Consider the space $\mathbb{C}^2$ with basis $B = \{[1,0]^T, [0,-i]^T \} = \{ \textbf{b}_1, \textbf{b}_2 \}$, where $i = \sqrt{-1}$ is the imaginary unit. The dual basis associated to $B$ is $B^\ast = \{\textbf{b}_1^\ast, \textbf{b}_2^\ast \}$. Let $\mathcal{T}:\mathbb{C}^2 \to \mathbb{C}^2$ such that $\mathcal{T}(z, w) = [i z, z +w]^T$. As a tensor, note that $\mathcal{T} \in \left( \mathbb{C}^2 \right)^\ast \otimes \mathbb{C}^2$ is a mixed tensor of type $(1,1)$. 
			
			To compute $\mathcal{T}$ in coordinates, first note that $\mathcal{T}(\textbf{b}_1) = [i,1]^T$ and $\mathcal{T}(\textbf{b}_2) = [0,-i]^T$. On the other hand, by interpreting $\mathcal{T}$ as a tensor product we know that
			
			$$\mathcal{T} = \sum_{j,k=1}^2 t_{jk} \ \textbf{b}_j^\ast \otimes \textbf{b}_k.$$

			Using this formula and the identification given in~\ref{product_isomorphism} we have that 
			
			$$\mathcal{T}(\textbf{b}_1) = t_{11} \textbf{b}_1^\ast(\textbf{b}_1) \cdot \textbf{b}_1 + t_{12} \textbf{b}_1^\ast( \textbf{b}_1) \cdot \textbf{b}_2 + t_{21} \textbf{b}_1^\ast(\textbf{b}_2) \cdot \textbf{b}_1 + t_{22} \textbf{b}_1^\ast(\textbf{b}_2) \cdot \textbf{b}_2 =
			\left[
			\begin{array}{c}
				t_{11}\\
				-i t_{12}
			\end{array}
			\right]$$ 
			
and  

			$$\mathcal{T}(\textbf{b}_2) = t_{11} \textbf{b}_1^\ast(\textbf{b}_2) \cdot \textbf{b}_1 + t_{12} \textbf{b}_2^\ast(\textbf{b}_1) \cdot \textbf{b}_2 + t_{21} \textbf{b}_2^\ast(\textbf{b}_2) \cdot \textbf{b}_1 + t_{22} \textbf{b}_2^\ast(\textbf{b}_2) \cdot \textbf{b}_2 = 
			\left[
			\begin{array}{c}
				-t_{21}\\
				-i t_{22}
			\end{array}
			\right].$$ 
			
			With this we conclude that $t_{11} = i, t_{12} = i, t_{21} = 0, t_{22} = 1$. Therefore we have that			
			$$\mathcal{T} = i \textbf{b}_1^\ast \otimes \textbf{b}_1 + i \textbf{b}_1^\ast \otimes \textbf{b}_2 + \textbf{b}_2^\ast \otimes \textbf{b}_2.$$
			
			There is a little subtlety here. Remember that, by remark~\ref{matrix_tensor_relation}, it is necessary to transpose the coordinate representation of $\mathcal{T}$. After transposing we get the matrix of $\mathcal{T}$ in basis $B$,
			
			$$\mathcal{T} = \left[
			\begin{array}{cc}
				i & 0\\
				i & 1
			\end{array}
			\right].$$
		
			The procedure described here is generalizable to any kind of linear map $\mathbb{K}^m \to \mathbb{K}^n$. It permits one to compute the associated tensor and its associated matrix.		
		\end{example}
		
		It is customary to use upper and lower index to distinguish between contravariant and covariant terms, but this won't be necessary here. In this work we will be only interested in studying tensors in the Euclidean space $\mathbb{K}^{I_1} \otimes \ldots \otimes \mathbb{K}^{I_L}$, and for this reason we will leave aside that index convention. 
		
		Let $\mathcal{T} = \textbf{v}^{(1)} \otimes \ldots \otimes \textbf{v}^{(L)} \in \mathbb{K}^{I_1} \otimes \ldots \otimes \mathbb{K}^{I_L}$. Remember we can consider $\mathcal{T}$ as the map $\textbf{v}^{(1)} \otimes \ldots \otimes \textbf{v}^{(L)}: (\mathbb{K}^{I_1})^\ast \otimes \ldots \otimes (\mathbb{K}^{I_L})^\ast \to \mathbb{K}$ given by
		
		$$\textbf{v}^{(1)} \otimes \ldots \otimes \textbf{v}^{(L)} \big( (\textbf{x}^{(1)})^\ast, \ldots, (\textbf{x}^{(L)})^\ast \big) = \langle \textbf{v}^{(1)}, \textbf{x}^{(1)} \rangle \cdot \ldots \cdot \langle \textbf{v}^{(L)}, \textbf{x}^{(L)} \rangle = $$
		\begin{equation}\label{rank1_map}
			= \sum_{i_1=1}^{I_1} \ldots \sum_{i_L=1}^{I_L} \overline{ x_{i_1}^{(1)} \ldots x_{i_L}^{(L)} } \cdot v_{i_1}^{(1)} \ldots v_{i_L}^{(L)}.
		\end{equation}
from which we conclude that $t_{i_1 \ldots i_L} = v_{i_1}^{(1)} \ldots v_{i_L}^{(L)}$. The values $t_{i_1 \ldots i_L}$ are the coordinates of the $\mathcal{T}$ as a multilinear map and as a tensor. It is also possible to use theorem~\ref{theorem_coordinates} to compute these coordinates. Although we are mainly concerned with tensors in $\mathbb{K}^{I_1} \otimes \ldots \otimes \mathbb{K}^{I_L}$, there are some classic examples of different types of tensors we want to show. We consider the canonical bases for all examples below.
		
		\begin{example}[Rank one matrix]
			Given two vectors $\textbf{v} \in\mathbb{C}^m, \textbf{u} \in \mathbb{C}^n$, consider the linear map with matrix $\textbf{u} \textbf{v}^\ast$. In this example we will see that the tensor associated to this map is $\mathcal{T} = \textbf{v}^\ast \otimes \textbf{u} \in (\mathbb{C}^m)^\ast \otimes \mathbb{C}^n$. First note that 
			$$\mathcal{T}(\textbf{x}, \textbf{y}^\ast) = \textbf{v}^\ast(\textbf{x}) \cdot \textbf{y}^\ast(\textbf{u}) = \langle \textbf{x}, \textbf{v} \rangle \cdot \langle \textbf{u}, \textbf{y} \rangle = \sum_{i=1}^n \sum_{j=1}^m x_i \overline{y}_j \overline{v}_i u_j$$
for all $\textbf{x} \in \mathbb{C}^m, \textbf{y} \in \mathbb{C}^n$. Considering $\mathcal{T}$ in coordinates, we can write $\mathcal{T} = \overline{\textbf{v}} \textbf{u}^T$, where $t_{ij} = \overline{v}_i u_j$. The matrix of the corresponding linear map is the transpose of this matrix (see~\ref{matrix_tensor_relation}), that is, $\textbf{A} = ( \overline{\textbf{v}} \textbf{u}^T)^T = \textbf{u} \textbf{v}^\ast$, as desired. Note that now we can write $\mathcal{T}(\textbf{x}, \textbf{y}^\ast) = \textbf{y}^\ast \textbf{A} \textbf{x}$. We may use the isomorphism $(\mathbb{C}^m)^\ast \otimes \mathbb{C}^n \cong L(\mathbb{C}^m; \mathbb{C}^n)$ from~\ref{product_isomorphism} and reinterpret $\mathcal{T}$ as the map $\mathcal{T}(\textbf{x}) = \langle \textbf{x}, \textbf{v} \rangle \cdot \textbf{u} = \textbf{A} \textbf{x}$.
		\end{example}  
		
		\begin{example}[SVD] \label{tensor_SVD}
			Again, let $\mathcal{T} \in (\mathbb{C}^m)^\ast \otimes \mathbb{C}^n$, but this time suppose there are vectors $\textbf{v}_1, \ldots, \textbf{v}_R \in \mathbb{C}^m$, $\textbf{u}_1, \ldots, \textbf{u}_R \in \mathbb{C}^n$, and scalars $\sigma_1, \ldots, \sigma_R \in \mathbb{C}$ such that $\mathcal{T} = \displaystyle \sum_{r=1}^R \sigma_r \textbf{v}_r^\ast \otimes \textbf{u}_r$ and let $\textbf{A} \in \mathbb{C}^{n \times m}$ be the matrix of the corresponding linear map, as before. Additionally, suppose this is the least $R$ with the property that such decomposition exists. As a consequence we have a SVD for $\textbf{A}$ given by $\displaystyle \textbf{A} = \textbf{U} \Sigma \textbf{V}^\ast = \sum_{r=1}^R \sigma_r \textbf{u}_r \textbf{v}_r^\ast$, where $\textbf{u}_r$ is the $r$-th column of $\textbf{U}$ (left singular vector), $\textbf{v}_r$ is the $r$-th column  $\textbf{V}$ (right singular vector), and $\Sigma = \text{diag}(\sigma_1, \ldots, \sigma_R)$.
		\end{example}
		
		\begin{example}[Matrix multiplication] \label{matrix_mult}
			For this last example, consider the tensor space isomorphism $(\mathbb{C}^{mn})^\ast \otimes (\mathbb{C}^{nl})^\ast \otimes \mathbb{C}^{ml} \cong \mathcal{L}\big( \mathbb{C}^{mn}, \mathbb{C}^{nl}; \mathbb{C}^{ml} \big)$. Each element of $\mathbb{C}^{mn}$ can be thought as a $m \times n$ matrix or a vector of size $mn$. Given a $m \times n$ matrix $\textbf{X}$, let $vec(\textbf{X})$ be the vector obtained by vertically stacking all columns of $\textbf{X}$. We will be identifying $\textbf{X}$ and $vec(\textbf{X})$ when it is convenient. The same considerations goes for the spaces $\mathbb{C}^{nl}$ and $\mathbb{C}^{ml}$. Let $\textbf{e}_{ij} \in \mathbb{C}^{m \times n}$ be the matrix with entry $(i,j)$ equal to 1 and all other entries equal to 0. Note that $\{ \textbf{e}_{11}, \textbf{e}_{12}, \ldots, \textbf{e}_{mn} \}$ is a basis for $\mathbb{C}^{m \times n}$, while $\{ vec(\textbf{e}_{11}), vec(\textbf{e}_{12}), \ldots, vec(\textbf{e}_{mn}) \}$ is the canonical basis of $\mathbb{C}^{mn}$. We will commit a little abuse of notation and use the same notation for the basis vectors of $\mathbb{C}^{nl}$ and $\mathbb{C}^{ml}$.  
	
		Now, define the tensor $\mathcal{T} \in (\mathbb{C}^{mn})^\ast \otimes (\mathbb{C}^{nl})^\ast \otimes \mathbb{C}^{ml}$ by

		$$\mathcal{T} = \sum_{i=1}^m \sum_{j=1}^n \sum_{k=1}^l \textbf{e}_{ij}^\ast \otimes \textbf{e}_{jk}^\ast \otimes \textbf{e}_{ik} = \sum_{i=1}^m \sum_{j=1}^n \sum_{k=1}^l \textbf{e}_{ij}^T \otimes \textbf{e}_{jk}^T \otimes \textbf{e}_{ik}.$$

		Given any matrices $\textbf{X} \in \mathbb{C}^{m \times n}, \textbf{Y} \in \mathbb{C}^{n \times l}$, and using the isomorphism~\ref{product_isomorphism}, we have that

		$$\mathcal{T}(\textbf{X}, \textbf{Y}) = \sum_{i=1}^m \sum_{j=1}^n \sum_{k=1}^l \langle \textbf{X}, \textbf{e}_{ij} \rangle \cdot \langle \textbf{Y}, \textbf{e}_{jk} \rangle \cdot \textbf{e}_{ik} = \sum_{i=1}^m \sum_{j=1}^n \sum_{k=1}^l x_{ij} \cdot y_{jk} \cdot \textbf{e}_{ik} = $$
		
		$$ =  \sum_{i=1}^m \sum_{k=1}^l \left( x_{i1} y_{1k} + \ldots + x_{in} y_{nk} \right) \cdot \textbf{e}_{ik} = \textbf{X} \cdot \textbf{Y}$$
	
		In short, $\mathcal{T}$ is a third order tensor describing the matrix multiplication. Note that we used $mnl$ terms in the summation defining $\mathcal{T}$. It is possible to use less terms, and the problem of finding the minimum number of terms is an open problem in mathematics. For instance, see chapter 1 of \cite{landsberg}.
		\end{example}	 
		
	\section{Canonical polyadic decomposition}
		When $\mathcal{T}$ is of the form $\mathcal{T} = \textbf{v}^{(1)} \otimes \ldots \otimes \textbf{v}^{(L)}$, we saw that $t_{i_1 \ldots i_L} = v_{i_1}^{(1)} \ldots v_{i_L}^{(L)}$. However, note that this formula does not apply for all the tensor space since not all tensors are of this form. An arbitrary tensor in $\mathbb{V}^{(1)} \otimes \ldots \otimes \mathbb{V}^{(L)}$ may be written as 
		
		\begin{equation} \label{cpd}
			\mathcal{T} = \sum_{r=1}^R \textbf{v}_r^{(1)} \otimes \ldots \otimes \textbf{v}_r^{(L)},
		\end{equation}		
where each $\textbf{v}_r^{(\ell)} \in \mathbb{V}^{(\ell)}$ is given by $\textbf{v}_r^{(\ell)} = \big[ v_{1r}^{(\ell)}, v_{2r}^{(\ell)}, \ldots, v_{I_\ell r}^{(\ell)} \big]^T$. It is of interest in applications to decompose $\mathcal{T}$ in this manner, such that $R$ is smallest as possible \cite{bro, cichoki2009, cichoki2014, savas, smilde, anandkumar}. The generalization coordinate representation of $\mathcal{T}$ is now given by

		\begin{equation} \label{tensor_coordinates_R}
			t_{i_1 \ldots i_L} = \sum_{r=1}^R v_{i_1 r}^{(1)} \ldots v_{i_L r}^{(L)}.
		\end{equation}
		
		Formula~\ref{cpd} realizes $\mathcal{T}$ as a sum of $R$ tensor products. For each $\ell=1 \ldots L$, the $\ell$-th \emph{factor matrix} associated to~\ref{cpd} is defined as $\textbf{V}^{(\ell)} = [ \textbf{v}_1^{(\ell)}, \ldots, \textbf{v}_R^{(\ell)} ] \in \mathbb{K}^{I_\ell \times R}$, and each one of its columns are called \emph{factors}. The space $\mathbb{V}^{(\ell)}$ sometimes is referred as the \emph{$\ell$-th mode}. When some definition depends on $\ell$ it is common to use some terminology which specify the current mode.
		
		\begin{definition}
			We say a tensor $\mathcal{T} \in \mathbb{V}^{(1)} \otimes \ldots \otimes \mathbb{V}^{(L)}$ has \emph{rank one} if there exists vectors $\textbf{v}^{(1)} \in \mathbb{V}^{(1)}, \ldots, \textbf{v}^{(L)} \in \mathbb{V}^{(L)}$ such that $\mathcal{T} = \textbf{v}^{(1)} \otimes \ldots \otimes \textbf{v}^{(L)}$.
		\end{definition}
		
		\begin{definition}
			We say a tensor $\mathcal{T} \in \mathbb{V}^{(1)} \otimes \ldots \otimes \mathbb{V}^{(L)}$ has \emph{rank $R$} if $R$ is the smallest number such that $\mathcal{T}$ can be written as a sum of $R$ rank one tensors. In this case we denote $rank(\mathcal{T}) = R$. 
		\end{definition}
		
		Suppose $rank(\mathcal{T}) = R$. Then the decomposition~\ref{cpd} is called a \emph{canonical polyadic decomposition (CPD)} for $\mathcal{T}$. Other known names for the CPD are \emph{PARAFAC} (parallel factors), \emph{CANDECOMP} (canonical decomposition) and \emph{CP decomposition}. In the case $R$ is not the rank of $\mathcal{T}$ we call this decomposition a \emph{rank-$R$ CPD for $\mathcal{T}$}. The first one to propose the notion of rank and this decomposition was Hitchcock \cite{hitchcock} in a work of 1927.
		
		In example~\ref{tensor_SVD} we discussed a connection between the SVD for matrices and tensors. If $\mathcal{T} = \displaystyle \sum_{r=1}^R \sigma_r \textbf{u}_r \textbf{v}_r^\ast$ is a SVD for $\mathcal{T}$, then we can write $\mathcal{T} = \displaystyle \sum_{r=1}^R \sigma_r \textbf{v}_r \otimes \textbf{u}_r^\ast$. In particular, this implies that $rank(\mathcal{T}) = R$. Conversely, if $rank(\mathcal{T}) = R$ and $\mathcal{T} = \displaystyle \sum_{r=1}^R \sigma_r \textbf{v}_r \otimes \textbf{u}_r^\ast$, then $\mathcal{T} = \displaystyle \sum_{r=1}^R \sigma_r \textbf{u}_r \textbf{v}_r^\ast$ is a SVD for $\mathcal{T}$ and, as a matrix, $\mathcal{T}$ has rank $R$. In conclusion, if $\mathcal{T}$ is a second order tensor, then its rank as a tensor equals its rank as a matrix. Furthermore, its CPD and its SVD coincide.  
		
		Given $\mathcal{T} \in \mathbb{V}^{(1)} \otimes \ldots \otimes \mathbb{V}^{(L)}$, we know in advance that $rank(\mathcal{T}) \leq \displaystyle \prod_{\ell=1}^L I_\ell$, since we always can write 
		
		$$\mathcal{T} = \sum_{i_1=1}^{I_1} \ldots \sum_{i_L=1}^{I_L} t_{i_1 \ldots i_L} \ \textbf{e}_{i_1}^{(1)} \otimes \ldots \otimes \textbf{e}_{i_L}^{(L)}.$$
		
		In fact, it is possible to write $\mathcal{T}$ with less terms as the next results shows. From theorem~\ref{isomorphism} we know it is possible to permute the spaces $\mathbb{V}^{(\ell)}$ without problems. Hence there is no loss of generality in considering $I_1 \geq I_2 \geq \ldots \geq I_L$.
		
		\begin{theorem}[Landsberg, \cite{landsberg}] \label{rank-bound}
			Let $\mathbb{V}^{(1)} \otimes \ldots \otimes \mathbb{V}^{(L)}$ be such that $I_1 \geq I_2 \geq \ldots \geq I_L$. Then $rank(\mathcal{T}) \leq \displaystyle \prod_{\ell=2}^{L} I_\ell$ for all $\mathcal{T} \in \mathbb{V}^{(1)} \otimes \ldots \otimes \mathbb{V}^{(L)}$.
		\end{theorem}
		
		\begin{corollary}
			$rank(\mathcal{T}) \leq \min \{ I_1 I_2,\ I_1 I_3,\ I_2 I_3 \}$ for all $\mathcal{T} \in \mathbb{V}^{(1)} \otimes \mathbb{V}^{(2)} \otimes \mathbb{V}^{(3)}$.
		\end{corollary}
		
		Remember the slices of third order tensors we showed in~\ref{slices}. The next result gives formulas for each one of these slices.
		
		\begin{theorem}
			Let $\mathcal{T} = \displaystyle\sum_{r=1}^R \textbf{x}_r \otimes \textbf{y}_r \otimes \textbf{z}_r \in \mathbb{K}^{I_1} \otimes \mathbb{K}^{I_2} \otimes \mathbb{K}^{I_3}$ be a third order tensor with rank $\leq R$, where each $\lambda_r$ is a scalar. Then the $i$-th horizontal slice of $\mathcal{T}$ is given by 
			$$\sum_{r=1}^R x_{ir} \ \textbf{y}_r \otimes \textbf{z}_r = \textbf{Y} \cdot \text{diag}(x_{i1}, \ldots, x_{iR}) \cdot \textbf{Z}^T,$$
the $j$-th lateral slice is given by 
			$$\sum_{r=1}^R y_{jr} \ \textbf{x}_r \otimes \textbf{z}_r = \textbf{X} \cdot \text{diag}(y_{j1}, \ldots, y_{jR}) \cdot \textbf{Z}^T.$$
and the $k$-th frontal slice slice is given by 
			$$\sum_{r=1}^R z_{kr} \ \textbf{x}_r \otimes \textbf{y}_r = \textbf{X} \cdot \text{diag}(z_{k1}, \ldots, x_{kR}) \cdot \textbf{Y}^T.$$
		\end{theorem}	
		
		Now suppose we have a tensor $\mathcal{T}$ with rank $R$ and we want to compute a CPD for $\mathcal{T}$. In practical applications obtaining equality as in formula~\ref{cpd} is not realistic. Usually one is content with an approximation
		
		\begin{equation} \label{cpd_approx}
			\mathcal{T} \approx \sum_{r=1}^R \textbf{v}_r^{(1)} \otimes \ldots \otimes \textbf{v}_r^{(L)}.
		\end{equation}	
		
		There are several algorithms to accomplish this goal and we will discuss some of them in chapter 3. For now we discuss rank properties. For instance, how should one proceed when $rank(\mathcal{T})$ is not known?	In order to obtain a CPD for $\mathcal{T}$ one would want to know its rank in the first place. Unfortunately this problem is known to be NP-hard \cite{lim2}. Another possibility would be to choose a large value $R$, an upper bound for $rank(\mathcal{T})$, and compute a rank-$R$ CPD for $\mathcal{T}$. As we will see soon this is not a good idea because a high rank CPD suffer from lack of uniqueness. In particular, this kind of CPD can overfit the data we are trying to model. The best choice here is to choose a low rank CPD approximation for $\mathcal{T}$. Caution is necessary to not take $R$ too small, because in that case our model will suffer from underfitting (high bias) and accuracy is lost. 
		
		A relevant property of higher order tensors is that their CPD are often unique (in a sense we will make clear soon). This property fail for matrices. For instance, consider a matrix $\mathcal{T} \in \mathbb{K}^{n \times m}$ together with a SVD given by $\mathcal{T} = \textbf{U} \cdot \Sigma \cdot \textbf{V}^\ast$, where $\Sigma = \text{diag}(\sigma_1, \ldots, \sigma_R)$. Making $\textbf{A} = \textbf{U} \cdot \Sigma$ and $\textbf{B} = \textbf{V}$ we have 
		
		$$\mathcal{T} = \textbf{AB}^\ast = \sum_{r=1}^R \textbf{A}_{:r} \textbf{B}_{:r}^\ast.$$ 		
		Notice this is a CPD for $\mathcal{T}$ since it is a sum of $R$ rank one terms. Now let $\textbf{W} \in \mathbb{K}^{R \times R}$ be unitary. Then we have 
		
		$$\mathcal{T} = \textbf{AW} (\textbf{BW})^\ast = \tilde{\textbf{A}}\tilde{\textbf{B}}^\ast = \sum_{r=1}^R \tilde{\textbf{A}}_{:r} \tilde{\textbf{B}}_{:r}^\ast.$$
		Varying $\textbf{W}$ we can obtain infinitely many different CPD's for $\mathcal{T}$. 
		
		Let $\mathcal{T} = \displaystyle \sum_{r=1}^R \mathcal{T}_r$ be a CPD for $\mathcal{T}$, where each $\mathcal{T}_r$ is a rank one term. Also, suppose $\mathcal{T}$ is a higher order (bigger than 2) tensor. The uniqueness of the CPD is up to the following trivial modifications:
		\begin{enumerate}
			\item Permutation of the ordering of the rank one terms. $\mathcal{T} = \displaystyle \sum_{r=1}^R \mathcal{T}_{\sigma(r)}$ is the same CPD, where $\sigma \in S_R$ is any permutation.
			\item Scaling indeterminacy. $\mathcal{T} = \displaystyle \sum_{r=1}^R \frac{1}{\lambda_r} \left( \lambda_r \mathcal{T}_r \right)$ is the same CPD, where $\lambda_r \neq 0$ is arbitrary for all $r = 1 \ldots R$.
		\end{enumerate}
		
		Sometimes one say the CPD is \emph{essentially unique}. Is this uniqueness what makes the CPD so attractive to applications. Most of the time the CPD is unique, but sometimes this may not be the case. For this reason we want to stablish uniqueness conditions. Additionally, we should clarify what means when we say the CPD is unique ``most of the time''. 
		
		The most well known result on uniqueness of tensors is due to J. B. Kruskal \cite{kruskal1, kruskal3} although it is limited to third order tensors. Posteriorly this result was extended to arbitrary higher order tensors by N. D. Sidiropoulos and R. Bro \cite{sidiropoulos}. Below we show this extended result.
		
		\begin{definition}
			Let $\textbf{X} \in \mathbb{K}^{m \times n}$ be a matrix. The \emph{k-rank} of $\textbf{X}$ is the maximum value $k$ such that any $k$ columns of $\textbf{X}$ are linearly independent. We denote this value by $k_\textbf{X}$.
		\end{definition}
		
		\begin{theorem}[N. D. Sidiropoulos and R. Bro]
			Let $\mathcal{T} \in \mathbb{V}^{(1)} \otimes \ldots \otimes \mathbb{V}^{(L)}$ be a tensor of rank $R$ with CPD given by formula~\ref{cpd} and let $\textbf{V}^{(\ell)}$ be the $\ell$-th factor matrix of this CPD, for $\ell = 1 \ldots L$. If
			
			$$\sum_{\ell=1}^L k_{\textbf{V}^{(\ell)}} \geq 2R + L - 1,$$
then this CPD is unique.			 
		\end{theorem}
		
		In the matrix case ($L = 2$), suppose that both factors, $\textbf{V}^{(1)}$ and $\textbf{V}^{(2)}$, have all columns linearly independent. Then we have that $k_{\textbf{V}^{(1)}} + k_{\textbf{V}^{(2)}} = R + R < 2R + 1$. Therefore the CPD is never unique in the matrix case, a fact we had already observed. With this we have a condition for uniqueness.  	
		
	\section{Tensor geometry} \label{topology}
		Given a tensor space $\mathbb{V}^{(1)} \otimes \ldots \otimes \mathbb{V}^{(L)}$ with a basis $\{ \textbf{e}_{i_1}^{(1)} \otimes \ldots \otimes \textbf{e}_{i_L}^{(L)} \}$, we can consider each tensor as a element of $\mathbb{K}^{I_1 \times \ldots \times I_L}$, that is, each tensor
		
		$$\mathcal{T} = \sum_{i_1=1}^{I_1} \ldots \sum_{i_L=1}^{I_L} t_{i_1 \ldots i_L}\ \textbf{e}_{i_1}^{(1)} \otimes \ldots \otimes \textbf{e}_{i_L}^{(L)}$$		
is identified with the multidimensional array with entries $t_{i_1 \ldots i_L}$. In this case we can consider $\mathbb{V}^{(1)} \otimes \ldots \otimes \mathbb{V}^{(L)}$ a space with inner product defined by

		$$\langle \mathcal{T}, \mathcal{S} \rangle = \sum_{i_1=1}^{I_1} \ldots \sum_{i_L=1}^{I_L} t_{i_1 \ldots i_L} \overline{s_{i_1 \ldots i_L}}.$$
		
		This induces the norm

		$$\| \mathcal{T} \| = \sqrt{ \sum_{i_1=1}^{I_1} \ldots \sum_{i_L=1}^{I_L} | t_{i_1 \ldots i_L} |^2 }.$$		 
		 
		This allow us to talk about proximity of tensors. This is relevant because usually one is interested is solving~\ref{cpd_approx} in the best way possible, that is, to obtain the rank-$R$ tensor closest to $\mathcal{T}$ between all rank-$R$ approximations. Unfortunately, computing the best rank-$R$ approximation of a tensor is, in general, a \emph{ill-posed}\footnote{We call a problem \emph{well-posed} if a solution exists, is unique, and is stable in the sense it depends continuously int the input data. A problem is ill-posed if it is not well-posed.} problem \cite{lim}. In particular we have the following result, first observed in \cite{bini} and then deeply explored in \cite{landsberg}.
	
		\begin{theorem}
			The limit of a sequence of rank tensors $R$ is not necessarily a rank-$R$ tensor.
		\end{theorem}
		
		Denote $\sigma_R \left( \mathbb{V}^{(1)} \otimes \ldots \otimes \mathbb{V}^{(L)} \right) = \left\{ \mathcal{T} \in \mathbb{V}^{(1)} \otimes \ldots \otimes \mathbb{V}^{(L)}: \ rank(\mathcal{T}) \leq R \right\}$ for the set of tensors with rank $\leq R$. With respect relation to the theorem above, the rank of the limit of tensors can give a ``jump''. Because of this, the set $\sigma_R$ is not necessarily closed in the norm topology. This motivates the following definition.
		
		\begin{definition}
			We say a tensor $\mathcal{T} \in \mathbb{V}^{(1)} \otimes \ldots \otimes \mathbb{V}^{(L)}$ has \emph{border rank-$R$} if $R$ is the smallest number such that $\mathcal{T} \in \overline{\sigma_R}$, where $\overline{\sigma_R}$ is the closure of $\sigma_R$. In this case we denote $\underline{rank}(\mathcal{T}) = R$. 
		\end{definition}	
		
		The term ``border rank'' first appeared in the paper \cite{bini} in the context of matrix multiplication. There is an interesting about the story of the border rank at the beginning of chapter 2 of \cite{landsberg2}. We have the following result as a direct consequence of the definition.	
	
		\begin{theorem}
			If $\underline{rank}(\mathcal{T}) = R$, then there exists a sequence of rank-$R$ tensors converging to $\mathcal{T}$ and there is not a sequence of tensors with rank $< R$ converging to $\mathcal{T}$.
		\end{theorem}	
	
		\begin{corollary}
			$\underline{rank}(\mathcal{T}) \leq rank(\mathcal{T})$.
		\end{corollary}	
		
		 As we observed, computing the best rank-$R$ approximation of a tensor is, in general, a ill-posed problem. However, this is not the case when $L=2$, that is, the matrix case. This result is known since 1936 with Eckart and Young \cite{young}.
		  
		 \begin{theorem}[Eckart-Young, 1936]
		 	Let $\textbf{M} = \displaystyle \sum_{r=1}^R \sigma_r \textbf{u}_r \textbf{v}_r^\ast$ be a SVD of a rank-$R$ matrix in $\mathbb{K}^{n \times m}$. For any $1 \leq \tilde{R} \leq R$, the best rank-$R$ approximation of $\textbf{M}$ is given by $\tilde{\textbf{M}} = \displaystyle \sum_{r=1}^{\tilde{R}} \sigma_r \textbf{u}_r \textbf{v}_r^\ast$.		 	 		
		 \end{theorem}
		 
		 The phenomenon of border rank is the one responsible for the ill-posedness of the approximation problem. If a tensor $\mathcal{T}$ has rank $R$ and border rank $\tilde{R} < R$, then there is a sequence of rank-$\tilde{R}$ tensors converging to $\mathcal{T}$. This implies, in particular, that it does not exists a tensor of rank $\tilde{R}$ closest to $\mathcal{T}$, so finding the best rank-$\tilde{R}$ approximation of $\mathcal{T}$ is a ill-posed problem. The first report of a such phenomenon was in \cite{bini2}, where they gave an explicit example of a sequence of rank-5 tensors converging to a rank 6 tensor in 1979. In \cite{lim} there is a simple example of a tensor of rank 3 and border rank 2 which we reproduce below for illustration purposes.
	
		\begin{theorem}[V. de Silva and L. H. Lim, 2008] \label{border-rank-lim}
			Let $I_1, I_2, I_3 \geq 2$. Let $\mathcal{T} \in \mathbb{R}^{I_1} \otimes \mathbb{R}^{I_2} \otimes \mathbb{R}^{I_3}$ be a tensor such that
			
			$$\mathcal{T} = \textbf{x}^{(1)} \otimes \textbf{x}^{(2)} \otimes \textbf{y}^{(3)} + 	\textbf{x}^{(1)} \otimes \textbf{y}^{(2)} \otimes \textbf{x}^{(3)} + \textbf{y}^{(1)} \otimes \textbf{x}^{(2)} \otimes \textbf{x}^{(3)}$$
where each pair $\textbf{x}^{(\ell)}, \textbf{y}^{(\ell)} \in \mathbb{R}^{I_\ell}$ is linearly independent. Then $rank(\mathcal{T}) = 3$ and 
			$$\mathcal{T}^{(n)} = n \left( \textbf{x}^{(1)} + \frac{1}{n} \textbf{y}^{(1)} \right) \otimes \left( \textbf{x}^{(2)} + \frac{1}{n} \textbf{y}^{(2)} \right) \otimes \left( \textbf{x}^{(3)} + \frac{1}{n} \textbf{y}^{(3)} \right) - n \textbf{x}^{(1)} \otimes \textbf{x}^{(2)} \otimes \textbf{x}^{(3)}$$
is a sequence of rank-2 tensors converging to $\mathcal{T}$. In particular, $\underline{rank}(\mathcal{T}) \leq 2$.			
		\end{theorem}
		
		The tensor $\mathcal{T}$ of the theorem is an example of a tensor that has no best rank-2 approximation. It is interesting to note that the limit expression for $\mathcal{T}^{(n)}$ may be regarded as a derivative. In fact, define the function $f:\mathbb{R} \to \mathbb{R}^{I_1} \otimes \mathbb{R}^{I_2} \otimes \mathbb{R}^{I_3}$ by
		
		$$f(t) = (\textbf{x}^{(1)} + t \textbf{y}^{(1)}) \otimes (\textbf{x}^{(2)} + t \textbf{y}^{(2)}) \otimes (\textbf{x}^{(3)} + t \textbf{y}^{(3)}) = $$	
		$$ = \textbf{x}^{(1)} \otimes \textbf{x}^{(2)} \otimes \textbf{x}^{(3)} + t \textbf{x}^{(1)} \otimes \textbf{x}^{(2)} \otimes \textbf{y}^{(3)} + t \textbf{x}^{(1)} \otimes \textbf{y}^{(2)} \otimes \textbf{x}^{(3)} + t^2 \textbf{x}^{(1)} \otimes \textbf{y}^{(2)} \otimes \textbf{y}^{(3)} + $$
		$$ + t \textbf{y}^{(1)} \otimes \textbf{x}^{(2)} \otimes \textbf{x}^{(3)} + t^2 \textbf{y}^{(1)} \otimes \textbf{x}^{(2)} \otimes \textbf{y}^{(3)} + t^2 \textbf{y}^{(1)} \otimes \textbf{y}^{(2)} \otimes \textbf{x}^{(3)} + t^3 \textbf{y}^{(1)} \otimes \textbf{y}^{(2)} \otimes \textbf{y}^{(3)}.$$
	
		On one hand, using the derivative rules and making $t = 0$ we obtain 
		
		$$f'(0) = \textbf{x}^{(1)} \otimes \textbf{x}^{(2)} \otimes \textbf{x}^{(3)} + \textbf{x}^{(1)} \otimes \textbf{x}^{(2)} \otimes \textbf{y}^{(3)} + \textbf{x}^{(1)} \otimes \textbf{y}^{(2)} \otimes \textbf{x}^{(3)}.$$	
On the other hand, using the limit definition for the derivative we obtain

		$$f'(0) = \lim_{t \to 0} \frac{ (\textbf{x}^{(1)} + t \textbf{y}^{(1)}) \otimes (\textbf{x}^{(2)} + t \textbf{y}^{(2)}) \otimes (\textbf{x}^{(3)} + t \textbf{y}^{(3)}) - \textbf{x}^{(1)}  \otimes \textbf{x}^{(2)}  \otimes \textbf{x}^{(3)} }{t}.$$
	
\begin{figure}[h] 
			\centering
			\includegraphics[scale=.4]{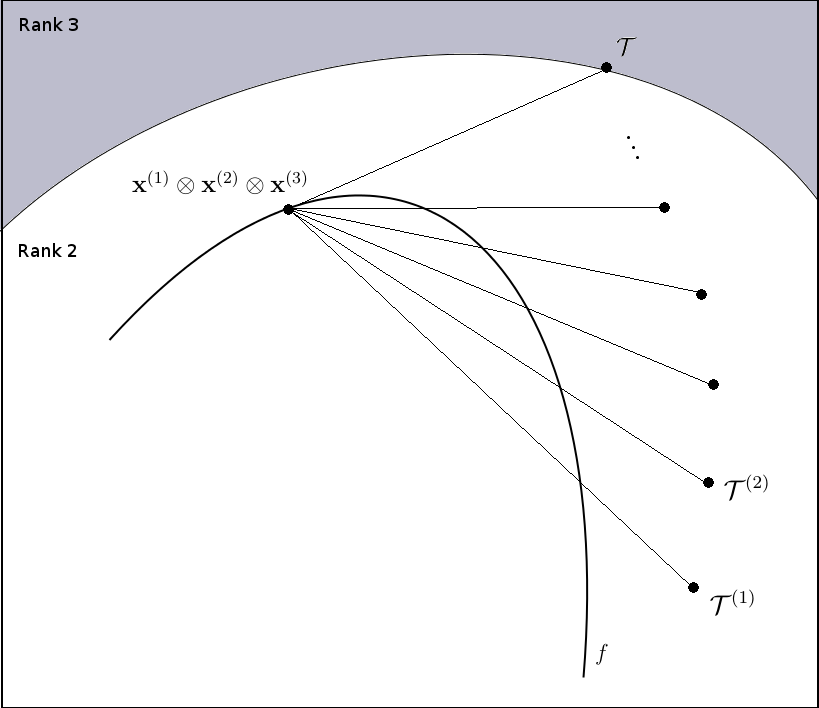}
			\caption{\footnotesize{Geometry of border rank.}}
			\label{border-rank}
		\end{figure}

		Making $t = 1/n$ and simplifying we obtain the expression of the theorem, that is, we have that $f'(0) = \lim_{n \to \infty} \mathcal{T}^{(n)}$. Notice that $f$ represents a curve in $\sigma_2$. From the point $f(0) = \textbf{x}^{(1)} \otimes \textbf{x}^{(2)} \otimes \textbf{x}^{(3)}$ we can draw secant lines in order to approximate the derivative $f'(0)$. Each secant line gives us a tensor in $\sigma_2$. At the limit we have the tangent tensor which, because of the theorem, will be outside $\sigma_2$. This is illustrated in figure~\ref{border-rank}.
		
		With respect to the topology of $\sigma_R$ in $\mathbb{V}^{(1)} \otimes \ldots \otimes \mathbb{V}^{(L)}$, we have already seen that $\sigma_R$ is not closed. It is also true that $\sigma_R$ is not open under certain conditions, as the next result shows.
		
		\begin{theorem}
			If $R < \dim\left( \mathbb{V}^{(1)} \otimes \ldots \otimes \mathbb{V}^{(L)} \right)$, then $\sigma_R$ is not open.
		\end{theorem}	
	
		\textbf{Proof:} Let $\mathcal{T} = \displaystyle \sum_{r=1}^R \textbf{v}_r^{(1)} \otimes \ldots \otimes \textbf{v}_r^{(L)} \in \sigma_R$ be a rank-$R$ tensor and let $\textbf{u}^{(1)} \otimes \ldots \otimes \textbf{u}^{(L)}$ be a tensor which it is not a linear combination of the tensor products $\textbf{v}_r^{(1)} \otimes \ldots \otimes \textbf{v}_r^{(L)}$. Then the sequence 
		
		$$\mathcal{T}^{(n)} = \frac{1}{n+1} \textbf{u}^{(1)} \otimes \ldots \otimes \textbf{u}^{(L)} + \mathcal{T}$$ 
converges to $\mathcal{T}$ and it is constituted of tensors with rank greater than $R$. Thus, every ball around $\mathcal{T}$ contains tensors outside $\sigma_R$. Therefore $\sigma_R$ is not open.$\hspace{4cm}\square$\bigskip

		Note that this theorem also holds for the matrix space. More precisely, there are sequences of rank-$R$ matrices converging to matrices with rank smaller than $R$. What does not occur with matrices is to have a sequence of rank-$R$ matrices converging to a matrix with rank greater than $R$. This last phenomenon is unique to tensors of order $ > 2 $, and this is where we see the issue of border rank.	
	
		The argument used in the previous theorem also shows that every rank-$R$ is an adherent point of $\left( \mathbb{V}^{(1)} \otimes \ldots \otimes \mathbb{V}^{(L)} \right) \backslash \sigma_R$, hence the set of rank-$R$ tensor is contained in $\overline{ \left( \mathbb{V}^{(1)} \otimes \ldots \otimes \mathbb{V}^{(L)} \right) \backslash \sigma_R }$. Let  $\partial(\sigma_R) = \overline{\sigma_R} \cap \overline{ \mathbb{V}^{(1)} \otimes \ldots \otimes \mathbb{V}^{(L)} \backslash \sigma_R }$ be the boundary of $\sigma_R$. Then it follows that the set of rank-$R$ tensors is contained in $\partial(\sigma_R)$. Next we give some results about norm invariance.
		
		\begin{theorem}[V. de Silva and L. H. Lim, 2008]
			Let $\mathcal{T} \in \mathbb{K}^{I_1} \otimes \ldots \otimes \mathbb{K}^{I_L}, \mathcal{S} \in \mathbb{K}^{I_1'} \otimes \ldots \otimes \mathbb{K}^{I_L'}$ and $\textbf{v}^{(1)} \in \mathbb{K}^{I_1}$, $\ldots$, $\textbf{v}^{(L)} \in \mathbb{K}^{I_L}$. Then the following statements holds.
			\begin{enumerate}
			\item $\|\textbf{v}^{(1)} \otimes \ldots \otimes \textbf{v}^{(L)} \| = \| \textbf{v}^{(1)} \| \cdot \ldots \cdot \| \textbf{v}^{(L)} \|$
			\item $\|\mathcal{T} \otimes \mathcal{S} \| = \| \mathcal{T} \| \cdot \| \mathcal{S} \|$
			\item If $\textbf{U}^{(1)} \in \mathbb{K}^{I_1 \times I_1}, \ldots, \textbf{U}^{(L)} \in \mathbb{K}^{I_L \times I_L}$ are unitary (orthogonal) matrices, then\\ $\| (\textbf{U}^{(1)}, \ldots, \textbf{U}^{(L)}) \cdot \mathcal{T} \| = \| \mathcal{T} \|$.
			\end{enumerate}
		\end{theorem}
		
		As already noted, $\sigma_R$ is not closed (except in the case of matrices). The next result shows other equivalent statements.
	
		\begin{theorem}[V. de Silva and L. H. Lim, 2008]
			Consider the space $\mathbb{K}^{I_1} \otimes \ldots \otimes \mathbb{K}^{I_L}$ with $L > 2$ and let $R \geq 2$. Then the following statements are equivalent.
			\begin{enumerate}
				\item $\sigma_R$ is not closed.
				\item There exists a sequence of tensors in $\sigma_R$ converging to a tensor with rank greater than $R$.
				\item There exists a tensor $\mathcal{S}$ of rank greater than $R$ such that $\displaystyle\inf_{\mathcal{T} \in \sigma_R } \| \mathcal{T} - \mathcal{S} \| = 0$.
				\item There exists a tensor $\mathcal{S}$ of rank greater than $R$ which does note have a best rank-$R$ approximation, that is, $\displaystyle\inf_{\mathcal{T} \in \sigma_R } \| \mathcal{T} - \mathcal{S} \|$ is not attained in $\sigma_R$.
			\end{enumerate}
		\end{theorem}
	
\begin{figure}[h] 
			\centering
			\includegraphics[scale=.4]{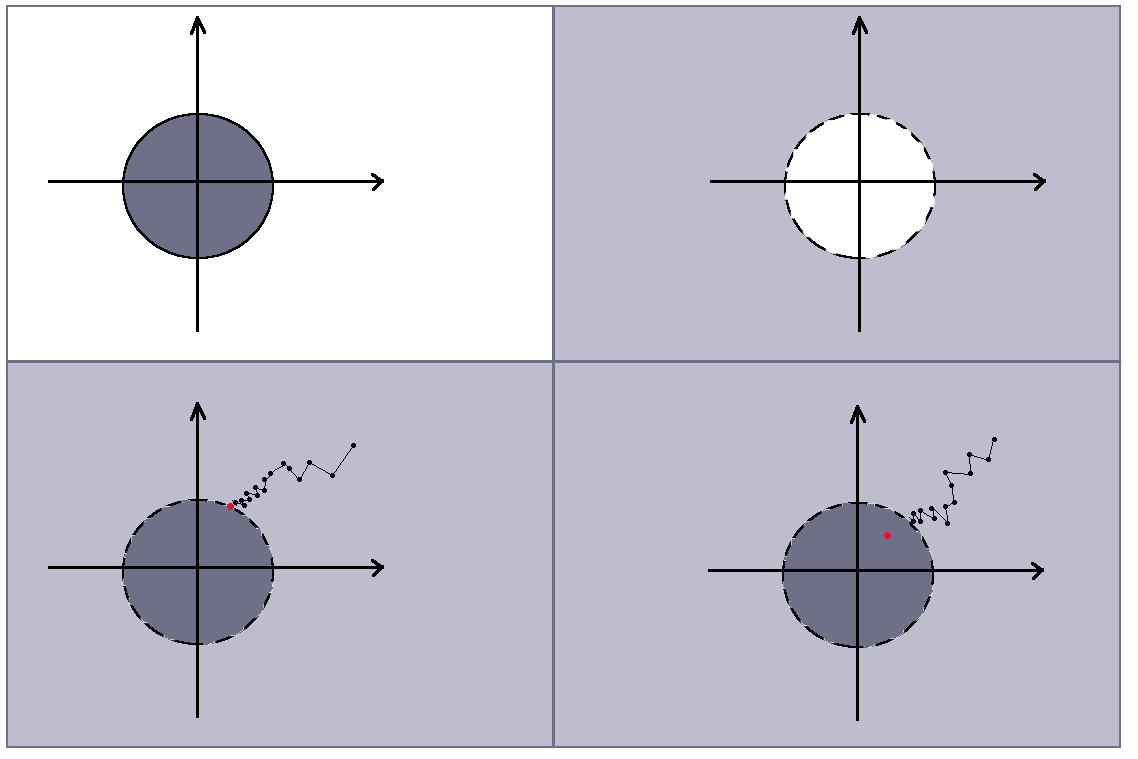}
			\caption{\footnotesize{Convergence issues.}}
			\label{approx-issues}
		\end{figure}

		It is important to emphasize that there are tensors of rank greater than $R$ which also can't be arbitrarily approximated by rank-$R$ tensors. Figure~\ref{approx-issues} illustrates the possible situations one can encounter. In the figure on the top left, the dark region represents a certain subset of tensors with rank greater than $R$. In the top right figure, the light region represents the set $\sigma_R$. The dotted line indicates that those border points are not in $\sigma_R$, they are part of the dark region. In the figure on the bottom left we have a sequence of points in $\sigma_R$ converging to the red dot, which is at the border between the dark and the light regions. This point is at the closure of $\sigma_R$, so that it is a tensor with rank greater than $R$ which can be approximated arbitrarily well by points in $\sigma_R$. This tensor has border rank equal to $R$. In the figure on the bottom right the sequence of points converges to the point of the border closest to the red point, but the limit of that convergence is not the point desired. In this case we have a tensor with rank greater than $R$ that does not have a best rank-$R$ approximation.
	
		Although everything we saw up to this point may indicate that the best rank-$R$ approximation problem is something to be avoided, there are two positive results showed below.
	
		\begin{theorem}[V. de Silva and L. H. Lim, 2008]
			Every tensor has a best rank 1 approximation.
		\end{theorem}
	
		\begin{theorem}[V. de Silva and L. H. Lim, 2008]
			Let $\mathcal{T} \in \mathbb{K}^{I_1} \otimes \ldots \otimes \mathbb{K}^{I_L}$ and integers $R_1, \ldots, R_L > 0$. Then $\mathcal{T}$ does have a best approximation of multilinear rank $ \leq (R_1, \ldots, R_L)$.
		\end{theorem}
	
		The next theorem shows one reason that can make tensor low rank approximations to fail.
	
		\begin{theorem}[V. de Silva and L. H. Lim, 2008]
			Let $\mathcal{T} \in \mathbb{K}^{I_1} \otimes \ldots \otimes \mathbb{K}^{I_L}$ a tensor with rank greater than $R$ and let $(\mathcal{T}^{(n)})$ be a sequence in $\sigma_R$ converging to $\mathcal{T}$. Furthermore, write
			
			$$\mathcal{T}^{(n)} = \sum_{r=1}^R \lambda_r^{(n)}\ \textbf{v}_r^{(n,1)} \otimes \ldots \otimes \textbf{v}_r^{(n,L)},$$
where each $\textbf{v}_r^{(n,\ell)} \in \mathbb{K}^{I_\ell}$ is a unitary vector and $\lambda_r^{(n)} \in \mathbb{K}$ a scalar. Then there exists two distinct numbers $1 \leq r_1, r_2, \leq R$ such that $\displaystyle \lim_{n \to \infty} |\lambda_{r_1}^{(n)}| = \lim_{n \to \infty} |\lambda_{r_2}^{(n)}| = \infty$.
		\end{theorem}
	
		Although we have factors diverging, the sequence still converges. What happens is that these divergent factors cause cancellations as $n$ increases. If $\mathcal{T}$ does not have a best rank-$R$ approximation, the process of computing more approximations $\mathcal{T}^{(n)}$ can continue indefinitely, with some coefficients $\lambda_r^{(n)}$ diverging, which is a problem when making computations with finite precision. Taking away the condition of the vectors being unitary, we will have vectors diverging, which will change nothing. This phenomenon of diverging terms has been observed in practical applications of multilinear models and is referred as ``degeneracy'' \cite{stegeman1, kruskal4, paatero, stegeman2, stegeman3}.	

%% file: CHAPTER_2.tex
\chapter{Tensor compression}\label{cap-2} 
	Tensor compression is an important tool to compute a CPD. It reduces the problem size, hence the computational and memory size. It relies on the computation of some SVDs of matrices, which is a familiar decomposition. Nowadays there are very fast implementations for the SVD, such as the randomized truncated SVD~\cite{rand_svd}. First we will see how to make unfoldings from tensors, then we go to multilinear rank and some related results, and finally we finish with the compression of tensors. Algorithms and their costs are shown along the way.
	
	\section{Multilinear multiplication}
		For each $\ell = 1 \ldots L$, let $B^{(\ell)} = \{ \textbf{e}_1^{(\ell)}, \ldots, \textbf{e}_{I_\ell}^{(\ell)} \}$ be a basis for each space $\mathbb{K}^{I_\ell}$, and consider a tensor $\mathcal{T} \in \mathbb{K}^{I_1} \otimes \ldots \otimes \mathbb{K}^{I_L}$ in coordinates $t_{i_1 \ldots i_L}$. As already observed, we can use these coordinates to interpret $\mathcal{T}$ as a multidimensional array in $\mathbb{K}^{I_1 \times \ldots \times I_L}$. Every time we refer to $\mathcal{T}$ as a element of $\mathbb{K}^{I_1 \times \ldots \times I_L}$ it will be implicit that there are fixed bases. Matrices can act on $\mathcal{T}$ by $L$ ``distinct directions'' through the usual matrix multiplication. Let $\textbf{M}^{(1)} \in \mathbb{K}^{I_1' \times I_1}, \ldots, \textbf{M}^{(L)} \in \mathbb{K}^{I_L' \times I_L}$ be any matrices, then we denote by $(\textbf{M}^{(1)}, \ldots, \textbf{M}^{(L)}) \cdot \mathcal{T}$ the ``multiplication'' between the $L$-tuple $(\textbf{M}^{(1)}, \ldots, \textbf{M}^{(L)}) $ and $\mathcal{T}$. The result of this multiplication is the tensor $\mathcal{S} = (\textbf{M}^{(1)}, \ldots, \textbf{M}^{(L)}) \cdot \mathcal{T} \in \mathbb{K}^{I_1'} \otimes \ldots \otimes \mathbb{K}^{I_L'}$ defined as
		
		$$s_{j_1 \ldots j_L} = \sum_{i_1=1}^{I_1} \ldots \sum_{i_L=1}^{I_L} m^{(1)}_{j_1 i_1} \ldots m^{(L)}_{j_L i_L} t_{i_1 \ldots i_L}.$$ 		
This operation is called \emph{multilinear multiplication}. In the particular case $\mathcal{T}$ is a matrix, we have that $(\textbf{M}, \textbf{N}) \cdot \mathcal{T} = \textbf{M} \cdot \mathcal{T} \cdot \textbf{N}^T$, and in the case $\mathcal{T}$ is a vector we have $(\textbf{M}) \cdot \mathcal{T} = \textbf{M} \cdot \mathcal{T}$.
		
		As a direct consequence of the definition, the tensor $\mathcal{S}$ is given in coordinates, although we do not have defined any basis for the spaces $\mathbb{K}^{I_1'}, \ldots, \mathbb{K}^{I_L'}$. The choice of these bases depends on each situation and in principle is arbitrary. The definition of the multilinear multiplication is motivated by the idea of making change of basis. 
		
		Consider the bases $B^{(\ell)} = \{ \textbf{e}_1^{(\ell)}, \ldots, \textbf{e}_{I_\ell}^{(\ell)} \}, \tilde{B}^{(\ell)} = \{ \tilde{\textbf{e}}_1^{(\ell)}, \ldots, \tilde{\textbf{e}}_{I_\ell}^{(\ell)} \}$ for each space $\mathbb{K}^{I_\ell}$ and let $\textbf{M}^{(\ell)} \in \mathbb{K}^{I_\ell \times I_\ell}$ be the change of basis matrix from $\tilde{B}^{(\ell)}$ to $B^{(\ell)}$, that is, we have that 
		
		$$\textbf{e}_i^{(\ell)} = \sum_{j=1}^{I_\ell} m_{ji}^{(\ell)} \tilde{\textbf{e}}_j^{(\ell)}.$$
		
		It follows that
		
		$$\mathcal{T} = \sum_{i_1=1}^{I_1} \ldots \sum_{i_L=1}^{I_L} t_{i_1 \ldots i_L} \ \textbf{e}_{i_1}^{(1)} \otimes \ldots \otimes \textbf{e}_{i_L}^{(L)} = $$
	
		$$ = \sum_{i_1=1}^{I_1} \ldots \sum_{i_L=1}^{I_L} t_{i_1 \ldots i_L} \left( \sum_{j_1=1}^{I_1} m_{j_1 i_1}^{(1)} \  \tilde{\textbf{e}}_{j_1}^{(1)} \right) \otimes \ldots \otimes \left( \sum_{j_L=1}^{I_L} m_{j_L i_L}^{(L)} \  \tilde{\textbf{e}}_{j_L}^{(L)} \right) = $$
	
		$$ = \sum_{i_1=1}^{I_1} \ldots \sum_{i_L=1}^{I_L} \sum_{j_1=1}^{I_1} \ldots \sum_{j_L=1}^{I_L} m_{j_1 i_1}^{(1)} \ldots m_{j_L i_L}^{(L)} \ t_{i_1 \ldots i_L} \ \tilde{\textbf{e}}_{j_1}^{(1)} \otimes \ldots \otimes \tilde{\textbf{e}}_{j_L}^{(L)} = $$
	
		$$ = \sum_{j_1=1}^{I_1} \ldots \sum_{j_L=1}^{I_L} s_{j_1 \ldots j_L} \ \tilde{\textbf{e}}_{j_1}^{(1)} \otimes \ldots \otimes \tilde{\textbf{e}}_{j_L}^{(L)},$$
where
		$$s_{j_1 \ldots j_L} = \sum_{i_1=1}^{I_1} \ldots \sum_{i_L=1}^{I_L} m_{j_1 i_1}^{(1)} \ldots m_{j_L i_L}^{(L)} \ t_{i_1 \ldots i_L}.$$
		
		\begin{remark}
			The change of basis is given by the formula $\mathcal{S} = (\textbf{M}^{(1)}, \ldots, \textbf{M}^{(L)}) \cdot \mathcal{T}$. The tensor $\mathcal{S}$ in coordinates represents $\mathcal{T}$ after this change of basis. The difference between $\mathcal{S}$ and $\mathcal{T}$ is only in its representation as a multidimensional array since it depends on coordinates. But as tensors they are the same object, which we refer as the \emph{abstract tensor}. One could be more precise and use some notation like $\mathcal{T}_{\tilde{B}^{(1)}, \ldots, \tilde{B}^{(L)}} = (\textbf{M}^{(1)}, \ldots, \textbf{M}^{(L)}) \cdot \mathcal{T}_{B^{(1)}, \ldots, B^{(L)}} $. This notation is more cumbersome and for this reason we will avoid it. Furthermore, it is not always the case that the matrices involved are change of basis matrices.
		\end{remark}
		
		Let $\mathcal{S} = (\textbf{M}^{(1)}, \ldots, \textbf{M}^{(L)}) \cdot \mathcal{T}$. Sometimes the following equality is useful:
			$$\mathcal{S} = \sum_{i_1=1}^{I_1} \ldots \sum_{i_L=1}^{I_L} t_{i_1 \ldots i_L} \ \textbf{M}_{: i_1}^{(1)} \otimes \ldots \otimes \textbf{M}_{: i_L}^{(L)}$$
		
		Denote by $GL(d, \mathbb{K})$ the linear group of matrices in $\mathbb{K}^{d \times d}$. If the field $\mathbb{K}$ is clear from the context, we just denote $GL(d)$. As we know, all change of basis matrices are invertible and every invertible matrix can be interpreted as a change of basis. This gives us a criterion of equivalence between tensors.
		
		\begin{definition} \label{equivalence}
			Let two tensors $\mathcal{T}, \mathcal{S} \in \mathbb{K}^{I_1 \times \ldots \times I_L}$. We say they are \emph{equivalent} if there are matrices $\textbf{M}^{(1)} \in GL(I_1), \ldots, \textbf{M}^{(L)} \in GL(I_L)$ such that $\mathcal{S} = (\textbf{M}^{(1)}, \ldots, \textbf{M}^{(L)}) \cdot \mathcal{T}$.
		\end{definition}
		
		\begin{theorem}
			Let two tensors $\mathcal{T}, \mathcal{S} \in \mathbb{K}^{I_1 \times \ldots \times I_L}$. They represent the same abstract tensor if, and only if, they are equivalent.
		\end{theorem}
		
		In the case of the theorem being valid, it is possible to have $\mathcal{S} \neq \mathcal{T}$ as multidimensional arrays and $\mathcal{S} = \mathcal{T}$ as abstract tensors. This is the same situation when we have distinct matrices representing the same linear map, the difference is only due to the choice of basis. Below there are some basic properties of the multilinear multiplication.
		
		\begin{theorem}[V. de Silva and L. H. Lim, 2008]
			Let two tensors $\mathcal{T}, \mathcal{S} \in \mathbb{K}^{I_1 \times \ldots \times I_L}$ and the matrices $\textbf{M}^{(1)} \in \mathbb{K}^{I_1' \times I_1}$, $\ldots$, $\textbf{M}^{(L)} \in \mathbb{K}^{I_L' \times I_L}$. Then
			\begin{enumerate}
				\item For all $\alpha, \beta \in \mathbb{K}$, we have 
				$$(\textbf{M}^{(1)}, \ldots, \textbf{M}^{(L)}) \cdot (\alpha \mathcal{T} + \beta \mathcal{S}) = \alpha ( (\textbf{M}^{(1)}, \ldots, \textbf{M}^{(L)}) \cdot \mathcal{T} ) + \beta ( (\textbf{M}^{(1)}, \ldots, \textbf{M}^{(L)}) \cdot \mathcal{S} ).$$
				\item For all $\textbf{N}^{(1)} \in \mathbb{K}^{I_1'' \times I_1'}, \ldots, \textbf{N}^{(L)} \in \mathbb{K}^{I_L'' \times I_L'}$, we have 
				$$(\textbf{N}^{(1)}, \ldots, \textbf{N}^{(L)}) \cdot ((\textbf{M}^{(1)}, \ldots, \textbf{M}^{(L)}) \cdot \mathcal{T}) = (\textbf{N}^{(1)} \textbf{M}^{(1)}, \ldots, \textbf{N}^{(L)} \textbf{M}^{(L)}) \cdot \mathcal{T}.$$
				\item For all $\alpha, \beta \in \mathbb{C}$ and all $\textbf{A}, \textbf{B} \in \mathbb{K}^{I_j' \times I_j}$, we have  
				$$(\textbf{M}^{(1)}, \ldots, \textbf{M}^{(j-1)}, \alpha \textbf{A} + \beta \textbf{B}, \textbf{M}^{(j+1)}, \ldots, \textbf{M}^{(L)}) \cdot \mathcal{T} = $$
				$$\hspace{-1cm} = \alpha ( (\textbf{M}^{(1)}, \ldots, \textbf{M}^{(j-1)}, \textbf{A}, \textbf{M}^{(j+1)}, \ldots, \textbf{M}^{(L)}) \cdot \mathcal{T} ) + \beta ( (\textbf{M}^{(1)}, \ldots, \textbf{M}^{(j-1)}, \textbf{B}, \textbf{M}^{(j+1)}, \ldots, \textbf{M}^{(L)}) \cdot \mathcal{T} ).$$
			\end{enumerate}
		\end{theorem}
		
		\begin{theorem}[V. de Silva and L. H. Lim, 2008] \label{equivalence-properties}
			Let $T \in \mathbb{K}^{I_1} \otimes \ldots \otimes \mathbb{K}^{I_K}$ and $\textbf{M}^{(1)} \in \mathbb{K}^{I_1' \times I_1}$, $\ldots$, $\textbf{M}^{(L)} \in \mathbb{K}^{I_L' \times I_L}$. Then
			\begin{enumerate}
				\item $rank((\textbf{M}^{(1)}, \ldots, \textbf{M}^{(L)})\cdot \mathcal{T}) \leq rank(\mathcal{T})$.
				\item If $\textbf{M}^{(1)} \in GL(I_1'), \ldots, \textbf{M}^{(L)} \in GL(I_L')$, then $rank((\textbf{M}^{(1)}, \ldots, \textbf{M}^{(L)})\cdot \mathcal{T}) = rank(\mathcal{T})$.
			\end{enumerate}
		\end{theorem}

		Now let's see how the multilinear multiplication and tensor product are related.
		
		\begin{theorem}[V. de Silva and L. H. Lim, 2008]
			Let $\textbf{v}^{(1)} \otimes \ldots \otimes \textbf{v}^{(L)} \in \mathbb{K}^{I_1} \otimes \ldots \otimes \mathbb{K}^{I_L}$ be a rank one tensor and let the matrices $\textbf{M}^{(1)} \in \mathbb{K}^{I_1' \times I_1}, \ldots, \textbf{M}^{(L)} \in \mathbb{K}^{I_L' \times I_L}$. Then
		$$(\textbf{M}^{(1)}, \ldots, \textbf{M}^{(L)}) \cdot \textbf{v}^{(1)} \otimes \ldots \otimes \textbf{v}^{(L)} = (\textbf{M}^{(1)} \textbf{v}^{(1)}) \otimes \ldots \otimes (\textbf{M}^{(L)} \textbf{v}^{(L)}).$$
		\end{theorem}

		\begin{corollary}
			Let $\mathcal{T} = \displaystyle\sum_{r=1}^R \textbf{v}_r^{(1)} \otimes \ldots \otimes \textbf{v}_r^{(L)} \in \mathbb{K}^{I_1} \otimes \ldots \otimes \mathbb{K}^{I_L}$ be a tensor with rank $\leq R$ and let the matrices $\textbf{M}^{(1)} \in \mathbb{K}^{I_1' \times I_1}, \ldots, \textbf{M}^{(L)} \in \mathbb{K}^{I_L' \times I_L}$. Then
			$$(\textbf{M}^{(1)}, \ldots, \textbf{M}^{(L)}) \cdot \mathcal{T} = \sum_{r=1}^R (\textbf{M}^{(1)} \textbf{v}_r^{(1)}) \otimes \ldots \otimes (\textbf{M}^{(L)} \textbf{v}_r^{(L)}).$$
		\end{corollary}
		
		\begin{theorem}
			Let $\mathcal{T} = \displaystyle\sum_{r=1}^R \lambda_r \ \textbf{v}_r^{(1)} \otimes \ldots \otimes \textbf{v}_r^{(L)} \in \mathbb{K}^{I_1} \otimes \ldots \otimes \mathbb{K}^{I_L}$ be a tensor with rank $\leq R$, where each $\lambda_r$ is a scalar. Then
			$$\mathcal{T} = (\textbf{V}^{(1)}, \ldots, \textbf{V}^{(L)}) \cdot \Lambda,$$
where $\textbf{V}^{(\ell)} = [\textbf{v}_1^{(\ell)}, \ldots, \textbf{v}_R^{(\ell)}] \in \mathbb{K}^{I_\ell \times R}$ for each $\ell = 1 \ldots R$, and $\Lambda = \text{diag}(\lambda_r) \in \mathbb{K}^{R \times \ldots \times R}$ is a diagonal tensor of order $L$.\footnote{These kind of tensor as sometimes called \emph{superdiagonal}. Denote by $\lambda_{i_1 \ldots i_L}$ the entries of $\Lambda$. Then we have that $\lambda_{i_1 \ldots i_L} = \lambda_r$ if $i_1 = \ldots = i_L = r$, and $\lambda_{i_1 \ldots i_L} = 0$ otherwise.} 
		\end{theorem}
	
		\begin{corollary}
			Let $\mathcal{T} = \displaystyle\sum_{r=1}^R \lambda_r \ \textbf{x}_r \otimes \textbf{y}_r \otimes \textbf{z}_r \in \mathbb{K}^{I_1} \otimes \mathbb{K}^{I_2} \otimes \mathbb{K}^{I_3}$ be a third order tensor with rank $\leq R$, where each $\lambda_r$ is a scalar. Then
			$$\mathcal{T} = (\textbf{X}, \textbf{Y}, \textbf{Z}) \cdot \Lambda,$$
where $\textbf{X} = [\textbf{x}_1, \ldots, \textbf{x}_R] \in \mathbb{K}^{I_1 \times R}$, $\textbf{Y} = [\textbf{y}_1, \ldots, \textbf{y}_R] \in \mathbb{K}^{I_2 \times R}$, $\textbf{Z} = [\textbf{z}_1, \ldots, \textbf{z}_R] \in \mathbb{K}^{I_3 \times R}$ and $\Lambda = \text{diag}(\lambda_r) \in \mathbb{K}^{R \times R \times R}$ is a diagonal tensor.
		\end{corollary}
		
		\subsection{Unfoldings}
			Suppose we have the mode-$j$ fibers of a tensor $\mathcal{T} \in \mathbb{V}^{(1)} \otimes \ldots \otimes \mathbb{V}^{(L)}$. These fibers are several vectors, as we've already seen. We can put them side by side and concatenate them to form a matrix. The ordering is not that important, we choose the ordering according to the order of the indexes. With this we have a matrix of shape $\displaystyle I_\ell \times \prod_{j \neq \ell} I_j$ called a \emph{unfolding of $\mathcal{T}$}. We denote this unfolding by $\mathcal{T}_{(\ell)}$. Other common names are \emph{matricization} and \emph{flattening}. The construction of $\mathcal{T}_{(\ell)}$ with the ordering we are using can be described by the following pseudo-code. Denote by $[ \ ]$ an ``empty matrix'' which will be filled column by column, where $\textbf{M} \leftarrow [ \textbf{M} | \textbf{v} ]$ means to add a column vector $\textbf{v}$ at the right of $\textbf{M}$ and then substitute $\textbf{M}$ by this new matrix.
			
			\begin{algorithm}[Unfolding]
				$ $\\
				\textbf{Input:} $\mathcal{T}, \ell$\vspace{.4cm}\\
				$\mathcal{T}_{(\ell)} = [ \ ]$\\
				$\verb|for | i_L = 1 \ldots I_L$\\	
				$\hspace{.5cm}\ddots$\\	
				$\hspace{1.5cm}\verb|for | i_{\ell+1} = 1 \ldots I_{\ell+1}$\\	
				$\hspace{2cm}\verb|for | i_{\ell-1} = 1 \ldots I_{\ell-1}$\\	
				$\hspace{2.7cm}\ddots$\\	
				$\hspace{3.5cm}\verb|for | i_1 = 1 \ldots I_1$\\	
				$\hspace{4cm} \mathcal{T}_{(\ell)} \leftarrow \left[ \mathcal{T}_{(\ell)} \ | \ \mathcal{T}_{i_1 \ldots i_{\ell-1}\ :\ i_{\ell+1} \ldots i_L} \right]$\vspace{.4cm}\\
				\textbf{Output:} $\mathcal{T}_{(\ell)}$
			\end{algorithm}
			
			\begin{example}
				Consider the tensor $\mathcal{T} \in \mathbb{R}^{3 \times 4 \times 2}$ given by
				
				$$\mathcal{T} = \left\{ 
				\left[
				\begin{array}{cccc}
					1 & 4 & 7 & 10\\
					2 & 5 & 8 & 11\\
					3 & 6 & 9 & 12
				\end{array}
				\right],
				\left[
				\begin{array}{cccc}
					13 & 16 & 19 & 22\\
					14 & 17 & 20 & 23\\
					15 & 18 & 21 & 24
				\end{array}
				\right]
				\right\},$$		
				where this is the representation of $\mathcal{T}$ through its frontal slices. The mode-1 fibers of $\mathcal{T}$ are
				
				$$\left[
				\begin{array}{c}
					1\\
					2\\
					3
				\end{array}
				\right],		
				\left[
				\begin{array}{c}
					4\\
					5\\
					6
				\end{array}
				\right],
				\left[
				\begin{array}{c}
					7\\
					8\\
					9
				\end{array}
				\right],
				\left[
				\begin{array}{c}
					10\\
					11\\
					12
				\end{array}
				\right],
				\left[
				\begin{array}{c}
					13\\
					14\\
					15
				\end{array}
				\right],
				\left[
				\begin{array}{c}
					16\\
					17\\
					18
				\end{array}
				\right],
				\left[
				\begin{array}{c}
					19\\
					20\\
					21
				\end{array}
				\right],
				\left[
				\begin{array}{c}
					22\\
					23\\
					24
				\end{array}
				\right].$$
				
				Note that we already ordered the vectors in accord with our convention. It follows that
				
				$$\mathcal{T}_{(1)} = \left[
				\begin{array}{cccccccc}
					1 & 4 & 7 & 10 & 13 & 16 & 19 & 22\\
					2 & 5 & 8 & 11 & 14 & 17 & 20 & 23\\
					3 & 6 & 9 & 12 & 15 & 18 & 21 & 24
				\end{array}
				\right].$$
	
				The mode-2 fibers are 
				
				$$\left[
				\begin{array}{c}
					1\\
					4\\
					7\\
					10
				\end{array}
				\right],	
				\left[
				\begin{array}{c}
					2\\
					5\\
					8\\
					11
				\end{array}
				\right],
				\left[
				\begin{array}{c}
					3\\
					6\\
					9\\
					12
				\end{array}
				\right],
				\left[
				\begin{array}{c}
					13\\
					16\\
					19\\
					22
				\end{array}
				\right],
				\left[
				\begin{array}{c}
					14\\
					17\\
					20\\
					23
				\end{array}
				\right],
				\left[
				\begin{array}{c}
					15\\
					18\\
					21\\
					24
				\end{array}
				\right].$$
				
				It follows that
				
				$$\mathcal{T}_{(2)} = \left[
				\begin{array}{cccccc}
					1 & 2 & 3 & 13 & 14 & 15\\
					4 & 5 & 6 & 16 & 17 & 18\\
					7 & 8 & 9 & 19 & 20 & 21\\
					10 & 11 & 12 & 22 & 23 & 24
				\end{array}
				\right].$$
	
				Finally, the mode-3 fibers are 
				
				$$\hspace{-.7cm}\left[
				\begin{array}{c}
					1\\
					13
				\end{array}
				\right],
				\left[
				\begin{array}{c}
					2\\
					14
				\end{array}
				\right],
				\left[
				\begin{array}{c}
					3\\
					15
				\end{array}
				\right],
				\left[
				\begin{array}{c}
					4\\
					16
				\end{array}
				\right],
				\left[
				\begin{array}{c}
					5\\
					17
				\end{array}
				\right],
				\left[
				\begin{array}{c}
					6\\
					18
				\end{array}
				\right],
				\left[
				\begin{array}{c}
					7\\
					19
				\end{array}
				\right],
				\left[
				\begin{array}{c}
					8\\
					20
				\end{array}
				\right],
				\left[
				\begin{array}{c}
					9\\
					21
				\end{array}
				\right],
				\left[
				\begin{array}{c}
					10\\
					22
				\end{array}
				\right],
				\left[
				\begin{array}{c}
					11\\
					23
				\end{array}
				\right],
				\left[
				\begin{array}{c}
					12\\
					24
				\end{array}
				\right].$$
				
				It follows that

				$$\mathcal{T}_{(3)} = \left[
				\begin{array}{cccccccccccc}
					1 & 2 & 3 & 4 & 5 & 6 & 7 & 8 & 9 & 10 & 11 & 12\\
					13 & 14 & 15 & 16 & 17 & 18 & 19 & 20 & 21 & 22 & 23 & 24
				\end{array}
				\right].$$
	
				We can note that $\mathcal{T}_{(1)} \in \mathbb{R}^{3 \times 4 \cdot 2},\ \mathcal{T}_{(2)} \in \mathbb{R}^{4 \times 3 \cdot 2},\ \mathcal{T}_{(3)} \in \mathbb{R}^{2 \times 3 \cdot 4}$, as expected.	
			\end{example}
			
			Unfoldings and multilinear multiplication are highly connected. The idea is that we can multiply an unfolding by some matrix and consider the result as the unfolding of a new tensor.
			
			\begin{definition}
				Let a tensor $\mathcal{T} \in \mathbb{K}^{I_1 \times \ldots \times I_L}$ and a matrix $\textbf{M} \in \mathbb{K}^{I_\ell' \times I_\ell}$. The \emph{product mode-$\ell$} between $\mathcal{T}$ and $\textbf{M}$ is the tensor $\mathcal{S} \in \mathbb{K}^{I_1 \times \ldots \times I_{\ell-1} \times I_\ell' \times I_{\ell+1} \times \ldots \times I_L}$ such that $\mathcal{S}_{(\ell)} = \textbf{M} \cdot \mathcal{T}_{(\ell)}$. We denote $\mathcal{S} = \mathcal{T} \times_\ell \textbf{M}$.
			\end{definition}
			
			Although the symbol $\times_\ell$ is at the right of the tensor, in the actual multiplication it comes at the left, that is, this is a left action on the tensor space. This is just a notational convention and should not cause confusion. Furthermore, we will omit parenthesis when making more than one of these products. More precisely, we will write $\mathcal{T} \times_\ell \textbf{M} \times_{\ell'} \textbf{N}$ instead of $(\mathcal{T} \times_\ell \textbf{M}) \times_{\ell'} \textbf{N}$. It is possible to obtain $\mathcal{S}$ explicitly in coordinates, this gives			
			$$\mathcal{S}_{i_1 \ldots i_{\ell-1} \hspace{.07cm} i_\ell' \hspace{.07cm} i_{\ell+1} \ldots i_L} = \sum_{k = 1}^{I_\ell'} m_{i_\ell' \hspace{.05cm} k} \cdot t_{i_1 \ldots i_{\ell-1} \hspace{.07cm} k \hspace{.07cm} i_{\ell+1} \ldots i_L}$$
for all $i_\ell' = 1 \ldots I_\ell'$. This formula remind us the formula of the multilinear multiplication. Indeed there is a connection. Let a tensor $T \in \mathbb{K}^{I_1 \times \ldots \times I_L}$ and matrices $\textbf{M}^{(1)} \in \mathbb{K}^{I_1' \times I_1}$, $\ldots$, $\textbf{M}^{(L)} \in \mathbb{K}^{I_L' \times I_L}$. Then 
			$$(\textbf{M}^{(1)}, \ldots, \textbf{M}^{(L)}) \cdot \mathcal{T} = \mathcal{T} \times_1 \textbf{M}^{(1)} \times_2 \textbf{M}^{(2)} \ldots \times_L \textbf{M}^{(L)}.$$
			
			This relationship gives us a clearer picture of how the multilinear multiplication acts on tensors. Each $\textbf{M}^{(\ell)}$ multiplies all the mode-$j$ fibers of $\mathcal{T}$ (which in the end is equivalent to multiply $\textbf{M}^{(\ell)}$ by $\mathcal{T}_{(\ell)}$), and thus we obtain a new tensor. The next example clarifies how works this relationship more concretely.
			
			\begin{example}
				Let $\mathcal{T}$ be the tensor of the previous example and let $\textbf{M} = \left[
				\begin{array}{ccc}
					1 & 3 & 5\\
					2 & 4 & 6
				\end{array}
				\right]$. By definition, we have that $\mathcal{T} \times_1 \textbf{M} \in \mathbb{R}^{2 \times 4 \times 2}$, with
				$$(\mathcal{T} \times_1 \textbf{M})_{(1)} = \textbf{M} \cdot \mathcal{T}_{(1)} = $$
		
				$$ = \left[
				\begin{array}{ccc}
					1 & 3 & 5\\
					2 & 4 & 6
				\end{array}
				\right] 
				\left[
			\begin{array}{cccccccc}
				1 & 4 & 7 & 10 & 13 & 16 & 19 & 22\\
				2 & 5 & 8 & 11 & 14 & 17 & 20 & 23\\
				3 & 6 & 9 & 12 & 15 & 18 & 21 & 24
			\end{array}
			\right] = 
				$$
		
				$$ = \left[
				\begin{array}{cccccccc}
					22 & 49 & 76 & 103 & 130 & 157 & 184 & 211\\
					28 & 64 & 100 & 136 & 172 & 208 & 244 & 280
				\end{array}
				\right].$$
		
				This is matrix of order $2 \times (4\cdot 2) = 2 \times 8$, so that the first half of the matrix is the first frontal slice of $\mathcal{T} \times_1 \textbf{M}$. From this we conclude that
				
				$$\mathcal{T} \times_1 \textbf{M} = \left\{ 
				\left[
				\begin{array}{cccc}
					22 & 49 & 76 & 103\\
					28 & 64 & 100 & 136
				\end{array}
				\right],
				\left[
				\begin{array}{cccc}
					130 & 157 & 184 & 211\\
					172 & 208 & 244 & 280
				\end{array}
				\right]
				\right\}.$$
		
				To make this example simple, we will only consider the action of the matrix $\textbf{M}$. Note that $\mathcal{T} \times_1 \textbf{M} = \mathcal{T} \times_1 \textbf{M} \times_2 \textbf{I}_4 \times_3 \textbf{I}_2 = (\textbf{M}, \textbf{I}_4, \textbf{I}_2) \cdot \mathcal{T}$. Now let $\mathcal{S} = (\textbf{M}, \textbf{I}_4, \textbf{I}_2) \cdot \mathcal{T}$. From the multilinear multiplication definition we have that  				
				$$\mathcal{S}_{i'j'k'} = \sum_{i=1}^3 \sum_{j=1}^4 \sum_{k=1}^2 m_{i' i} \cdot (\textbf{I}_4)_{j' j} \cdot  (\textbf{I}_2)_{k' k} \cdot t_{ijk} = \sum_{i=1}^3 m_{i' i} \cdot t_{i j' k'} = \textbf{M}_{i':} \cdot T_{:j' k'}.$$
				
				This last expression is the product of the $i'$-th row of $\textbf{M}$ by the column vector of $\mathcal{T}$ obtained by fixing $j',k'$ and varying the rows, that is, a mode-1 fiber. By varying $i'$ to form a mode-1 fiber of $\mathcal{S}$ we obtain				
				$$\mathcal{S}_{:j' k'} = \left[
				\begin{array}{c}
					t_{1 j' k'}\\
					t_{2 j' k'}
				\end{array}
				\right]	=	
				\left[
				\begin{array}{c}
					\textbf{M}_{1:} \cdot \mathcal{T}_{: j' k'}\\
					\textbf{M}_{2:} \cdot \mathcal{T}_{: j' k'}
				\end{array}
				\right] = 
				\left[
				\begin{array}{c}
					\textbf{M}_{1:}\\
					\textbf{M}_{2:}
				\end{array}
				\right] \cdot \mathcal{T}_{: j' k'} = \textbf{M} \cdot \mathcal{T}_{: j' k'}$$
	
				From this we conclude that

				$$\mathcal{S}_{(1)} = \left[\mathcal{S}_{:11}, \mathcal{S}_{:21}, \mathcal{S}_{:31}, \mathcal{S}_{:41}, \mathcal{S}_{:12}, \mathcal{S}_{:22}, \mathcal{S}_{:32}, \mathcal{S}_{:42} \right] = $$
		
				$$ = \left[ \textbf{M} \mathcal{T}_{:11}, \textbf{M} \mathcal{T}_{:21}, \textbf{M} \mathcal{T}_{:31}, \textbf{M} \mathcal{T}_{:41}, \textbf{M} \mathcal{T}_{:12}, \textbf{M} \mathcal{T}_{:22}, \textbf{M} \mathcal{T}_{:32}, \textbf{M} \mathcal{T}_{:42} \right] = $$
		
				$$ = \textbf{M} \cdot \left[ \mathcal{T}_{:11}, \mathcal{T}_{:21}, \mathcal{T}_{:31}, \mathcal{T}_{:41}, \mathcal{T}_{:12}, \mathcal{T}_{:22}, \mathcal{T}_{:32}, \mathcal{T}_{:42} \right] = $$
		
				$$ = \textbf{M} \cdot \mathcal{T}_{(1)}.$$
			\end{example}			
		
		\subsection{Multilinear rank}	 
			A natural thing to do is to consider the rank of unfoldings.
	
			\begin{definition}
				For each $\ell = 1 \ldots L$, the \emph {mode-$\ell$ rank of $\mathcal{T}$} is the rank of $\mathcal{T}_{(\ell)}$. We denote this rank by $rank_{(\ell)}(\mathcal{T})$. 
			\end{definition} 
	
			\begin{definition}
				The \emph{multilinear rank of $\mathcal{T}$} is the $L$-tuple $\left( rank_{(1)}(\mathcal{T}), \ldots, rank_{(L)}(\mathcal{T}) \right)$. We denote this rank by $rank_\boxplus(\mathcal{T})$.
			\end{definition}
	
			The multilinear rank is also often called the \emph{Tucker rank} of $\mathcal{T}$. Now we will see some important results regarding this rank.
			
			\begin{theorem}[V. de Silva and L. H. Lim, 2008] \label{multirank-properties}
				Let $\mathcal{T} \in \mathbb{K}^{I_1} \otimes \ldots \otimes \mathbb{K}^{I_L}$ and $\textbf{M}^{(1)} \in \mathbb{K}^{I_1' \times I_1}, \ldots, \textbf{M}^{(L)} \in \mathbb{K}^{I_L' \times I_L}$. Then the following statements holds.
				\begin{enumerate}
					\item The multilinear rank doesn't depend on the field being real or complex. 
					\item $rank_\boxplus(\mathcal{T}) = (1,1, \ldots, 1)$ if, and only if, $rank(\mathcal{T}) = 1$.
					\item $rank_{(\ell)}(\mathcal{T}) \leq \min\{ rank(\mathcal{T}), I_1, \ldots, I_L\}$ for all $\ell = 1 \ldots L$.
					\item $\| rank_\boxplus(\mathcal{T}) \|_\infty \leq rank(T)$.
					\item $rank_\boxplus((\textbf{M}^{(1)}, \ldots, \textbf{M}^{(L)})\cdot \mathcal{T}) \leq rank_\boxplus(\mathcal{T})$.
					\item If $\textbf{M}^{(1)} \in GL(I_1'), \ldots, \textbf{M}^{(L)} \in GL(I_L')$, then $rank_\boxplus((\textbf{M}^{(1)}, \ldots, \textbf{M}^{(L)}) \cdot \mathcal{T}) = rank_\boxplus(\mathcal{T})$.
				\end{enumerate}
			\end{theorem}
			
			\begin{remark} \label{smaller-tensor}
				Consider a tensor in coordinates $\mathcal{T} \in \mathbb{K}^{I_1 \times \ldots \times I_L}$ such that $rank_\boxplus(\mathcal{T}) = (R_1, \ldots, R_L) $. Theorem~\ref{multirank-properties}-6 together with theorem~\ref{equivalence-properties}-2 makes it possible to study rank properties of $\mathcal{T}$ as a tensor in $\mathbb{K}^{R_1 \times \ldots \times R_L}$. 
			\end{remark}
			
			\begin{theorem}[V. de Silva and L. H. Lim, 2008]
				Let $\mathcal{T} \in \mathbb{K}^{I_1 \times \ldots \times I_L}$ be a tensor such that $rank_\boxplus(\mathcal{T}) \leq (R_1, \ldots, R_L)$. Then there there exists full rank matrices $\textbf{M}^{(1)} \in \mathbb{K}^{I_1 \times R_1}$, $\ldots$, $\textbf{M}^{(L)} \in \mathbb{K}^{I_L \times R_L}$ and a tensor $\mathcal{S} \in \mathbb{K}^{R_1 \times \ldots \times R_L}$ such that $\mathcal{T} = (\textbf{M}^{(1)}, \ldots, \textbf{M}^{(L)}) \cdot \mathcal{S}$ and $rank(\mathcal{T}) = rank(\mathcal{S})$.
			\end{theorem}			
	
	\section{Compressing with the multilinear singular value decomposition}
		\subsection{Tucker decomposition}
			\begin{definition}
				Let $\mathcal{T} \in \mathbb{K}^{I_1 \ldots I_L}$ be a tensor. Then a \emph{Tucker decomposition of $\mathcal{T}$} is a decomposition of the form $\mathcal{T} = (\textbf{M}^{(1)}, \ldots, \textbf{M}^{(L)}) \cdot \mathcal{S}$, where $\textbf{M}^{(\ell)} \in \mathbb{K}^{I_\ell' \times I_\ell}$ for $\ell = 1 \ldots L$, and $\mathcal{S} \in \mathbb{K}^{I_1' \times \ldots \times I_L'}$. 
			\end{definition}	
	
			Each matrix $\textbf{M}^{(\ell)}$ is called a \emph{factor matrix} and the tensor $\mathcal{S}$ is called the \emph{core tensor}. Usually one defines this decomposition assuming all factor matrices are unitary (orthogonal), but we will prefer the more general definition. This way the CPD may be seen as a particular case of the Tucker decomposition. 
		
			The Tucker decomposition was first introduced in 1963 and refined subsequently \cite{tucker1, tucker2}. This decomposition can be considered as a form of high-order PCA (Principal Component Analysis)\cite{kolda} when the core tensor has lower dimensions than the original one. In this case we consider $\mathcal{S}$ as a compressed form of $\mathcal{T}$. Usually the CPD is indicated for latent parameter estimation and the Tucker decomposition is indicated for compression and dimensionality reduction.  
		
			From what we've seen on multilinear multiplication, we can write a Tucker decomposition $\mathcal{T} = (\textbf{M}^{(1)}, \ldots, \textbf{M}^{(L)}) \cdot \mathcal{S}$ as 
		
			$$\mathcal{T} = \sum_{i_1=1}^{I_1} \ldots \sum_{i_L=1}^{I_L} s_{i_1 \ldots i_L} \ \textbf{M}_{: i_1}^{(1)} \otimes \ldots \otimes \textbf{M}_{: i_L}^{(L)}.$$
If $\mathcal{S} = \text{diag}(s_r) \in \mathbb{K}^{R \times \ldots \times R}$ is diagonal and $\textbf{M}^{(\ell)} \in \mathbb{K}^{I_\ell \times R}$, then we have a rank-$R$ CPD of $\mathcal{T}$ and we can write

		$$\mathcal{T} = \sum_{r=1}^R s_r \ \textbf{M}_{: r}^{(1)} \otimes \ldots \otimes \textbf{M}_{: r}^{(L)}.$$
We can see that the factor matrices of the Tucker decomposition agrees with the factor matrices of the CPD. On the other extreme, let $\{ \textbf{e}_1^{(\ell)}, \ldots, \textbf{e}_{I_\ell}^{(\ell)} \}$ be a basis of $\mathbb{K}^{I_\ell}$ for each $\ell = 1 \ldots L$. In this case we are able to write $\mathcal{T}$ as 

			$$\mathcal{T} = \sum_{i_1=1}^{I_1} \ldots \sum_{i_L=1}^{I_L} t_{i_1 \ldots i_L} \ \textbf{e}_{i_1}^{(1)} \otimes \ldots \otimes \textbf{e}_{i_L}^{(L)}.$$
Denoting $\textbf{E}^{(\ell)} = [\textbf{e}_1^{(\ell)}, \ldots, \textbf{e}_L^{(\ell)}] \in \mathbb{K}^{I_\ell \times I_\ell}$, we can write $\mathcal{T} = (\textbf{E}^{(1)}, \ldots, \textbf{E}^{(L)}) \cdot \mathcal{T}$. This is a trivial Tucker decomposition, whereas the CPD can be seen as the ``ultimate'' Tucker decomposition. Between these two there are other useful decompositions to consider.

			\begin{theorem}[T. G. Kolda, 2006] \label{unfolding-formula}
				Let $\mathcal{T} \in \mathbb{K}^{I_1 \ldots I_L}$ be a tensor with a Tucker decomposition given by $\mathcal{T} = (\textbf{M}^{(1)}, \ldots, \textbf{M}^{(L)}) \cdot \mathcal{S}$. Then, for each $\ell = 1 \ldots L$, the unfolding $\mathcal{T}_{(\ell)}$ is given by 
				$$\mathcal{T}_{(\ell)} = \textbf{M}^{(\ell)} \cdot \mathcal{S}_{(\ell)} \cdot (\textbf{M}^{(L)} \tilde{\otimes} \ldots \tilde{\otimes} \textbf{M}^{(\ell+1)} \tilde{\otimes} \textbf{M}^{(\ell-1)} \tilde{\otimes} \ldots \tilde{\otimes} \textbf{M}^{(1)})^T,$$
where $\tilde{\otimes}$ is the \emph{Kronecker product}, see appendix~\ref{appenB}. 
			\end{theorem}
		
			\begin{corollary}
				Let $\mathcal{T} \in \mathbb{K}^{I_1 \times I_2 \times I_3}$ be a third order tensor with a CPD given by $\mathcal{T} = \displaystyle \sum_{r=1}^R \lambda_r \ \textbf{x}_r \otimes \textbf{y}_r \otimes \textbf{z}_r = (\textbf{X}, \textbf{Y}, \textbf{Z}) \cdot \Lambda$, where $\Lambda = \text{diag}(\lambda_r) \in \mathbb{K}^{R \times R \times R}$. Then
				$$\Lambda_{(1)} = \Lambda_{(2)} = \Lambda_{(3)} = \left[
				\begin{array}{cccc}
					\lambda_1 \ 0 \ \ldots \ 0 & 0 \ \ 0 \ \ldots \ \ 0 & \ldots & 0 \ \ 0 \ \ldots \ \ 0\\
					0 \ \ 0 \ \ldots \ \ 0 & 0 \ \lambda_2 \ \ldots \ 0 & \ldots & 0 \ \ 0 \ \ldots \ \ 0\\
					\vdots & \vdots & \ldots & \vdots\\
					\underbrace{0 \ \ 0 \ \ldots \ 0}_{R \text{ columns}} & \underbrace{0 \ \ 0 \ \ldots \ 0}_{R \text{ columns}} & \ldots & \underbrace{0 \ \ 0 \ \ldots \ \lambda_R}_{R \text{ columns}}
				\end{array}
				\right]$$
and, denoting this matrix by $[\Lambda]$, we also have that
				$$\mathcal{T}_{(1)} = \textbf{X} \cdot [\Lambda] \cdot (\textbf{Z} \tilde{\otimes} \textbf{Y})^T,$$
				$$\mathcal{T}_{(2)} = \textbf{Y} \cdot [\Lambda] \cdot (\textbf{Z} \tilde{\otimes} \textbf{X})^T,$$
				$$\mathcal{T}_{(3)} = \textbf{Z} \cdot [\Lambda] \cdot (\textbf{Y} \tilde{\otimes} \textbf{X})^T.$$
			\end{corollary}
			
			\begin{theorem}[T. G. Kolda, 2006] \label{unfolding-formula2}
				Let $\mathcal{T} \in \mathbb{K}^{I_1 \ldots I_L}$ be a tensor with a Tucker decomposition given by $\mathcal{T} = (\textbf{M}^{(1)}, \ldots, \textbf{M}^{(L)}) \cdot \mathcal{I}_{R \times \ldots \times R}$. Then, for each $\ell = 1 \ldots L$, the unfolding $\mathcal{T}_{(\ell)}$ is given by 
				$$\mathcal{T}_{(\ell)} = \textbf{M}^{(\ell)} \cdot (\textbf{M}^{(L)} \odot \ldots \odot \textbf{M}^{(\ell+1)} \odot \textbf{M}^{(\ell-1)} \odot \ldots \odot \textbf{M}^{(1)})^T.$$ 
			\end{theorem}
			
			\subsection{Multilinear singular value decomposition}
				Now we will see how to generalize the SVD to tensors. This generalization is called \emph{multilinear singular value decomposition (MLSVD)}, but sometimes it is also called \emph{High order singular value decomposition (HOSVD)}. Some texts even call this as being the Tucker decomposition. This generalization has been investigated in psychometrics \cite{tucker2} as the \emph{Tucker model}, which basically was a special Tucker decomposition of third order tensors. The first work to formalize this as a high order singular value decomposition is \cite{mlsvd}.  
				
				Let $\textbf{M}$ be a matrix with SVD given by $\textbf{M} = \textbf{U} \Sigma \textbf{V}^\ast$. We can rewrite this equation using the multilinear multiplication notation, then we obtain $\textbf{M} = (\textbf{U}, \overline{\textbf{V}}) \cdot \Sigma$, where $\overline{\textbf{V}}$ is the conjugate of $\textbf{V}$ coordinatewise. Denoting $\textbf{U}^{(1)} = \textbf{U}$ and $\textbf{U}^{(2)} = \overline{\textbf{V}}$, we can write $\textbf{M} = (\textbf{U}^{(1)}, \textbf{U}^{(2)}) \cdot \Sigma$. The MLSVD is a generalization of this observation.  
				
				Given a tensor $\mathcal{S} \in \mathbb{K}^{I_1 \times \ldots \times I_L}$, let $\mathcal{S}_{i_\ell = k} \in \mathbb{K}^{I_1 \times \ldots \times I_{\ell-1} \times I_{\ell+1} \times \ldots \times I_L}$ be the subtensor of $\mathcal{S}$ obtained by fixing the $\ell$-th index of $\mathcal{S}$ with value equal to $k$ and varying all the other indexes. More precisely, $\mathcal{S}_{i_\ell = k} = \mathcal{S}_{: \ldots : k : \ldots :}$, where the value $k$ is at the $\ell$-th index. We call these subtensors by \emph{hyperslices}. In the case of a third order tensors, these subtensors are the slices described in~\ref{slices}. $\mathcal{S}_{i_1 = k}$ are the horizontal slices, $\mathcal{S}_ {i_2 = k}$ the lateral slices and $\mathcal{S}_ {i_3 = k}$ the frontal slices.
				
				\begin{theorem}[L. De Lathauwer, B. De Moor, J. Vandewalle, 2000] \label{MLSVD}
					Let $\mathcal{T} \in \mathbb{K}^{I_1 \times \ldots \times I_L}$ be arbitrary. Then there exists unitary (orthogonal) matrices $\textbf{U}^{(1)} \in \mathbb{K}^{I_1 \times I_1}, \ldots$, $\textbf{U}^{(L)} \in \mathbb{K}^{I_L \times I_L}$ and a tensor $\mathcal{S} \in \mathbb{K}^{I_1 \times \ldots \times I_L}$ such that 
					\begin{enumerate}
						\item $\mathcal{T} = (\textbf{U}^{(1)}, \ldots, \textbf{U}^{(L)}) \cdot \mathcal{S}$.
						\item For all $\ell = 1 \ldots L$, the subtensors $\mathcal{S}_{i_\ell = 1}, \ldots, \mathcal{S}_{i_\ell = I_\ell}$ are orthogonal with respect to each other.
						\item For all $\ell = 1 \ldots L$, $\|\mathcal{S}_{i_\ell = 1}\| \geq \ldots \geq \|\mathcal{S}_{i_\ell = I_\ell}\|$.
					\end{enumerate}
				\end{theorem}
				
				\textbf{Proof:} For each unfolding $\mathcal{T}_{(\ell)}$, consider the corresponding reduced SVD 
				
				$$\mathcal{T}_{(\ell)} = \textbf{U}^{(\ell)} \cdot \Sigma^{(\ell)} \cdot (\textbf{V}^{(\ell)})^\ast$$
where $\Sigma^{(\ell)} = \text{diag}\big( \sigma_1^{(\ell)}, \ldots, \sigma_{I_\ell}^{(\ell)} \big) \in \mathbb{R}^{I_\ell \times I_\ell}$ and $\textbf{U} \in \mathbb{K}^{I_\ell \times I_\ell}, \textbf{V} \in \mathbb{K}^{ \left( \prod_{j \neq \ell} I_j \right) \times I_\ell}$ are unitary (orthogonal) matrices. Let $\mathcal{S} \in \mathbb{K}^{I_1 \times \ldots I_L}$ be the tensor defined as

				$$\mathcal{S} = \big( \big( \textbf{U}^{(1)} \big)^\ast, \ldots, \big( \textbf{U}^{(L)} \big)^\ast \big) \cdot \mathcal{T}.$$

				By theorem~\ref{unfolding-formula} we know that				
				$$\mathcal{S}_{(\ell)} = (\textbf{U}^{(\ell)})^\ast \cdot \mathcal{T}_{(\ell)} \cdot \big( \big( \textbf{U}^{(L)} \big)^\ast \tilde{\otimes} \ldots \tilde{\otimes} \big( \textbf{U}^{(\ell+1)} \big)^\ast \tilde{\otimes} \big( \textbf{U}^{(\ell-1)} \big)^\ast \tilde{\otimes} \ldots \tilde{\otimes} \big( \textbf{U}^{(1)} \big)^\ast \big)^T = $$
	
				$$ = (\textbf{U}^{(\ell)})^\ast \cdot \mathcal{T}_{(\ell)} \cdot \left( \overline{\textbf{U}^{(L)}} \tilde{\otimes} \ldots \tilde{\otimes} \overline{\textbf{U}^{(\ell+1)}} \tilde{\otimes} \overline{\textbf{U}^{(\ell-1)}} \tilde{\otimes} \ldots \tilde{\otimes} \overline{\textbf{U}^{(1)}} \right) = $$
	
				$$ = \Sigma^{(\ell)} \cdot (\textbf{V}^{(\ell)})^\ast \cdot \left( \overline{\textbf{U}^{(L)}} \tilde{\otimes} \ldots \tilde{\otimes} \overline{\textbf{U}^{(\ell+1)}} \tilde{\otimes} \overline{\textbf{U}^{(\ell-1)}} \tilde{\otimes} \ldots \tilde{\otimes} \overline{\textbf{U}^{(1)}} \right).$$

				Denote by $(\mathcal{S}_{(\ell)})_{k :}$ the $k$-th row of $\mathcal{S}_{(\ell)}$. From the definition we have that $(\mathcal{S}_{(\ell)})_{k :} = \mathcal{S}_{i_\ell = k}$. Now let $1 \leq k, k' \leq I_\ell$ distinct, then we have $ \langle \mathcal{S}_{i_\ell = k}, \mathcal{S}_{i_\ell = k'} \rangle = 0$ because $\mathcal{S}_{(\ell)}$ has orthonormal rows. This last assertion follows from the fact that $\mathcal{S}_{(\ell)}$ is the product of a diagonal matrix by two of unitary (orthogonal) matrices. Finally, note that $\| \mathcal{S}_{i_\ell = k} \| = | \sigma_k^{(\ell)} | = \sigma_k^{(\ell)}$.
Hence it follows that $\| \mathcal{S}_{i_\ell = 1} \| \geq \ldots \geq \| \mathcal{S}_{i_\ell = I_\ell} \|$.$\hspace{6.5cm}\square$\bigskip
				
				Property 2 is called \emph{all-orthogonality}. It means that $\langle \mathcal{S}_{i_\ell=k}, \mathcal{S}_{i_\ell=k'} \rangle = 0$ for all $k, k' = 1 \ldots I_\ell$ and $k \neq k'$. The inner product considered is the one mentioned at the beginning of the tensor geometry section. The ordering given in property 3 can be considered as a convention, made to fix a particular ordering of the columns of the matrices $\textbf{U}^{(\ell)}$. Additionally, this ordering gives us a more precise parallel to the SVD. We will denote $\sigma_k^{(\ell)} = \| \mathcal{S}_{i_\ell = k} \| $ and call this a \emph{singular value mode-$\ell$ of $\mathcal{T} $}, while each column vector $\textbf{U}_{:j}^{(\ell)}$ is called a \emph{singular vector mode-$\ell$ of $\mathcal{T}$}. An algorithm for computing the MLSVD is presented below.
				
				\begin{algorithm}[MLSVD] \label{MLSVD-alg}
					$ $\\
					\textbf{Input:} $\mathcal{T} \in \mathbb{R}^{I_1 \times \ldots \times I_L}$\vspace{.4cm}\\
					\verb|for | $\ell=1 \ldots L$\\
						$\hspace{1cm} \mathcal{T}_{(\ell)} \leftarrow $ \verb|Unfolding| $(\mathcal{T}, \ell)$\\
						$\hspace{1cm} \textbf{U}^{(\ell)}, \Sigma^{(\ell)}, \big(\textbf{V}^{(\ell)}\big)^\ast \leftarrow $ \verb|SVD| $(\mathcal{T}_{(\ell)})$\\
					$\mathcal{S} \leftarrow \big( \big( \textbf{U}^{(1)} \big)^\ast, \ldots, \big( \textbf{U}^{(L)} \big)^\ast \big) \cdot \mathcal{T}$\vspace{.4cm}\\
					\textbf{Output:} $\textbf{U}^{(\ell)}, \Sigma^{(\ell)} \in \mathbb{R}^{I_\ell \times I_\ell}$ for $\ell = 1 \ldots L$ and $\mathcal{S} \in \mathbb{R}^{I_1 \times \ldots \times I_L}$
				\end{algorithm}
				
\begin{figure}[h] 
					\centering
					\includegraphics[scale=.15]{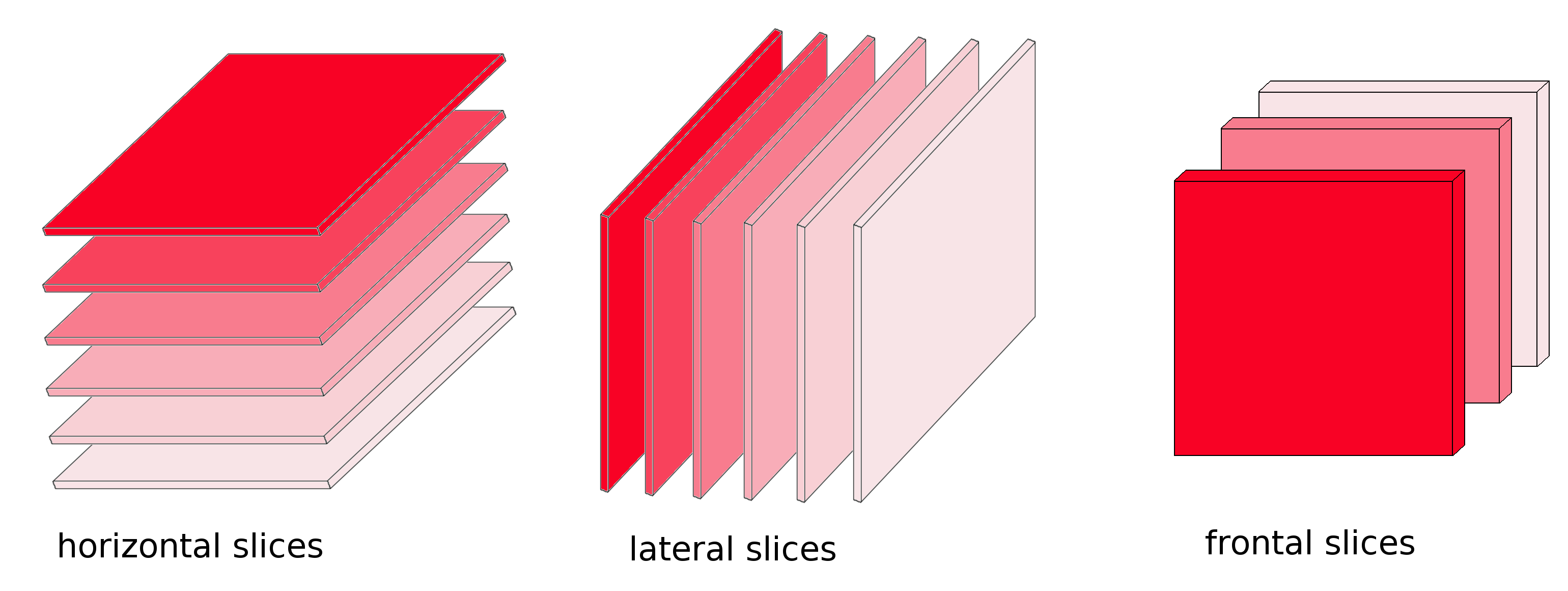}
					\caption{\footnotesize{Energy of the slices of a core third order tensor $\mathcal{S}$ obtained after a MLSVD.}}
					\label{slices-energy}
				\end{figure}

				The core tensor $\mathcal{S}$ of the MLSVD distributes the ``energy'' (i.e., the magnitude of its entries) in such a way so that it concentrates more energy around the first entry $s_{11 \ldots 1}$ and disperses as we move along each dimension. Figure~\ref{slices-energy} illustrates the energy distribution when $\mathcal{S}$ is a third order tensor. The red slices contains more energy and it changes to white when the slice contains less energy. Note that the energy of the slices are given precisely by the singular values.
				
				\begin{theorem}[L. De Lathauwer, B. De Moor, J. Vandewalle, 2000]
					Let $\mathcal{S}_{(\ell)}$ be an unfolding of $\mathcal{S}$. If $(\mathcal{S}_{(\ell)})_{k :}$ is the $k$-th row of $\mathcal{S}_{(\ell)}$, then $\| (\mathcal{S}_{(\ell)})_{k :} \| = \sigma_k^{(\ell)}$.
				\end{theorem}
				
				Let $\mathcal{T} = (\textbf{U}^{(1)}, \ldots, \textbf{U}^{(L)}) \cdot \mathcal{S}$ be the MLSVD of $\mathcal{T}$. Define the matrix $\Sigma^{(\ell)} \in \mathbb{K}^{I_\ell \times I_\ell}$ given by $\Sigma^{(\ell)} = \text{diag}\big( \sigma_1^{(\ell)}, \ldots, \sigma_{I_\ell}^{(\ell)} \big)$, and let $\tilde{\mathcal{S}}_{(\ell)} \in \mathbb{K}^{I_\ell \times I_\ell}$ be defined by the relation $\mathcal{S}_{(\ell)} = \Sigma^{(\ell)} \cdot \tilde{\mathcal{S}}_{(\ell)}$. Notice that $\tilde{\mathcal{S}}_{(\ell)}$ is the normalized version of $\mathcal{S}_{(\ell)}$ (it has unit-length rows). In the case none of the $\sigma_k^{(\ell)}$ is null, we can write $\tilde{\mathcal{S}}_{(\ell)} = \left( \Sigma^{(\ell)} \right)^{-1} \cdot \mathcal{S}_{(\ell)}$. Finally, define the matrix $\textbf{V}^{(\ell)} \in \mathbb{K}^{ \left( \prod_{j \neq \ell} I_j \right) \times I_\ell}$ given by the relation 
				$$(\textbf{V}^{(\ell)})^\ast = \tilde{\mathcal{S}}_{(\ell)} \cdot \left( \overline{\textbf{U}^{(L)}} \tilde{\otimes} \ldots \tilde{\otimes} \overline{\textbf{U}^{(\ell+1)}} \tilde{\otimes} \overline{\textbf{U}^{(\ell-1)}} \tilde{\otimes} \ldots \tilde{\otimes} \overline{\textbf{U}^{(1)}} \right).$$
With these notations and the formula for $\mathcal{T}_{(\ell)}$ (theorem~\ref{unfolding-formula}) we can write				
				$$\mathcal{T}_{(\ell)} = \textbf{U}^{(\ell)} \cdot \Sigma^{(\ell)} \cdot (\textbf{V}^{(\ell)})^\ast.$$
				
				\begin{theorem}[L. De Lathauwer, B. De Moor, J. Vandewalle, 2000]
					With the notations above, $\mathcal{T}_{(\ell)} = \textbf{U}^{(\ell)} \cdot \Sigma^{(\ell)} \cdot (\textbf{V}^{(\ell)})^\ast$ is a SVD for $\mathcal{T}_{(\ell)}$, for each $\ell = 1 \ldots L$.
				\end{theorem}
				
				With this theorem we can see that the MLSVD is a reasonable extension of the SVD, because one property desired for such an extension is to be able to use the MLSVD to obtain the SVD of each unfolding of $\mathcal{T}$. Soon we will give more reasons to consider this as a good extension of the SVD.    
				
				In the MLSVD, just as in the SVD, it is possible to have the last singular values mode-$\ell$ equal to zero, that is, it can exist a number $1 \leq R_\ell \leq I_\ell$ such that $\sigma_k^{(\ell)} = 0$ for $k = R_\ell + 1 \ldots I_\ell $. Therefore all hyperslices $\mathcal{S}_{i_\ell = k}$ are null (all of its entries are zero) for $k = R_\ell + 1 \ldots I_\ell$. Consequently, if we have $i_\ell > R_\ell$ for a multi-index $(i_1, \ldots, i_L)$, then $s_{i_1 \ldots i_L} = 0$. With this, we can write the tensor $\mathcal{T}$ as
								
	$$\mathcal{T} = \sum_{i_1=1}^{R_1} \ldots \sum_{i_L=1}^{R_L} s_{i_1 \ldots i_L} \ \textbf{U}_{: i_1}^{(1)} \otimes \ldots \otimes \textbf{U}_{: i_L}^{(L)},$$
which shows $\mathcal{T}$ as a sum of $\displaystyle\prod_{\ell=1}^L R_\ell$ rank one terms.	We illustrate the form of $\mathcal{S}$ in the figure below as a third order core tensor. The gray part corresponds to the part of nonzero values, while the white part consists only of zeros.
				
				\begin{figure}[h] 
					\centering
					\includegraphics[scale=.4]{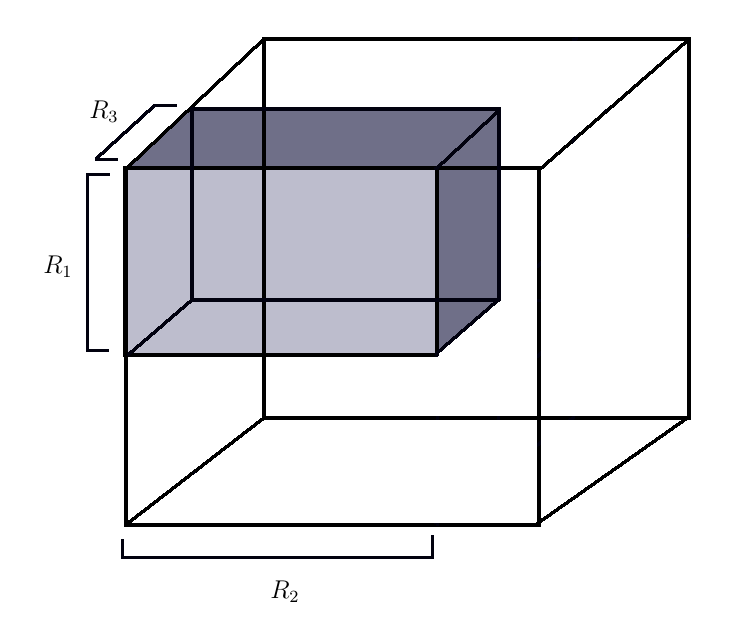}
					\caption{\footnotesize{Representation of the null slices slices of core third order tensor $\mathcal{S}$ obtained after a MLSVD.}}
				\end{figure}
				
				As a consequence of the above observation and theorem~\ref{multirank-properties}-4, we obtain
				
				$$\max_\ell R_\ell \leq rank(\mathcal{T}) \leq \prod_{\ell=1}^L R_\ell,$$ 
which gives us an upper and lower bound for the rank of $\mathcal{T}$. 

				Instead of considering $\mathcal{S}$ as a tensor in $\mathbb{K}^{I_1 \times \ldots \times I_L}$ with lots of zeros, we can focus only on the dense part, hence we consider $\mathcal{S}$ as a tensor in $\mathbb{K}^{R_1 \times \ldots \times R_L}$. This possibility had already been discussed shortly in~\ref{smaller-tensor}. In doing so, we can also discard the last $I_\ell - R_\ell$ columns of each $\textbf{U}^{(\ell)}$, hence we have $\textbf{U}^{(\ell)} \in \mathbb{K}^{I_\ell \times R_\ell} $. The equality $\mathcal{T} = (\textbf{U}^{(1)}, \ldots, \textbf{U}^{(L)}) \cdot \mathcal{S}$ remains intact after these changes, we only removed the unnecessary terms. Note that we passed from $\displaystyle\prod_{\ell=1}^L I_\ell$ terms to just $\displaystyle\prod_{\ell=1}^L R_\ell$.  The procedure of obtaining such decomposition clearly shows that $\mathcal{S}$, after deleting the unnecessary zeros, is a \emph{compressed version of $\mathcal{T}$}. The MLSVD in compressed form can also be called a \emph{reduced MLSVD}, in contrast to the full MLSVD of theorem~\ref{MLSVD}.
				
				\begin{remark} \label{mlsvd-cpd}
					After computing the MLSVD of $\mathcal{T}$, notice that computing a CPD for $\mathcal{S}$ is equivalent to computing a CPD for $\mathcal{T}$. Indeed, if $\mathcal{S} = \displaystyle \sum_{r=1}^R \textbf{w}_r^{(1)} \otimes \ldots \otimes \textbf{w}_r^{(L)}$, then we can write $\mathcal{S} = (\textbf{W}^{(1)}, \ldots, \textbf{W}^{(L)}) \cdot \mathcal{I}_{R \times \ldots \times R}$, which implies that
					$$\mathcal{T} = (\textbf{U}^{(1)}, \ldots, \textbf{U}^{(L)}) \cdot \left( (\textbf{W}^{(1)}, \ldots, \textbf{W}^{(L)}) \cdot \mathcal{I}_{R \times \ldots \times R} \right) = $$
					$$ = (\textbf{U}^{(1)} \textbf{W}^{(1)}, \ldots, \textbf{U}^{(L)} \textbf{W}^{(L)}) \cdot \mathcal{I}_{R \times \ldots \times R}.$$
				\end{remark}
				
				Now we go back to our claim that the MLSVD is indeed a good generalization of the SVD. The next results demonstrate such claim.
				
				\begin{theorem}[L. De Lathauwer, B. De Moor, J. Vandewalle, 2000]
					In the case $L = 2$, the MLSVD is equal to the SVD of matrices.
				\end{theorem}
				
				\begin{theorem}[L. De Lathauwer, B. De Moor, J. Vandewalle, 2000] \label{MLSVD-properties}
					Let $\mathcal{T} = (\textbf{U}^{(1)}, \ldots, \textbf{U}^{(L)}) \cdot \mathcal{S}$ be  MLSVD of $\mathcal{T}$ and, for each $\ell = 1 \ldots L$, let $R_\ell$ be the largest index such that $\sigma_{R_\ell}^{(\ell)} > 0$. Then the following holds.
					\begin{enumerate}
						\item $rank_\boxplus(\mathcal{T}) = (R_1, \ldots, R_L)$.
						\item $\text{Im}(\mathcal{T}_{(\ell)}) = \text{span}(\textbf{U}_1^{(\ell)}, \ldots, \textbf{U}_{R_\ell}^{(\ell)})$ for each $\ell = 1 \ldots L$.
						\item $\ker(\mathcal{T}_{(\ell)}^\ast) = \text{span}(\textbf{U}_{R_\ell+1}^{(\ell)}, \ldots, \textbf{U}_{I_\ell}^{(\ell)})$ for each $\ell = 1 \ldots L$.
						\item $\| \mathcal{T} \|^2 = \| \mathcal{S} \|^2 = \displaystyle \sum_{r=1}^{R_1} \left( \sigma_r^{(1)} \right)^2 = \ldots = \sum_{r=1}^{R_L} \left( \sigma_r^{(L)} \right)^2$.
					\end{enumerate}
				\end{theorem}
				
				Item 1 of this theorem tells us that the values $R_\ell$ coincide with the rank mode-$\ell$ of $\mathcal{T}$. This is in agreement with the relation between the singular values of matrices and rank. All other items also have their immediate version for the classic SVD. A relevant SVD property that unfortunately does not extend to the MLSVD is that of the tensor closest to $\mathcal{T}$ with a fixed lower multilinear rank. In the matrix case we have the Eckart-Young theorem, which states that the rank-$\tilde{R}$ matrix closest to $\textbf{M} = \displaystyle \sum_{r=1}^R \sigma_r \textbf{u}_r \textbf{v}_r^\ast$ is  
$\tilde{\textbf{M}} = \displaystyle \sum_{r=1}^{\tilde{R}} \sigma_r \textbf{u}_r \textbf{v}_r^\ast$, where $\tilde{R} \leq R$ and the decomposition for $\textbf{M}$ is a SVD. $\tilde{\textbf{M}}$ is a low rank approximation for $\textbf{M}$ (the best one with rank $\tilde{R}$) and in this case we have the following theorem about the error of this approximation.	
			
				\begin{theorem} \label{SVD-truncation}
					With the notations above, we have that\footnote{Remember that we are using the Frobenius norm.} 
					$$\|\textbf{M} - \tilde{\textbf{M}}\|^2 = \sum_{r=\tilde{R}+1}^R \sigma_r^2.$$
				\end{theorem}
				
				We can obtain a truncation $\tilde{\mathcal{T}}$ for $\mathcal{T}$ in a totally analogous way. Let $rank_\boxplus(\mathcal{T}) = (R_1, \ldots, R_L) $. Given a lower multilinear rank $(\tilde{R}_1, \ldots, \tilde{R}_L) < (R_1, \ldots, R_L)$, we obtain a truncation by zeroing all coefficients $\tilde{t}_{i_1 \ldots i_L} $ when some of its indexes satisfies $i_\ell > \tilde{R}_\ell$. Then we have the a tensor of multilinear rank $ (\tilde{R}_1, \ldots, \tilde{R}_L) $ given by
				
				$$\tilde{\mathcal{T}} = \sum_{i_1=1}^{\tilde{R}_1} \ldots \sum_{i_L=1}^{\tilde{R}_L} t_{i_1 \ldots i_L} \ \textbf{U}_{: i_1}^{(1)} \otimes \ldots \otimes \textbf{U}_{: i_L}^{(L)}.$$
				
				The following result first appeared in \cite{mlsvd} and the first correct proof appeared in \cite{sthosvd}. 
				
				\begin{theorem} \label{energy-mlsvd}
					With the notation above, we have that
					$$\|\mathcal{T} - \tilde{\mathcal{T}}\|^2 \leq \sum_{r = \tilde{R}_1+1}^{R_1} \left( \sigma_r^{(1)} \right)^2 + \ldots + \sum_{r = \tilde{R}_L+1} ^{R_L} \left( \sigma_r^{(L)} \right)^2.$$		
				\end{theorem}
	
				In general, $\tilde{\mathcal{T}}$ is a very good low multilinear rank approximation for $\mathcal{T}$, but unlike the matrix case, it is not necessarily the best approximation with low multilinear rank $(\tilde{R}_1, \ldots, \tilde{R}_L)$. At the time we know that \cite{grasedyck}
				
				$$\| \mathcal{T} - \tilde{\mathcal{T}} \| \leq \sqrt{L} \min_{rank_\boxplus(\mathcal{X}) = (\tilde{R}_1, \ldots, \tilde{R}_L)} \| \mathcal{T} - \mathcal{X} \|.$$  
				
				The current algorithms to obtain the best approximation usually use the MLSVD to obtain $\tilde{\mathcal{T}}$ and then use this tensor (which is already very close to $\mathcal{T}$) as the starting point for some iterative algorithm \cite{grassmann}.
				
				A few words about the cost of algorithm~\ref{MLSVD-alg} are necessary. Start noting that the construction of the unfoldings of $\mathcal{T}$ requires noncontiguous memory accesses. Although the computational time to achieve this may be noticeable, we disregard this cost because the cost to compute the SVDs dominates the algorithm. In general, the unfolding $\mathcal{T}_{(\ell)}$ probably will have much more columns than rows.

				The procedure described in~\ref{MLSVD-SVD} permits one to compute its singular values $\sigma_k^{(\ell)}$ and its left singular vectors $\textbf{U}_{:k}^{(\ell)}$ with low cost since the left singular vectors are the eigenvectors obtained in lemma~\ref{svd-transpose-lemma}. With this approach we don't compute the right singular vectors, but that is fine since they are not used for the MLSVD. In fact, not computing them is a good thing, we avoid unnecessary computations (note that this wouldn't be possible if we computed the SVD of $\mathcal{T}_{(\ell)}$ directly). The cost of this approach is of $\mathcal{O}\left( I_\ell \displaystyle\prod_{\ell'=1}^L I_{\ell'} + I_\ell^ 3 \right)$ flops. The case where $\mathcal{T}_{(\ell)}$ have less columns than rows is a bit more problematic, we have a unbalanced tensor (definition~\ref{unbalanced}) to deal with. We could take a similar route and compute the eigenvalues of $\mathcal{T}_{(\ell)}^T \mathcal{T}_{(\ell)}$ but now the eigenvectors obtained are the right singular vectors of $\mathcal{T}_{(\ell)}$, which we don't need to use. This is a situation where we have to use the SVD and lose some performance. Computing the reduced SVD in this case has a cost of, at least, $\mathcal{O}\left( 2 I_\ell \displaystyle\prod_{\ell' \neq \ell} I_{\ell'}^2 + 2 \displaystyle\prod_{\ell' \neq \ell} I_{\ell'}^3 \right)$ flops. Instead of computing the full SVD and deal with these high costs, we can take another route and compute a truncated SVD. Since each $\mathcal{T}_{(\ell)}$ has rank not bigger than $\min\{I_\ell, R\}$, we define $P_\ell = \min\{I_\ell, R\}$ and compute the truncated SVD with $P_\ell$ singular values. This approach has a cost of $\mathcal{O}\left( \displaystyle P_\ell \prod_{\ell'=1}^L I_{\ell'} \right)$ flops, which is much better than the previous cost. If we are willing to lose some precision to gain speed, it is possible to use randomized algorithms for the truncated SVD \cite{rand_svd}. This algorithm has a cost of $\mathcal{O}\left( \displaystyle \log(P_\ell) \prod_{\ell'=1}^L I_{\ell'} \right)$ flops and the loss in precision is negligible. Since we need to perform this computation for each $\ell=1 \ldots L$, the overall cost is of $\mathcal{O}\left( \displaystyle\sum_{\ell=1}^L \log(P_\ell) \prod_{\ell'=1}^L I_{\ell'} \right)$ flops.
				
				The approach just described to compute the MLSVD can be referred as the \emph{classic truncated MLSVD}, where the idea is to compute the truncated SVD of each unfolding of $\mathcal{T}$. A much faster and slightly less precise approach is the \emph{sequentially truncated MLSVD} algorithm introduced by N. Vannieuwenhove, R. Vandebril and K. Meerberge \cite{sthosvd}, which consists of interlacing the computation of the core tensor and the factor matrices. The cost is of this algorithm is of $\mathcal{O}\left( \displaystyle\sum_{\ell=1}^L \log(P_\ell) \displaystyle\prod_{\ell'=1}^{\ell-1} R_{\ell'} \prod_{\ell'=\ell}^L I_{\ell'} + R_\ell^2 \displaystyle\prod_{\ell'=1}^{\ell-1} R_{\ell'} \prod_{\ell'=\ell+1}^L I_{\ell'} \right)$ flops. The first is cost is due to the truncated SVD and the second one is due to the matrix-matrix multiplication, see algorithm below.\bigskip
				
				\begin{algorithm}[ST-MLSVD] \label{ST-MLSVD-alg}
					$ $\\
					\textbf{Input:} $\mathcal{T} \in \mathbb{R}^{I_1 \times \ldots \times I_L}$\vspace{.4cm}\\
					$\mathcal{S}^{(0)} = \mathcal{T}$\\
					\verb|for | $\ell=1 \ldots L$\\
						$\hspace{1cm} \mathcal{S}_{(\ell)}^{(\ell-1)} \leftarrow $ \verb|Unfolding| $(\mathcal{S}^{(\ell-1)}, \ell)$\\
						$\hspace{1cm} \textbf{U}^{(\ell)}, \Sigma^{(\ell)}, \textbf{V}^{(\ell)} \leftarrow $ \verb|SVD| $(\mathcal{S}_{(\ell)})$\\
						$\hspace{1cm} \text{Set } \mathcal{S}^{(\ell)}, \text{ where } \mathcal{S}_{(\ell)}^{(\ell)} = \Sigma^{(\ell)} \big( \textbf{V}^{(\ell)} \big)^T$\vspace{.4cm}\\
					\textbf{Output:} $\textbf{U}^{(\ell)} \in \mathbb{R}^{I_\ell \times R_\ell}, \Sigma^{(\ell)} \in \mathbb{R}^{R_\ell \times R_\ell}$ for $\ell = 1 \ldots L$ and $\mathcal{S} \in \mathbb{R}^{R_1 \times \ldots \times R_L}$
				\end{algorithm}
				
				The first version of Tensor Fox only used the classic truncated algorithm to compute the MLSVD, but now it has the possibility to switch between this one and the sequentially truncated algorithm.   

%% file: CHAPTER_3.tex
\chapter{Gauss-Newton algorithm}\label{cap-3} 
	In this chapter we really begin to touch the computational aspects of the CPD. Some of the history and known algorithms are discussed, and then we go to nonlinear least squares algorithms, which includes the Gauss-Newton method. This is the method used by Tensor Fox. The main challenge is the approximated Hessian, which is singular by construction. To overcome this issue we introduce a regularization with a new kind of matrix not seen in the literature. In the end of the chapter we show one of the main contributions of this work: a set of formulas for matrix-vector products that appear in DGN applied to CPD which result in a reduction of complexity with respect to a straightforward calculation.
	
	We also give the constructive proofs of some already known results from tensor algebra. The reason of this is to give the reader an understanding of the computational aspects related of our problems, this makes a nice parallel with actual coding necessary to make Tensor Fox. Therefore all proofs are given in coordinates. 

	\section{Preliminaries}
		In 1927, Hitchcock \cite{hitchcock1, hitchcock2} proposed the idea of the CPD. This idea only became popular in 1970 inside the psychometrics community, called CANDECOMP (canonical decomposition) by Carrol and Chang \cite{carrol}, and PARAFAC (parallel factors) by Harshman \cite{harshman}. Gradually the CPD began to be applied in more areas and today is a successful tool for multi-dimensional data analysis and prediction, such as blind source separation, food industry, dimensionality reduction, pattern/image recognition, machine learning and data mining \cite{bro, cichoki2009, cichoki2014, savas, smilde, anandkumar}.  
	
		As we mentioned in section 1.4, computing the rank of a tensor is a NP-hard problem. Given a tensor $\mathcal{T}$ one may try to find its rank together with its CPD by searching for a best fit, that is, for each $R = 1, 2, \ldots$, compute a rank-$R$ CPD for $\mathcal{T}$ until some criteria of good fit is reached. This is a possible but risky procedure because of the border rank phenomenon. More precisely, if $\underline{rank}(\mathcal{T}) < rank(\mathcal{T})$, then it is possible to obtain arbitrarily good fits with low rank approximations for $\mathcal{T}$. This may cause problems in practice \cite{mitchell, paatero, mitchell2}. 
	
		Fix a tensor $\mathcal{T} \in \mathbb{K}^{I_1 \times \ldots \times I_L}$ and a value $R$. Our problem is to compute a rank-$R$ approximation $\tilde{\mathcal{T}} = (\textbf{W}^{(1)}, \ldots, \textbf{W}^{(L)}) \cdot \mathcal{I}_{R \times \ldots \times R} \in \mathbb{K}^{I_1 \times \ldots \times I_L}$. More precisely, $\tilde{\mathcal{T}}$ is given by
	
		\begin{equation} \label{compute_cpd}		
			\tilde{\mathcal{T}} = \sum_{r=1}^R \textbf{w}_r^{(1)} \otimes \ldots \otimes \textbf{w}_r^{(L)},
		\end{equation}
where $\textbf{W}^{(\ell)} = [\textbf{w}_1^{(\ell)}, \ldots, \textbf{w}_R^{(\ell)}]^T$ are the factor matrices, for $\ell = 1, \ldots L$, and $\mathcal{I}_{R \times \ldots \times R} \in \mathbb{K}^{R \times \ldots \times R}$ is the diagonal tensor with $L$ dimensions.

		Note that we don't assume $rank(\mathcal{T}) = R$, we only assume that such $R$ is given. Frequently the data at hand is noisy, because of this the rank of $\mathcal{T}$ may be different from the noiseless data, and this is the one we are interested in. Usually is a good idea to compute low rank approximations since the noise can cause the rank to be typical or generic, that is, the noisy tensors are not special (because they are ``everywhere'') but the objective tensor is. We want this approximation to be such that $\| \mathcal{T} - \tilde{\mathcal{T}} \|$ is smallest as possible. The most common approach is to formulate the problem as an unconstrained minimization problem
		
		\begin{equation} \label{min}
			\min_{\tilde{\mathcal{T}}} \frac{1}{2} \| \mathcal{T} - \tilde{\mathcal{T}} \|^2.
		\end{equation}	
		
		Define the function $F:\mathbb{K}^{R \sum_{\ell=1}^L I_\ell} \to \mathbb{R}$ given by
		
		$$F(\textbf{w}) = \frac{1}{2} \| \mathcal{T} - (\textbf{W}^{(1)}, \ldots, \textbf{W}^{(L)}) \cdot \mathcal{I}_{R \times \ldots \times R} \|^2,$$
where $\textbf{w}$ is the just the vector obtained by vertically stacking\footnote{Remember, from example~\ref{matrix_mult}, that $vec(\textbf{W}^{(\ell)})$ is the vector obtained by vertically stacking all columns of $\textbf{W}^{(\ell)}$.} all columns of all factor matrices. More precisely,
		$$\textbf{w} = \left[
		\begin{array}{c}
			vec(\textbf{W}^{(1)})\\
			\vdots\\
			vec(\textbf{W}^{(L)})
		\end{array}		
		\right].$$		
		
		Note that solving~\ref{min} is equivalent to minimizing $F$. In this work we will be considering only real tensors, but the complex case is considered in \cite{tensorlab2} for instance. The first thing we should note is that any solution of~\ref{min} is a critical point of $F$ so it is of interest to derive an expression for the derivative of $F$ and its critical points. 
		
		\begin{lemma}
			For any vectors $\textbf{x}^{(1)}, \textbf{y}^{(1)} \in \mathbb{R}^{I_1}, \ldots, \textbf{x}^{(L)}, \textbf{y}^{(L)} \in \mathbb{R}^{I_L}$ we have that
			$$\prod_{\ell=1}^L \langle \textbf{x}^{(\ell)}, \textbf{y}^{(\ell)} \rangle = \sum_{i_1=1}^{I_1} \ldots \sum_{i_L=1}^{I_L} x_{i_1}^{(1)} \ldots x_{L}^{(L)} y_{i_1}^{(1)} \ldots y_{i_L}^{(L)}.$$
		\end{lemma}
		
		\begin{theorem} \label{partialF}
			For all $\ell' = 1 \ldots L$, $i' = 1 \ldots I_\ell$ and $r' = 1 \ldots R$, we have			
			$$\frac{\partial F}{w_{i' r'}^{(\ell')}}(\textbf{w}) = - \mathcal{T}(\textbf{w}_{r'}^{(1)}, \ldots, \textbf{e}_{i'}^{(\ell')}, \ldots, \textbf{w}_{r'}^{(L)}) + \sum_{r=1}^R w_{i' r}^{(\ell')} \prod_{\ell \neq \ell'} \langle \textbf{w}_r^{(\ell)}, \textbf{w}_{r'}^{(\ell)} \rangle.$$
		\end{theorem}
		
		\textbf{Proof:} $$\frac{\partial F}{w_{i' r'}^{(\ell')}}(\textbf{w}) = \frac{1}{2} \frac{\partial}{\partial w_{i' r'}^{(\ell')}} \left( \sum_{i_1=1}^{I_1} \ldots \sum_{i_L=1}^{I_L} \left( t_{i_1 \ldots i_L} - \sum_{r=1}^R w_{i_1 r}^{(1)} \ldots w_{i_L r}^{(L)} \right)^2 \right) = $$
		
		$$ = \frac{1}{2} \sum_{i_1=1}^{I_1} \ldots \sum_{i_L=1}^{I_L} \frac{\partial}{\partial w_{i' r'}^{(\ell')}} \left( t_{i_1 \ldots i_L} - \sum_{r=1}^R w_{i_1 r}^{(1)} \ldots w_{i_L r}^{(L)} \right)^2 = $$
		
		$$ = \frac{1}{2} \sum_{i_1=1}^{I_1} \ldots \sum_{i_L=1}^{I_L} 2 \left( t_{i_1 \ldots i_L} - \sum_{r=1}^R w_{i_1 r}^{(1)} \ldots w_{i_L r}^{(L)} \right) \cdot \frac{\partial}{\partial w_{i' r'}^{(\ell')}} \left( t_{i_1 \ldots i_L} - \sum_{r=1}^R w_{i_1 r}^{(1)} \ldots w_{i_L r}^{(L)} \right) = $$
		
		$$ = \sum_{i_1=1}^{I_1} \ldots \sum_{i_L=1}^{I_L} \left( t_{i_1 \ldots i_L} - \sum_{r=1}^R w_{i_1 r}^{(1)} \ldots w_{i_L r}^{(L)} \right) \cdot \left( - \frac{\partial}{\partial w_{i' r'}^{(\ell')}} w_{i_1 r'}^{(1)} \ldots w_{i_L r'}^{(L)} \right) = $$
		
		$$ = \sum_{i_1=1}^{I_1} \ldots \sum_{i_L=1}^{I_L} \left( t_{i_1 \ldots i_L} - \sum_{r=1}^R w_{i_1 r}^{(1)} \ldots w_{i_L r}^{(L)} \right) \cdot \left( - \prod_{\ell \neq \ell'} w_{i_\ell r'}^{(\ell)} \cdot \frac{\partial}{\partial w_{i' r'}^{(\ell')}} w_{i_{\ell'} r'}^{(\ell')} \right) = $$
		
		$$ = \sum_{i_1=1}^{I_1} \ldots\sum_{i_{\ell'-1}=1}^{I_{\ell'-1}} \ \sum_{i_{\ell'+1}=1}^{I_{\ell'+1}} \ldots \sum_{i_L=1}^{I_L} \left( t_{i_1 \ldots i_L} - \sum_{r=1}^R w_{i' r}^{(\ell')} \prod_{\ell \neq \ell'} w_{i_\ell r}^{(\ell)} \right) \cdot \left( - \prod_{\ell \neq \ell'} w_{i_\ell r'}^{(\ell)} \cdot \frac{\partial}{\partial w_{i' r'}^{(\ell')}} w_{i' r'}^{(\ell')} \right) = $$
		
		$$ = \sum_{i_1=1}^{I_1} \ldots\sum_{i_{\ell'-1}=1}^{I_{\ell'-1}} \ \sum_{i_{\ell'+1}=1}^{I_{\ell'+1}} \ldots \sum_{i_L=1}^{I_L} \left( t_{i_1 \ldots i_L} - \sum_{r=1}^R w_{i' r}^{(\ell')} \prod_{\ell \neq \ell'} w_{i_\ell r}^{(\ell)} \right) \cdot \left( - \prod_{\ell \neq \ell'} w_{i_\ell r'}^{(\ell)} \right) = $$
		
		$$\hspace{-1cm} = - \sum_{i_1=1}^{I_1} \ldots\sum_{i_{\ell'-1}=1}^{I_{\ell'-1}} \ \sum_{i_{\ell'+1}=1}^{I_{\ell'+1}} \ldots \sum_{i_L=1}^{I_L}  t_{i_1 \ldots i_L} \prod_{\ell \neq \ell'} w_{i_\ell r'}^{(\ell)} \ + \ \sum_{r=1}^R w_{i' r}^{(\ell')} \sum_{i_1=1}^{I_1} \ldots\sum_{i_{\ell'-1}=1}^{I_{\ell'-1}} \ \sum_{i_{\ell'+1}=1}^{I_{\ell'+1}} \ldots \sum_{i_L=1}^{I_L} \prod_{\ell \neq \ell'} w_{i_\ell r}^{(\ell)} \prod_{\ell \neq \ell'} w_{i_\ell r'}^{(\ell)}.$$
		
		$$ = - \mathcal{T}(\textbf{w}_{r'}^{(1)}, \ldots, \textbf{e}_{i'}^{(\ell')}, \ldots, \textbf{w}_{r'}^{(L)}) + \sum_{r=1}^R w_{i' r}^{(\ell')} \prod_{\ell \neq \ell'} \langle \textbf{w}_r^{(\ell)}, \textbf{w}_{r'}^{(\ell)} \rangle,$$		
where the last term was calculated using the previous lemma.$\hspace{4.5cm}\square$\bigskip	

		\begin{corollary}
			Let $\displaystyle\frac{\partial F}{\textbf{w}_{r'}^{(\ell')}}(\textbf{w}) = \left[ \frac{\partial F}{w_{1 r}^{(\ell')}}, \ldots, \frac{\partial F}{w_{I_{\ell'} r}^{(\ell')}} \right]^T$, then			
			$$\hspace{-.5cm}\frac{\partial F}{\textbf{w}_{r'}^{(\ell')}}(\textbf{w}) = 
			- \left[ \begin{array}{c} 
				\mathcal{T}(\textbf{w}_{r'}^{(1)}, \ldots, \textbf{e}_1^{(\ell')}, \ldots, \textbf{w}_{r'}^{(L)})\\ 
				\vdots\\ 
				\mathcal{T}(\textbf{w}_{r'}^{(1)}, \ldots, \textbf{e}_{I_{\ell'}}^{(\ell')}, \ldots, \textbf{w}_{r'}^{(L)})
			\end{array} \right] 
			+ \sum_{r=1}^R \prod_{\ell \neq \ell'} \langle \textbf{w}_r^{(\ell)}, \textbf{w}_{r'}^{(\ell)} \rangle \textbf{w}_r^{(\ell')}.$$
		\end{corollary}
		
		Before we start talking about algorithms, it is relevant to mention what are the biggest challenges we will be facing when designing an algorithm to solve~\ref{min}. In the following, let $\mathcal{T} = \displaystyle \sum_{r=1}^R \textbf{v}_r^{(1)} \otimes \ldots \otimes \textbf{v}_r^{(L)}$. 
		
		\begin{enumerate}
			\item As mentioned in section 1.5, finding the best rank-$R$ approximation may be a ill-posed problem.
			\item Ill-posed problems are not ``rare'' in the sense that the set of ill-posed problems (i.e., the set of tensors such that~\ref{min} is ill-posed) usually have positive volume \cite{lim}.
			\item We say $\mathcal{T}$ has a \emph{bottleneck} if there is a mode $1 \leq \ell \leq L$ such that at least two vectors $\textbf{v}_{r'}^{(\ell)}$ and $\textbf{v}_{r''}^{(\ell)}$ are almost collinear. This means $\mathcal{T}$ has two ``problematic'' vectors in $\mathbb{R}^{I_\ell}$, and they are considered problematic because two different rank one terms will have almost collinear vectors at the $\ell$-th position.
			\item We say $\mathcal{T}$ has a \emph{swamp} if there are bottlenecks in all modes. 
			\item We say $\mathcal{T}$ has a \emph{degeneracy} if some factors diverge but cancel out when the process of computing a best fit is performed. This is the same degeneracy mentioned at the end of section 1.5. 
		\end{enumerate}	
		
		The last three items in this list constitutes a classification of possible causes of bad numerical performance in numerical algorithms to compute the CPD. This classification was proposed already in 1989 \cite{kruskal4} and still is used today.  To get a better understand of bottlenecks and swamps, we illustrate these occurrences in the next two examples.
		
		\begin{example}[Bottleneck]
			Consider the tensor $\mathcal{T} \in \mathbb{R}^{2 \times 2 \times 2}$ given by			
			$$\mathcal{T} = 
			\left[ \begin{array}{c}
				1\\
				2
			\end{array} \right]
			\otimes \left[ \begin{array}{c}
				0\\
				1
			\end{array} \right]
			\otimes \left[ \begin{array}{c}
				3\\
				4
			\end{array} \right] \ + \ 
			\left[ \begin{array}{c}
				1.01\\
				1.9
			\end{array} \right]
			\otimes \left[ \begin{array}{c}
				1\\
				1
			\end{array} \right]
			\otimes \left[ \begin{array}{c}
				2\\
				1
			\end{array} \right] \ + \ 
			\left[ \begin{array}{c}
				0\\
				2
			\end{array} \right]
			\otimes \left[ \begin{array}{c}
				2\\
				2
			\end{array} \right]
			\otimes \left[ \begin{array}{c}
				4\\
				4
			\end{array} \right].$$
			
			The vectors of the first mode are
			
			$$\left[ \begin{array}{c}
				1\\
				2
			\end{array} \right], \quad
			\left[ \begin{array}{c}
				1.01\\
				1.9
			\end{array} \right], \quad
			\left[ \begin{array}{c}
				0\\
				2
			\end{array} \right].$$
			
			The vectors of the second mode are
			
			$$\left[ \begin{array}{c}
				0\\
				1
			\end{array} \right],
			\left[ \begin{array}{c}
				1\\
				1
			\end{array} \right],
			\left[ \begin{array}{c}
				2\\
				2
			\end{array} \right].$$
			
			The vectors of the third mode are
			
			$$\left[ \begin{array}{c}
				3\\
				4
			\end{array} \right],
			 \left[ \begin{array}{c}
				2\\
				1
			\end{array} \right],
			\left[ \begin{array}{c}
				4\\
				4
			\end{array} \right].$$			
			
			We can note that the vectors of the first mode introduces a possible bottleneck in $\mathcal{T}$. More precisely, consider the two vectors
			$\left[ \begin{array}{c}
				1\\
				2
			\end{array} \right]$
and
		$\left[ \begin{array}{c}
				1.01\\
				1.9
			\end{array} \right]$
of the first mode. Compared to the vectors of the other modes, these two vectors indeed are very close to be equal (hence collinear) so we can consider that $\mathcal{T}$ has a bottleneck. This is a single bottleneck, but it is possible to have multiple bottlenecks. The extreme case is when we have bottlenecks at all modes, which represents a swamp.			
		\end{example}
		
		\begin{example}[Swamp]
			Consider the tensor $\mathcal{T} \in \mathbb{R}^{2 \times 2 \times 2}$ given by
			
			$$\mathcal{T} = 
			\left[ \begin{array}{c}
				1\\
				2
			\end{array} \right]
			\otimes \left[ \begin{array}{c}
				0\\
				1
			\end{array} \right]
			\otimes \left[ \begin{array}{c}
				3\\
				4
			\end{array} \right] \ + \ 
			\left[ \begin{array}{c}
				1.01\\
				1.9
			\end{array} \right]
			\otimes \left[ \begin{array}{c}
				0\\
				1
			\end{array} \right]
			\otimes \left[ \begin{array}{c}
				1\\
				1
			\end{array} \right] \ + \ 
			\left[ \begin{array}{c}
				0.09\\
				2.01
			\end{array} \right]
			\otimes \left[ \begin{array}{c}
				10^{-6}\\
				1
			\end{array} \right]
			\otimes \left[ \begin{array}{c}
				2.99\\
				3.99
			\end{array} \right].$$
			
			The vectors of the first mode are
			
			$$\left[ \begin{array}{c}
				1\\
				2
			\end{array} \right], \quad
			\left[ \begin{array}{c}
				1.01\\
				1.9
			\end{array} \right], \quad
			\left[ \begin{array}{c}
				0.09\\
				2.01
			\end{array} \right].$$
			
			The vectors of the second mode are
			
			$$\left[ \begin{array}{c}
				0\\
				1
			\end{array} \right],
			\left[ \begin{array}{c}
				0\\
				1
			\end{array} \right],
			\left[ \begin{array}{c}
				10^{-6}\\
				1
			\end{array} \right].$$
			
			The vectors of the third mode are
			
			$$\left[ \begin{array}{c}
				3\\
				4
			\end{array} \right],
			 \left[ \begin{array}{c}
				1\\
				1
			\end{array} \right],
			\left[ \begin{array}{c}
				2.99\\
				3.99
			\end{array} \right].$$			
			
			In all modes there are at least two almost collinear vectors, therefore $\mathcal{T}$ has a swamp.
		\end{example}	
		
		Now we briefly describe the most used approaches to compute a CPD.

	\section{Alternating least squares} 
		Several different approaches to solving~\ref{min} were proposed in the past years, but before that, a single algorithm were used from decades: the \emph{alternating least squares} (ALS). For decades this was considered to be the ``workhorse'' algorithm to compute the CPD. We present the algorithm for third order tensor and generalize after that. From this point, we write $m \times n \times p$ instead $I_1 \times I_2 \times I_3$ for third order tensor.  
		
		Consider a tensor $\mathcal{T} \in \mathbb{R}^{m \times n \times p}$ and suppose we want to compute a rank-$R$ CPD for $\mathcal{T}$. Let $\tilde{\mathcal{T}} = \displaystyle \sum_{r=1}^R \textbf{x}_r \otimes \textbf{y}_r \times \textbf{z}_r$ be the approximating tensor, and $\textbf{X} = [\textbf{x}_1, \ldots, \textbf{x}_R] \in \mathbb{R}^{m \times R}$, $\textbf{Y} = [\textbf{y}_1, \ldots, \textbf{y}_R] \in \mathbb{R}^{n \times R}$, $\textbf{Z} = [\textbf{z}_1, \ldots, \textbf{z}_R] \in \mathbb{R}^{p \times R}$ the factor matrices.
		
		The first step of the ALS algorithm consists in generating a initial tensor $\tilde{\mathcal{T}}$ to start the iterations. The method of initialization is not part of the algorithm so we won't discuss it here. Now fix $\textbf{Y}, \textbf{Z}$ and solve the minimization problem		
		$$\min_\textbf{X} \| \mathcal{T} - (\textbf{X}, \textbf{Y}, \textbf{Z}) \cdot \mathcal{I}_{R \times \times R} \|.$$
Note that we solve it just for $\textbf{X}$. After this is done, fix $\textbf{X}, \textbf{Z}$ and solve the minimization problem 		$$\min_\textbf{Y} \| \mathcal{T} - (\textbf{X}, \textbf{Y}, \textbf{Z}) \cdot \mathcal{I}_{R \times \times R} \|$$
for $\textbf{Y}$. After that, fix $\textbf{X}, \textbf{Y}$ and solve the minimization problem 		
		$$\min_\textbf{Z} \| \mathcal{T} - (\textbf{X}, \textbf{Y}, \textbf{Z}) \cdot \mathcal{I}_{R \times \times R} \|$$
for $\textbf{Z}$. Once the three factor matrices are updated we repeat this procedure all over again, updating as many times as necessary. The name comes from the fact we are alternating the factor matrices to be solved at each linear least squares problem. To see that these minimization problems indeed are linear least squares problems, note that we can rewrite the first one as	
		$$\min_\textbf{X} \| \mathcal{T}_{(1)} - \tilde{\mathcal{T}}_{(1)} \| = \min_\textbf{X} \| \mathcal{T}_{(1)} - \textbf{X} \cdot (\textbf{Z} \odot \textbf{Y})^T \|,$$
where the equality $\tilde{\mathcal{T}}_{(1)} = \textbf{X} \cdot (\textbf{Z} \odot \textbf{Y})^T$ comes from theorem~\ref{unfolding-formula2}. The solution of this problem is given explicitly by $\textbf{X} = \mathcal{T}_{(1)} \cdot ((\textbf{Z} \odot \textbf{Y})^T)^\dagger$. Theorem~\ref{special-products} gives a formula for the pseudoinverse of a Khatri-Rao product, so we have 

		$$((\textbf{Z} \odot \textbf{Y})^T)^\dagger = ((\textbf{Z} \odot \textbf{Y})^\dagger)^T = \left( ((\textbf{Z}^T\textbf{Z}) \ast (\textbf{Y}^T\textbf{Y}))^\dagger (\textbf{Z} \odot \textbf{Y})^T \right)^T = $$
	
		$$ = (\textbf{Z} \odot \textbf{Y}) \big( \underbrace{((\textbf{Z}^T\textbf{Z}) \ast (\textbf{Y}^T\textbf{Y}))^T}_{\text{symmetric}} \big)^\dagger = (\textbf{Z} \odot \textbf{Y}) \left( (\textbf{Z}^T\textbf{Z}) \ast (\textbf{Y}^T\textbf{Y}) \right)^\dagger.$$
	
		This formulation only requires computing the pseudoinverse of a $R \times R$ matrix rather than a $np \times R$, which is the size of $(\textbf{Z} \odot \textbf{Y})^T$. Below we give the algorithm for the general case.
		
		\begin{algorithm}[ALS]
			$ $\\
			\textbf{Input:} $\mathcal{T}, R$\vspace{.4cm}\\
			\verb|Initialize | $\textbf{W}^{(1)}, \ldots, \textbf{W}^{(L)}$\\
			\verb|repeat|\\
				$\hspace{1cm}$ \verb|for | $\ell = 1 \ldots L$\\
					$\hspace{2cm} \textbf{W}^{(\ell)} \leftarrow \verb|argmin|_{\textbf{W}^{(\ell)}} \| \mathcal{T}_{(\ell)} - \textbf{W}^{(\ell)} \cdot (\textbf{W}^{(L)} \odot \ldots \odot \textbf{W}^{(\ell+1)} \odot \textbf{W}^{(\ell-1)} \odot \ldots \odot \textbf{W}^{(1)})^T \|$\\
			\verb|until stopping criteria is met|\vspace{.4cm}\\		
			\textbf{Output:} $\textbf{W}^{(1)}, \ldots, \textbf{W}^{(L)}$
		\end{algorithm}
		
		The ALS algorithm is simple to understand and to implement, but it has some serious drawbacks. Usually it take many iterations to converge, it is not guaranteed to converge to a global minimum or even to a stationary point, and the final solution can depend heavily on the initialization. Furthermore, the ALS algorithm is known from its convergence problems in the presence of bottlenecks and swamps. Several approaches were implemented in order to improve the ALS performance \cite{comon2, lathauwer, tensorlab2, kolda, kolda3} and still today there are people trying improve ALS. 		
	
	\section{Optimization methods}
		Because of the ALS limitations, researches tried to investigate different approaches to~\ref{min}. A natural idea is to consider it just as an unconstrained minimization problem and apply the known algorithms. Gradient-based algorithms are the main choice in this case. In the literature one can find theoretical and experimental work with \emph{nonlinear conjugate gradient method} and \emph{limited-memory BFGS method} \cite{tensorlab2, cp_opt}. Improvements to these methods were tried, such as adding line search\footnote{It seems that the \emph{Mor\' e -Thuente} line search is the most popular choice of line search for tensor CPD, see \cite{more} for more information.}, regularization and dogleg trust region \cite{madsen}.  
	
	\section{Nonlinear least squares} 
		The method of \emph{least squares} is a standard approach to approximate the solution of overdetermined systems, that is, sets of equations in which there are more equations than unknowns. The term  ``least squares'' comes from the fact that we want to minimize the overall sum of the squares of the residuals made in the results of every single equation. Let $(\textbf{x}^{(i)}, y^{(i)})_{i=1 \ldots n}$ be a dataset of $n$ points, where $\textbf{x}^{(i)} \in \mathbb{R}^p$ is the independent variable and $y^{(i)} \in \mathbb{R}$ is the dependent variable. The \emph{model} is a function $h_\textbf{w}:\mathbb{R}^p \to \mathbb{R}$ (also called the \emph{hypothesis}) in the variable $\textbf{x}$ with parameters $\textbf{w} = [ w_1, \ldots, w_m]^T$ to be adjusted. Our objective is to find the parameters that best fits the data. The ``best fit'' is in the sense that the choice of $\textbf{w}$ should minimize the value
		
		$$\sum_{i=1}^n \left( y^{(i)} - h_\textbf{w}(\textbf{x}^{(i)}) \right)^2.$$ 
		
		Denote $f_i(\textbf{w}) = y^{(i)} - h_\textbf{w}(\textbf{x}^{(i)})$. We have that each $f_i$ is a \emph{residual} of the model. The \emph{error function} (or \emph{cost function}) is the function $F: \mathbb{R}^m \to \mathbb{R}$ given by
		
		$$F(\textbf{w}) = \frac{1}{2} \sum_{i=1}^n \left( y^{(i)} - h_\textbf{w}(\textbf{x}^{(i)}) \right)^2 = \frac{1}{2} \sum_{i=1}^n f_i(\textbf{w})^2 = \frac{1}{2} \| f(\textbf{w}) \|^2,$$
where $f:\mathbb{R}^m \to \mathbb{R}^n$ is given by $f = (f_1, \ldots, f_n)$. Our objective is to minimize $F$. 

		\begin{remark}
			The ideal situation occurs when $\textbf{w}$ is such that $F(\textbf{w}) = 0$, so $h_\textbf{w}(\textbf{x}^{(i)}) = y^{(i)}$ for all $i=1 \ldots n$. This means $\textbf{w}$ is a solution of the system
			
			$$\left\{ \begin{array}{c}
				h_\textbf{w}(\textbf{x}^{(1)}) = y^{(1)}\\
				\vdots\\
				h_\textbf{w}(\textbf{x}^{(n)}) = y^{(n)}
			\end{array}\right..$$
			
			Note that this is a system of $n$ equations and $m$ variables. Since we always expect to have more data then parameters (i.e., $m < n$), this is a overdetermined system. In practice $\textbf{w}$ will not be an exact solution but an approximated solution for this system.
		\end{remark}  
		
		We call the model \emph{linear} when $h_\textbf{w}$ is a linear combination of the parameters, that is, $h_\textbf{w}$ is of the form $h_\textbf{w}(\textbf{x}) = \displaystyle\sum_{j=1}^m w_j \ \phi_j(\textbf{x})$, where each $\phi_j$ is any function of $\textbf{x}$. In this case one can write
		
		$$\min_\textbf{w} F(\textbf{w}) = \frac{1}{2} \min_\textbf{w} \left\| \left[ 
		\begin{array}{ccc} 
			\phi_1(\textbf{x}^{(1)}) & \ldots & \phi_m(\textbf{x}^{(1)})\\
			\vdots & & \vdots\\
			\phi_1(\textbf{x}^{(n)}) & \ldots & \phi_m(\textbf{x}^{(n)})
		\end{array} \right]
		\left[ 
		\begin{array}{c}
			w_1\\
			\vdots\\
			w_m
		\end{array} \right]
		- \left[ 
		\begin{array}{c}
			y^{(1)}\\
			\vdots\\
			y^{(n)}
		\end{array} \right] \right\|^2 = \frac{1}{2} \| \textbf{A}\textbf{w} - \textbf{y} \|^2,$$
where $a_{ij} = \phi_j(\textbf{x}^{(i)})$ and $\textbf{y} = [y^{(1)}, \ldots, y^{(n)}]^T$. The minimizer is the solution of the normal equation $\textbf{A}^T \textbf{A} \textbf{w} = \textbf{A}^T \textbf{y}$, and is given explicitly by $\textbf{w} = (\textbf{A}^T \textbf{A})^{-1} \textbf{A}^T \textbf{y}$. 

		We call the model \emph{nonlinear} when it is not linear. In this case there is no guarantee of a closed-form solution as in the linear case. 
		
		\begin{lemma} \label{gradF}
			Let $\textbf{J}_f(\textbf{w})$ be the Jacobian matrix of $f$ at $\textbf{w}$. Then
			
			$$\nabla F(\textbf{w}) = \textbf{J}_f^T(\textbf{w}) \cdot f(\textbf{w}).$$
		\end{lemma}	
			
		\textbf{Proof}: First note that $\displaystyle\frac{\partial F}{\partial w_j}(\textbf{w}) = \sum_{i=1}^n f_i(\textbf{w}) \frac{\partial f_i}{\partial w_j}(\textbf{w})$ for all $j = 1 \ldots m$. In particular, we have 

		$$\nabla F(\textbf{w}) = 
		\left[ \begin{array}{c}
			\sum_{i=1}^n f_i(\textbf{w}) \frac{\partial f_i}{\partial w_1}(\textbf{w})\\ 
			\vdots\\ 
			\sum_{i=1}^n f_i(\textbf{w}) \frac{\partial f_i}{\partial w_m}(\textbf{w})
		\end{array} \right]  = 
		\left[ \begin{array}{ccc}
			\frac{\partial f_1}{\partial w_1}(\textbf{w}) & \ldots & \frac{\partial f_n}{\partial w_1}(\textbf{w})\\
			\vdots & & \vdots\\
			\frac{\partial f_1}{\partial w_m}(\textbf{w}) & \ldots & \frac{\partial f_n}{\partial w_m}(\textbf{w})
		\end{array} \right]
		\left[ \begin{array}{c}
			f_1(\textbf{w})\\
			\vdots\\
			f_n(\textbf{w})
		\end{array} \right] = $$
		
		$$\hspace{2.5cm} = \left[ \begin{array}{ccc}
			\frac{\partial f_1}{\partial w_1}(\textbf{w}) & \ldots & \frac{\partial f_1}{\partial w_m}(\textbf{w})\\
			\vdots & & \vdots\\
			\frac{\partial f_n}{\partial w_1}(\textbf{w}) & \ldots & \frac{\partial f_n}{\partial w_m}(\textbf{w})
		\end{array} \right]^T
		\left[ \begin{array}{c}
			f_1(\textbf{w})\\
			\vdots\\
			f_n(\textbf{w})
		\end{array} \right] = 
		\textbf{J}_f^T(\textbf{w}) \cdot f(\textbf{w}). \hspace{2.5cm}\square$$

		Now we bring all this to the context of tensors. Let a tensor $\mathcal{T} \in \mathbb{R}^{I_1 \times \ldots \times I_L}$ and a rank $R$. The dataset points are indexed by multi-indexes $(i_1, \ldots, i_L) \in I_1 \times \ldots \times I_L$. Each observation in this case are the coordinates of $\mathcal{T}$, that is, we have $y^{(i_1, \ldots, i_L)} = t_{i_1 \ldots i_L}$. The independent variables are just the multi-indexes\footnote{Our ``population'' is composed by the indexes of the tensor, while the variables are the actual values of the tensor at such coordinates.}, to indicate what coordinate of $\mathcal{T}$ we are considering. The model $h_\textbf{w}: I_1 \times \ldots \times I_L \to \mathbb{R}$ tries to approximate each coordinate of $\mathcal{T}$ with the parameters in $\textbf{w} = [vec(\textbf{W}^{(1)})^T, \ldots, vec(\textbf{W}^{(L)})^T]^T$ through the CPD formulation given at the beginning of this chapter. More precisely, we have that 
		
		$$h_\textbf{w}(i_1, \ldots, i_L) = \tilde{t}_{i_1 \ldots i_L} = \sum_{r=1}^R w_{i_1 r}^{(1)} \ldots w_{i_L r}^{(L)}.$$   
		
		At this point it should be clear that 
		
		\begin{equation} \label{error_function}
			F(\textbf{w}) = \frac{1}{2} \sum_{i_1=1}^{I_1} \ldots \sum_{i_L=1}^{I_L} \left( t_{i_1 \ldots i_L} - \sum_{r=1}^R w_{i_1 r}^{(1)} \ldots w_{i_L r}^{(L)} \right)^2 = \frac{1}{2} \| \mathcal{T} - \tilde{\mathcal{T}} \|^2
		\end{equation}
and 
		
		\begin{equation} \label{residual_function}
			f_{i_1 \ldots i_L}(\textbf{w}) = t_{i_1 \ldots i_L} - \sum_{r=1}^R w_{i_1 r}^{(1)} \ldots w_{i_L r}^{(L)}.
		\end{equation}
		
		The following result is immediate.
		
		\begin{lemma} \label{partialf}	
			For all $\ell' = 1 \ldots L$, $i' = 1 \ldots I_\ell$ and $r' = 1 \ldots R$, we have 		
			
			$$\frac{\partial f_{i_1 \ldots i_{\ell'} \ldots i_L}}{\partial w_{i' r'}^{(\ell')}}(\textbf{w}) = 
			\left\{ \begin{array}{c}
				\displaystyle - \prod_{\ell \neq \ell'} w_{i_\ell r'}^{(\ell)}, \quad \text{ if } i' = i_{\ell'}\\
				0, \quad \text{ otherwise }
			\end{array} \right..$$		
		\end{lemma}	
			
		From this lemma it follows that
		
		$$\frac{\partial F}{\partial w_{i' r'}^{(\ell')}}(\textbf{w}) = \sum_{i_1=1}^{I_1} \ldots \sum_{i_L=1}^{I_L} f_{i_1 \ldots i_L}(\textbf{w}) \cdot \frac{\partial f_{i_1 \ldots i_L}}{\partial w_{i' r'}^{(\ell')}}(\textbf{w}) = $$
		
		$$ = \sum_{i_1=1}^{I_1} \ldots \sum_{i_{\ell'-1}=1}^{I_{\ell'-1}} \ \sum_{i_{\ell'+1}=1}^{I_{\ell'+1}} \ldots \sum_{i_L=1}^{I_L} f_{i_1 \ldots i' \ldots i_L}(\textbf{w}) \cdot \frac{\partial f_{i_1 \ldots i' \ldots i_L}}{\partial w_{i' r'}^{(\ell')}}(\textbf{w}) = $$
		
		$$ =  \sum_{i_1=1}^{I_1} \ldots \sum_{i_{\ell'-1}=1}^{I_{\ell'-1}} \ \sum_{i_{\ell'+1}=1}^{I_{\ell'+1}} \ldots \sum_{i_L=1}^{I_L} \left( t_{i_1 \ldots i_L} - \sum_{r=1}^R w_{i' r}^{(\ell')} \prod_{\ell \neq \ell'} w_{i_\ell r}^{(\ell)} \right) \left( - \prod_{\ell \neq \ell'} w_{i_\ell r'}^{(\ell)} \right).$$
		
		Observe how these manipulations leads to a faster and clean proof of theorem~\ref{partialF}.
		
	\section{Gauss-Newton} \label{gn_section}
		There are several ways to work a nonlinear least squares problem, and our chosen method is based on the \emph{Gauss-Newton} algorithm. This algorithm can only be used to minimize a sum of squared function values, but it has the advantage that second derivatives, which can be challenging to compute, are not required. Consider the same notations used in the previous section.
		
		The first step of the Gauss-Newton algorithm is to consider a first order approximation of $f$ at a point $\textbf{w}^{(0)}$, that is,
		
		\begin{equation} \label{linear_approx}
		f(\textbf{w}^{(0)} + \underbrace{(\textbf{w} - \textbf{w}^{(0)})}_{step} ) = f(\textbf{w}) \approx f(\textbf{w}^{(0)}) + \textbf{J}_f(\textbf{w}^{(0)}) \cdot (\textbf{w} - \textbf{w}^{(0)}). 
		\end{equation}
	
		In order to minimize~\ref{linear_approx} at the neighborhood of $\textbf{w}^{(0)}$ we can compute the minimum of $\| f(\textbf{w}^{(0)}) + \textbf{J}_f(\textbf{w}^{(0)}) \cdot (\textbf{w} - \textbf{w}^{(0)}) \|$ for $\textbf{w}$. Note that minimizing~\ref{linear_approx} is a least squares problem since we can rewrite this problem as $\displaystyle \min_\textbf{x} \| \textbf{A}\textbf{x} - \textbf{b}\|$ for $\textbf{A} = \textbf{J}_f(\textbf{w}^{(0)})$, $\textbf{x} = \textbf{w} - \textbf{w}^{(0)}$, $\textbf{b} = - f(\textbf{w}^{(0)})$.
	
		The solution of~\ref{linear_approx} gives the next iterate $\textbf{w}^{(1)}$. More generally, we obtain $\textbf{w}^{(k+1)}$ from $\textbf{w}^{(k)}$ by defining $\textbf{w}^{(k+1)} = \textbf{x}_\ast + \textbf{w}^{(k)}$, where $\textbf{x}_\ast$ is the solution of the normal equations
		
		\begin{equation} \label{normal_eq}
			\textbf{A}^T \textbf{A} \textbf{x} = \textbf{A}^T \textbf{b} 
		\end{equation}
for $\textbf{A} = \textbf{J}_f(\textbf{w}^{(k)})$, $\textbf{x} = \textbf{w} - \textbf{w}^{(k)}$, $\textbf{b} = - f(\textbf{w}^{(k)})$. The explicit formula for $\textbf{w}^{(k+1)}$ is

		$$\textbf{w}^{(k+1)} = \textbf{w}^{(k)} - \left( \textbf{J}_f\left( \textbf{w}^{(k)} \right)^T \textbf{J}_f\left( \textbf{w}^{(k)} \right) \right)^{-1} \cdot \textbf{J}_f\left( \textbf{w}^{(k)} \right)^T \cdot f(\textbf{w}^{(k)})$$ 
although we should clarify that the inverse above never will be explicitly computed. Usually one uses an iterative algorithm to solve~\ref{normal_eq}. This process of obtaining successive $\textbf{w}^{(k)}$ is the Gauss-Newton algorithm, and it is guaranteed to converge to a local minima of $F$. For more details about this algorithm and its properties I recommend \cite{madsen}. We show the most relevant results here.

		\begin{theorem}[K. Madsen, H. B. Nielsen, and O. Tingleff, 2004] \label{descent}
			$\textbf{w}^{(k+1)} - \textbf{w}^{(k)}$ is a descent direction for $F$ at $\textbf{w}^{(k)}$.
		\end{theorem}
		
		\textbf{Proof:} From the formula
		
		$$\textbf{w}^{(k+1)} = \textbf{w}^{(k)} - \left( \textbf{J}_f( \textbf{w}^{(k)} )^T \textbf{J}_f( \textbf{w}^{(k)} ) \right)^{-1} \cdot \textbf{J}_f( \textbf{w}^{(k)} )^T \cdot f(\textbf{w}^{(k)})$$ 
we can conclude that

		$$ - \textbf{J}_f( \textbf{w}^{(k)} )^T \textbf{J}_f( \textbf{w}^{(k)} ) \cdot \left( \textbf{w}^{(k+1)} - \textbf{w}^{(k)} \right) = \textbf{J}_f( \textbf{w}^{(k)} )^T \cdot f(\textbf{w}^{(k)}).$$
	
		With the above observation and lemma~\ref{gradF} we have that
		
		$$ \langle \nabla F(\textbf{w}^{(k)}), \textbf{w}^{(k+1)} - \textbf{w}^{(k)} \rangle = \langle \textbf{J}_f(\textbf{w}^{(k)})^T f(\textbf{w}^{(k)}), \textbf{w}^{(k+1)} - \textbf{w}^{(k)} \rangle = $$
		
		$$ = \langle - \textbf{J}_f(\textbf{w}^{(k)})^T \textbf{J}_f(\textbf{w}^{(k)}) \left( \textbf{w}^{(k+1)} - \textbf{w}^{(k)} \right), \textbf{w}^{(k+1)} - \textbf{w}^{(k)} \rangle = $$
		
		$$ = - \langle \textbf{J}_f(\textbf{w}^{(k)}) \left( \textbf{w}^{(k+1)} - \textbf{w}^{(k)} \right), \textbf{J}_f(\textbf{w}^{(k)}) \left( \textbf{w}^{(k+1)} - \textbf{w}^{(k)} \right) \rangle = $$
		
		$$ = - \| \textbf{J}_f(\textbf{w}^{(k)}) \left( \textbf{w}^{(k+1)} - \textbf{w}^{(k)} \right) \|^2 \leq 0.$$
	
		Since the derivative of $F$ at $\textbf{w}^{(k)}$ in the direction $\textbf{w}^{(k+1)} - \textbf{w}^{(k)}$ is negative, it follows that this is a descent direction.$\hspace{11cm}\square$\bigskip 
		
		\subsection{Deriving the Gauss-Newton method from the Newton's method}
			Denote the Hessian matrix of $F$ at $\textbf{w}$ by $\textbf{H}_F(\textbf{w})$. We can relate the derivatives of $f$ and $F$ through the following result.
	
			\begin{lemma} 
				We have that			
				$$\textbf{H}_F(\textbf{w}) = \textbf{J}_f^T(\textbf{w}) \textbf{J}_f(\textbf{w}) + \sum_{i_1=1}^{I_1} \ldots \sum_{i_L=1}^{I_L} f_{i_1 \ldots i_L}(\textbf{w}) \cdot \textbf{H}_{f_{i_1 \ldots i_L}}(\textbf{w}),$$ 
where $\textbf{H}_{f_{i_1 \ldots i_L}}$ is the Hessian matrix of $f_{i_1 \ldots i_L}$.
			\end{lemma}
	
	As the algorithm converges we expect to have $f_{i_1 \ldots i_L} \approx 0$ for all $i_1, \ldots, i_L$. Therefore this lemma implies that $\textbf{H}_F \approx \textbf{J}_f^T \textbf{J}_f$ when close to an optimal point.		
		
			We can apply Newton's method to solve the system $\nabla F = 0$. Then the iteration formula becomes
			
			$$\textbf{w}^{(k+1)} = \textbf{w}^{(k)} - \left( \textbf{H}_F(\textbf{w}^{(k)}) \right)^{-1} \cdot \nabla F(\textbf{w}^{(k)}) \approx $$
			
			$$\approx \textbf{w}^{(k)} - \left( \textbf{J}_f^T(\textbf{w}^{(k)}) \textbf{J}_f(\textbf{w}^{(k)}) \right)^{-1} \cdot \nabla F(\textbf{w}^{(k)}) = $$
			
			$$ = \textbf{w}^{(k)} - \left( \textbf{J}_f^T(\textbf{w}^{(k)}) \textbf{J}_f(\textbf{w}^{(k)}) \right)^{-1} \cdot \textbf{J}_f^T(\textbf{w}^{(k)}) f(\textbf{w}^{(k)}).$$
			
			Here we recover the Gauss-Newton iteration formula for $\textbf{w}$. This approximation makes sense if the Hessian is nonsingular, otherwise, the approximation will be poor, which will cause the Gauss-Newton algorithm to make very small steps, slowing down convergence. We can also have problems if the function $f$ is highly nonlinear, in which case the approximation of the Hessian will again be poor. 			
		
		\subsection{Damped Gauss-Newton} \label{damped-gauss-newton}	
			Everything of this section until this point was very general, without specifying tensors. For all this section we consider $f$ and $F$ in the tensor context. In this context, we don't consider $f$ highly nonlinear but, unfortunately, the approximate Hessian converges to a singular matrix as the Gauss-Newton iterations converges. We make this statement more precise soon. For this part we will denote $\displaystyle\frac{\partial f_{i_1 \ldots i_L}}{\partial w_{i' r'}^{(\ell')}} = \frac{\partial f_{i_1 \ldots i_L}}{\partial w_{i' r'}^{(\ell')}}(\textbf{w})$, that is, we suppress the point $\textbf{w}$ where the derivative is being evaluated. This will make notation simpler and won't cause any confusion. Also let 
			
			$$\frac{\partial f_{i_1 \ldots i_L}}{\partial \textbf{w}_{r'}^{(\ell')}} = \left[ \frac{\partial f_{i_1 \ldots i_L}}{\partial w_{1 r'}^{(\ell')}}, \ldots, \frac{\partial f_{i_1 \ldots i_L}}{\partial w_{I_{\ell'} r'}^{(\ell')}} \right] \in \mathbb{R}^{I_{\ell'}}$$  
and 
			
			$$\frac{\partial f}{\partial \textbf{w}_{r'}^{(\ell')}} = 
			\left[ \begin{array}{c}
				\displaystyle\frac{\partial f_{1 \ldots 1}}{\partial \textbf{w}_{r'}^{(\ell')}}\\\\
				\displaystyle\frac{\partial f_{1 \ldots 2}}{\partial \textbf{w}_{r'}^{(\ell')}}\\\\
				\vdots\\\\
				\displaystyle\frac{\partial f_{I_1 \ldots I_L-1}}{\partial \textbf{w}_{r'}^{(\ell')}}\\\\
				\displaystyle\frac{\partial f_{I_1 \ldots I_L}}{\partial \textbf{w}_{r'}^{(\ell')}}
			\end{array} \right] \in \mathbb{R}^{\prod_{\ell=1}^L I_\ell \times I_{\ell'}}.$$
			
			Notice that the ordering of the rows follows the indexes $i_1 \ldots i_L$ as numbers in increasing order. We can consider just one more level of compact notation and define 
			
			$$\frac{\partial f}{\partial \textbf{W}^{(\ell')}} = \left[ \frac{\partial f}{\partial \textbf{w}_1^{(\ell')}}, \ldots, \frac{\partial f}{\partial \textbf{w}_R^{(\ell')}} \right] \in \mathbb{R}^{ \prod_{\ell=1}^L I_\ell \times I_{\ell'} R},$$
so we have that

			$$\hspace{2.1cm} \displaystyle J_f = \left[ \frac{\partial f}{\partial \textbf{W}^{(1)}}, \ldots, \frac{\partial f}{\partial \textbf{W}^{(L)}} \right] \in \mathbb{R}^{ \prod_{\ell=1}^L I_\ell \times R \sum_{\ell=1}^L I_\ell}.$$
			
			In order to understand the structure of $\textbf{J}_f$ first we have to understand the structure of each block $\displaystyle\frac{\partial f}{\partial \textbf{w}_{r'}^{(\ell')}}$. Remember that the entries of this matrix were computed in lemma~\ref{partialf}, now we just need to find some structure on it. Given a multi-index $i_1 \ldots i_{\ell'} \ldots i_L$, the associated row has only one non zero entry, which is the column $i_{\ell'}$, because
			
			$$\frac{\partial f_{i_1 \ldots i_{\ell'} \ldots i_L}}{\partial \textbf{w}_{r'}^{(\ell')}} = \left[ \frac{\partial f_{i_1 \ldots i_{\ell'} \ldots i_L}}{\partial w_{1 r'}^{(\ell')}}, \ldots, \frac{\partial f_{i_1 \ldots i_{\ell'} \ldots i_L}}{\partial w_{i_{\ell'} r'}^{(\ell')}}, \ldots,  \frac{\partial f_{i_1 \ldots i_{\ell'} \ldots i_L}}{\partial w_{I_{\ell'} r'}^{(\ell')}} \right] = \big[ 0, \ldots, \underbrace{- \prod_{\ell \neq \ell'} w_{i_\ell r'}^{(\ell)}}_{i_{\ell'} \text{ column}}, \ldots, 0 \big].$$ 
			
\begin{figure}[H]
				\centering
				\includegraphics[scale=.45]{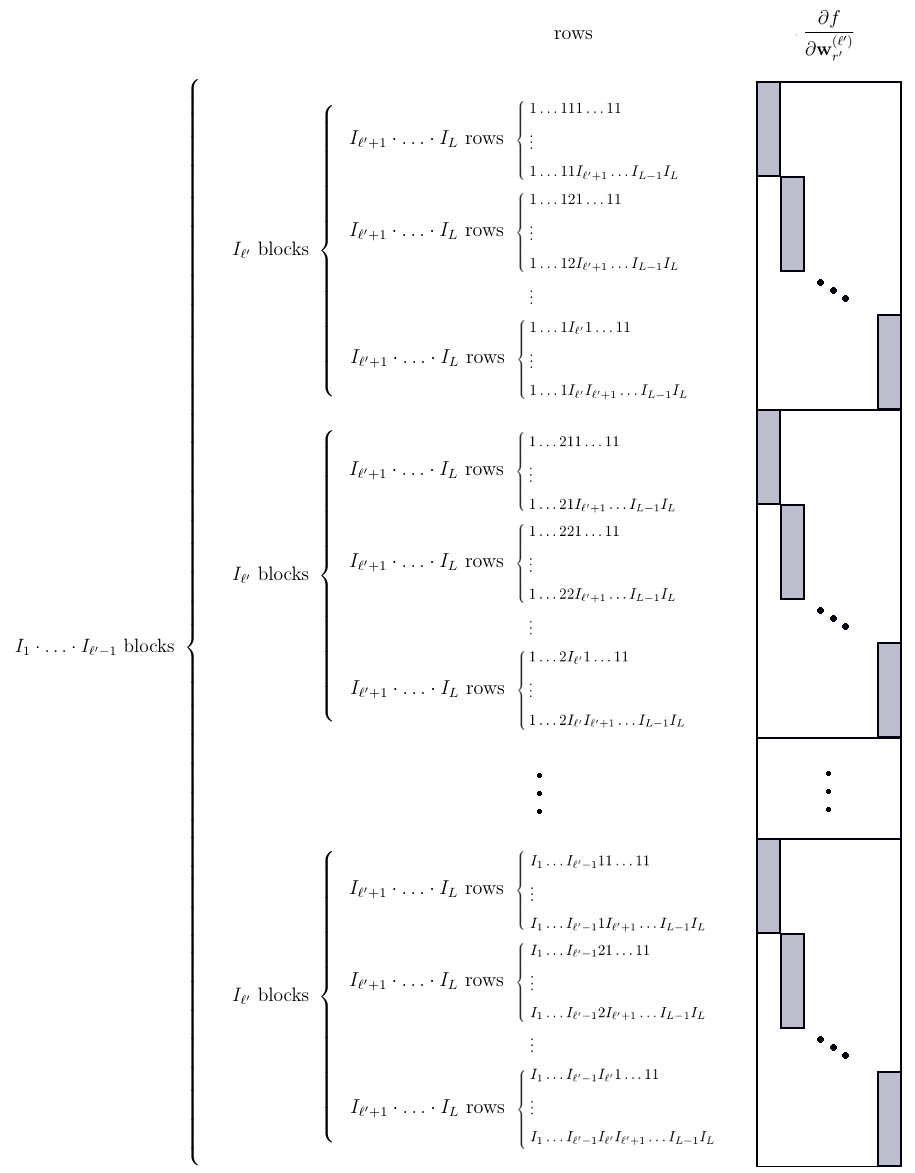}
				\caption{\footnotesize{Sparse structure of $\frac{\partial f}{\partial \textbf{w}_r^{(\ell')}}$. The gray part correspond to the non zero entries and the rest are full of zeros, and each gray column is a vector of size $\displaystyle \prod_{\ell=\ell'+1}^L I_\ell$.}}
				\label{sparse}
			\end{figure}

			Now we can see that $\displaystyle\frac{\partial f}{\partial \textbf{w}_{r'}^{(\ell')}}$ will have a sparse structure, and this structure is dictated by the multi-indexes: at row $i_1 \ldots i_{\ell'} \ldots i_L$, the only non zero entry is given by the value $i_{\ell'}$. Since this structure does not depend on $r'$, the block matrix $\displaystyle \frac{\partial f}{\partial \textbf{W}^{(\ell')}}$ consists in a repetition of $R$ blocks with the same sparse structure (but not necessarily the same values). Because of the ordering of $i_1 \ldots i_{\ell'} \ldots i_L$, each change of non zero column is done periodically as the scheme in figure~\ref{sparse} illustrates.			
			
\begin{figure}[H] 
				\centering
				\includegraphics[scale=.42]{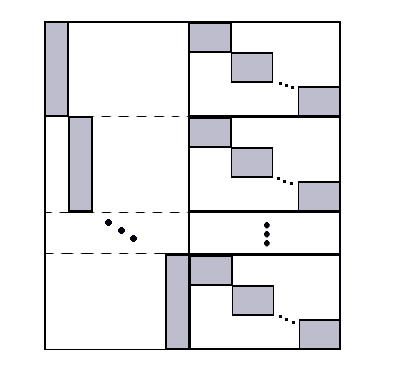}
				\caption{\footnotesize{Consider the top block of the previous figure, relative to $\frac{\partial f}{\partial w_{r'}^{(\ell')}}$. While run through the rows of this block, at the same time there will be $\displaystyle \prod_{\ell=1}^{\ell'} I_\ell$ blocks relative to $\frac{\partial f}{\partial w_{r'}^{(\ell'+1)}}$.}}
				\label{columns}
			\end{figure}

			Note that when we go to row $1 \ldots 1 1 1 \ldots 1 1$ to row $1 \ldots 1 2 1 \ldots 1 1$, the gray block moves one column to the right. When this happens, the next index, $i_{\ell'+1}$, already changed all its values exactly one time. It means all this moves to the right occurred with respect to $i_{\ell'+1}$ when it happened just one with respect to $i_{\ell'}$. This is due to the ordering we used, which functions just as natural numbers in increasing order. This explanation gets clear with a picture, so figure~\ref{columns} illustrates better what is going on.  			
			
\begin{figure} 
			\centering
			\includegraphics[scale=.5]{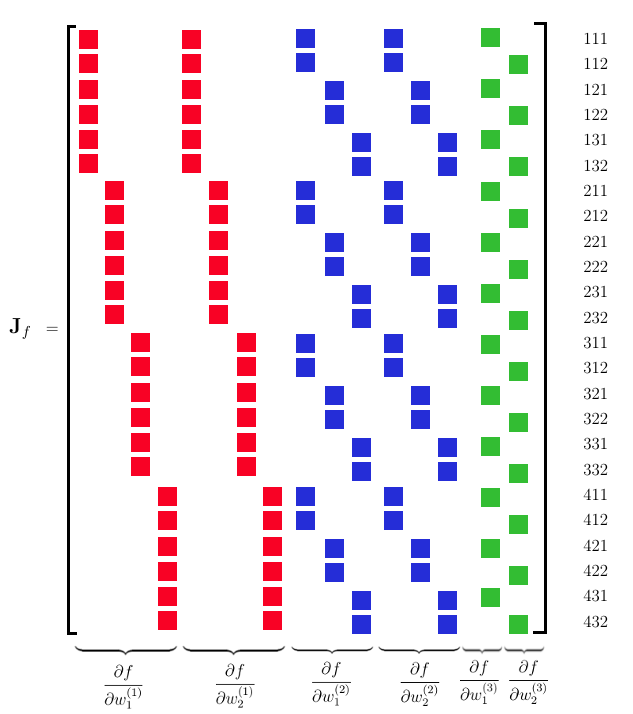}
			\caption{\footnotesize{Sparse structure of $\textbf{J}_f$ when $\mathcal{T} \in \mathbb{R}^{4 \times 3 \times 2}$ and $R = 2$.}}
			\label{structure-concrete}
		\end{figure}

			We finish this discussion of the sparse structure displaying the structure of a concrete tensor shape. If still there are any doubts at this point, we invite the reader to reproduce the structure of $\textbf{J}_f$ when $\mathcal{T} \in \mathbb{R}^{4 \times 3 \times 2}$ and $R = 2$. This structure is given in figure~\ref{structure-concrete}. We separate the modes by color to facilitate understanding. At the right to the matrix we also added the respective multi-index associated to each row. The green color correspond to the right index, the blue correspond to the middle index and the red correspond to the left index.  			
			
			\begin{lemma}[T. G. Kolda, E. Acar, D. M. Dunlavy, 2011] \label{kronecker}
				For all $\ell' = 1 \ldots L$ and $r' = 1 \ldots R$, we have				
				$$\frac{\partial f}{\partial \textbf{w}_{r'}^{(\ell')}} = - \textbf{w}_{r'}^{(1)} \tilde{\otimes} \ldots \tilde{\otimes} \textbf{w}_{r'}^{(\ell'-1)} \tilde{\otimes} \textbf{I}_{I_{\ell'}} \tilde{\otimes} \textbf{w}_{r'}^{(\ell'+1)} \tilde{\otimes} \ldots \tilde{\otimes} \textbf{w}_{r'}^{(L)}.$$
			\end{lemma} 
			
			\textbf{Proof:} First note that			
			$$\textbf{I}_{I_{\ell'}} \tilde{\otimes}\textbf{w}_{r'}^{(\ell'+1)} \tilde{\otimes} \ldots \tilde{\otimes} \textbf{w}_{r'}^{(L)} = 
			\left[ \begin{array}{ccc}
				\textbf{w}_{r'}^{(\ell'+1)} \tilde{\otimes} \ldots \tilde{\otimes} \textbf{w}_{r'}^{(L)} & & \\
				& \ddots & \\
				& & \textbf{w}_{r'}^{(\ell'+1)} \tilde{\otimes} \ldots \tilde{\otimes} \textbf{w}_{r'}^{(L)} 
			\end{array} \right]$$
is a block diagonal matrix, with $I_{\ell'} \times I_{\ell'}$ block shape where each entry is a vector of size $\displaystyle \prod_{\ell=\ell'+1}^L I_\ell$. With this we have
			$$\textbf{w}_{r'}^{(1)} \tilde{\otimes} \ldots \tilde{\otimes} \textbf{w}_{r'}^{(\ell'-1)} \tilde{\otimes} 
			\left[ \begin{array}{ccc}
				\textbf{w}_{r'}^{(\ell'+1)} \tilde{\otimes} \ldots \tilde{\otimes} \textbf{w}_{r'}^{(L)} & & \\
				& \ddots & \\
				& & \textbf{w}_{r'}^{(\ell'+1)} \tilde{\otimes} \ldots \tilde{\otimes} \textbf{w}_{r'}^{(L)} 
			\end{array} \right] = $$
			
			$$ = 
			\left[ \begin{array}{c}
				w_{1 r'}^{(1)} \cdot \ldots \cdot w_{1 r'}^{(\ell'-1)}
				\left[ \begin{array}{ccc}
					\textbf{w}_{r'}^{(\ell'+1)} \tilde{\otimes} \ldots \tilde{\otimes} \textbf{w}_{r'}^{(L)} & & \\
					& \ddots & \\
					& & \textbf{w}_{r'}^{(\ell'+1)} \tilde{\otimes} \ldots \tilde{\otimes} \textbf{w}_{r'}^{(L)} 
				\end{array} \right]\\\\
				\vdots\\\\
				w_{I_1 r'}^{(1)} \cdot \ldots \cdot w_{I_{\ell'-1} \ r'}^{(\ell'-1)}
				\left[ \begin{array}{ccc}
					\textbf{w}_{r'}^{(\ell'+1)} \tilde{\otimes} \ldots \tilde{\otimes} \textbf{w}_{r'}^{(L)} & & \\
					& \ddots & \\
					& & \textbf{w}_{r'}^{(\ell'+1)} \tilde{\otimes} \ldots \tilde{\otimes} \textbf{w}_{r'}^{(L)} 
				\end{array} \right]
			\end{array} \right].$$
			
			By lemma~\ref{kronecker-cpd} we have that			
			$$\textbf{w}_{r'}^{(\ell'+1)} \tilde{\otimes} \ldots \tilde{\otimes} \textbf{w}_{r'}^{(L)} = 
			\left[ \begin{array}{c}
				w_{1 r'}^{(\ell'+1)} \cdot \ldots \cdot w_{1 r'}^{(L)}\\
				\vdots\\
				w_{I_{\ell'+1} \ r'}^{(\ell'+1)} \cdot \ldots \cdot w_{I_L r'}^{(L)}
			\end{array} \right],$$
hence we can write more explicitly the last expression as \footnotesize

			$$ 
			\left[ \begin{array}{c}
				w_{1 r'}^{(1)} \cdot \ldots \cdot w_{1 r'}^{(\ell'-1)}
				\left[ \begin{array}{ccc}
					\left[ \begin{array}{c}
						w_{1 r'}^{(\ell'+1)} \cdot \ldots \cdot w_{1 r'}^{(L)}\\
						\vdots\\
						w_{I_{\ell'+1} \ r'}^{(\ell'+1)} \cdot \ldots \cdot w_{I_L r'}^{(L)}
					\end{array} \right] & & \\
					& \ddots & \\
					& & \left[ \begin{array}{c}
						w_{1 r'}^{(\ell'+1)} \cdot \ldots \cdot w_{1 r'}^{(L)}\\
						\vdots\\
						w_{I_{\ell'+1} \ r'}^{(\ell'+1)} \cdot \ldots \cdot w_{I_L r'}^{(L)}
					\end{array} \right] 
				\end{array} \right]\\\\
				\vdots\\\\
				w_{I_1 r'}^{(1)} \cdot \ldots \cdot w_{I_{\ell'-1} \ r'}^{(\ell'-1)}
				\left[ \begin{array}{ccc}
					\left[ \begin{array}{c}
						w_{1 r'}^{(\ell'+1)} \cdot \ldots \cdot w_{1 r'}^{(L)}\\
						\vdots\\
						w_{I_{\ell'+1} \ r'}^{(\ell'+1)} \cdot \ldots \cdot w_{I_L r'}^{(L)}
					\end{array} \right] & & \\
					& \ddots & \\
					& & \left[ \begin{array}{c}
						w_{1 r'}^{(\ell'+1)} \cdot \ldots \cdot w_{1 r'}^{(L)}\\
						\vdots\\
						w_{I_{\ell'+1} \ r'}^{(\ell'+1)} \cdot \ldots \cdot w_{I_L r'}^{(L)}
					\end{array} \right] 
				\end{array} \right]
			\end{array} \right] = $$
			
			$$\hspace{-1cm} = 
			\left[ \begin{array}{c}
				\left[ \begin{array}{ccc}
					\left[ \begin{array}{c}
						w_{1 r'}^{(1)} \cdot \ldots \cdot w_{1 r'}^{(\ell'-1)} \cdot w_{1 r'}^{(\ell'+1)} \cdot \ldots \cdot w_{1 r'}^{(L)}\\
						\vdots\\
						w_{1 r'}^{(1)} \cdot \ldots \cdot w_{1 r'}^{(\ell'-1)} \cdot w_{I_{\ell'+1} \ r'}^{(\ell'+1)} \cdot \ldots \cdot w_{I_L r'}^{(L)}
					\end{array} \right] & & \\
					& \ddots & \\
					& & \left[ \begin{array}{c}
						w_{1 r'}^{(1)} \cdot \ldots \cdot w_{1 r'}^{(\ell'-1)} \cdot w_{1 r'}^{(\ell'+1)} \cdot \ldots \cdot w_{1 r'}^{(L)}\\
						\vdots\\
						w_{1 r'}^{(1)} \cdot \ldots \cdot w_{1 r'}^{(\ell'-1)} \cdot w_{I_{\ell'+1} \ r'}^{(\ell'+1)} \cdot \ldots \cdot w_{I_L r'}^{(L)}
					\end{array} \right] 
				\end{array} \right]\\\\
				\vdots\\\\
				\left[ \begin{array}{ccc}
					\left[ \begin{array}{c}
						w_{I_1 r'}^{(1)} \cdot \ldots \cdot w_{I_{\ell'-1} \ r'}^{(\ell'-1)} \cdot w_{1 r'}^{(\ell'+1)} \cdot \ldots \cdot w_{1 r'}^{(L)}\\
						\vdots\\
						w_{I_1 r'}^{(1)} \cdot \ldots \cdot w_{I_{\ell'-1} \ r'}^{(\ell'-1)} \cdot w_{I_{\ell'+1} \ r'}^{(\ell'+1)} \cdot \ldots \cdot w_{I_L r'}^{(L)}
					\end{array} \right] & & \\
					& \ddots & \\
					& & \left[ \begin{array}{c}
						w_{I_1 r'}^{(1)} \cdot \ldots \cdot w_{I_{\ell'-1} \ r'}^{(\ell'-1)} \cdot w_{1 r'}^{(\ell'+1)} \cdot \ldots \cdot w_{1 r'}^{(L)}\\
						\vdots\\
						w_{I_1 r'}^{(1)} \cdot \ldots \cdot w_{I_{\ell'-1} \ r'}^{(\ell'-1)} \cdot w_{I_{\ell'+1} \ r'}^{(\ell'+1)} \cdot \ldots \cdot w_{I_L r'}^{(L)}
					\end{array} \right] 
				\end{array} \right]
			\end{array} \right].$$\bigskip
			
			\normalsize
			This reveals the block structure of $\textbf{J}_f$ already observed and illustrated in figures~\ref{sparse}, ~\ref{columns}, ~\ref{structure-concrete}. Now, given any multi-index $i_1 \ldots i_{\ell'} \ldots i_L$ we just need to show that $\displaystyle\frac{\partial f_{i_1 \ldots i_{\ell'} \ldots i_L}}{\partial \textbf{w}_{r'}^{(\ell')}}$ is equal to the negative of row $i_1 \ldots i_{\ell'} \ldots i_L$ of $- \textbf{w}_{r'}^{(1)} \tilde{\otimes} \ldots \tilde{\otimes} \textbf{w}_{r'}^{(\ell'-1)} \tilde{\otimes} \textbf{I}_{I_{\ell'}} \tilde{\otimes} \textbf{w}_{r'}^{(\ell'+1)} \tilde{\otimes} \ldots \tilde{\otimes} \textbf{w}_{r'}^{(L)}$. Note that this row is the $i_{\ell'+1} \ldots i_L$ row of the $\ell'$-th diagonal term of 
			
			$$- w_{i_1 r'}^{(1)} \cdot \ldots \cdot w_{i_{\ell'-1} \ r'}^{(\ell'-1)}
				\left[ \begin{array}{ccc}
					\textbf{w}_{r'}^{(\ell'+1)} \tilde{\otimes} \ldots \tilde{\otimes} \textbf{w}_{r'}^{(L)} & & \\
					& \ddots & \\
					& & \textbf{w}_{r'}^{(\ell'+1)} \tilde{\otimes} \ldots \tilde{\otimes} \textbf{w}_{r'}^{(L)} 
			\end{array} \right].$$
This row is given by

			$$ - w_{i_1 r'}^{(1)} \cdot \ldots \cdot w_{i_{\ell'-1} \ r'}^{(\ell'-1)} \cdot \big[ 0, \ldots, 0, \underbrace{w_{i_{\ell'+1} r'}^{(\ell'+1)} \cdot \ldots \cdot w_{i_L \ r'}^{(L)}}_{\ell'-\text{ column }}, 0 \ldots, 0 \big] = $$
			
			$$\hspace{3cm} = \big[ 0, \ldots, 0, - \underbrace{\prod_{\ell \neq \ell'} w_{i_\ell r'}^{(\ell)}}_{\ell'-\text{ column }}, 0 \ldots, 0 \big] = \frac{\partial f_{i_1 \ldots i_{\ell'} \ldots i_L}}{\partial \textbf{w}_{r'}^{(\ell')}}. \hspace{4.4cm} \square$$
			
			\begin{theorem}
				$\textbf{J}_f^T \textbf{J}_f$ is singular.
			\end{theorem}
			
			\textbf{Proof:} Notice that suffices to prove that $\textbf{J}_f$ is not a full rank matrix. Remember that $f$ maps possible rank-$R$ CPDs $(\textbf{W}^{(1)}, \ldots, \textbf{W}^{(L)})$ to its individual residuals, a vector with size $\displaystyle \prod_{\ell=1}^L I_\ell$. Because of the scale indeterminacy, we know that the inputs $(\lambda \textbf{W}^{(1)}, \ldots, \textbf{W}^{(L)})$ and $(\textbf{W}^{(1)}, \ldots, \lambda \textbf{W}^{(\ell)}, \ldots, \textbf{W}^{(L)})$ correspond to the same tensor for any $\ell \neq 1$. This is because the former correspond to the tensor

			$$(\lambda \textbf{w}_1^{(1)}) \otimes \ldots \otimes \textbf{w}_1^{(L)} + \ldots + (\lambda \textbf{w}_R^{(1)}) \otimes \ldots \otimes \textbf{w}_R^{(L)},$$
while the latter correspond to  	

			$$\textbf{w}_1^{(1)} \otimes \ldots \otimes (\lambda \textbf{w}_1^{(\ell)})\otimes \ldots \otimes \textbf{w}_1^{(L)} + \ldots + \textbf{w}_R^{(1)} \otimes \ldots \otimes (\lambda \textbf{w}_R^{(\ell)})\otimes \ldots \otimes \textbf{w}_R^{(L)}.$$	
	
			There are several ways of using the phenomenon of scale indeterminacy to generate different inputs to $f$ which correspond to the same tensor. In particular, if some tensor is a local minima for $F$, it is also a critical point, and because of the scale indeterminacy there will be infinitely many inputs associated to the same critical point. This means the critical point is singular. Since $\nabla F = \textbf{J}_f^T \cdot f$, we conclude that $\textbf{J}_f$ is not of full rank at critical points.$\hspace{13.1cm}\square$\bigskip

			Since $\textbf{J}_f^T \textbf{J}_f$ is singular, the Gauss-Newton update rule given by~\ref{normal_eq} is an ill-posed problem. We can address this issue by introducing a regularization term, so we can avoid singularity and improve convergence. A common approach is to introduce a suitable regularization matrix $\textbf{L} \in \mathbb{R}^{R \sum_{\ell=1}^L I_\ell \times R \sum_{\ell=1}^L I_\ell}$ called \emph{Tikhonov matrix}. Instead of solving equation~\ref{normal_eq} we are now solving
		\begin{equation}
			\left( \textbf{A}^T \textbf{A} + \textbf{L}^T \textbf{L} \right) \textbf{x} = \textbf{A}^T \textbf{b}, \label{dGN_eq}
		\end{equation}
where $\textbf{A} = \textbf{J}_f(\textbf{w}^{(k)}), \textbf{x} = \textbf{w} - \textbf{w}^{(k)}, b = -f(\textbf{w}^{(k)})$. The \emph{damped Gauss-Newton} (dGN)\footnote{This algorithm is also called \emph{Levenberg-Marquardt} in the literature. It is possible that this second name is used more often but still we prefer the other. The name comes from the fact that K. Levenberg \cite{levenberg} and D. Marquardt \cite{marquardt} are the ones responsible for introducing the damping parameter in the Gauss-Newton method. } algorithm is the case when $\textbf{L}^T \textbf{L} = \mu \textbf{I}_{R \sum_{\ell=1}^L I_\ell}$, where $\mu > 0$ is the \emph{damping parameter}.  Usually this parameter is updated at each iteration. These updates are very important since $\mu$ influences both the direction and the size of the step at each iteration. Instead of working with the regularization matrix $\mu \textbf{I}_{R \sum_{\ell=1}^L I_\ell}$ as is usual, we consider the more general matrix $\mu \textbf{D}$, where $\textbf{D}$ is a positive diagonal $\left( R \displaystyle \sum_{\ell=1}^L I_\ell \right) \times \left( R \displaystyle \sum_{\ell=1}^L I_\ell \right)$ matrix.

		\begin{theorem} The following holds.
			\begin{enumerate}
				\item $\textbf{J}_f^T \textbf{J}_f + \mu \textbf{D}$ is a positive definite matrix for all $\mu > 0$.
				\item $\textbf{w}^{(k+1)} - \textbf{w}^{(k)}$ is a descent direction for $F$ at $\textbf{w}^{(k)}$.
				\item If $\mu$ is large enough, then $\displaystyle \textbf{w}^{(k+1)} - \textbf{w}^{(k)} \approx - \frac{1}{\mu} \textbf{D}^{-1} \nabla F(\textbf{w}^{(k)})$.
				\item If $\mu$ is small enough, then $\textbf{w}^{(k+1)} - \textbf{w}^{(k)} \approx \textbf{w}_{GN}^{(k+1)} - \textbf{w}^{(k)}$, where $\textbf{w}_{GN}^{(k+1)}$ is the point we would obtain using classic Gauss-Newton iteration (i.e.,  without regularization). 
			\end{enumerate}
		\end{theorem} 
		
		\textbf{Proof:}
		To prove 1, just take any $\textbf{w} \in \mathbb{R}^{R \sum_{\ell=1}^L I_\ell}$ and note that
		
		$$\left\langle \left( \textbf{J}_f^T \textbf{J}_f + \mu \textbf{D} \right)\textbf{w}, \textbf{w} \right\rangle = $$
		$$ = \left\langle \textbf{J}_f^T \textbf{J}_f \textbf{w}, \textbf{w} \right\rangle + \left\langle \mu \textbf{D} \textbf{w}, \textbf{w} \right\rangle = $$
		$$ = \left\langle \textbf{J}_f \textbf{w}, \textbf{J}_f \textbf{w} \right\rangle + \left\langle \sqrt{\mu} \sqrt{\textbf{D}} \textbf{w}, \sqrt{\mu} \sqrt{\textbf{D}} \textbf{w} \right\rangle = $$
		$$ = \| \textbf{J}_f \textbf{w} \|^2 + \| \sqrt{\mu} \sqrt{\textbf{D}} \textbf{w} \|^2 > 0.$$
	
		The proof of 2 is very similar to the previous proof in the classic Gauss-Newton. From the iteration formula
		
		$$\textbf{w}^{(k+1)} = \textbf{w}^{(k)} - \left( \textbf{J}_f(\textbf{w}^{(k)})^T \textbf{J}_f(\textbf{w}^{(k)}) + \mu \textbf{D} \right)^{-1} \textbf{J}_f(\textbf{w}^{(k)})^T \cdot f(\textbf{w}^{(k)})$$
we can conclude that

		$$- \left( \textbf{J}_f(\textbf{w}^{(k)})^T \textbf{J}_f(\textbf{w}^{(k)}) + \mu \textbf{D} \right) \cdot \left( \textbf{w}^{(k+1)} - \textbf{w}^{(k)} \right) = \textbf{J}_f(\textbf{w}^{(k)})^T \cdot f(\textbf{w}^{(k)}).$$
		
		Now, with this identity, note that
		
		$$\left\langle \nabla F(\textbf{w}^{(k)}),  \textbf{w}^{(k+1)} - \textbf{w}^{(k)} \right\rangle = $$
		$$ = \left\langle \textbf{J}_f(\textbf{w}^{(k)})^T f(\textbf{w}^{(k)}), \textbf{w}^{(k+1)} - \textbf{w}^{(k)} \right\rangle = $$
		$$ = - \left\langle \left( \textbf{J}_f(\textbf{w}^{(k)})^T \textbf{J}_f(\textbf{w}^{(k)}) + \mu \textbf{D} \right) \cdot \left( \textbf{w}^{(k+1)} - \textbf{w}^{(k)} \right), \textbf{w}^{(k+1)} - \textbf{w}^{(k)} \right\rangle < 0. $$
		
		The inequality above follows from the fact that $\textbf{J}_f(\textbf{w}^{(k)})^T \textbf{J}_f(\textbf{w}^{(k)}) + \mu \textbf{D}$ is positive definite.
	
		To prove 3, take $\mu$ such that $\| \textbf{D}^{-1} \textbf{J}_f(\textbf{w}^{(k)})^T \textbf{J}_f(\textbf{w}^{(k)}) \| \ll \mu$ (this is the ``large enough'' in our context). We know $\textbf{J}_f(\textbf{w}^{(k)})^T \textbf{J}_f(\textbf{w}^{(k)}) + \mu \textbf{D}$ is invertible since it is positive definite. Also, by the definition of $\mu$ we have that
		
		$$\left( \textbf{J}_f(\textbf{w}^{(k)})^T \textbf{J}_f(\textbf{w}^{(k)}) + \mu \textbf{D} \right)^{-1} = \left( \mu \textbf{D} \left( \displaystyle \frac{1}{\mu} \textbf{D}^{-1} \textbf{J}_f(\textbf{w}^{(k)})^T \textbf{J}_f(\textbf{w}^{(k)}) + \textbf{I} \right) \right)^{-1} \approx $$
		$$ \approx \left( \mu \textbf{D} \left( \textbf{0} + \textbf{I} \right) \right)^{-1} = \frac{1}{\mu} \textbf{D}^{-1}.$$
		
		Using the iteration formula with this approximation gives 
		$$\textbf{w}^{(k+1)} \approx \textbf{w}^{(k)} - \frac{1}{\mu} \textbf{D}^{-1} \textbf{J}_f(\textbf{w}^{(k)})^T f(\textbf{w}^{(k)}) = \textbf{w}^{(k)} - \frac{1}{\mu} \textbf{D}^{-1} \nabla F(\textbf{w}^{-1}).$$	

		Finally, to prove 4 just consider $\mu \approx 0$ and substitute in the iteration formula. Then we get the classical formula trivially.$\hspace{10cm}\square$\bigskip		
	
		\begin{remark}
			If $\textbf{D} \approx \textbf{I}$, item 3 can be used when the current iteration is far from the solution, since $\displaystyle - \frac{1}{\mu} \nabla F(\textbf{w}^{(k)})$ is a short step in the descent direction (we want to be careful when distant to the solution). This shows that dGN behaves as the gradient descent algorithm when distant to the solution. On the other hand, item 4 is to be used when the current iteration is close to the solution, since the step is closer to the classic Gauss-Newton, we may attain quadratic convergence at the final iterations. 
			
			It is also important to note that the approximation in item 3 comes from the fact that 
			$$\displaystyle \left\| \frac{1}{\mu} \textbf{D}^{-1} \textbf{J}_f(\textbf{w}^{(k)})^T \textbf{J}_f(\textbf{w}^{(k)}) \right\| = \displaystyle \frac{1}{\mu} \left\| \textbf{D}^{-1} \textbf{J}_f(\textbf{w}^{(k)})^T \textbf{J}_f(\textbf{w}^{(k)}) \right\| \ll 1,$$
hence $\displaystyle \frac{1}{\mu} \textbf{D}^{-1} \textbf{J}_f(\textbf{w}^{(k)})^T \textbf{J}_f(\textbf{w}^{(k)}) \approx 0$. With this approximation we have the approximation in item 4 too.
		\end{remark} 
		
		The damping parameter rule update is a very important issue in this algorithm and we should comment a few things about it (the interested reader may check \cite{madsen} for more about this subject). Let $\mu^{(0)}$ be the initial damping parameter and $\mu^{(k)}$ be the damping parameter used at the $k$-th iteration. A rule of thumb for $\mu^{(0)}$ is to set $\mu^{(0)} = \tau \cdot \max_i a_{ii}$, where $\textbf{A} = \textbf{J}_f^T(\textbf{w}^{(0)}) \textbf{J}_f(\textbf{w}^{(0)})$ and $\tau \in (0, 1]$. We should choose $\tau$ smaller when $\textbf{w}^{(0)}$ is closer to the solution. The algorithm usually is not so sensitive to the choice of $\tau$ so there is no problem in using other initial values. More important is the update rule for $\mu^{(k)}$.
		
		After computing $\textbf{w}^{(k+1)}$, the \emph{gain ratio} is defined as
		
		$$g = \frac{ F(\textbf{w}^{(k)}) - F(\textbf{w}^{(k+1)}) }{ F(\textbf{w}^{(k)}) - \frac{1}{2}\| \tilde{f}(\textbf{w}^{(k+1)}) \|^2 } = \frac{ \| \mathcal{T} - \tilde{\mathcal{T}}^{(k)} \|^2 - \| \mathcal{T} - \tilde{\mathcal{T}}^{(k+1)} \|^2 }{ \| \mathcal{T} - \tilde{\mathcal{T}}^{(k)} \|^2 - \| \tilde{f}(\textbf{w}^{(k+1)}) \|^2 },$$
where $\tilde{\mathcal{T}}^{(k)}$ is the approximating tensor computed at the $k$-th iteration and $\tilde{f}(\textbf{w}^{(k+1)}) = f(\textbf{w}^{(k)}) + \textbf{J}_f(\textbf{w}^{(k)}) \cdot (\textbf{w}^{(k+1)} - \textbf{w}^{(k)})$ is the first order approximation of $f$ at $\textbf{w}^{(k+1)}$. The denominator is the predicted improvement\footnote{Remember that we obtain $\textbf{w}^{(k+1)}$ as a solution to the minimization problem $\displaystyle \min_\textbf{w} \| \tilde{f}(\textbf{w}) \|^2$ so we expect it to be close to $\displaystyle \min_\textbf{w} \| f(\textbf{w}) \|^2 = \displaystyle \min_\textbf{w} F(\textbf{w})$. Therefore $\| \tilde{f}(\textbf{w}^{(k+1)})\|^2$ gives the expected error while $F(\textbf{w}^{(k+1)})$ gives the actual error.} from the $k$-th to $(k+1)$-th iteration, whereas the numerator measures the actual improvement.

		A large value of $g$ indicates that $\| \tilde{f}(\textbf{w}^{(k+1)})\|^2$ is a good approximation to $\|\mathcal{T} - \tilde{\mathcal{T}}^{(k+1)}\|^2$, in other words, the linear model $\tilde{f}$ is making good predictions about the actual errors. A small or negative $g$ indicates that this approximation is poor. In the former case we decrease the damping parameter so the next iteration is more like a Gauss-Newton iteration, and in the latter case we increase the damping parameter to have more regularization and make smaller steps, this way the steps are more guaranteed to be in the steepest direction. This is the general idea, but in reality these updates depends a lot on the problem at hand and there is not a single procedure which works for all problems. We will talk more about it when introducing our own algorithm to compute the CPD.
		
		One thing to notice is that the denominator is always positive by construction since
		
		$$F(\textbf{w}^{(k)})- \| \tilde{f}(\textbf{w}^{(k+1)}) \|^2 = F(\textbf{w}^{(k)}) - \| f(\textbf{w}^{(k)}) + \textbf{J}_f(\textbf{w}^{(k)}) \cdot (\textbf{w}^{(k+1)} - \textbf{w}^{(k)}) \|^2 = $$		
		$$ = F(\textbf{w}^{(k)}) - \| f(\textbf{w}^{(k)}) \|^2 - 2 \langle f(\textbf{w}^{(k)}), \textbf{J}_f(\textbf{w}^{(k)})(\textbf{w}^{(k+1)} - \textbf{w}^{(k)}) \rangle - \| \textbf{J}_f(\textbf{w}^{(k)})(\textbf{w}^{(k+1)} - \textbf{w}^{(k)}) \|^2 = $$	
		$$ = -2 \langle f(\textbf{w}^{(k)}), \textbf{J}_f(\textbf{w}^{(k)})(\textbf{w}^{(k+1)} - \textbf{w}^{(k)}) \rangle - \| \textbf{J}_f(\textbf{w}^{(k)})(\textbf{w}^{(k+1)} - \textbf{w}^{(k)}) \|^2 =$$		
		$$ = -2 \langle \nabla F(\textbf{w}^{(k)}), \textbf{w}^{(k+1)} - \textbf{w}^{(k)} \rangle + \langle \nabla F(\textbf{w}^{(k)}), \textbf{w}^{(k+1)} - \textbf{w}^{(k)} \rangle = $$
		$$ = - \langle \nabla F(\textbf{w}^{(k)}), \textbf{w}^{(k+1)} - \textbf{w}^{(k)} \rangle = $$
		$$ = \| \textbf{J}_f(\textbf{w}^{(k)})(\textbf{w}^{(k+1)} - \textbf{w}^{(k)}) \|^2 > 0.$$
		
		Notice we used a few identities from theorem~\ref{descent} to make some of the manipulations above. The important thing here is to realize that $g$ will be negative only when the numerator is negative, and this means the error of the approximation increases at iteration $k+1$. In section~\ref{not_mon} we will see that it is expected that $g < 0$ in a few iterations, and we argue that this apparent drawback can actually be beneficial. 
		
		Different versions of the dGN algorithm were implemented and tested \cite{tensorlab, tensorlab2, tomasi, tomasi2, cp_opt, cichoki2013}. It is already known that dGN have faster local convergence when close to the optimal point since it uses information from the (approximated) Hessian. In the presence of bottlenecks or swamps the ALS presents problems, and even with improvements such as regularization, line search, rotation, etc, these problems doesn't disappear completely. A severe limitation of the ALS is the fact that it updates only one factor per iteration. The dGN uses all information to make the updates, and this leads to a more robust algorithm which is insensible to bottlenecks and swamps. It is also observed that dGN is less sensitive to the initialization, whereas ALS is highly sensitive, and dGN is more accurate than ALS in general \cite{tomasi2, tensorlab2}. This means dGN is a more robust and reliable algorithm in every sense. The drawback of the dGN algorithm is the computation of a solution to~\ref{dGN_eq} at each iteration. To solve these normal equations we have to deal with a big and dense matrix, which is computationally costly. In the next section we will present our approach to this problem.	
		
		\subsection{Dealing with the Hessian}		
			To work the normal equations~\ref{dGN_eq} we need to exploit some structure of $\textbf{J}_f^T \textbf{J}_f$ in order to make fast computations with low memory cost. The next theorem is a first step in this direction.   	
	
			\begin{theorem} \label{JfT_Jf}
				Denote $\omega_{r' r''}^{(\ell)} = \langle \textbf{w}_{r'}^{(\ell)}, \textbf{w}_{r''}^{(\ell)} \rangle$. Then we have that			
			
				$$\textbf{J}_f^T \textbf{J}_f = 
					\left[ \begin{array}{ccc}
						\textbf{H}_{11} & \ldots & \textbf{H}_{1L}\\
						\vdots & & \vdots\\
						\textbf{H}_{L1} & \ldots & \textbf{H}_{LL}
					\end{array} \right],$$
where

				$$\textbf{H}_{\ell' \ell''} = 
					\left[ \begin{array}{ccc}
					\displaystyle \prod_{\ell \neq \ell', \ell''} \omega_{11}^{(\ell)} \cdot \textbf{w}_1^{(\ell')} \textbf{w}_1^{(\ell'')^T} & \ldots & \displaystyle \prod_{\ell \neq \ell', \ell''} \omega_{1R}^{(\ell)} \cdot \textbf{w}_R^{(\ell')} \textbf{w}_1^{(\ell'')^T}\\
					\vdots & & \vdots\\
					\displaystyle \prod_{\ell \neq \ell', \ell''} \omega_{R1}^{(\ell)} \cdot \textbf{w}_1^{(\ell')} \textbf{w}_R^{(\ell'')^T} & \ldots & \displaystyle \prod_{\ell \neq \ell', \ell''} \omega_{RR}^{(\ell)} \cdot \textbf{w}_R^{(\ell')} \textbf{w}_R^{(\ell'')^T}
				\end{array} \right]$$
for $\ell' \neq \ell''$, and

				$$\textbf{H}_{\ell' \ell'} = 
					\left[ \begin{array}{ccc}
					\displaystyle \prod_{\ell \neq \ell'} \omega_{11}^{(\ell)} \cdot \textbf{I}_{I_{\ell'}} & \ldots & \displaystyle \prod_{\ell \neq \ell'} \omega_{1R}^{(\ell)} \cdot \textbf{I}_{I_{\ell'}}\\
					\vdots & & \vdots\\
					\displaystyle \prod_{\ell \neq \ell'} \omega_{R1}^{(\ell)} \cdot \textbf{I}_{I_{\ell'}} & \ldots & \displaystyle \prod_{\ell \neq \ell'} \omega_{RR}^{(\ell)} \cdot \textbf{I}_{I_{\ell'}}
				\end{array} \right].$$
				
			\end{theorem}
			
			\textbf{Proof:} First notice that
			
				$$\textbf{J}_f^T \textbf{J}_f = 
			 	\left[ \begin{array}{c}
			 		\displaystyle\frac{\partial f}{\partial \textbf{W}^{(1)}}^T\\ 
			 		\vdots\\ 
			 		\displaystyle\frac{\partial f}{\partial \textbf{W}^{(L)}}^T 
				 \end{array} \right]
			 	 \left[ \frac{\partial f}{\partial \textbf{W}^{(1)}}, \ldots, \frac{\partial f}{\partial \textbf{W}^{(L)}} \right] = \left[ \begin{array}{ccc}
			 	 	\displaystyle\frac{\partial f}{\partial \textbf{W}^{(1)}}^T \frac{\partial f}{\partial \textbf{W}^{(1)}} & \ldots & \displaystyle\frac{\partial f}{\partial \textbf{W}^{(1)}}^T \frac{\partial f}{\partial \textbf{W}^{(L)}}\\
			  		\vdots & & \vdots\\
			  		\displaystyle\frac{\partial f}{\partial \textbf{W}^{(L)}}^T \frac{\partial f}{\partial \textbf{W}^{(1)}} & \ldots & \displaystyle\frac{\partial f}{\partial \textbf{W}^{(L)}}^T \frac{\partial f}{\partial \textbf{W}^{(L)}}\\
				 \end{array} \right],$$
where

				$$\hspace{-.5cm} \frac{\partial f}{\partial \textbf{W}^{(\ell')}}^T \frac{\partial f}{\partial \textbf{W}^{(\ell'')}} = 
				\left[ \begin{array}{c}
					\displaystyle\frac{\partial f}{\partial \textbf{w}_1^{(\ell')}}^T\\
					\vdots\\ 
					\displaystyle\frac{\partial f}{\partial \textbf{w}_R^{(\ell')}}^T 
				\end{array} \right]
				\left[ \frac{\partial f}{\partial \textbf{w}_1^{(\ell'')}}, \ldots, \frac{\partial f}{\partial \textbf{w}_R^{(\ell'')}} \right] = 
				\left[ \begin{array}{ccc}
					\displaystyle\frac{\partial f}{\partial \textbf{w}_1^{(\ell')}}^T \frac{\partial f}{\partial \textbf{w}_1^{(\ell'')}} & \ldots & \displaystyle\frac{\partial f}{\partial \textbf{w}_1^{(\ell')}}^T \frac{\partial f}{\partial \textbf{w}_R^{(\ell'')}}\\
					\vdots & & \vdots\\
					\displaystyle\frac{\partial f}{\partial \textbf{w}_R^{(\ell')}}^T \frac{\partial f}{\partial \textbf{w}_1^{(\ell'')}} & \ldots & \displaystyle\frac{\partial f}{\partial \textbf{w}_R^{(\ell')}}^T \frac{\partial f}{\partial \textbf{w}_R^{(\ell'')}}
				\end{array} \right].$$
			
				Let $\omega_{r' r''}^{(\ell)} = \langle \textbf{w}_{r'}^{(\ell)}, \textbf{w}_{r''}^{(\ell)} \rangle$ and assume, without loss of generality, that $1 \leq \ell' < \ell'' \leq L$. We can use lemma~\ref{kronecker} and theorem~\ref{special-products} to write the entries of this matrix as
			
				$$\frac{\partial f}{\partial \textbf{w}_{r'}^{(\ell')}}^T \frac{\partial f}{\partial \textbf{w}_{r''}^{(\ell'')}} = $$
			
				$$\hspace{-.5cm} = \left( \textbf{w}_{r'}^{(1)} \tilde{\otimes} \ldots \tilde{\otimes} \textbf{w}_{r'}^{(\ell'-1)} \tilde{\otimes} \textbf{I}_{I_{\ell'}} \tilde{\otimes} \textbf{w}_{r'}^{(\ell'+1)} \tilde{\otimes} \ldots \tilde{\otimes} \textbf{w}_{r'}^{(L)} \right)^T \left( \textbf{w}_{r''}^{(1)} \tilde{\otimes} \ldots \tilde{\otimes} \textbf{w}_{r''}^{(\ell''-1)} \tilde{\otimes} \textbf{I}_{I_{\ell''}} \tilde{\otimes} \textbf{w}_{r''}^{(\ell''+1)} \tilde{\otimes} \ldots \tilde{\otimes} \textbf{w}_{r''}^{(L)} \right) = $$
			
				$$ = \omega_{r' r''}^{(1)} \tilde{\otimes} \ldots \tilde{\otimes} \omega_{r' r''}^{(\ell'-1)} \tilde{\otimes} \left( \textbf{I}_{\ell'} \textbf{w}_{r''}^{(\ell')} \right) \tilde{\otimes} \omega_{r' r''}^{(\ell'+1)} \tilde{\otimes} \ldots \tilde{\otimes} \omega_{r' r''}^{(\ell''-1)} \tilde{\otimes} \left( \textbf{w}_{r'}^{(\ell'')^T} \textbf{I}_{\ell''} \right) \tilde{\otimes} \omega_{r' r''}^{(\ell''+1)} \tilde{\otimes} \ldots \tilde{\otimes} \omega_{r' r''}^{(L)} = $$ 
			
				$$ = \prod_{\ell \neq \ell', \ell''} \omega_{r' r''}^{(\ell)} \cdot \textbf{w}_{r''}^{(\ell')} \textbf{w}_{r'}^{(\ell'')^T}.$$
			
				In the case $\ell' = \ell''$ we have
			
				$$\frac{\partial f}{\partial \textbf{w}_{r'}^{(\ell')}}^T \frac{\partial f}{\partial \textbf{w}_{r''}^{(\ell')}} = $$
			
				$$ = \omega_{r' r''}^{(1)} \tilde{\otimes} \ldots \tilde{\otimes} \omega_{r' r''}^{(\ell'-1)} \tilde{\otimes} \  \textbf{I}_{I_{\ell'}}^2 \tilde{\otimes} \ \omega_{r' r''}^{(\ell'+1)} \tilde{\otimes} \ldots \tilde{\otimes} \omega_{r' r''}^{(L)} = $$ 
			
				$$ = \prod_{\ell \neq \ell'} \omega_{r' r''}^{(\ell)} \cdot \textbf{I}_{I_{\ell'}}.$$
			
				Therefore,
			
				$$\frac{\partial f}{\partial \textbf{W}^{(\ell')}}^T \frac{\partial f}{\partial \textbf{W}^{(\ell'')}} = 
				\left[ \begin{array}{ccc}
					\displaystyle \prod_{\ell \neq \ell', \ell''} \omega_{11}^{(\ell)} \cdot \textbf{w}_1^{(\ell')} \textbf{w}_1^{(\ell'')^T} & \ldots & \displaystyle \prod_{\ell \neq \ell', \ell''} \omega_{1R}^{(\ell)} \cdot \textbf{w}_R^{(\ell')} \textbf{w}_1^{(\ell'')^T}\\
					\vdots & & \vdots\\
					\displaystyle \prod_{\ell \neq \ell', \ell''} \omega_{R1}^{(\ell)} \cdot \textbf{w}_1^{(\ell')} \textbf{w}_R^{(\ell'')^T} & \ldots & \displaystyle \prod_{\ell \neq \ell', \ell''} \omega_{RR}^{(\ell)} \cdot \textbf{w}_R^{(\ell')} \textbf{w}_R^{(\ell'')^T}
				\end{array} \right]$$
when $\ell' \neq \ell''$. Finally, we have that

				$$\hspace{2.8cm} \frac{\partial f}{\partial \textbf{W}^{(\ell')}}^T \frac{\partial f}{\partial \textbf{W}^{(\ell')}} = 
				\left[ \begin{array}{ccc}
					\displaystyle \prod_{\ell \neq \ell'} \omega_{11}^{(\ell)} \cdot \textbf{I}_{I_{\ell'}} & \ldots & \displaystyle \prod_{\ell \neq \ell'} \omega_{1R}^{(\ell)} \cdot \textbf{I}_{I_{\ell'}}\\
					\vdots & & \vdots\\
					\displaystyle \prod_{\ell \neq \ell'} \omega_{R1}^{(\ell)} \cdot \textbf{I}_{I_{\ell'}} & \ldots & \displaystyle \prod_{\ell \neq \ell'} \omega_{RR}^{(\ell)} \cdot \textbf{I}_{I_{\ell'}}
				\end{array} \right]. \hspace{2.8cm}\square$$

			\begin{remark}
				Each $\textbf{H}_{\ell' \ell''}$ is block $R \times R$ matrix, where each of its blocks is a $I_{\ell'} \times I_{\ell''}$ matrix, so $\textbf{H}_{\ell' \ell''}$ has shape $I_{\ell'}R \times I_{\ell''}R$. Since $\textbf{J}_f^T \textbf{J}_f$ is a block $L \times L$ matrix, with the $(\ell', \ell'')$ block being $\textbf{H}_{\ell' \ell''}$, we conclude that $\textbf{J}_f^T \textbf{J}_f$ has shape $R \displaystyle \sum_{\ell=1}^L I_\ell \times R \displaystyle \sum_{\ell=1}^L I_\ell$. This shape is much smaller than the shape of $\textbf{J}_f$, which is of $\displaystyle \prod_{\ell=1}^L I_\ell \times R \displaystyle \sum_{\ell=1}^L I_\ell$, hence we avoid the curse of dimensionality with this approach. It should be noted that this is only true if $R < \frac{ \prod_{\ell=1}^L I_\ell}{ \sum_{\ell=1}^L I_\ell}$. This will be almost always the case since this is the same as saying $R$ is smaller than the generic rank, and in fact almost always we will choose $R$ to satisfy this property. Finally, we want to remark that the notation $\textbf{H}_{\ell' \ell''}$ comes from the fact that $\textbf{J}_f^T \textbf{J}_f \approx \textbf{H}_F$ as we converges to the solution. 
			\end{remark}
			
			It is convenient to store the products of the terms $\omega_{r' r''}^{(\ell)}$ in matrix form, so we define
			
			$$\Pi^{(\ell', \ell'')} = 
			\left[ \begin{array}{ccc}
				\displaystyle \prod_{\ell \neq \ell', \ell''} \omega_{1 1}^{(\ell)} & \ldots & \displaystyle \prod_{\ell \neq \ell', \ell''} \omega_{1 R}^{(\ell)}\\
				\vdots & & \vdots\\
				\displaystyle \prod_{\ell \neq \ell', \ell''} \omega_{R 1}^{(\ell)} & \ldots & \displaystyle \prod_{\ell \neq \ell', \ell''} \omega_{R R}^{(\ell)}
			\end{array} \right]$$
for $\ell' \neq \ell''$, and

			$$\Pi^{(\ell')} = 
			\left[ \begin{array}{ccc}
				\displaystyle \prod_{\ell \neq \ell'} \omega_{1 1}^{(\ell)} & \ldots & \displaystyle \prod_{\ell \neq \ell'} \omega_{1 R}^{(\ell)}\\
				\vdots & & \vdots\\
				\displaystyle \prod_{\ell \neq \ell'} \omega_{R 1}^{(\ell)} & \ldots & \displaystyle \prod_{\ell \neq \ell'} \omega_{R R}^{(\ell)}
			\end{array} \right].$$
			
			Notice that these (symmetric) matrices are the Hadamard product of some Gramian matrices. Define $\pi^{(\ell)} = \textbf{W}^{(\ell)^T} \textbf{W}^{(\ell)}$. Then we have that
		$$\Pi^{(\ell', \ell'')} = \pi^{(1)} \ast \ldots \ast \pi^{(\ell'-1)} \ast \pi^{(\ell'+1)} \ast \ldots \ast \pi^{(\ell''-1)} \ast \pi^{(\ell''+1)} \ast \ldots \ast \pi^{(L)}$$		
and
			$$\Pi^{(\ell')} = \pi^{(1)} \ast \ldots \ast \pi^{(\ell'-1)} \ast \pi^{(\ell'+1)} \ast \ldots \ast \pi^{(L)}.$$	
		
			Now we can see that is possible to retrieve $\textbf{J}_f^T \textbf{J}_f$ from the factor matrices with few computations and low memory cost. First we compute and store all the products $\textbf{W}^{(\ell)^T} \textbf{W}^{(\ell)}$ and  this has a total computational cost of $\mathcal{O}\left(R^2 \displaystyle\sum_{\ell=1}^L I_\ell \right)$. To compute each $\Pi^{(\ell', \ell'')}$ we have to perform $L-2$ Hadamard products between $R \times R$ matrices, and this has a cost of $\mathcal{O}\left( (L-2)R^2 \right)$ flops. The cost to store all the matrices $\Pi^{(\ell', \ell')}$ is of $L^2R^2$ floats. We remark that it is possible to cut all these costs by half since $\Pi^{(\ell', \ell'')} = \Pi^{(\ell'', \ell')}$ for all $\ell' \neq \ell''$.  Finally, note that $\Pi^{(\ell')} = \Pi^{(\ell', \ell'')} \ast \big(\textbf{W}^{(\ell'')^T} \textbf{W}^{(\ell'')}\big)$, so we can construct $\Pi^{(\ell')}$ with only one Hadamard product, which has a cost of $\mathcal{O}(R^2)$ flops. Doing this for all $\ell$ amounts to $\mathcal{O}(LR^2)$ flops. The memory cost to store each $\Pi^{(\ell')}$ is the same of $\Pi^{(\ell', \ell')}$, so we have a total memory cost of $\mathcal{O}(LR^2)$ floats. Overall, the computational cost to compute all $\Pi^{(\ell', \ell'')}$ and all $\Pi^{(\ell')}$ is of $\mathcal{O}\left( R^2 \left( L + \displaystyle\sum_{\ell=1}^L I_\ell \right) \right)$ flops, and the total memory cost is of $\mathcal{O}\left( R^2(L + L^2) \right)$ floats. 
			
			By writing the matrices $\Pi^{(\ell', \ell'')}$ and $\Pi^{(\ell')}$ as the result of many Hadamard products we have a way to simplify the expressions for the blocks $\textbf{H}_{\ell' \ell''}$.
			
			\begin{corollary}
				Let $\textbf{1}_{m \times n}$ be the $m \times n$ matrix constituted only by ones and define
				$$\textbf{K}^{(\ell', \ell'')} = 
				\left[ \begin{array}{ccc}
					\textbf{w}_1^{(\ell')} \textbf{w}_1^{(\ell'')^T} & \ldots & \textbf{w}_R^{(\ell')} \textbf{w}_1^{(\ell'')^T}\\
					\vdots & & \vdots\\
					\textbf{w}_1^{(\ell')} \textbf{w}_R^{(\ell'')^T} & \ldots & \textbf{w}_R^{(\ell')} \textbf{w}_R^{(\ell'')^T}
				\end{array} \right]$$ 
for $\ell' \neq \ell''$. Then 
				$$\textbf{H}_{\ell', \ell''} = \left( \Pi^{(\ell', \ell'')} \tilde{\otimes} \textbf{1}_{I_{\ell'} \times I_{\ell''}} \right) \ast \textbf{K}^{(\ell', \ell'')}$$
for $\ell' \neq \ell''$, and
				$$\textbf{H}_{\ell', \ell'} = \Pi^{(\ell')} \tilde{\otimes} \textbf{I}_{\ell'}.$$
			\end{corollary}
			
			As already mentioned, $\textbf{J}_f^T \textbf{J}_f$ will be used to solve the normal equations~\ref{dGN_eq}. The algorithm of choice to accomplish this is the conjugate gradient method (see appendix~\ref{appenA}). This classical algorithm is particularly efficient to solve normal equations where the matrix is positive definite, which is our case. Furthermore, our version of the conjugate gradient is \emph{matrix-free}, that is, we are able to compute matrix-vector products $\textbf{J}_f^T \textbf{J}_f \cdot \textbf{v}$ without actually constructing $\textbf{J}_f^T \textbf{J}_f$. By exploiting the block structure of $\textbf{J}_f^T \textbf{J}_f$ we can save memory and the computational cost still is lower than the naive cost of $R^2 \left( \displaystyle \sum_{\ell=1}^L I_\ell \right)^2$ flops.  
			
			The next theorem is a new contribution of this work. With this result we are able to perform fast matrix-vector products with low memory cost. 
		
			\begin{theorem} \label{Hv}
				Given any vector $\textbf{v} \in \mathbb{R}^{R \sum_{\ell=1}^L I_\ell}$, write 			
			$$\textbf{v} = 
			\left[ \begin{array}{c}
				vec(\textbf{V}^{(1)})\\
				\vdots\\
				vec(\textbf{V}^{(L)})
			\end{array} \right]$$
where $\textbf{V}^{(\ell)} = [\textbf{v}_1^{(\ell)}, \ldots, \textbf{v}_R^{(\ell)}] \in \mathbb{R}^{I_\ell \times R}$ and each $\textbf{v}_r^{(\ell)}$ is a column of $\textbf{V}^{(\ell)}$. Then  
			$$\textbf{J}_f^T \textbf{J}_f \cdot \textbf{v} = \left[ \begin{array}{c}
					\displaystyle \sum_{\ell=1}^L \frac{\partial f}{\partial \textbf{W}^{(1)}}^T \frac{\partial f}{\partial \textbf{W}^{(\ell)}} \cdot 	vec(\textbf{V}^{(\ell)})\\
					\vdots\\	
					\displaystyle \sum_{\ell=1}^L \frac{\partial f}{\partial \textbf{W}^{(L)}}^T \frac{\partial f}{\partial \textbf{W}^{(\ell)}} \cdot 	vec(\textbf{V}^{(\ell)})
				\end{array} \right]$$
where 
			$$\frac{\partial f}{\partial \textbf{W}^{(\ell')}}^T \frac{\partial f}{\partial \textbf{W}^{(\ell'')}} \cdot 	vec(\textbf{V}^{(\ell'')}) = vec\left( \textbf{W}^{(\ell')} \cdot \Big( \Pi^{(\ell', \ell'')} \ast \big( \textbf{V}^{(\ell'')^T} \cdot \textbf{W}^{(\ell'')} \big) \Big) \right)$$
for $\ell' \neq \ell''$ and
			$$\frac{\partial f}{\partial \textbf{W}^{(\ell')}}^T \frac{\partial f}{\partial \textbf{W}^{(\ell')}} \cdot 	vec(\textbf{V}^{(\ell')}) = vec\left( \textbf{V}^{(\ell')} \cdot \Pi^{(\ell')} \right).$$			
			\end{theorem}
			
			\textbf{Proof:} First notice that

			$$\textbf{J}_f^T \textbf{J}_f \cdot \textbf{v} = 
			\left[ \begin{array}{ccc}
			 	\displaystyle\frac{\partial f}{\partial \textbf{W}^{(1)}}^T \frac{\partial f}{\partial \textbf{W}^{(1)}} & \ldots & \displaystyle\frac{\partial f}{\partial \textbf{W}^{(1)}}^T \frac{\partial f}{\partial \textbf{W}^{(L)}}\\
			  	\vdots & & \vdots\\
			   \displaystyle\frac{\partial f}{\partial \textbf{W}^{(L)}}^T \frac{\partial f}{\partial \textbf{W}^{(1)}} & \ldots & \displaystyle\frac{\partial f}{\partial \textbf{W}^{(L)}}^T \frac{\partial f}{\partial \textbf{W}^{(L)}}\\
			\end{array} \right]
			\left[ \begin{array}{c}
				vec(\textbf{V}^{(1)})\\
				\vdots\\
				vec(\textbf{V}^{(L)})
			\end{array} \right] = $$
			
			$$ = 
				\left[ \begin{array}{c}
					\displaystyle \sum_{\ell=1}^L \frac{\partial f}{\partial \textbf{W}^{(1)}}^T \frac{\partial f}{\partial \textbf{W}^{(\ell)}} \cdot 	vec(\textbf{V}^{(\ell)})\\
					\vdots\\	
					\displaystyle \sum_{\ell=1}^L \frac{\partial f}{\partial \textbf{W}^{(L)}}^T \frac{\partial f}{\partial \textbf{W}^{(\ell)}} \cdot 	vec(\textbf{V}^{(\ell)})
				\end{array} \right].$$
			
			Now we simplify each term in the summation above. It is necessary to consider two separate cases.
			
			\textbf{Case 1 (different modes):} If $\ell' \neq \ell''$, then 
			
			$$\frac{\partial f}{\partial \textbf{W}^{(\ell')}}^T \frac{\partial f}{\partial \textbf{W}^{(\ell'')}} \cdot 	vec(\textbf{V}^{(\ell'')}) = $$
			
			$$ = 
			\left[ \begin{array}{ccc}
				\displaystyle \prod_{\ell \neq \ell', \ell''} \omega_{11}^{(\ell)} \cdot \textbf{w}_1^{(\ell')} \textbf{w}_1^{(\ell'')^T} & \ldots & \displaystyle \prod_{\ell \neq \ell', \ell''} \omega_{1R}^{(\ell)} \cdot \textbf{w}_R^{(\ell')} \textbf{w}_1^{(\ell'')^T}\\
				\vdots & & \vdots\\
				\displaystyle \prod_{\ell \neq \ell', \ell''} \omega_{R1}^{(\ell)} \cdot \textbf{w}_1^{(\ell')} \textbf{w}_R^{(\ell'')^T} & \ldots & \displaystyle \prod_{\ell \neq \ell', \ell''} \omega_{RR}^{(\ell)} \cdot \textbf{w}_R^{(\ell')} \textbf{w}_R^{(\ell'')^T}
			\end{array} \right] 
			\left[ \begin{array}{c}
				\textbf{v}_1^{(\ell'')}\\
				\vdots\\
				\textbf{v}_R^{(\ell'')}
			\end{array} \right] = $$
				
			$$ = 
			\left[ \begin{array}{c}
				\displaystyle \sum_{r=1}^R \prod_{\ell \neq \ell', \ell''} \omega_{1r}^{(\ell)} \cdot \textbf{w}_r^{(\ell')} \textbf{w}_1^{(\ell'')^T} \cdot \textbf{v}_r^{(\ell'')}\\
				\vdots\\
				\displaystyle \sum_{r=1}^R \prod_{\ell \neq \ell', \ell''} \omega_{Rr}^{(\ell)} \cdot \textbf{w}_r^{(\ell')} \textbf{w}_R^{(\ell'')^T} \cdot \textbf{v}_r^{(\ell'')}\\
			\end{array} \right] =  
			\left[ \begin{array}{c}
				\displaystyle \sum_{r=1}^R \prod_{\ell \neq \ell', \ell''} \omega_{1r}^{(\ell)} \cdot \textbf{w}_r^{(\ell')} \langle \textbf{w}_1^{(\ell'')}, \textbf{v}_r^{(\ell'')} \rangle\\
				\vdots\\
				\displaystyle \sum_{r=1}^R \prod_{\ell \neq \ell', \ell''} \omega_{Rr}^{(\ell)} \cdot \textbf{w}_r^{(\ell')} \langle \textbf{w}_R^{(\ell'')}, \textbf{v}_r^{(\ell'')} \rangle
			\end{array} \right] = $$
				
			$$ = 
			\left[ \begin{array}{c}
				\left[ \textbf{w}_1^{(\ell')}, \ldots, \textbf{w}_R^{(\ell')} \right]
				\left[ \begin{array}{c}
					\displaystyle \prod_{\ell \neq \ell', \ell''} \omega_{11}^{(\ell)} \langle \textbf{w}_1^{(\ell'')}, \textbf{v}_1^{(\ell'')} \rangle\\
					\vdots\\
					\displaystyle \prod_{\ell \neq \ell', \ell''} \omega_{1R}^{(\ell)} \langle \textbf{w}_1^{(\ell'')}, \textbf{v}_R^{(\ell'')} \rangle
				\end{array} \right]\\
				\vdots\\
				\left[ \textbf{w}_1^{(\ell')}, \ldots, \textbf{w}_R^{(\ell')} \right]
				\left[ \begin{array}{c}
					\displaystyle \prod_{\ell \neq \ell', \ell''} \omega_{R1}^{(\ell)} \langle \textbf{w}_R^{(\ell'')}, \textbf{v}_1^{(\ell'')} \rangle\\
					\vdots\\
					\displaystyle \prod_{\ell \neq \ell', \ell''} \omega_{RR}^{(\ell)} \langle \textbf{w}_R^{(\ell'')}, \textbf{v}_R^{(\ell'')} \rangle
				\end{array} \right]
			\end{array} \right] = 
			\left[ \begin{array}{c}
				\textbf{W}^{(\ell')} \cdot \Big( \Pi_1^{(\ell', \ell'')} \ast \big( \textbf{V}^{(\ell'')^T} \cdot \textbf{w}_1^{(\ell'')} \big) \Big)\\
				\vdots\\
				\textbf{W}^{(\ell')} \cdot \Big( \Pi_R^{(\ell', \ell'')} \ast \big( \textbf{V}^{(\ell'')^T} \cdot \textbf{w}_R^{(\ell'')} \big) \Big)
			\end{array} \right] = $$
			
			$$ = 
			vec\left( \left[ \textbf{W}^{(\ell')} \cdot \Big( \Pi_1^{(\ell', \ell'')} \ast \big( \textbf{V}^{(\ell'')^T} \cdot \textbf{w}_1^{(\ell'')} \big) \Big), \ldots, \textbf{W}^{(\ell')} \cdot \Big( \Pi_R^{(\ell', \ell'')} \ast \big( \textbf{V}^{(\ell'')^T} \cdot \textbf{w}_R^{(\ell'')} \big) \Big) \right] \right) = $$
			
			$$ =
			vec\left( \textbf{W}^{(\ell')} \cdot \Big( \Pi^{(\ell', \ell'')} \ast \big( \textbf{V}^{(\ell'')^T} \cdot \textbf{W}^{(\ell'')} \big) \Big) \right)$$
where each $\Pi_{r'}^{(\ell', \ell'')}$ is the $r'$-th column of $\Pi^{(\ell', \ell'')}$. Despite the notation refers to the rows of $\Pi^{(\ell', \ell'')}$, this is not a problem since this matrix is symmetric.\\
				
			\textbf{Case 2 (equal modes):} For a mode $\ell'$ we have 
			
			$$\frac{\partial f}{\partial \textbf{W}^{(\ell')}}^T \frac{\partial f}{\partial \textbf{W}^{(\ell')}} \cdot 	vec(\textbf{V}^{(\ell')}) = $$	
			
			$$ = 
			\left[ \begin{array}{ccc}
				\displaystyle \prod_{\ell \neq \ell'} \omega_{11}^{(\ell)} \cdot \textbf{I}_{I_\ell'} & \ldots & \displaystyle \prod_{\ell \neq \ell'} \omega_{1R}^{(\ell)} \cdot \textbf{I}_{I_\ell'}\\
				\vdots & & \vdots\\
				\displaystyle \prod_{\ell \neq \ell'} \omega_{R1}^{(\ell)} \cdot \textbf{I}_{I_\ell'} & \ldots & \displaystyle \prod_{\ell \neq \ell'} \omega_{RR}^{(\ell)} \cdot \textbf{I}_{I_\ell'}
			\end{array} \right]
			\left[ \begin{array}{c}
				\textbf{v}_1^{(\ell')}\\
				\vdots\\
				\textbf{v}_R^{(\ell')}
			\end{array}	\right] =  
			\left[ \begin{array}{ccc}
				\displaystyle \sum_{r=1}^R \prod_{\ell \neq \ell'} \omega_{1r}^{(\ell)} \cdot \textbf{v}_r^{(\ell')}\\
				\vdots\\
				\displaystyle \sum_{r=1}^R \prod_{\ell \neq \ell'} \omega_{Rr}^{(\ell)} \cdot \textbf{v}_r^{(\ell')}
			\end{array} \right] = $$
			
			$$\hspace{1cm} = 
			\left[ \begin{array}{c}
				\textbf{V}^{(\ell')} \cdot \Pi_1^{(\ell')}\\ 
				\vdots\\ 
				\textbf{V}^{(\ell')} \cdot \Pi_R^{(\ell')}
			\end{array} \right] = 
			vec\left( \textbf{V}^{(\ell')} \cdot \Pi_1^{(\ell')}, \ldots, \textbf{V}^{(\ell')} \cdot \Pi_R^{(\ell')} \right) = vec\left( \textbf{V}^{(\ell')} \cdot \Pi^{(\ell')} \right). \hspace{1cm}\square$$
				
			For $\ell' \neq \ell''$, the computation of $vec\left( \textbf{W}^{(\ell')} \cdot \Big( \Pi^{(\ell', \ell')} \ast \big( \textbf{V}^{(\ell'')^T} \cdot \textbf{W}^{(\ell'')} \big) \Big) \right)$ requires two matrix-matrix multiplications and one Hadamard product. This approach benefits from the BLAS-3 efficiency while the Hadamard product can efficiently be done in parallel. Now let's we make some complexity analysis. First, since $\textbf{V}^{(\ell'')^T}$ is of shape $R \times I_{\ell''}$ and $\textbf{W}^{(\ell'')}$ is $I_{\ell''} \times R$, multiplying them requires $I_{\ell''} R^2$ flops. The resulting matrix is of shape $R \times R$, so the Hadamard product by $\Pi^{(\ell', \ell')}$ requires more $R^2$ flops. Finally, since $\textbf{W}^{(\ell')}$ is of shape $I_{\ell'} \times R$ and $\Big( \Pi^{(\ell', \ell')} \ast \big( \textbf{V}^{(\ell'')^T} \cdot \textbf{W}^{(\ell'')} \big) \Big)$ is of shape $R \times R$, multiplying them requires $I_{\ell'} R^2$ flops. Overall the cost is of $(1 + I_{\ell'} + I_{\ell''})R^2$ flops.
			
			The cost of computing $\textbf{V}^{(\ell')} \ast \Pi^{(\ell')}$ is more straightforward, being of $I_{\ell'} R^2$ flops. For each $\ell' = 1 \ldots L$, the compute the sum $\displaystyle \sum_{\ell=1}^L \frac{\partial f}{\partial \textbf{W}^{(\ell')}}^T \frac{\partial f}{\partial \textbf{W}^{(\ell)}} \cdot vec(\textbf{V}^{(\ell)})$ requires to perform $(1 + I_{\ell'} + I_{\ell}) R^2$ flops for all $\ell \neq \ell'$ and them more $I_{\ell'} R^2$ flops for $\ell'$. The total cost of doing this is of $I_{\ell'} R^2 + \displaystyle \sum_{\ell \neq \ell'} (1 + I_{\ell'} + I_{\ell})R^2 = \left( (L-1) + (L-1) I_{\ell'} + \left( \displaystyle \sum_{\ell=1}^L I_\ell \right) \right) R^2$ flops. Since we have to do it for each $\ell'$, the total cost is of $\displaystyle \sum_{\ell'=1}^L \left( (L-1) + (L-1) I_{\ell'} + \left( \displaystyle \sum_{\ell=1}^L I_\ell \right) \right) R^2 = \left( L(L-1) + (2L-1)\left( \displaystyle \sum_{\ell=1}^L I_\ell \right) \right) R^2$ flops. We can summarize al this analysis by saying that the cost of the matrix-vector multiplication $\textbf{J}_f^T \textbf{J}_f \cdot \textbf{v}$ is of $\mathcal{O}\left( L R^2 \displaystyle \sum_{\ell=1}^L I_\ell \right)$ flops. Recall that the naive multiplication costs $R^2 \left( \displaystyle \sum_{\ell=1}^L I_\ell \right)^2$ flops. Hence the gain in performance comes from the exploitation of the BLAS-3 implementation and the low complexity cost which was possible because of the block structure of $\textbf{J}_f^T \textbf{J}_f$. 
			
			Table 3.1 gather the information of all costs we obtained until now. From this point we will always summarize all costs analysis in tables at the final of each section. This will help the reader to make an overview of the costs when necessary.\\
			
			\begin{table}
				\centering
				\begin{tabular}{|c|c|c|}
					\hline
					\textbf{Task} & \textbf{Memory} & \textbf{Computational time}\\
					\hline
					Computing all $\Pi^{(\ell', \ell'')}$ and $\Pi^{(\ell')}$ & $R^2(L + L^2)$ & $\mathcal{O}\left( R^2 \left( L + \displaystyle \sum_{\ell=1}^L I_\ell \right) \right)$\\
					\hline
					Computing $\textbf{J}_f^T \textbf{J}_f \cdot \textbf{v}$ & $R^2 \displaystyle \sum_{\ell=1}^L I_\ell$ & $\mathcal{O}\left( L R^2 \displaystyle \sum_{\ell=1}^L I_\ell \right)$\\
					\hline
				\end{tabular}
				\caption{\footnotesize{Memory and computational costs - I.}} \bigskip
			\end{table}	

%% file: CHAPTER_4.tex
\chapter{Computational experiments}\label{cap-4} 
	This chapter starts presenting Tensor Fox, with a detailed discussion of main subroutines and their costs. Comprehensive experiments follows, first with the introduction of the tensors used for the benchmarks. We tried to use a wide distinct choice of tensors: positive tensors, tensors from machine learning, ill-conditioned tensors, and so on. The parameters of Tensor Fox are fine tuned against this set of tensors, so they are general enough. Then we conduct the benchmarks and discuss the results. The difference between Tensor Fox and other packages are discussed too. 
	
	The last five sections discuss are a discussion of the main aspects of Tensor Fox. In section~\ref{not_mon} we connect the gain ratio and the fact that Tensor Fox is not monotonic, and how this is a good thing. In sections~\ref{diag_reg},~\ref{cond_section} and~\ref{par} we discuss the diagonal regularization, conditioning and parallelism, respectively, with lots of more experiments validating our claims. Finally, in section~\ref{main_feat} we will see more details about the features of Tensor Fox. Which of them are more relevant, and which are less relevant.    

	\section{Tensor Fox}
		One of the main contributions of this work is the algorithm described in this section. An implementation of this algorithm is available (open and free) for Python, by the name of \emph{Tensor Fox}, check the link \url{https://github.com/felipebottega/Tensor-Fox}. In this section we present and explain it in details, together with computational and memory costs. The algorithm will be presented in parts, following the computational flow as showed below in figure~\ref{tensorfox-flow}. Each box represents a major part of the algorithm, and each one is an algorithm by itself. 
		
		The \emph{Compress} box is the computation of the MLSVD of the tensor. Most of the time it is unnecessary to work in the original space since the compressed tensor is equivalent to the original one in the sense of definition~\ref{equivalence}. In theory, if we find and exactly CPD and uncompress it, then this uncompressed CPD is an exact CPD for the tensor at the original space. After compressing, we must generate an initial tensor to start the dGN iterations. This is done at the \emph{Initialization} part. More ahead we will show an original way to generating good initializations based on the MLSVD. The \emph{dGN} algorithm is the most costly part and is there which lies some new ideas obtained after a lot of experimentation. These ideas are the result of much ``mathematical alchemy''. Finally, after computing a CPD we just need to \emph{Uncompress} the solution.  
		
		\begin{figure} 
		\begin{center}		
		\begin{tikzpicture}[>=stealth, thick]
			\node (A) at (0,0) [draw, terminal, text width=6cm, minimum height=0.5cm, align=flush center] 
			{\textbf{Input:} $R \in \mathbb{N}, \ \mathcal{T} \in \mathbb{R}^{I_1 \times \ldots \times I_L}$};
			
			\node (B) at (0,-2) [draw, process, minimum height=0.5cm, align=flush center] 
			{Compression};

			\node (C) at (0,-4) [draw, process, minimum height=0.5cm, align=flush center] 
			{Initialization};
			
			\node (D) at (0,-6) [draw, process, minimum height=0.5cm, align=flush center] 
			{dGN};
			
			\node (E) at (0,-8) [draw, process, minimum height=0.5cm, align=flush center] 
			{Uncompression};
			
			\node (F) at (0,-10) [draw, terminal, text width=10cm, minimum height=0.5cm, align=flush center] 
			{\textbf{Output:} $\textbf{W}^{(\ell)} \in \mathbb{R}^{I_\ell \times R} \text{ for } \ell = 1 \ldots L$ and $\Lambda \in \mathbb{R}^{R \times \ldots \times R}$ diagonal, such that each $\textbf{W}^{(\ell)}$ has unit columns and $\mathcal{T} \approx \sum_{r=1}^R \lambda_r \ \textbf{w}_r^{(1)} \otimes \ldots \otimes \textbf{w}_r^{(L)}$};

			\draw[->] (A) -- (B);
			\draw[->] (B) -- (C);
			\draw[->] (C) -- (D);
			\draw[->] (D) -- (E);
			\draw[->] (E) -- (F);
		\end{tikzpicture}
		\end{center}
		\caption{\footnotesize{Flow chart of the main parts of Tensor Fox.}}
		\label{tensorfox-flow}
		\end{figure}
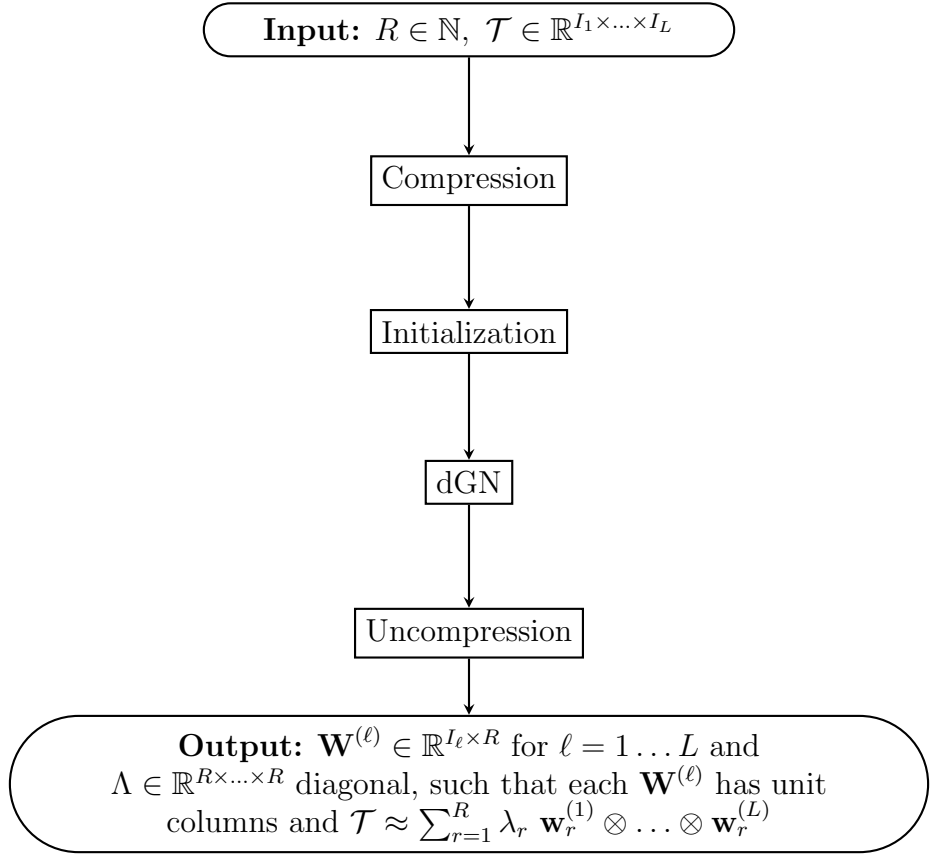	
			
		\subsection{Compression} \label{compression}		 
			Working with ``raw'' data is, usually, not advised because of the typical large data size. A standard approach is to compress the data before starting the actual work, and this is not different in the context of tensors. This is an efficient way of reducing the effects of the curse of dimensionality. The main tool to compress a tensor is the MLSVD, and it is interesting to use it because remark~\ref{mlsvd-cpd} guarantees one can use the core tensor of the MLSVD to compute a CPD for the original tensor. The idea is very similar to the procedure of compressing matrix data with the SVD: we compute the SVD of the matrix and truncate it to a matrix with lower rank. The choice of this rank is such that the truncated version is small enough to be tractable but close enough to the original matrix. Now we use this same idea on the tensor context. 
			
			First, compute the MLSVD of $\mathcal{T}$, obtaining a decomposition $\mathcal{T} = (\textbf{U}^{(1)}, \ldots, \textbf{U}^{(L)}) \cdot \mathcal{S}$ where $\textbf{U}^{(\ell)} \in \mathbb{R}^{I_\ell \times I_\ell}$ are orthogonal for $\ell = 1 \ldots L$, and $\mathcal{S} \in \mathbb{R}^{I_1 \times \ldots \times I_L}$ is the core tensor (see theorem~\ref{MLSVD}). If $rank_\boxplus(\mathcal{T}) = (R_1, \ldots, R_L)$, then we can discard the last $I_\ell - R_\ell$ columns of each $\textbf{U}^{(\ell)}$ and consider $\textbf{U}^{(\ell)} \in \mathbb{R}^{I_\ell \times R_\ell}$, also we discard all hyperslices $\mathcal{S}_{i_\ell = k}$ for $k = R_\ell+1 \ldots I_\ell$ and consider $\mathcal{S} \in \mathbb{R}^{R_1 \times \ldots \times R_L}$. After these truncations, the equality $\mathcal{T} = (\textbf{U}^{(1)}, \ldots, \textbf{U}^{(L)}) \cdot \mathcal{S}$ remains intact. Actually, these truncations are not ``real truncations'' because we only transformed the decomposition from the full format to the reduced format, in the same way there is the full SVD and the reduced SVD for matrices. To get a real truncation with low multilinear rank we need to keep deleting columns of $\textbf{U}^{(\ell)}$ and hyperslices of $\mathcal{S}$. 
			
\begin{figure}
				\centering
				\includegraphics[scale=.5]{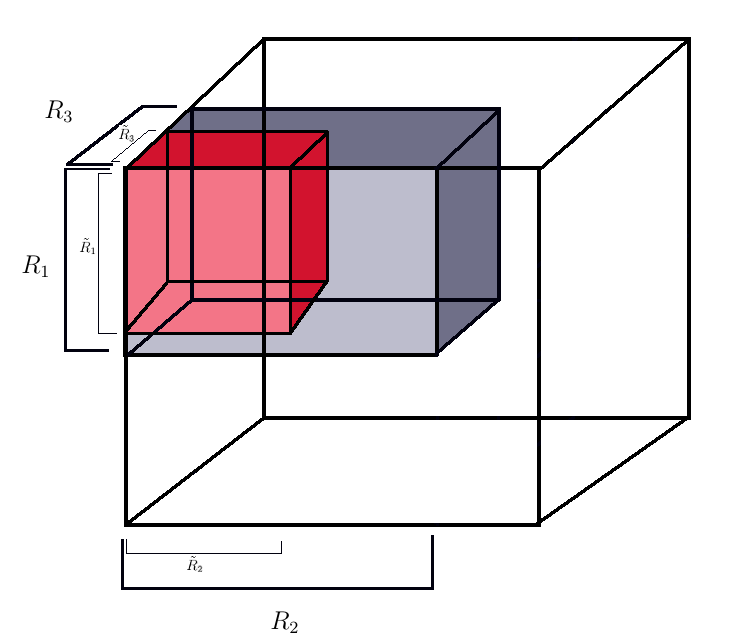}
				\caption{\footnotesize{Truncated tensor $\tilde{\mathcal{S}}$.}}
				\label{trunc}
			\end{figure}

			Let $(\tilde{R}_1, \ldots, \tilde{R}_L) \leq (R_1, \ldots, R_L)$ be a lower multilinear rank. We define $\tilde{\textbf{U}}^{(\ell)} = \left[ \textbf{U}_{:1}^{(\ell)}, \ldots, \textbf{U}_{:\tilde{R}_\ell}^{(\ell)} \right] \in \mathbb{R}^{I_\ell \times \tilde{R}_\ell}$ to be the matrix composed by the first columns of $\textbf{U}^{(\ell)}$, and $\tilde{\mathcal{S}} \in \mathbb{R}^{\tilde{R}_1 \times \ldots \times \tilde{R}_L}$ is such that $\tilde{s}_{i_1 \ldots i_L} = s_{i_1 \ldots i_L}$ for $1 \leq i_1 \leq \tilde{R}_1, \ldots, 1 \leq i_L \leq \tilde{R}_L$. Figure~\ref{trunc} illustrates such a truncation in the case of a third order tensor. The white part correspond to $\mathcal{S}$ after we computed the full MLSVD (the one of theorem~\ref{MLSVD}), the gray tensor is the reduced format of $\mathcal{S}$, and the red tensor is the truncated tensor $\tilde{\mathcal{S}}$. Our goal is to find the smallest $(\tilde{R}_1, \ldots, \tilde{R}_L)$ such that $\| \mathcal{S} - \tilde{\mathcal{S}} \|$ is not so large\footnote{Actually, $\mathcal{S}$ and $\tilde{\mathcal{S}}$ belongs to different spaces. We committed an abuse of notation and wrote $\| \mathcal{S} - \tilde{\mathcal{S}} \|$ considering the projection of $\tilde{\mathcal{S}}$ over the space of $\mathcal{S}$, that is, enlarge $\tilde{\mathcal{S}}$ so it has the same size of $\mathcal{S}$ and consider these new entries as zeros. This is how we are projecting.}. Since reducing $(\tilde{R}_1, \ldots, \tilde{R}_L)$ too much causes $\| \mathcal{S} - \tilde{\mathcal{S}} \|$ to increase, there is a trade off we have to manage in the best way possible. 
			
			Ideally, we want that each estimate $\tilde{R}_\ell$ to be exactly equal to $R_\ell$, but that is not always possible. Furthermore, if $\mathcal{T}$ comes from data with noise, we are willing to truncate in order to work only with the relevant information. This relevant information is concentrated after the MLSVD is performed, and we can measure it in terms of the energy distribution, a topic already discussed in chapter 2. Because of theorem~\ref{MLSVD-properties}-1, we have that $R_\ell$ is the rank of the unfolding $\mathcal{T}_{(\ell)}$. Since the computation of the MLSVD requires the computation of the SVD of each unfolding, we can truncate these SVDs right after their computation. In fact, we can already start computing a truncated SVD with $P_\ell = \min(I_\ell, R)$ singular values since it is guaranteed that $R_\ell \leq \min(I_\ell, R)$. Let \verb|mlsvd_tol| be a small positive tolerance value. For each truncated SVD of $\mathcal{T}_{(\ell)}$ with $\tilde{R}_\ell$ singular values, we want to verify its relative error with respect to $\mathcal{T}_{(\ell)}$. We start this process with $\tilde{R}_\ell = 1$ and keep increasing $\tilde{R}_\ell$ until the relative error is smaller than $\frac{1}{L} \cdot \verb|mlsvd_tol|$. Because of theorem~\ref{SVD-truncation} we don't need to actually construct these truncations and compute their respective relative errors. Denote by $\tilde{\mathcal{T}}_{(\ell)}$ the truncation with $\tilde{R}_\ell$ singular values, then we have that
			$$ \displaystyle \frac{\| \mathcal{T}_{(\ell)} - \tilde{\mathcal{T}}_{(\ell)} \|^2}{\| \mathcal{T}_{(\ell)} \|^2} = \frac{ \sum_{r=\tilde{R}_\ell+1}^{P_\ell} \big(\sigma_r^{(\ell)}\big)^2 }{\sum_{r=1}^{P_\ell} \big(\sigma_r^{(\ell)}\big)^2}.$$ 		
			
			After obtaining each $\tilde{R}_\ell$ we are able construct the truncated core tensor $\tilde{\mathcal{S}}$, which has shape $\tilde{R}_1 \times \ldots \times \tilde{R}_L$. It is not hard to see that the error between this truncation and the original core tensor is bounded as 
			
			$$\frac{ \| \mathcal{S} - \tilde{\mathcal{S}} \|^2}{\| \mathcal{S} \|^2} \leq \frac{ \sum_{i_1=\tilde{R}_1+1}^{P_1} \big( \sigma_{i_1}^{(1)} \big)^2 + \ldots + \sum_{i_L=\tilde{R}_L+1}^{P_L} \big( \sigma_{i_1}^{(1)} \big)^2 }{\| \mathcal{S} \|^ 2} < L \cdot \frac{1}{L} \verb|mlsvd_tol| = \verb|mlsvd_tol|,$$
where we use the fact that $\| \mathcal{S} \| = \| \mathcal{T} \| = \| \mathcal{T}_{(\ell)} \|$. Note that these computations only requires $\displaystyle \sum_{\ell=1}^L (R_\ell - \tilde{R}_\ell)$ flops. The following algorithm summarizes all we have described so far. We are using the notation $\Sigma^{(\ell)} = \{ \sigma_1^{(\ell)}, \ldots, \sigma_{P_\ell}^{(\ell)} \}$ for the set of the singular values of $\mathcal{T}_{(\ell)}$.\newpage
			
			\begin{algorithm}[Compression] \label{compression-alg}
				$ $\\
				\textbf{Input:} $\mathcal{T} \in \mathbb{R}^{I_1 \times \ldots \times I_L}, \ \verb|mlsvd_tol| > 0$\vspace{.4cm}\\
				$\big\{ \textbf{U}^{(1)}, \ldots, \textbf{U}^{(L)}, \Sigma^{(1)}, \ldots, \Sigma^{(L)}, \mathcal{S} \big\} \leftarrow $ \verb|MLSVD|$(\mathcal{T})$\\
				\verb|for | $\ell = 1 \ldots L$\\
					$\hspace{1cm}$ \verb|for | $i = 1 \ldots P_\ell$\\
						$\hspace{2cm}$ \verb|rel_error| $ = \frac{ \sum_{r=\tilde{R}_\ell+1}^{P_\ell} \big(\sigma_r^{(\ell)}\big)^2 }{\sum_{r=1}^{P_\ell} \big(\sigma_r^{(\ell)}\big)^2}$\\
						$\hspace{2cm}$ \verb|if rel_error| $ < \frac{1}{L}$ \verb|mlsvd_tol|\\
							$\hspace{3cm} \tilde{R}_\ell = i$\\ 
							$\hspace{3cm}$ \verb|break|\\
				\verb|for | $\ell=1 \ldots L$\\
					$\hspace{1cm} \tilde{\textbf{U}}^{(\ell)} \leftarrow$ \verb|truncation of | $\textbf{U}^{(\ell)}$ \verb| to have | $\tilde{R}_\ell$ \verb| columns|\\
				$\tilde{\mathcal{S}} \leftarrow$ \verb|truncation of | $\mathcal{S}$ \verb| to have shape | $\tilde{R}_1 \times \ldots \times \tilde{R}_L$\vspace{.4cm}\\
				\textbf{Output:} $\tilde{\textbf{U}^{(\ell)}} \in \mathbb{R}^{I_\ell \times \tilde{R}_\ell}$ for $\ell = 1 \ldots L$ and $\tilde{\mathcal{S}} \in \mathbb{R}^{\tilde{R}_1 \times \ldots \times \tilde{R}_L}$
			\end{algorithm}
		
			The function \verb|MLSVD| is described in~\ref{MLSVD-alg}, which has a cost of $\mathcal{O}\left( \displaystyle \log(P_\ell) \prod_{\ell'=1}^L I_{\ell'} \right)$ flops if we compute the truncated SVD with $P_\ell$ singular values for each unfolding. This is the dominant cost in this algorithm. With respect to the memory storage, each $\Sigma^{(\ell)}$ is stored as a vector of size $P_\ell$, it is necessary to store $I_\ell P_\ell$ floats for each $\textbf{U}^{(\ell)}$, and not more than $\displaystyle \prod_{\ell=1}^L P_\ell$ floats for $\mathcal{S}$. After computing the MLSVD, the worse case of the truncation stage has a cost of $\mathcal{O}\left( \displaystyle \sum_{\ell=1}^L P_\ell^2 \right)$ flops, which is cheap. Table 4.1 summarizes all main costs necessary to compress a tensor. Note that we are considering only the classic truncated MLSVD costs here. It is possible to obtain even lower costs with the sequentially truncated algorithm.\\
		
			\begin{table}
				\centering
		 		\begin{tabular}{|c|c|c|}
					\hline
					\textbf{Task} & \textbf{Memory} & \textbf{Computational cost}\\
					\hline
					Computing the MLSVD & $\displaystyle R \sum_{\ell=1}^L P_\ell + \prod_{\ell=1}^L P_\ell$ & $\mathcal{O}\left( \displaystyle \sum_{\ell=1}^L \log(P_\ell) \prod_{\ell'=1}^L I_{\ell'} \right)$\\
					\hline
					Truncating the MLSVD & $\displaystyle R \sum_{\ell=1}^L \tilde{R}_\ell + \prod_{\ell=1}^L \tilde{R}_\ell$ & $\mathcal{O}\left( \displaystyle \sum_{\ell=1}^L P_\ell^2 \right)$\\
					\hline
				\end{tabular}
				\caption{\footnotesize{Memory and computational costs - II.}} \bigskip
			\end{table}  
			
		\subsection{Initialization}
			As we mentioned before, the dGN method is not so sensitive to the initialization, but it still can benefits of good initializations. Furthermore, some initializations actually can lead to local minimum, and when this happens it is necessary to try again. In order to prevent local minimum it is interesting to have a method which generates initialization sufficiently close to the global minimum. Note that our objective is to compute a rank-$R$ CPD for $\mathcal{S} \in \mathbb{R}^{R_1 \times \ldots \times R_L}$, which is already truncated at this moment.  
			
			The first approach we consider is very usual in other implementations. The program may just generate random initial random matrices $\textbf{W}^{(\ell)} \in \mathbb{R}^{R_\ell \times R}$ with the entries draw from the standard normal distribution (that is, mean 0 and variance 1). With this we have an initial approximated CPD $\mathcal{S} = (\textbf{W}^{(1)}, \ldots, \textbf{W}^{(L)}) \cdot \mathcal{I}_{R \times \ldots \times R}$ which we use to start iterating. For general purposes this approach is good enough, but sometimes one can gain performance when the initialization is close enough to the objective tensor. With this in mind we propose another approach.
			
			Our second approach is a simple method based on the truncated MLSVD $\mathcal{T} \approx (\textbf{U}^{(1)}, \ldots, \textbf{U}^{(L)}) \cdot \mathcal{S}$. Remember the discussion about energy we had in chapter 2. The energy gets higher as the coordinates of $\mathcal{S}$ are close to $11 \ldots 1$, the very first entry of $\mathcal{S}$. In this case the energy of a coordinate can be seen as its magnitude, regardless it is positive or negative. Since we want a rank-$R$ CPD, one idea is to take $R$ entries of $\mathcal{S}$ with high energy and construct an approximated CPD from them. Since the energy is highly concentrated around $s_{11 \ldots 1}$, this approximated CPD may already be very close to the objective tensor. Now, choose $R$ multi-indexes $(j_1^{(1)}, \ldots, j_1^{(L)}), \ldots, (j_R^{(1)}, \ldots, j_R^{(L)}) \in R_1 \times \ldots \times R_L$ close to $(11 \ldots 1)$. These choices gives a rank-$R$ approximation
			
			$$\mathcal{S} \approx \sum_{r=1}^R s_{j_r^{(1)}, \ldots, j_r^{(L)}} \ \textbf{e}_{j_r^{(1)}}^{(1)} \otimes \ldots \otimes \textbf{e}_{j_r^{(L)}}^{(L)},$$
where each $\textbf{e}_j^{(\ell)}$ is the $j$-th basis canonical vector of $\mathbb{R}^{R_\ell}$. We denote this approximation by $\tilde{\mathcal{S}}$. We also denote $J = \{ (j_1^{(1)}, \ldots, j_1^{(L)}), \ldots, (j_R^{(1)}, \ldots, j_R^{(L)}) \}$ the set of multi-indexes we use to construct $\tilde{\mathcal{S}}$. Then the corresponding error is given by
			
			$$ \| \mathcal{S} - \tilde{\mathcal{S}} \|^2 = \sum_{(i_1, \ldots, i_L) \in R_1 \times \ldots \times R_L} s_{i_1 \ldots i_L}^2 - \sum_{(i_1, \ldots, i_L) \in J} s_{i_1 \ldots i_L}^2 = \sum_{(i_1, \ldots, i_L) \in R_1 \times \ldots \times R_L \backslash J} s_{i_1 \ldots i_L}^2.$$ 
			
			If our choice is adequate we can minimize this error. Furthermore, depending on the distribution energy of $\mathcal{S}$ this error may actually be very small. In my experience with tensors, that kind of initialization works well for several problems, and sometimes it works much better than the random initialization. The computational cost of the initialization is of $\mathcal{O}\left( \displaystyle R \sum_{\ell=1}^L R_\ell \right)$ flops, the size of the starting point. 
		
		\subsection{dGN} 
			Let $\mathcal{T} \approx (\textbf{U}^{(1)}, \ldots, \textbf{U}^{(L)}) \cdot \mathcal{S}$ be a truncated MLSVD for $\mathcal{T}$ which we assume to be close enough. The tensor obtained by the initialization will be denoted by $\mathcal{S}^{(0)}$. We denote by $\mathcal{S}^{(k)}$ the tensor obtained at the iteration $k$ of the dGN algorithm. Our goal is to iteratively produce successive approximations $\mathcal{S}^{(0)}, \mathcal{S}^{(1)}, \ldots, \mathcal{S}^{(k)}, \mathcal{S}^{(k+1)}, \ldots$ converging to $\mathcal{S}$. Associated with each $\mathcal{S}^{(k)}$ there will be the factor matrices $\textbf{W}^{(1,k)}, \ldots, \textbf{W}^{(L,k)}$ where $\textbf{W}^{(\ell, k)} \in \mathbb{R}^{R_\ell \times R}$ for each $\ell = 1 \ldots R$. As have been done before, we denote $\textbf{w}^{(k)} = \left[ vec(\textbf{W}^{(1,k)})^T, \ldots, vec(\textbf{W}^{(L,k)})^T \right]^T$.
			
			Let $K$ be the last iteration of the dGN. Then we expect to have
			
			$$\mathcal{S} \approx \underbrace{(\textbf{W}^{(1,K)}, \ldots, \textbf{W}^{(L,K)}) \cdot \mathcal{I}_{R \times R}}_{\mathcal{S}^{(K)}},$$
which leads to

			$$\mathcal{T} \approx \big( (\textbf{U}^{(1)})^\ast \textbf{W}^{(1,K)}, (\textbf{U}^{(L)})^\ast \textbf{W}^{(L,K)} \big) \cdot \mathcal{I}_{R \times \ldots \times R} = $$
			
			$$ = \big( \tilde{\textbf{W}}^{(1)}, \ldots, \tilde{\textbf{W}}^{(L)} \big) \cdot \mathcal{I}_{R \times \ldots \times R}.$$
			
			$$ = \sum_{r=1}^R \tilde{\textbf{W}_{:r}}^{(1)} \otimes \ldots \otimes \tilde{\textbf{W}_{:r}}^{(L)}.$$
			
			\subsubsection{Main parameters}
				First we need to establish the values of three important parameters: the initial damping parameter $\mu^{(0)}$, the maximum number of iterations \verb|maxiter| and the tolerance \verb|tol|. As mentioned in the ``damped Gauss-Newton'' section, the initial damping parameter is of the form $\mu^{(0)} = \tau \cdot \max_i a_{ii}$ (see \cite{madsen} to know more about this choice), where $\textbf{A} = \textbf{J}_f^T(\textbf{w}^{(0)}) \textbf{J}_f(\textbf{w}^{(0)})$, but instead of this we use the similar definition $\mu^{(0)} = \tau \cdot \mathbb{E}( |\mathcal{S}| )$, where $\mathbb{E}( |\mathcal{S}| )$ is the average of the entries of $\mathcal{S}$ in absolute value. This initial value was observed to perform well in practice. Also, since we want to reinforce regularization at the first iterations, we use $\tau = 1$. With regard to the maximum number of iterations, nonlinear least squares methods usually converges within a few hundreds iterations, whereas ALS and general unconstrained optimization methods need thousands iterations to converge. A reasonable choice in this case is \verb|maxiter|$= 200$. Usually there are many different tolerance parameters, one for each stopping condition. Although Tensor Fox does have the possibility to choose different values for each tolerance parameter, but for simplicity we consider that all tolerances are equal to \verb|tol|. At the moment we are using \verb|tol|$= 10^{-6}$. All these choices are the default values in TensorFox. They were obtained after several tests and experiments which will be showed soon. 
			
			\subsubsection{Computing the residual}
				At each iteration, the first task is to compute the residual function $f = (f_{11 \ldots 1}, \ldots, f_{R_1 R_2 \ldots R_L})$, where 
			
				$$f_{i_1 \ldots i_L}(\textbf{w}^{(k)}) = s_{i_1 \ldots i_L} - \sum_{r=1}^R w_{i_1 r}^{(1,k)} \ldots w_{i_L r}^{(L, k)}.$$
The cost of this computation is of $L-1$ flops, but since this has to be done for all the residuals, the total cost is of $(L-1) \displaystyle \prod_{\ell=1}^L R_\ell$ flops. 

			\subsubsection{Computing the gradient}
				Recall~\ref{gradF} from that $\nabla F(\textbf{w}) = \textbf{J}_f^T(\textbf{w}) \cdot f(\textbf{w})$, where $F$ is the error function defined in~\ref{error_function}. In section~\ref{damped-gauss-newton} we observed that $\textbf{J}_f(\textbf{w})$ is sparse with a certain structure we can exploit. Usually $\textbf{J}_f$ has $LR$ nonzero entries in each row, these are the values $-\displaystyle \prod_{\ell \neq \ell'} w_{i_\ell r'}^{(\ell)}$ for $\ell' = 1 \ldots L$ and $r' = 1 \ldots R$. Since $\textbf{J}_f(\textbf{w})$ has $L R \displaystyle \prod_{\ell=1}^L R_\ell$ nonzero entries, we can expect to compute $\textbf{J}_f^T(\textbf{w}) \cdot f(\textbf{w})$ with $L R \displaystyle \prod_{\ell=1}^L R_\ell$ operations, which is still much better than $R \displaystyle \sum_{\ell=1}^L R_\ell \displaystyle \prod_{\ell=1}^L R_\ell$ if we considered $\textbf{J}_f(\textbf{w})$ as a dense matrix. Consider the products $\ast$ and $\odot$ defined in~\ref{appenB} and remember that
				$$\Pi^{(\ell')} = \pi^{(1)} \ast \ldots \ast \pi^{(\ell'-1)} \ast \pi^{(\ell'+1)} \ast \ldots \ast \pi^{(L)},$$
where $\pi^{(\ell)} = \textbf{W}^{(\ell)^T} \textbf{W}^{(\ell)}$.	Then we have the following result from \cite{cp_opt}.	
		
				\begin{theorem}[T. G. Kolda, E. Acar, D. M. Dunlavy, 2011]
					The the partial derivatives of $F$ with respect to $\textbf{W}^{(\ell)}$ are given by
						$$\frac{\partial F}{\partial \textbf{W}^{(\ell)}}(\textbf{w}) = \textbf{W}^{(\ell)} \Pi^{(\ell)} - \mathcal{T}_{(\ell)} \left( \textbf{W}^{(L)} \odot \ldots \odot \textbf{W}^{(\ell+1)} \odot \textbf{W}^{(\ell-1)} \odot \ldots \odot \textbf{W}^{(1)} \right).$$
				\end{theorem}
		
				Since $\nabla F(\textbf{w}) = \left[ vec\left( \frac{\partial F}{\partial \textbf{W}^{(1)}}(\textbf{w}) \right)^T, \ldots, vec\left( \frac{\partial F}{\partial \textbf{W}^{(L)}}(\textbf{w}) \right)^T \right]^T$, from lemma~\ref{gradF} we conclude that
				$$\textbf{J}_f^T(\textbf{w}) \textbf{f}(\textbf{w}) = - \left[ vec\left( \frac{\partial F}{\partial \textbf{W}^{(1)}}(\textbf{w}) \right)^T, \ldots, vec\left( \frac{\partial F}{\partial \textbf{W}^{(L)}}(\textbf{w}) \right)^T \right]^T.$$		
		
				The dominant cost is the computation of $\mathcal{T}_{(\ell)} \left( \textbf{W}^{(L)} \odot \ldots \odot \textbf{W}^{(\ell+1)} \odot \textbf{W}^{(\ell-1)} \odot \ldots \odot \textbf{W}^{(1)} \right)$, which is $\mathcal{O}\left( R \displaystyle \prod_{\ell=1}^L I_\ell \right)$ flops. Since we have to perform this $L$ times, the total cost is of $\mathcal{O}\left( L R \displaystyle \prod_{\ell=1}^L I_\ell \right)$ flops. Efficient ways to deal with these Khatri-Rao products can be found in \cite{fast_kr}.

			\subsubsection{Conjugate gradient} \label{dGN-cg}
				The next stage is the computation of the step to take. Remember the description given at the beginning of section~\ref{gn_section}. At the point $\textbf{w}^{(k)}$ we want to take a step in direction $\textbf{x}$ such that the new point $\textbf{w}^{(k+1)} = \textbf{w}^{(k)} + \textbf{x}$ minimizes the residual at the neighborhood of $\textbf{w}^{(k)}$. This leads to a normal equations which we regularize (see~\ref{normal_eq}), obtaining the system
			
				$$( \textbf{A}^T \textbf{A} + \mu^{(k)} \textbf{D} ) \textbf{x} = \textbf{A}^T \textbf{b},$$
where $\textbf{D}$ is a diagonal $\left( R \displaystyle \sum_{\ell=1}^L R_\ell \right) \times \left( R \displaystyle \sum_{\ell=1}^L R_\ell \right)$ matrix, $\textbf{A} = \textbf{J}_f(\textbf{w}^{(k)})$, $\textbf{x} = \textbf{w} - \textbf{w}^{(k)}$, $\textbf{b} = -f(\textbf{w}^{(k)})$. The optimal solution $\textbf{x}_\ast$ gives the optimal point $\textbf{w}^{(k+1)}$ by setting $\textbf{w}^{(k+1)} = \textbf{w}^{(k)} + \textbf{x}_\ast$. 

				The regularization matrix $\textbf{D}$ depends on $\textbf{w}^{(k)}$ and is chosen to make $\textbf{A}^T \textbf{A} + \textbf{D}$ diagonally dominant (notice the absence of the damping parameter). As the iteration goes we will have $\mu^{(k)} \to 0$ and the effect of the regularization decreases. When close to the optimal point the only effect of $\textbf{D}$ is to guarantee that the system is well-posed. Still, the system converges to an ill-posed system as $\mu^{(k)} \to 0$ so we will be working with ill-conditioned system when close to the objective point. One way to mitigate the effects of this ill-conditioning is preconditioning the system with a Jacobi preconditioner as explained in appendix~\ref{appenA}. After the preconditioning we have the new system
			
				\begin{equation} \label{dGN-eq}			
					\textbf{M}^{-1/2}( \textbf{A}^T \textbf{A} + \mu^{(k)} \textbf{D} )\textbf{M}^{-1/2} \textbf{x} = \textbf{M}^{-1/2} \textbf{A}^T \textbf{b},
				\end{equation}
where $\textbf{M} = \text{diag}(a_{ii} + \mu^{(k)} d_{ii})$. Notice this equation is the same as the one preconditioned in~\ref{precond4}, but here we have $\textbf{A}^T \textbf{A} + \mu^{(k)} \textbf{D}$ instead of just $\textbf{A}$, and $\textbf{A}^T \textbf{b}$ instead of just $\textbf{b}$. To solve this system we rely on the conjugate gradient, see appendix~\ref{cg-alg}. The algorithm should be adapted to our current situation. 

				\begin{algorithm}[Conjugate gradient - dGN] \label{cg-dGN-alg}
					$ $\\
					\textbf{Input:} $\textbf{A}, \textbf{M}, \textbf{D}, \textbf{b}, \mu^{(k)}$ \vspace{.4cm}\\
					$\textbf{B} \leftarrow \textbf{M}^{-1/2}( \textbf{A}^T \textbf{A} + \mu^{(k)} \textbf{D}) \textbf{M}^{-1/2}$\\
					$\textbf{x}^{(0)} \leftarrow 0 \in \mathbb{R}^n$\\
					$\textbf{r}^{(0)} \leftarrow \textbf{M}^{-1/2} \textbf{A}^T \textbf{b}$\\
					$\textbf{p}^{(0)} \leftarrow \textbf{r}^{(0)}$\\
					\verb|for | $i = 1 \ldots$ \verb|cg_maxiter| \\
						   $\hspace{1cm} \textbf{z} \leftarrow \textbf{B} \cdot p^{(i-1)}$\\
						$\hspace{1cm} \alpha^{(i)} \leftarrow \frac{ \| \textbf{r}^{(i-1)} \|^2 }{ \| \textbf{p}^{(i-1)} \|_{\textbf{B}}^2 }$\\
						$\hspace{1cm} \textbf{x}^{(i)} \leftarrow \textbf{x}^{(i-1)} + \alpha^{(i)} \textbf{p}^{(i-1)}$\\
						$\hspace{1cm} \textbf{r}^{(i)} \leftarrow \textbf{r}^{(i-1)} - \alpha^{(i)} \textbf{z}$\\
						$\hspace{1cm} \varepsilon \leftarrow \| \textbf{r}^{(i)} \|^2$\\
						$\hspace{1cm} \beta^{(i)} \leftarrow \frac{ \| \textbf{r}^{(i)} \|^2 }{ \| \textbf{r}^{(i-1)} \|^2 }$\\
						$\hspace{1cm} \textbf{p}^{(i)} \leftarrow \textbf{r}^{(i)} + \beta^{(i)} \textbf{p}^{(i-1)}$\\
						$\hspace{1cm}$\verb|if | $\varepsilon < $ \verb|tol|\\
							$\hspace{2cm}$\verb|break| \vspace{.4cm}\\
					\textbf{Output:} $\textbf{x}^{(i')}$, where $i'$ is the last index of the iterations
				\end{algorithm}
			
				Now we comment a few things about the algorithm, from top to bottom. The first thing we should remark is that the matrix $\textbf{B}$ is never computed explicitly. As already mentioned, this is a matrix-free algorithm because of theorem~\ref{Hv}. The vector $\textbf{r}^{(0)}$ demands some computational effort as observed before. Multiplying $\textbf{A}^T \textbf{b}$ by $\textbf{M}^{-1/2}$ amounts to just $R \displaystyle \sum_{\ell=1}^L R_\ell$ operations since $\textbf{M}^{-1/2}$ is diagonal. The result is a vector of size $R \displaystyle\sum_{\ell=1}^L R_\ell$. 
			
				We call the conjugate gradient algorithm by CG for short. The parameter \verb|cg_maxiter| is the maximum number of iterations permitted. When in situations like this, where the CG is used as a step for other iterative algorithm, it is usual to set \verb|cg_maxiter| to a low value. Sometimes \verb|cg_maxiter| $= 10$ is already enough to produce a reasonable step for the ``bigger'' algorithm (in this case, the dGN). It should be pointed that this is very case dependent and just there isn't a universal rule for how to set this parameter. We tested it for several fixed values and none was satisfactory. Low values made the dGN converge to local minimum most of the time, higher values improved the convergence but were very costly. It was noted a better performance when we increased \verb|cg_maxiter| a little for each iteration of the dGN. However this approach forced the program to make too much CG iterations at the end of the dGN, which was wasteful. Several attempts to solve this issue were tried. The problem in all of them was the fact that they all were deterministic. As soon as the value \verb|cg_maxiter| was set to be a random integer everything worked much better. At iteration $k$ of the dGN, we define \verb|cg_maxiter| as being a random integer draw from the uniform distribution in the interval $\left[ 1 + \lceil k^{0.4} \rceil, 2 + \lceil k^{0.9} \rceil \right]$. This strategy to obtain \verb|cg_maxiter| (instead of just fixing it) is the result of many experiments. 
				
				\begin{remark} \label{unusual_steps}
					Allowing the number of iterations to change dynamically allows the algorithm to make unusual steps sometimes, which proved to be a successful way to avoid local minima. This is a very specific format of interval so you can imagine how much of ``mathematical alchemy'' was necessary to create such solution.     
				\end{remark} 
			
				The product $\textbf{B} \cdot \textbf{p}^{(i-1)}$ is to be computed in three steps as showed below.			
				\begin{flalign*}
					& \textbf{z} \leftarrow \textbf{M}^{-1/2} \cdot \textbf{p}^{(i-1)}\\
					& \textbf{z} \leftarrow \textbf{A}^T \textbf{A} \cdot \textbf{z} + \mu^{(k)} \textbf{z}\\
					& \textbf{z} \leftarrow \textbf{M}^{-1/2} \cdot \textbf{z}
				\end{flalign*}			
	
				However before starting these computations we need to construct all the matrices $\Pi^{(\ell', \ell'')}$ and $\Pi^{(\ell')}$. We remark that the computation of the gradient $\textbf{A}^T \cdot \textbf{b}$ and these matrices are performed before the CG loop, so their costs are accounted only once. We already observed that, since $\textbf{M}^{-1/2}$ is diagonal, the first product can be computed with $R \displaystyle \sum_{\ell=1}^L R_\ell$ operations. The same goes for the multiplication by the scalar $\mu^{(k)}$. Theorem~\ref{Hv} and the discussion after it showed that $\textbf{A}^T \textbf{A} \cdot \textbf{z}$ can be computed with $\mathcal{O}\left( L R^2 \displaystyle \sum_{\ell=1}^L R_\ell \right)$ operations. Finally we have to make more $R \displaystyle \sum_{\ell=1}^L R_\ell$ operations in the last step $\textbf{z} \leftarrow \textbf{M}^{-1/2} \cdot \textbf{z}$.  The overall cost to obtain $\textbf{B} \cdot \textbf{p}^{(i-1)}$ is of 				
				$$\mathcal{O}\left( (3 R + L R^2) \displaystyle \sum_{\ell=1}^L R_\ell \right)$$ 
flops. All the remaining lines of the algorithm together have a cost of $\mathcal{O}\left( 6R \displaystyle \sum_{\ell=1}^L R_\ell \right)$ flops. Putting everything together we can see that the cost of the CG algorithm is of 
				$$\mathcal{O}\left( L R \displaystyle \prod_{\ell=1}^L R_\ell + R^2 \left( L + \displaystyle \sum_{\ell=1}^L R_\ell \right) + i_{CG} (9 R + L R^2) \displaystyle \sum_{\ell=1}^L R_\ell \right),\vspace{-0.5cm}$$
				$$\hspace{0.2cm}\underbrace{ \hspace{1.5cm} }_{ \text{gradient} } \hspace{0.5cm} 
				\underbrace{ \hspace{3cm} }_{\Pi^{(\ell', \ell'')} \text{ and } \Pi^{(\ell')}} \hspace{0.5cm} 
				\underbrace{ \hspace{4cm} }_{\text{CG loop}}$$
where $i_{CG}$ is the number of iterations of the CG. We know that $i_{CG}$ is a random integer draw from the interval $\left[ 1 + \lceil k^{0.4} \rceil, 2 + \lceil k^{0.9} \rceil \right]$, where $k$ is the current iteration of the dGN. Let's consider the worst case $k = \verb|maxiter| = 200$ just to have an idea. In this, case the expected $i_{CG}$ is
				
				$$\mathbb{E}(i_{CG}) = \left\lfloor \frac{2 + \lceil 200^{0.9} \rceil - 1 - \lceil 200^{0.4} \rceil}{2} \right\rfloor = 55,$$
which is a reasonable value.

			\subsubsection{Updates}	 \label{updates}
				After the CG is performed we have a point $\textbf{x}_\ast$, which is an approximated solution of~\ref{dGN-eq}. In turn, this equation implies that $\textbf{w}^{(k)} + \textbf{x}_\ast$ minimizes the residual at the neighborhood of $\textbf{w}^{(k)}$ (see the discussion about~\ref{linear_approx}). For this reason, we set $\textbf{w}^{(k+1)} = \textbf{w}^{(k)} + \textbf{x}_\ast$. This is the first update of the dGN algorithm.
			
				Depending on the size of the dimensions $R_1, \ldots, R_L$, this second step is likely to be the most costly part of all routines inside dGN. We are talking about the evaluation of the error function $F$. Let $\mathcal{S}^{(k)} = \displaystyle \sum_{r=1}^R \textbf{w}_r^{(1,k)} \otimes \ldots \otimes \textbf{w}_r^{(L,k)}$ be the approximating tensor computed at the $k$-th iteration. To compute the error at iteration $k$ we must compute
			
				$$F(\textbf{w}^{(k)}) = \frac{1}{2}\left\| \mathcal{S} - \mathcal{S}^{k} \right\|^2 = \frac{1}{2}\left\| \mathcal{S} - \sum_{r=1}^R \textbf{w}_r^{(1,k)} \otimes \ldots \otimes \textbf{w}_r^{(L,k)} \right\|^2 = $$
			
				$$ = \frac{1}{2} \sum_{i_1=1}^{R_1} \ldots \sum_{i_L=1}^{R_L} \left( s_{i_1 \ldots i_L} - \sum_{r=1}^R w_{i_1 r}^{(1, k)} \cdot \ldots \cdot w_{i_L r}^{(L, k)} \right)^2.$$     
			
				we don't take in account the cost to compute this error since it is just a matter of summing the squares of the residual whose cost is already considered. Before making this evaluation is it interesting to ``normalize'' the factors, that is, scale them so we have $\| \textbf{w}_r^{(1,k)} \| = \| \textbf{w}_r^{(2,k)} \| = \ldots = \| \textbf{w}_r^{(L,k)} \|$ for each $r = 1 \ldots R$. This is always possible, has a low cost of $R \displaystyle \sum_{\ell=1}^L R_\ell$ flops and can improve the conditioning of the problem, making it easier for the next CG iterations to find a good solution. When the factors are scaled this way they are called \emph{norm-balanced}. In section 5 of \cite{nick} the conditioning of norm-balanced representations of tensors is discussed at some detail. In particular, the lower bound of the condition number is minimized when we use balanced-norm representations.  
			
				The last update is the damping parameter update. Let $\mu^{(k)}$ be the damping parameter at iteration $k$ of dGN. Previously, in the ``Damped Gauss-Newton'' section, we saw that this update depends on the gain ratio
			
				$$g = \frac{ \| \mathcal{S} - \tilde{\mathcal{S}}^{(k-1)} \|^2 - \| \mathcal{S} - \tilde{\mathcal{S}}^{(k)} \|^2 }{ \| \mathcal{S} - \tilde{\mathcal{S}}^{(k-1)} \|^2 - \| \tilde{f}(\textbf{w}^{(k)}) \|^2 }.$$
A large value of $g$ indicates that $\| \tilde{f}(\textbf{w}^{(k)})\|^2$ is a good approximation to $\|\mathcal{S} - \tilde{\mathcal{S}}^{(k)}\|^2$, in other words, the linear model $\tilde{f}$ is making good predictions about the actual errors. A small or negative $g$ indicates that this approximation is poor (the reason for that is indicated in section~\ref{damped-gauss-newton}). When the approximation is good we decrease the damping parameter so the iterations becomes more like a Gauss-Newton iteration, which converges rapidly. Otherwise we increase the damping parameter, which adds regularization to the model and make it more alike the gradient descent method. 

				An update strategy widely used is the following.			
				\begin{flalign*}
					& \texttt{if } g < 0.25\\
					& \hspace{.8cm} \mu \leftarrow 2 \mu\\
					& \texttt{else if } g > 0.75\\
					& \hspace{.8cm} \mu \leftarrow \mu/3
				\end{flalign*}
			
				This strategy was originally proposed by Marquardt in \cite{marquardt}. This update strategy is not so sensitive to minor changes in the thresholds $0.25, 0.75$ or the update factors $2, 1/3$. Another strategy frequently used (which was also used to compute CPDs in \cite{cichoki2013}) is the following.
				\begin{flalign*}
					& \texttt{if } g > 0\\
					& \hspace{.8cm} \mu \leftarrow \mu \cdot \max\{ 1/3, 1 - (2g -1)^3 \}\\
					& \hspace{.8cm} \nu \leftarrow 2\\
					& \texttt{else} \\
					& \hspace{.8cm} \mu \leftarrow \mu \cdot \nu\\
					& \hspace{.8cm} \nu \leftarrow 2 \nu
				\end{flalign*}
			 			   
				In general this strategy outperforms the previous one, as demonstrated in \cite{nielsen}. The value $\nu$ usually is initiated to $\nu = 2$, but minor changes doesn't affect the performance significantly. We tested the second strategy in our problem for a vast of different updates and none was very effective. The first one seemed to perform better for some specific choices of thresholds and factors. In the end, we decided to use the following update strategy.			
				\begin{flalign*}
					& \texttt{if } g < 0.25\\
					& \hspace{.8cm} \mu \leftarrow 3 \mu\\
					& \texttt{else if } g > 0.75\\
					& \hspace{.8cm} \mu \leftarrow \mu/3
				\end{flalign*}
			
				These values were obtained experimentally as a result of several tests over several distinct tensors. Still, there is a reasoning to explain why these values works well. As we already know, the maximum number of iterations of the CG will be very low at the first iterations of the dGN (\verb|cg_maxiter| can be even 2 or 3). This means the solution found by the CG has a large residual, which means $\| \tilde{f}(\textbf{w}^{(k)})\|$ will be large when $k$ is small. Even in this situation the dGN is likely to decrease the errors monotonically, so we will have $\| \mathcal{S} - \mathcal{S}^{(k)} \| \leq \| \mathcal{S} - \mathcal{S}^{(k-1)} \| \leq \| \tilde{f}(\textbf{w}^{(k)})\|$ regardless the number of iterations performed by the CG. This implies that $g < 0$. If we apply the original Marquardt strategy, then the damping parameter will increase while also increasing the regularization of the problem, which is already very regularized at the first iterations. This excessive regularization almost always leads to solutions which are some local minimum of the problem. Even if we use some good initialization (which clearly means that decreasing the damping parameter is the right move to do), the fact that \verb|cg_maxiter| is low at the beginning will lead to more regularization, which leads to local minima. As we keep doing dGN iterations, it is expected to have $\| \tilde{f}(\textbf{w}^{(k)})\| \ll \| \mathcal{S} - \mathcal{S}^{(k)} \| \lessapprox \| \mathcal{S} - \mathcal{S}^{(k-1)} \|$ so $g$ is positive and close to 0. To reinforce the fact that we want to decrease the damping parameter, it is crucial to choose thresholds which favors this. Therefore, most of the time the algorithm will decrease the damping parameter and only sometimes it will increase it. 
			  	
			\subsubsection{Stopping conditions} \label{stop}		
				As in any iterative algorithm, good stopping conditions are crucial to improve performance. The conditions we used to stop iterating are based in other implementations, but they were tested and refined to lead to the best performance possible. We present our four stopping conditions below and in the same order they are implemented. The program stops at iteration $k$ if
			
				\begin{enumerate}
					\item \textbf{Relative error}: \ \ \ \
						$\displaystyle \frac{ \| \mathcal{S} - \mathcal{S}^{(k)} \| }{\| \mathcal{S} \|} < $ \verb|tol|\\
					\item \textbf{Step size}: \ \ \ \
						$\| \textbf{w}^{(k-1)} - \textbf{w}^{(k)} \| < $ \verb|tol|\\
					\item \textbf{Relative error improvement}: \ \ \ \
						$\displaystyle \left| \frac{ \| \mathcal{S} - \mathcal{S}^{(k-1)} \| }{\| \mathcal{S} \|} - \frac{ \| \mathcal{S} - \mathcal{S}^{(k)} \| }{\| \mathcal{S} \|} \right| < $ \verb|tol|\\
					\item \textbf{Gradient norm}: \ \ \ \
						$\| \nabla F(\textbf{w}^{(k)}) \|_\infty <  \verb|tol|$\\
					\item \textbf{Average error:} \ \ \ \
					 	$$\displaystyle \frac{1}{c} \sum_{k=k_0 - 2c}^{k_0-c}\frac{ \| \mathcal{S} - \mathcal{S}^{(k)} \| }{\| \mathcal{S} \|} - \frac{1}{c} \sum_{k=k_0-c}^{k_0}\frac{ \| \mathcal{S} - \mathcal{S}^{(k)} \| }{\| \mathcal{S} \|} <  \texttt{tol}$$
where $k_0 > 2c$ and $\displaystyle c = 1 + \left\lceil \frac{\texttt{maxiter}}{10} \right\rceil$\\
					 \item \textbf{Average improvement:} \ \ \ \
					 	$$\displaystyle \frac{1}{c} \sum_{k=k_0-c}^{k_0} \left| \frac{ \| \mathcal{S} - \mathcal{S}^{(k-1)} \| }{\| \mathcal{S} \|} - \frac{ \| \mathcal{S} - \mathcal{S}^{(k)} \| }{\| \mathcal{S} \|} \right| <  10^{-3} \cdot \frac{1}{c} \sum_{k=k_0-c}^{k_0}\frac{ \| \mathcal{S} - \mathcal{S}^{(k)} \| }{\| \mathcal{S} \|}$$
where $k_0 > 2c$ and $\displaystyle c = 1 + \left\lceil \frac{\texttt{maxiter}}{10} \right\rceil$\\
					\item \textbf{Divergence:} \ \ \ \
						$\displaystyle \frac{ \| \mathcal{S} - \mathcal{S}^{(k)} \| }{\| \mathcal{S} \|} > \frac{ \max\{ 1, \ \|\mathcal{T}\|^2 \} }{10^{-16} + \texttt{tol}}$
				\end{enumerate}
				
\begin{figure}[H]
					\centering
					\includegraphics[scale=.5]{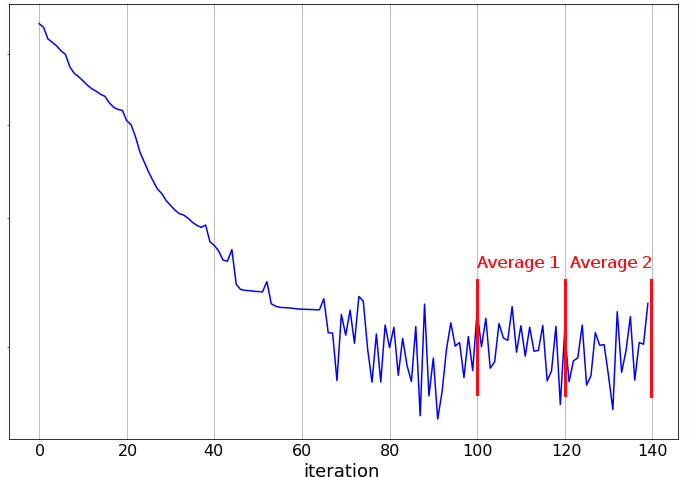}
					\caption{\footnotesize{The blue line represents the evolution of the relative error in a CPD computation. The program can stop because the average 2 is bigger than average 1.}}
					\label{Stopping condition: average error}
				\end{figure}

				The first three conditions are standard and will be present in almost any solver in one form or another. The gradient condition can be found in Tensor Toolbox. Since the algorithm converges to a critical point, this is a natural stopping condition to include. The purpose of the last condition is clear, and the term $10^{-16}$ in the denominator is there to prevent division by zero in the case of \verb|tol| $= 0$. Conditions 5 and 6 deserves a little more explanation. First of all, we remark that the constants $c$ and $10^{-3}$ were obtained empirically and there is nothing special about them. Secondly, these two conditions are only verified at each $c$ iterations so we are always comparing the average of a batch with the next one. Condition 5 avoids errors oscillating too much time without any overall decreasing, like illustrated in figure~\ref{Stopping condition: average error}. 
				
\begin{figure}[H]
					\centering
					\includegraphics[scale=.5]{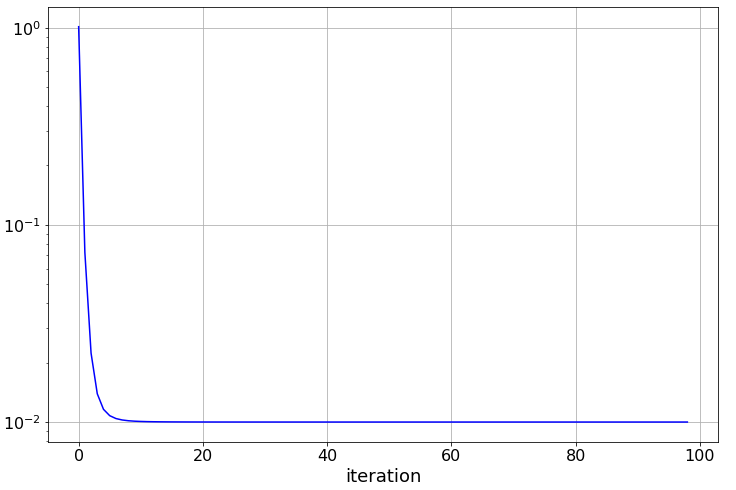}
					\caption{\footnotesize{Note that the program could have stopped much earlier. Even if the errors are strictly decreasing, the additional accuracy is irrelevant compared to the final error.}}
					\label{Stopping condition: average improvement}
				\end{figure}

				Condition 6 avoids too long periods of negligible improvements. For instance, if the error is of order $\mathcal{O}(10^{-2})$, it is not necessary the program to waste time making hundreds iterations of improvements of order $\mathcal{O}(10^{-6})$. Figure~\ref{Stopping condition: average improvement} illustrates this kind of situation. 
				
				The cost of the step size is of $R \displaystyle \sum_{\ell=1}^L R_\ell$ flops. Remember that $\nabla F(\textbf{w}^{(k)}) = \textbf{J}_f^T \cdot f(\textbf{w}^{(k)})$ (lemma~\ref{gradF}). With the notations used in the CG algorithm we have that $\nabla F(\textbf{w}^{(k)}) = - \textbf{A}^T \cdot \textbf{b}$. Since this vector had to be computed in the CG algorithm, we can just store it be used for the stopping condition. The computation of its norm amounts to $R \displaystyle \sum_{\ell=1}^L R_\ell$ flops. The other stopping conditions are cheaper to compute so we don't count them.
			
			\subsubsection{Overall cost of dGN}
				We already have done all the complexity analysis of each algorithm and each small routine inside the dGN algorithm. Here we just put all them together in order to facilitate the presentation. We assign some memory costs with the symbol ``$-$'' when the respective cost is irrelevant.
			
				Overall, the memory cost of one dGN iteration is of 
			
				$$\displaystyle R \sum_{\ell=1}^L R_\ell + R^2 (L + L^2) + R \displaystyle \sum_{\ell=1}^L R_\ell = \displaystyle 2R \sum_{\ell=1}^L R_\ell + R^2 (L + L^2) \text{ floats,} \vspace{-0.3cm}$$
				$$\hspace{-6cm} \underbrace{ \hspace{1.5cm} }_{ \text{gradient} } \hspace{0.3cm}
				\underbrace{ \hspace{2.2cm} }_{ \Pi^{(\ell', \ell'')} \text{ and } \Pi^{(\ell')} } \hspace{0.3cm}
				\underbrace{ \hspace{1.5cm} }_{ \text{vector result} }$$
while the computational cost is of 

				$$\mathcal{O}\left( (L-1) \displaystyle \prod_{\ell=1}^L R_\ell + L R \displaystyle \prod_{\ell=1}^L R_\ell + R^2 \left( L + \displaystyle \sum_{\ell=1}^L R_\ell \right) + i_{CG} (9 R + L R^2) \displaystyle \sum_{\ell=1}^L R_\ell + 3 R \displaystyle \sum_{\ell=1}^L R_\ell \right) =\vspace{-0.4cm}$$
				$$\hspace{0.9cm}
				\underbrace{ \hspace{2.2cm} }_{\text{residual}} \hspace{0.5cm} 
				\underbrace{ \hspace{1.8cm} }_{ \text{gradient} } \hspace{0.4cm} 
				\underbrace{ \hspace{3cm} }_{ \Pi^{(\ell', \ell'')} \text{ and } \Pi^{(\ell')} } \hspace{0.5cm} 
				\underbrace{ \hspace{3.9cm} }_{ \text{CG loop} } \hspace{0.5cm} 
				\underbrace{ \hspace{1.6cm} }_{ \text{stop conditions} }$$	
			
				$$ = \mathcal{O}\left( (LR+L-1) \displaystyle \prod_{\ell=1}^L R_\ell + R^2 \left( L + \displaystyle \sum_{\ell=1}^L R_\ell \right) + 3 R \displaystyle \sum_{\ell=1}^L R_\ell + i_{CG} (9 R + L R^2) \displaystyle \sum_{\ell=1}^L R_\ell \right) \text{ flops.}$$	 
				
				\begin{table}
					\centering
					\begin{tabular}{|c|c|c|}
						\hline
						\textbf{Task} & \textbf{Memory} & \textbf{Computational cost}\\
						\hline
						Computing the residual & $\displaystyle \prod_{\ell=1}^L R_\ell$ & $\mathcal{O}\left( (L-1) \displaystyle \prod_{\ell=1}^L R_\ell \right)$\\
						\hline
						Computing the gradient & $R \displaystyle \sum_{\ell=1}^L R_\ell$ & $\mathcal{O}\left( L R \displaystyle \prod_{\ell=1}^L R_\ell \right)$\\
						\hline
						Computing $\Pi^{(\ell', \ell'')}$ and $\Pi^{(\ell')}$ & $R^2 (L + L^2)$ & $\mathcal{O}\left( R^2 \left( L + \displaystyle \sum_{\ell=1}^L R_\ell \right) \right)$\\
						\hline
						CG method & $R \displaystyle \sum_{\ell=1}^L R_\ell $ & $\mathcal{O}\left( i_{CG} (9 R + L R^2) \displaystyle \sum_{\ell=1}^L R_\ell \right)$\\
						\hline
						Normalizing the factors & $-$ & $\mathcal{O}\left( R \displaystyle \sum_{\ell=1}^L R_\ell \right)$\\
						\hline
						Damping parameter update & $-$ & $\mathcal{O}(1)$\\
						\hline
						Stopping conditions & $-$ & $\mathcal{O}\left( 3R \displaystyle \sum_{\ell=1}^L R_\ell \right)$\\
						\hline
					\end{tabular}
					\caption{\footnotesize{Memory and computational costs - III.}} \bigskip
				\end{table}	
				
			\subsubsection{Comparison to other algorithms} \label{comparison}
				As mentioned in chapter 0, the solvers Tensorlab and Tensor Box also implements the Gauss-Newton method. It is of interest to compare the implementations and their complexities since they are based on the same method. To simplify out analysis, assume that the compression stage is already performed and limit the discussion here to cubic tensors of shape $R \times R \times \ldots \times R$.   
				
				Tensorlab is a tensor package with several algorithms to compute a more general decomposition, namely the \emph{Block Term Decomposition} (BTD). They propose to work with the more general approximate Hessian, which have structure to be exploited, just as in this work, see theorem 4.5 \cite{tensorlab2}. Since their formulation is more general, they end with more equations to perform matrix-vector multiplication, which then can be simplified in the particular case of the CPD. Preconditioning is employed too, but instead of a diagonal preconditioner they use a block diagonal preconditioner (also called a \emph{block Jacobi preconditioner}). This approach may improve the conditioning substantially, but it comes with the cost of having to solve a linear system at each CG iteration. The cost to compute the error at each iteration is of $\mathcal{O}(2R^{L+1})$ flops and the cost to compute the cogradient is of $\mathcal{O}(2LR^{L+1})$ flops, see appendix of \cite{tensorlab2} for more details. The faster algorithm is the inexact Gauss-Newton, which costs $\mathcal{O}\left(c_1\left( \frac{5}{2} L^2R^2 + 8LR^2 + \frac{1}{3}LR^3 \right)\right)$ flops to solve the inverse problem of each Gauss-Newton iteration. The value $c_1$ is the number of CG iterations performed, which is set to $15$ by default. For more details we refer to the appendix of the same article.
				
				Tensor Box starts decomposing the approximated Hessian as $\textbf{H} = \textbf{G} + \textbf{Z} \textbf{K} \textbf{Z}^T$, where all these matrices are sparse with some structure, see theorem 4.2 of \cite{cichoki2013}. This structure is exploited so then can write the inverse of the damped approximated as $(\textbf{H} + \mu \textbf{I})^{-1} = \tilde{\textbf{G}}_\mu - \textbf{L}_\mu \textbf{B}_\mu \textbf{L}_\mu^T$, where all these matrices have some precise structure which can be used to accelerate the iterations. Their algorithm depends on the construction of several intermediate matrices, so we will only mention their associate costs, but the reader more interested is encouraged to read section 4.3 of \cite{cichoki2013}. Building matrix $\textbf{K}$ costs $\mathcal{O}(L^3R^2)$ flops. Inverting all matrices $\Gamma_\mu^{(n)}$ costs $\mathcal{O}(LR^3)$ flops. The main cost of their algorithm is the computation of damped factors, which amounts to $\mathcal{O}(LR^{L+1})$ flops. These factor are used to construct a inverse linear problem with size $LR^2 \times LR^2$, which is solved with $\mathcal{O}(L R^3 + L^3R^6)$ flops. 
				
				For each solver, the costs to perform a single iteration is summarized below.
				
				$$\begin{array}{cc}
					\hspace{-0.98cm}\textbf{Tensor Fox:} & \mathcal{O}\left( LR^{L+1} + (L-1)R^L + 5 LR^2 + LR^3 + c_2 (9 LR^2 + L^2 R^3) \right) \text{ flops}\\
				    \hspace{-1.2cm}\textbf{Tensorlab:} & \hspace{-0.8cm} \mathcal{O}\left( 2R^{L+1} + 2LR^{L+1} + c_1\left( \frac{5}{2} L^2R^2 + 8LR^2 + \frac{1}{3}LR^3 \right) \right) \text{ flops}\\
				    \hspace{-.8cm}\textbf{Tensor Box:} & \hspace{-2.7cm}\mathcal{O}\left( L^3R^2 + LR^3 + LR^{L+1} + LR^3 + L^3R^6 \right) \text{ flops}
				\end{array}$$	
				
				At first, the analysis is not clear because $c_2$ is stochastic (it is the maximum number of CG iterations). Assuming that all $200$ dGN iterations of Tensor Fox will be performed, the average number of CG iterations is $c_2 = 28$. The term $R^{L+1}$ adds more cost to Tensorlab as the order increases. We note that Tensor Box is slower for low $L$ and faster for bigger $L$ (more precisely, for $L \geq 8$), with Tensor Fox being the second faster. For $L=3$ and small $R$, we note that Tensorlab is faster than the others, but Tensor Fox starts to be the faster as $R$ increases (around $R = 50$). For $L = 4, 5, 6, 7$ Tensor Fox is the faster algorithm for any $R \geq 10$. We remark that these complexity analysis depends on the constants associated to each cost. 					 
			
		\subsection{Uncompression}	
			After we have computed a CPD $(\textbf{W}^{(1)}, \ldots, \textbf{W}^{(L)})$ for $\mathcal{S}$, we expect that 
			
			$$\mathcal{S} \approx (\textbf{W}^{(1)}, \ldots, \textbf{W}^{(L)}) \cdot \mathcal{I}_{R \times \ldots \times R}.$$ On the other hand we know that $\mathcal{T} \approx (\textbf{U}^{(1)}, \ldots, \textbf{U}^{(L)}) \cdot \mathcal{S}$, hence 
			
			$$\mathcal{T} \approx (\tilde{\textbf{W}}^{(1)}, \ldots, \tilde{\textbf{W}}^{(L)}) \cdot \mathcal{I}_{R \times \ldots \times R},$$ 
where $\tilde{\textbf{W}}^{(\ell)} = \textbf{U}^{(\ell)} \textbf{W}^{(\ell)}$ for each $\ell = 1 \ldots R$. With this we can write
			
			$$\mathcal{T} \approx \sum_{r=1}^R \tilde{\textbf{W}}_{: r}^{(1)} \otimes \ldots \otimes \tilde{\textbf{W}}_{: r}^{(L)},$$
an approximated CPD for $\mathcal{T}$. 

			Note that $\tilde{\textbf{W}}^{(\ell)} \in \mathbb{R}^{R_\ell \times R}$, while $\textbf{W}^{(\ell)} \in \mathbb{R}^{I_\ell \times R}$. The process of transforming smaller dimensional factor matrices $\tilde{\textbf{W}}^{(\ell)}$ into the bigger factor matrices $\textbf{W}^{(\ell)}$ is what we call \emph{uncompression}. This is the last stage of all, after that we have a CPD for $\mathcal{T}$. The cost to produce each $\tilde{\textbf{W}}^{(\ell)}$ is the cost to multiply a $I_\ell \times R_\ell$ matrix by a $R_\ell \times R$ matrix, which has a cost of $\mathcal{O}(I_\ell R_\ell R)$ flops. Therefore, the cost to produce the uncompressed CPD is of $\mathcal{O}\left( R \displaystyle \sum_{\ell=1}^L I_\ell R_\ell \right)$ flops.
			
			We finish this section with a succinct summary of all costs to produce a CPD (tables 4.3 and 4.4) in Tensor Fox. Just as we denoted by $i_{CG}$ the number of iterations of the CG algorithm, we denote by $i_{dGN}$ the number of iterations of the dGN algorithm.\\
			
			\begin{table} 
				\centering
				\begin{tabular}{|c|c|}
					\hline
					\textbf{Task} & \textbf{Memory}\\
					\hline
					Compression & $(R+1) \displaystyle \sum_{\ell=1}^L I_\ell + \displaystyle \prod_{\ell=1}^L I_\ell + \displaystyle \prod_{\ell=1}^L R_\ell$ \\
					\hline
					Initialization & $\displaystyle R \sum_{\ell=1}^L R_\ell$ \\
					\hline
					dGN & $\displaystyle \prod_{\ell=1}^L R_\ell + R^2 L^2 + R \displaystyle \sum_{\ell=1}^L R_\ell$ \\
					\hline
					Uncompression & $R \displaystyle \sum_{\ell=1}^L I_\ell$ \\
					\hline
				\end{tabular}
				\label{total-memory}
				\caption{\footnotesize{Memory costs - IV.}} \bigskip
			\end{table} 
				
			\begin{table}
				\centering
				\begin{tabular}{|c|c|}
					\hline
					\textbf{Task} & \textbf{Computational cost}\\
					\hline
					Compression & $\mathcal{O}\left( \displaystyle \left( 1 + \prod_{\ell=1}^L I_\ell \right) \displaystyle \sum_{\ell=1}^L I_\ell + \displaystyle \sum_{\ell=1}^L I_\ell^3 \right)$\\
					\hline
					Initialization & $\mathcal{O}\left( \displaystyle R \sum_{\ell=1}^L R_\ell \right)$\\
					\hline
					dGN & \footnotesize $\mathcal{O}\left( i_{dGN} \left( (LR+R+L-1) \displaystyle \prod_{\ell=1}^L R_\ell + 3 R \displaystyle \sum_{\ell=1}^L R_\ell + R^2 \displaystyle \sum_{\ell=1}^L R_\ell + i_{CG} (8 R + L R^2) \displaystyle \sum_{\ell=1}^L R_\ell \right) \right)$ \normalsize\\
					\hline
					Uncompression & $\mathcal{O}\left( R \displaystyle \sum_{\ell=1}^L I_\ell R_\ell \right)$\\
					\hline
				\end{tabular}
				\label{total-cost}
				\caption{\footnotesize{Computational costs - IV.}} \bigskip
			\end{table}	

	\section{Warming up}
		Now we start to make computational experiments. Let's start with a simple example in order to show some computational aspects we are interested in. For this suppose a multivariate function $\varphi: \mathbb{R}^L \to \mathbb{R}$ as described in example~\ref{func-grid}. To make everything more concrete let's consider the highly nonlinear function		
		$$\varphi(x,y,z) = \cos(1-x) \sin(1+y^2) - \sin(1+x^2) e^{x^2+y^2+z^2} - \cos(x) \ln(1+z).$$ 
Now consider some uniform grids in the interval $[0,1]$, with six samples for $x$, five samples for $y$ and four samples for $z$. With these points we form the tensor $\mathcal{T} \in \mathbb{R}^{6 \times 5 \times 4}$ defined as $t_{ijk} = \varphi(x_i, y_j, z_k)$, where $x_i = 0, 0.2, 0.4, 0.6, 0.8, 1$, $y_j = 0, 0.25, 0.5, 0.75, 1$ and $z_k = 0, 0.33333333, 0.66666667, 1$. We can write $\mathcal{T} = \displaystyle \sum_{\ell=1}^3 \varphi_1^{(\ell)}(x) \varphi_2^{(\ell)}(y) \varphi_3^{(\ell)}(z)$, where
		\begin{flalign*}			
			& \varphi_1^{(1)}(x) = \cos(1-x),\\
			& \varphi_2^{(1)}(y) = \sin(1+y^2),\\
			& \varphi_3^{(1)}(z) = 1,\\
			& \varphi_1^{(2)}(x) = -\sin(1+x^2) e^{x^2},\\
			& \varphi_2^{(2)}(y) = e^{y^2},\\
			& \varphi_3^{(2)}(z) = e^{z^2},\\
			& \varphi_1^{(3)}(x) = -\cos(x),\\
			& \varphi_2^{(3)}(y) = 1,\\
			& \varphi_3^{(3)}(z) = \ln(1+z).\\
		\end{flalign*} 

\begin{figure}
			\centering
			\includegraphics[scale=.4]{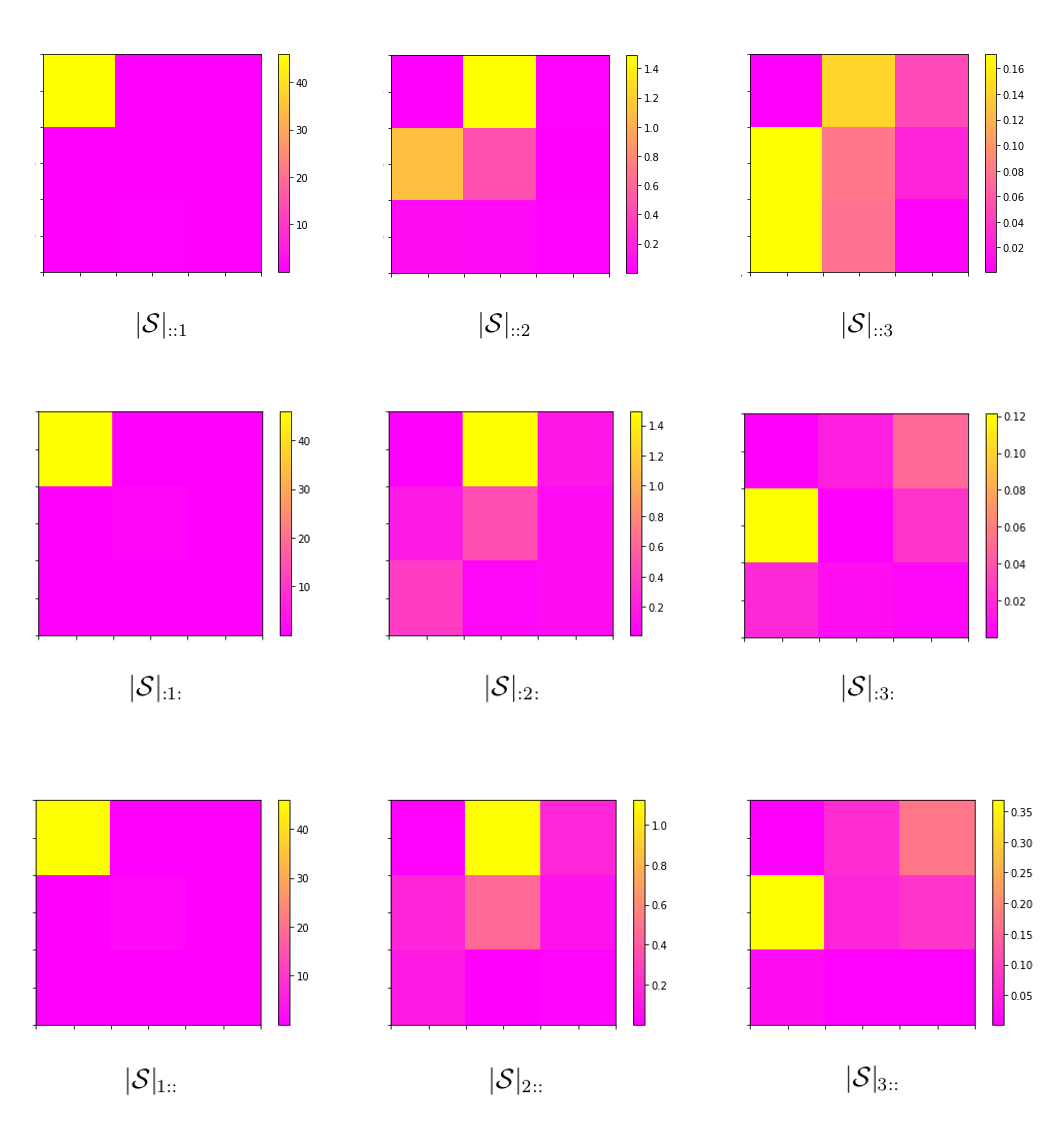}
			\caption{\footnotesize{Energy distribution of truncation $\tilde{\mathcal{S}}$.}}
			\label{energy-dist}
		\end{figure}

		Let $\mathcal{T} = (\textbf{U}^{(1)}, \textbf{U}^{(2)}, \textbf{U}^{(3)}) \cdot \mathcal{S}$ be the MLSVD of $\mathcal{T}$. After the truncation process (described in~\ref{compression}) we obtain a tensor $\tilde{\mathcal{S}} \in \mathbb{R}^{3 \times 3 \times 3}$ whose respective relative error is $\frac{\| \mathcal{S} - \tilde{\mathcal{S}} \|}{\| \mathcal{S} \|} = 3.85 \cdot 10^{-14}$. Since the multilinear rank is limited by the rank, we know in advance that $rank_{\boxplus}(\mathcal{T}) \leq (3, 3, 3)$. Since the error showed is practically is minimal, the truncated MLSVD obtained is the actual MLSVD of $\mathcal{T}$, and we have the equality $rank_{\boxplus}(\mathcal{T}) = (3, 3, 3)$.  We can visualize the energy distribution of $\mathcal{S}$ in figure~\ref{energy-dist}, which plots the energy of the slices respective to each mode. Each square is an entry of a slice in absolute value. The bars on the side of each image indicates the magnitudes of the entries.
		
		The tensor $\mathcal{S}$ can be showed explicitly, it is given by 
			
			$$\hspace{-1cm} \mathcal{S} = \left\{
			\left[ \begin{array}{ccc}
				45.84 & 0.012 & 0.00006\\
 				0.014 & 0.16 & 0.12\\
 				0.002 & -0.36 & -0.019
			\end{array} \right],
			\left[ \begin{array}{ccc}
				0.0005 & 1.48 & -0.015\\
 				1.12 & -0.46 & -0.00031\\
 				0.065 & -0.053 & 0.0069
			\end{array} \right],
			\left[ \begin{array}{ccc}
				0.00096 & -0.142 & 0.049\\
 				0.17 & -0.08 & -0.024\\
 				-0.17 & 0.076 & 0.0038
			\end{array} \right]
			\right\}.$$
This presentation of $\mathcal{S}$ is based on its frontal slices, that is, $\mathcal{S} = \{ \mathcal{S}_{::1}, \mathcal{S}_{::2}, \mathcal{S}_{::3} \}$.	
		
		Just as we did in~\ref{compression}, we could have obtained a smaller truncation $\tilde{\mathcal{S}}$ and then we would forget $\mathcal{S}$ by defining $\mathcal{S} = \tilde{\mathcal{S}}$. The matrices $\textbf{U}^{(\ell)}$ would have some of its last columns deleted and we would redefine these matrices too. The original MLSVD is forgotten and we use the truncated version in its place. Note that, by limiting ourselves to the truncated MLSVD, the best possible CPD we can achieve will have relative error equal\footnote{It is still possible to have smaller errors by pure luck, but we won't take this in account.} to the relative error between $\tilde{\mathcal{S}}$ and $\mathcal{S}$. By varying $1 \leq i,j,k \leq 3$ and checking all possible truncations $\mathcal{S}_{:i, :j, :k}$ we can get a better picture of the relation between truncation and error. Table 4.5 below shows how the truncation and errors are related in this example.\\
		
		\begin{table}
			\centering
			\footnotesize
			\begin{tabular}{|c|c|}
				\hline
				\textbf{Shape of truncation} & \textbf{Relative error}\\
				\hline
				$1\times 1\times 1$ & 0.04345277\\
				\hline
 				$1\times 1\times 2$ & 0.04345277\\
 				\hline
				$1\times 1\times 3$ & 0.04345276\\
				\hline
				$1\times 2\times 1$ & 0.04345198\\
				\hline
				$1\times 3\times 1$ & 0.04345198\\
				\hline
				$2\times 1\times 1$ & 0.04345156\\
				\hline
 				$3\times 1\times 1$ & 0.04345154\\
 				\hline
 				$2\times 2\times 1$ & 0.04330007\\
 				\hline
 				$2\times 3\times 1$ & 0.04321934\\
 				\hline
 				$3\times 2\times 1$ & 0.04255157\\
 				\hline
				$3\times 3\times 1$ & 0.04246718\\
				\hline
 				$2\times 1\times 2$ & 0.03591381\\
 				\hline
 				$3\times 1\times 2$ & 0.03588523\\
 				\hline
 				$2\times 1\times 3$ & 0.03572117\\
 				\hline
 				$3\times 1\times 3$ & 0.03549749\\
 				\hline
 				$1\times 2\times 2$ & 0.02888081\\
 				\hline
 				$1\times 3\times 2$ & 0.02887894\\
 				\hline
 				$1\times 2\times 3$ & 0.02871292\\
 				\hline
 				$1\times 3\times 3$ & 0.02869078\\
 				\hline
 				$2\times 2\times 2$ & 0.01093718\\
 				\hline
 				$2\times 3\times 2$ & 0.01060795\\
 				\hline
 				$2\times 2\times 3$ & 0.00964929\\
 				\hline
 				$2\times 3\times 3$ & 0.00919564\\
 				\hline
 				$3\times 2\times 2$ & 0.00720686\\
 				\hline
 				$3\times 3\times 2$ & 0.00668082\\
 				\hline
 				$3\times 2\times 3$ & 0.00296115\\
 				\hline
 				$3\times 3\times 3$ & 0\\
				\hline
			\end{tabular}
			\caption{\footnotesize{Error of all possible truncations.}} \bigskip
		\end{table} 
			
		Our method of initialization based on the MLSVD chooses the three brightest yellow squares (where the energy is most concentrated) of figure~\ref{energy-dist} to generate an initial approximation $\mathcal{S}^{(0)}$. This initial tensor is 		
		$$\mathcal{S}^{(0)} = s_{111}\ \textbf{e}_1 \otimes \textbf{e}_1 \otimes \textbf{e}_1 + s_{212}\ \textbf{e}_2 \otimes \textbf{e}_1
\otimes \textbf{e}_2 + s_{122}\ \textbf{e}_1 \otimes \textbf{e}_2 \otimes \textbf{e}_2 = $$
		$$ = 45.84\ \textbf{e}_1 \otimes \textbf{e}_1 \otimes \textbf{e}_1 + 1.12\ \textbf{e}_2 \otimes \textbf{e}_1
\otimes \textbf{e}_2 + 1.48\ \textbf{e}_1 \otimes \textbf{e}_2 \otimes \textbf{e}_2.$$	

\begin{figure} 
			\centering
			\includegraphics[scale=.5]{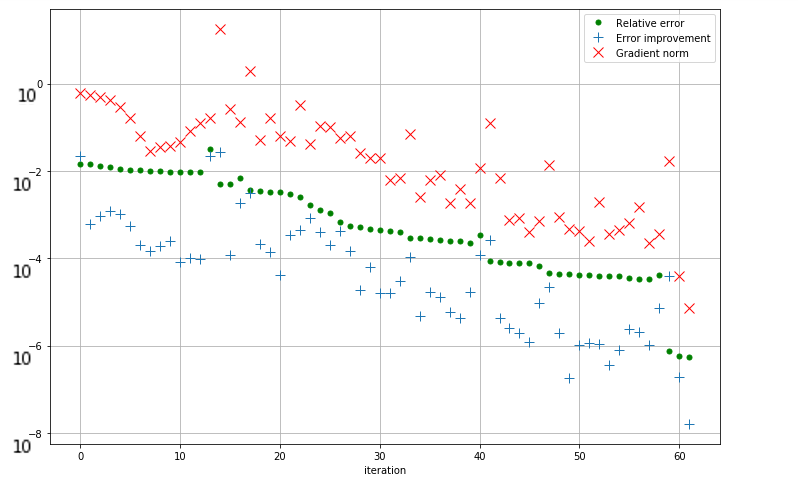}
			\caption{\footnotesize{Plot of the evolution of several measures made during the computation of a CPD.}}
			\label{warming_up-results}
		\end{figure}

		The respective error of this approximated rank-$3$ CPD is 
		$$\displaystyle \frac{\| \mathcal{T} - (\textbf{U}^{(1)}, \textbf{U}^{(2)}, \textbf{U}^{(3)}) \cdot \mathcal{S}^{(0)} \|}{\| \mathcal{T} \|} = 0.0153,$$
which is already small enough for a first iteration. From this tensor we start the dGN iterations in order to compute a better rank-3 CPD for $\mathcal{T}$. The results of the computations are summarized in figure~\ref{warming_up-results}. All the plots are in $\log_{10}$ scale. We can see the error decreasing but going up sometimes. The algorithm is not necessarily monotonic, but it always auto-correct itself. Adding line-search methods at each iteration were attempted, and this indeed made the algorithm monotonic (at least in practice). However this monotonicity seems to make the iterations to be drawn by local minima. Every attempt to decrease the error a little bit at each iteration led to some kind of local minima attraction. Our current algorithm is more ``erratic'' sometimes but is this behavior is what makes the steps ``run away'' from local minima. The \emph{error improvement} showed is the absolute value of the difference between two consecutive relative errors. Figure~\ref{warming_up-results} shows the evolution of the error, improvement and gradient as the iterative process progresses.
			
		 After all computations are finished the program outputs the approximated CPD $(\textbf{X}, \textbf{Y}, \textbf{Z}) \cdot \mathcal{I}_{3 \times 3 \times 3}$, where
		 
		 $$\textbf{X} = 
		 \left[ \begin{array}{ccc}
		 	-0.81709105 & 0.78059637 & -0.53290962\\
 			-0.87159249 & 0.76495679 & -0.6872218 \\
 			-1.04470608 & 0.71882424 & -0.81413498\\
 			-1.36099384 & 0.64403858 & -0.90858919\\
 			-1.83712919 & 0.54358195 & -0.96681853\\
 			-2.40011984 & 0.42145802 & -0.98650216
		 \end{array} \right],$$
		 
		 $$\textbf{Y} = 
		  \left[ \begin{array}{ccc}
		 	-0.97352551 & -0.72100749 & 0.83684308\\
 			-1.03631189 & -0.72099582 & 0.86877435\\
 			-1.25003023 & -0.72096898 & 0.94377806\\
 			-1.70859268 & -0.72095345 & 0.99448294\\
 			-2.64633452 & -0.72099573 & 0.90429673
		 \end{array} \right],$$
		 
		 $$\textbf{Z} = 
		  \left[ \begin{array}{ccc}
		 	-1.05782625 & -0000.246843654 & -1.01914978\\
 			-1.18213922 & 0.510766716 & -1.01896187\\
 			-1.64980857 & 0.907135088 & -1.01880278\\
 			-2.87546993 & 1.23098132 & -1.01864246
		 \end{array} \right].$$
		 
		 If everything worked correctly, the first column of $\textbf{X}$ should be some multiple\footnote{This issue comes from the scaling we are using. It is possible to re-scale the factors to fit the original functions without changing the CPD.} of the outputs of $\varphi_1^{(1)}(x) = \cos(1-x)$ for $x = 0, 0.2, 0.4, 0.6, 0.8, 1$. We can try scaling the columns of $\textbf{X}$ to fit the data. If after several trials none of them seems to work, it may be a good idea to compute another CPD, or change the initialization method, or introduce more restrictions to the model, etc. 
		 
	\section{Benchmark tensors}	
		As we already said sometimes, most parameters and subroutines of Tensor Fox were obtained by experience. The first attempts were based on other implementations, and from there we start to fine tuning to get the most of every single parameter. To make sure the choices are optimal and generic enough, it is necessary to test the performance for a diverse family of tensors, taking into account what are the known difficulties and what is applicable. In this section, we present all the tensors used for testing and benchmarking.
		
		\subsection{Swimmer} 
\begin{figure} 
				\centering
				\includegraphics[scale=.5]{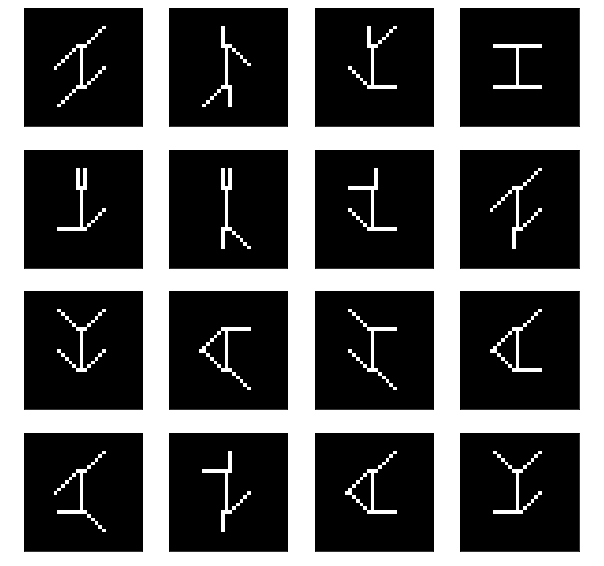}
				\caption{\footnotesize{Swimmer tensor.}}
				\label{swimmer}
			\end{figure}

			This tensor was constructed based on the paper \cite{shashua} as an example of a nonnegative tensor. It is a set of 256 images of dimensions $32 \times 32$ representing a swimmer. Each image contains a torso (the invariant part) of 12 pixels in the center and four limbs of 6 pixels that can be in one of 4 positions. In this work they proposed to use a rank $R = 50$ tensor to approximate, and we do the same for our test. In figure~\ref{swimmer} we can see some frontal slices of this tensor.
			
		\subsection{Handwritten digits} \label{hw}
\begin{figure} 
				\centering
				\includegraphics[scale=.5]{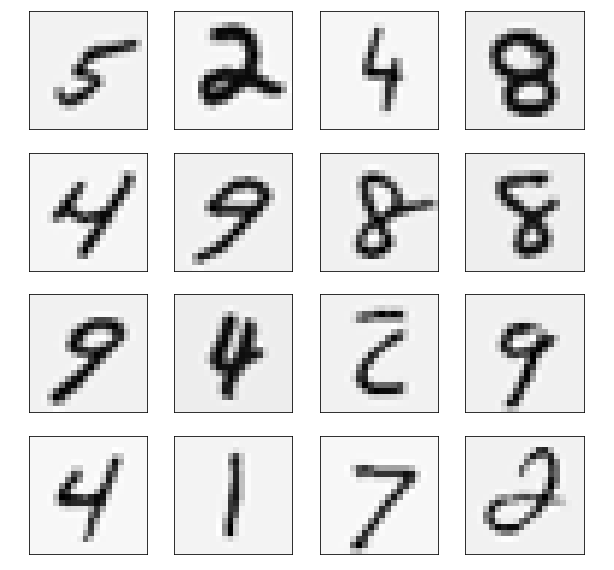}
				\caption{\footnotesize{Handwritten digits tensor.}}
				\label{handwritten}
			\end{figure}

			This is a classic tensor in machine learning, it is the MNIST\footnote{\url{http://yann.lecun.com/exdb/mnist/}} database of handwritten digits. Each slice is a image of dimensions $20 \times 20$ of a handwritten digit. Also, each 500 consecutive slices correspond to the same digit, so the first 500 slices correspond to the digit 0, the slices 501 to 1000 correspond to the digit 1, and so on. We choose $R = 150$ as a good rank to construct the approximating CPD to this tensor. In figure~\ref{handwritten} we can see some frontal slices of this tensor.
			
		\subsection{Border rank} 
			Remember the definition and discussion about border rank in~\ref{topology}. We observed that the strict inequality $\underline{rank}(\mathcal{T}) < rank(\mathcal{T}) = R$ can happen, and this means that the set of rank-$R$ tensors is not closed. This phenomenon makes the CPD computation a challenging problem. The paper \cite{lim} has a great discussion on this subject. In the same paper (and theorem~\ref{border-rank-lim}) they show that 
			$$\mathcal{T}^{(n)} = n \left( \textbf{x}^{(1)} + \frac{1}{n} \textbf{y}^{(1)} \right) \otimes \left( \textbf{x}^{(2)} + \frac{1}{n} \textbf{y}^{(2)} \right) \otimes \left( \textbf{x}^{(3)} + \frac{1}{n} \textbf{y}^{(3)} \right) - n \textbf{x}^{(1)} \otimes \textbf{x}^{(2)} \otimes \textbf{x}^{(3)}$$
is a sequence of rank sequence of rank 2 tensors converging to a tensor of rank 3, where each pair $\textbf{x}^{(\ell)}, \textbf{y}^{(\ell)} \in \mathbb{R}^{I_\ell}$ is linearly independent. Remember that the limit tensor is $\mathcal{T} = \textbf{x}^{(1)} \otimes \textbf{x}^{(2)} \otimes \textbf{y}^{(3)} + 	\textbf{x}^{(1)} \otimes \textbf{y}^{(2)} \otimes \textbf{x}^{(3)} + \textbf{y}^{(1)} \otimes \textbf{x}^{(2)} \otimes \textbf{x}^{(3)}$. We choose to compute a CPD of rank $R = 2$ to see how the algorithms behaves when we try to approximate a problematic tensor by tensor with low rank. In theory it is possible to have arbitrarily good approximations.

		\subsection{Matrix multiplication}
			Let $\mathcal{M}_N \in \mathbb{R}^{N^2 \times N^2 \times N^2}$ be the tensor associated with the multiplication between two matrices in $\mathbb{R}^{N \times N}$. The classic form of $\mathcal{M}_N$ is given by
			$$\mathcal{M}_N = \sum_{i=1}^m \sum_{j=1}^n \sum_{k=1}^l vec(\textbf{e}_{ij}) \otimes vec(\textbf{e}_{jk}) \otimes vec(\textbf{e}_{ik}),$$
where $\textbf{e}_{ij}$ is the matrix $N \times N$ with entry $(i,j)$ equal to 1 and the remaining entries equal to zero. Since Strassen \cite{strassen} it is known that matrix multiplication between matrices of dimensions $N \times N$ can be made with $\mathcal{O}(N^{\log_2 7})$ operations. Many improvements were made after Strassen but we won't enter in the details here. For the purpose of testing we choose the small value $N = 5$ and the rank $R = \lceil 5^{\log_2 7} \rceil = 92$ in honor of Strassen. However note that this is probably not the exact rank of $\mathcal{M}_5$, so this test is about a strict low rank approximation of a difficult tensor.			

		\subsection{Collinear factors}
			The phenomenon of swamps occurs when all factors in each mode are almost collinear. Their presence is a challenge for many algorithms because they can slow down convergence. Now we will create synthetic data to simulate various degrees of collinearity between the factors. We begin generating three random matrices $\textbf{M}_X \in \mathbb{R}^{m \times R}, \textbf{M}_Y \in \mathbb{R}^{n \times R}, \textbf{M}_Z \in \mathbb{R}^{p \times r}$, where each entry is drawn from the normal distribution with mean 0 and variance 1. After that we perform QR decomposition of each matrix, obtaining the decompositions $\textbf{M}_X = \textbf{Q}_X \textbf{R}_X, \textbf{M}_Y = \textbf{Q}_Y \textbf{R}_Y, \textbf{M}_Z = \textbf{Q}_Z \textbf{R}_Z$. The matrices $\textbf{Q}_X, \textbf{Q}_Y, \textbf{Q}_Z$ are orthogonal. Now fix three columns $\textbf{q}_X^{(i')}, \textbf{q}_Y^{(j')}, \textbf{q}_Z^{(k')}$ of each one of these matrices. The factors $\textbf{X} = [ \textbf{x}^{(1)}, \ldots, \textbf{x}^{(R)} ] \in \mathbb{R}^{m \times R}$, $\textbf{Y} = [ \textbf{y}^{(1)}, \ldots, \textbf{y}^{(R)} ] \in \mathbb{R}^{n \times R}$, $\textbf{Z} = [ \textbf{z}^{(1)}, \ldots, \textbf{z}^{(R)} ] \in \mathbb{R}^{p \times R}$ are generated by the equations below.
			
			$$\textbf{x}^{(i)} = \textbf{q}_X^{(i')} + c \cdot \textbf{q}_X^{(i)}, \quad i=1 \ldots R$$
			$$\textbf{y}^{(j)} = \textbf{q}_Y^{(j')} + c \cdot \textbf{q}_Y^{(j)}, \quad j=1 \ldots R$$
			$$\textbf{z}^{(k)} = \textbf{q}_Z^{(k')} + c \cdot \textbf{q}_Z^{(k)}, \quad k=1 \ldots R$$  
		
			The parameter $c \geq 0$ defines the degree of collinearity between the vectors of each factor. A value of $c$ close to 0 indicates high degree of collinearity, while a high value of $c$ indicates low degree of collinearity.
			
			Another phenomenon that occurs in practice is the presence of noise in the data. So we will treat these two phenomena at once in this benchmark. After generating the factors $\textbf{X}, \textbf{Y}, \textbf{Z}$ we have a tensor $\mathcal{T} = (\textbf{X}, \textbf{Y}, \textbf{Z}) \cdot \mathcal{I}_{R \times R \times R}$. That is, $\textbf{X}, \textbf{Y}, \textbf{Z} $ are the exact CPD of $\mathcal{T}$. Now consider a noise $\mathcal{N} \in \mathbb{R}^{m \times n \times p}$ such that each entry of $\mathcal{N} $ is obtained by the normal distribution with mean 0 and variance 1. Thus we form the tensor $\hat{\mathcal{T}} = \mathcal{T} + \nu \cdot \mathcal{N}$, where $ \nu > 0 $ defines the magnitude of the noises. The idea is to compute a CPD of $ \hat{\mathcal{T}}$ of rank $R$ and then evaluate the relative error between this tensor and $\mathcal{T}$. We expect the computed CPD to clear the noises and to be close to $\mathcal{T}$ (even if it is not close to $ \hat{\mathcal{T}}$). We will fix $\nu = 0.01$ and generate tensors for $c = 0.1, 0.5, 0.9$. In all cases we will be using $ m = n = p = 300 $ and $R = 15 $. This is a particularly difficult problem since we are considering swamps and noises at once. The same procedure to generate tensors were used for benchmarking in \cite{cichoki2013}.
			
		\subsection{Double bottlenecks}
			We proceed almost in the same as before for swamps, we used the same procedure to generate the first two columns of each factor matrix, then the remaining columns are equal to the columns of the QR decomposition. After generating the factors $\textbf{X}, \textbf{Y}, \textbf{Z}$ we consider a noise $\mathcal{N} \in \mathbb{R}^{m \times n \times p}$ such that each entry of $\mathcal{N} $ is obtained by the normal distribution with mean 0 and variance 1. Thus we form the tensor $\hat{\mathcal{T}} = \mathcal{T} + \nu \cdot \mathcal{N}$, where $ \nu > 0 $ defines the magnitude of the noises. The collinear parameter used for the tests are $c = 0.1, 0.5$, and we fix $\nu = 0.01$ as before. The procedure before the noise is presented below.
			
			$$\textbf{x}^{(i)} = \textbf{q}_X^{(i')} + c \cdot \textbf{q}_X^{(i)}, \quad i=1, 2$$
			$$\textbf{y}^{(j)} = \textbf{q}_Y^{(j')} + c \cdot \textbf{q}_Y^{(j)}, \quad j=1, 2$$
			$$\textbf{z}^{(k)} = \textbf{q}_Z^{(k')} + c \cdot \textbf{q}_Z^{(k)}, \quad k=1, 2$$ 
			
			$$\textbf{x}^{(i)} = \textbf{q}_X^{(i')}, \quad i=3 \ldots R$$
			$$\textbf{y}^{(j)} = \textbf{q}_Y^{(j')}, \quad j=3 \ldots R$$
			$$\textbf{z}^{(k)} = \textbf{q}_Z^{(k')}, \quad k=3 \ldots R$$ 
			
	\section{Fine tuning}
		We will conduct the tests using only third order tensors The dimensions will be written as $m, n, p$ so the tensors are elements of $\mathbb{R}^{m \times n \times p}$. As already mentioned, the parameters \verb|maxiter|, \verb|tol|, \verb|cg_maxiter| were found after several experimentations and tests. Here we consider varying each one of them around the default ones and test them with our test tensors. The results we are going to show here are just a fraction of the total work done in order to achieve these values. In particular, \verb|cg_maxiter| was a result of several different models and parameter fine tuning. We will conduct an \emph{hyperparameter} search on the parameters. This consists in considering a high dimensional grid of values, where each dimension corresponds to a certain parameter. Each point in the grid correspond to a model (a particular choice of parameters) on which we evaluate its performance. The grid we are considering here is given by $P_1 \times P_2 \times P_3 \times P_4$, where		
		 \begin{flalign*}
			& P_1 = \{ 100, 200, 400 \},\\
			& P_2 = \{ 10^{-4}, 10^{-6}, 10^{-8} \},\\
			& P_3 = \{ 0.2, 0.3, 0.4, 0.5 \},\\
			& P_4 = \{ 0.6, 0.7, 0.8, 0.9 \}.
		\end{flalign*}
			
\begin{figure} 
			\hspace*{-1cm}
			\includegraphics[scale=.42]{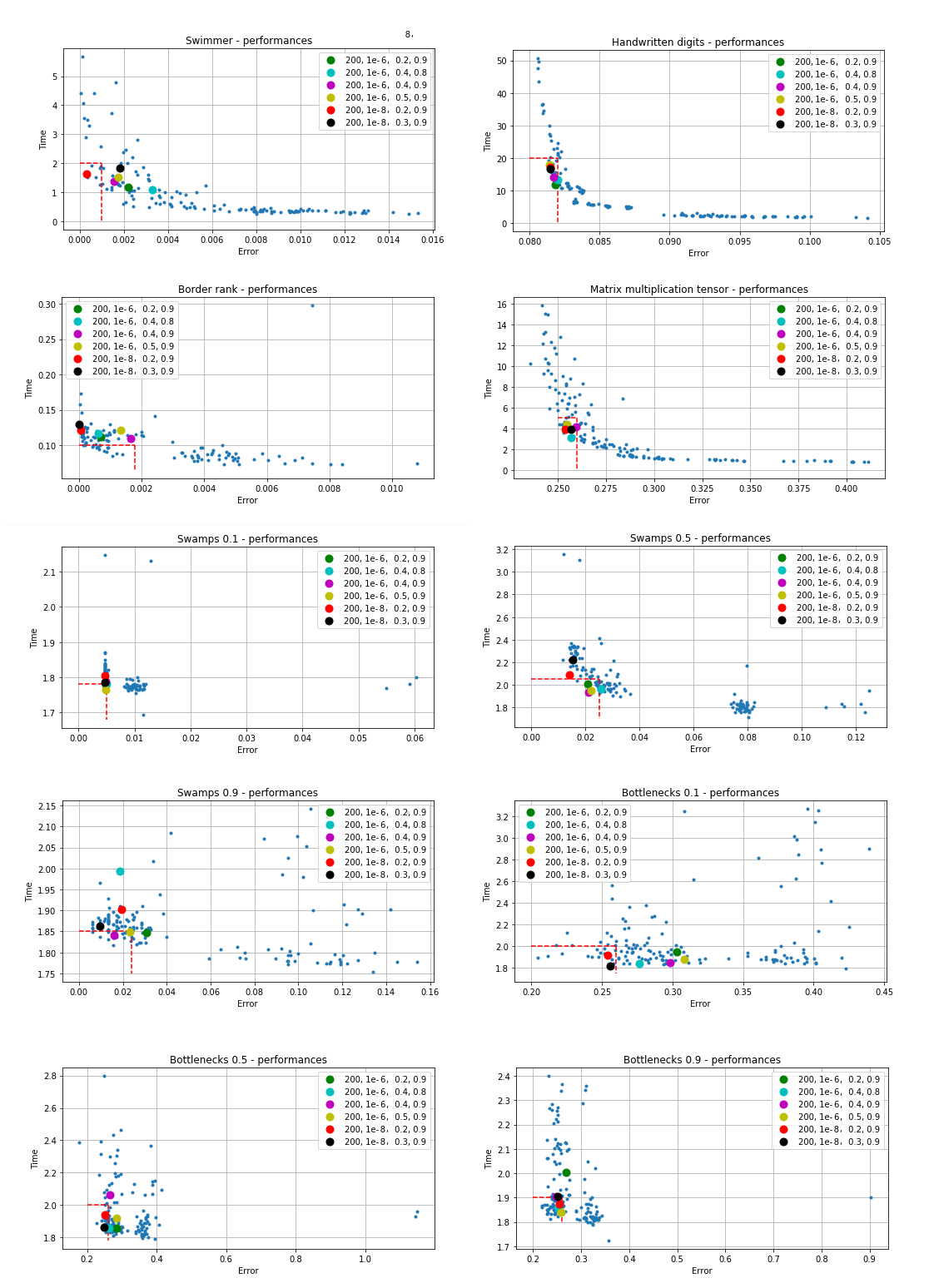}
			\caption{\footnotesize{The location of the best models may help to decide which one is better.}}
			\label{hyperparameter}
		\end{figure}

		The idea is that \verb|maxiter| $= p_1$, \verb|tol| $= p_2 \cdot mnp$, $a = p_3, b = p_4$, where $p_i \in P_i$ for all $i$. Consider that $\verb|cg_maxiter| \sim \left[ 1 + \lceil k^a \rceil, 2 + \lceil k^b \rceil \right]$. It is possible to use more values and larger intervals, but experience showed that good performance is attained around these values. We test each tensor 20 times and take the average time and average error to measure its performance. In the plots of figure~\ref{hyperparameter}, each small point correspond to a model (a particular choice of parameters). The better choices correspond to points closer to the origin, which translates to small error in a small time. We selected some regions of interest in order to discard the failed models. The sizes of these regions depend on the considered tensor. Sometimes the precision is more important that the time (at a certain degree) and sometimes we can afford to lose precision to gain in time. This is a subtle matter, but in general we will avoid models falling too far away from the region close to the origin. The goal is to obtain a well balanced set of parameters. The region between the red dotted lines and the axes indicates the region of interest for each test tensor. Every point inside these regions are registered. The idea is too see if there are models appearing inside these regions for different tensors. This would indicate that the model performs well in more than just one tensor. The ones appearing more frequently are more likely to be well balanced and generic. Table 4.6 shows the results of all models which appeared at least four times inside some of the defined regions. The models marked in red are the ones with five or more occurrences. These are the best ones by our criteria. We used the abbreviations Sw = Swimmer, Hw = Handwritten, Br = Border rank, Mm = Matrix multiplication, Swp 0.1 = Swamp with $c = 0.1$, and similar for the other swamp tensors, and Bn 0.1 = Bottleneck with $c = 0.1$, with the same considerations as the swamp tensors. To decide which model to choose we also plot the position of these four models. The first thing we note is that all models performs reasonably well and quite similar for the first four test tensors. The differences begin to appear when we introduce collinearity. In the case of the swamp tensors the models with tolerance $10^{-8}$ are able to attain the same error of the other models but in a bigger time. On the other hand, in the case of bottleneck tensors the models with tolerance $10^{-6}$ takes the same time of the other models but they get a noticeable bigger error for $c = 0.1$. 
		
		All these four models are good, but we have to decide between one of them. We can repeat the previous tests but with more repetitions to get more accurate performance statistics. Since there are few models to test this is a possible task. The tests are repeated with 50 repetitions instead of 20, and now we also consider the variance of the error and time.\footnote{Given random variables $X_1, \ldots,  X_N$ i.i.d. (independent and identically distributed), the associated \emph{empirical variance} is $\overline{\sigma}^2 = \frac{1}{N-1} \sum_{i=1}^N (X_i - \overline{X})^2$, where $\overline{X} = \frac{1}{N} \sum_{i=1}^N X_i$ is the \emph{empirical mean}} The results are showed in the next tables. We start discarding the model where \verb|tol| $= 10^{-8}$ since this one clearly is more demanding. At the end of the day, ``speed shall prevail''. The other three models are more well balanced so any choice we make is good enough. Our final choice is \verb|maxiter| $= 200$, \verb|tol| $= 10^{-6}$, $a = 0.4, b = 0.9$ since it has competitive errors measures with better timings in most of the tests. 
		
		\begin{table}
			\scriptsize
			\hspace*{-1cm}
  			\begin{tabular}{|c|c|c|c|c|c|c|c|c|c|c|c|c|c|} 
				\hline
				\verb|maxiter| & \verb|tol| & $a$ & $b$ & \textbf{Sw} & \textbf{Hw} & \textbf{Br} & \textbf{Mm} & \textbf{Swp 0.1} & \textbf{Swp 0.5} & \textbf{Swp 0.9} & \textbf{Bn 0.1} & \textbf{Bn 0.5} & \textbf{Bn 0.9}\\
				\hline
				100  &   $10^{-8}$  &  0.4  &  0.7 & & & & & & x & x & x & x & \\   
				\hline
				200  &   $10^{-6}$  &  0.4  &  0.9 & \textcolor{red}{\textbf{x}} & \textcolor{red}{\textbf{x}} & & & \textcolor{red}{\textbf{x}} & & \textcolor{red}{\textbf{x}} & \textcolor{red}{\textbf{x}} & & \\
				\hline
				200  &   $10^{-6}$  &  0.5  &  0.9 & & \textcolor{red}{\textbf{x}} & & \textcolor{red}{\textbf{x}} & & & \textcolor{red}{\textbf{x}} & \textcolor{red}{\textbf{x}} & & \textcolor{red}{\textbf{x}} \\   
				\hline			
				200  &   $10^{-8}$  &  0.5  &  0.9 & \textcolor{red}{\textbf{x}} & \textcolor{red}{\textbf{x}} & \textcolor{red}{\textbf{x}} & \textcolor{red}{\textbf{x}} & & & & & \textcolor{red}{\textbf{x}} & \\
				\hline
				400  &   $10^{-6}$  &  0.2  &  0.8 & & & & & x & x & & x & x & \\   
				\hline
				400  &   $10^{-6}$  &  0.2  &  0.9 & x & x & & & & x & & & & x \\   
				\hline
				400  &   $10^{-6}$  &  0.3  &  0.9 & & x & & & x & x & & x & & \\   
				\hline
				400  &   $10^{-6}$  &  0.4  &  0.7 & & & & & x & x & & x & x & \\   
				\hline
				400  &   $10^{-6}$  &  0.4  &  0.8 & \textcolor{red}{\textbf{x}} & \textcolor{red}{\textbf{x}} & & & & \textcolor{red}{\textbf{x}} & \textcolor{red}{\textbf{x}} & & & \textcolor{red}{\textbf{x}} \\ 
				\hline 
				400  &   $10^{-6}$  &  0.4  &  0.9 & & x & & & & x & & x & x & \\   
				\hline			
				400  &   $10^{-6}$  &  0.5  &  0.9 & & x & & & x & x & & & & x \\   
				\hline				
			\end{tabular}
			\caption{\footnotesize{Hyperparameter grid search - best models.}} 
		\end{table}	
		
		\begin{table}
			\scriptsize
			\centering
  			\begin{tabular}{|c|c|c|c|c|} 
				\hline
				\textbf{Model} & \textbf{Average error} & \textbf{Variance error} & \textbf{Average time} & \textbf{Variance time} \\
				\hline
				$200,\ 10^{-6},\ 0.4,\ 0.9$ & 0.0011059 & 0.0000071 & 1.27 sec & 0.13 sec \\
				\hline
				$200,\ 10^{-6},\ 0.5,\ 0.9$ & 0.0018093 & 0.0000155 & 1.43 sec & 0.19 sec \\   
				\hline
				$200,\ 10^{-8},\ 0.5,\ 0.9$ & 0.0013749 & 0.0000139 & 1.78 sec & 0.02 sec \\
				\hline
				$400,\ 10^{-6},\ 0.4,\ 0.8$ & 0.0029579 & 0.0000282 & 1.39 sec & 0.25 sec \\
				\hline			
			\end{tabular}
			\caption{\footnotesize{Swimmer tensor - final best performances.}} 
		\end{table}	
		
		\begin{table}
			\scriptsize
			\centering
  			\begin{tabular}{|c|c|c|c|c|} 
				\hline
				\textbf{Model} & \textbf{Average error} & \textbf{Variance error} & \textbf{Average time} & \textbf{Variance time} \\
				\hline
				$200,\ 10^{-6},\ 0.4,\ 0.9$ & 0.0816178 & 0.0000001 & 15.32 sec & 7.49 sec \\
				\hline
				$200,\ 10^{-6},\ 0.5,\ 0.9$ & 0.0815007 & 0.0000000 & 15.58 sec & 4.09 sec \\   
				\hline
				$200,\ 10^{-8},\ 0.5,\ 0.9$ & 0.0814086 & 0.00000001 & 17.50 sec & 1.74 sec \\
				\hline
				$400,\ 10^{-6},\ 0.4,\ 0.8$ & 0.0817002 & 0.00000001 & 16.69 sec & 7.57 sec \\
				\hline			
			\end{tabular}
			\caption{\footnotesize{Handwritten tensor - final best performances.}} 
		\end{table}	
		
		\begin{table}
			\scriptsize
			\centering
  			\begin{tabular}{|c|c|c|c|c|} 
				\hline
				\textbf{Model} & \textbf{Average error} & \textbf{Variance error} & \textbf{Average time} & \textbf{Variance time} \\
				\hline
				$200,\ 10^{-6},\ 0.4,\ 0.9$ & 0.0016909 & 0.0000044 & 0.11 sec & 0.0003 sec \\
				\hline
				$200,\ 10^{-6},\ 0.5,\ 0.9$ & 0.0027261 & 0.0000046 & 0.11 sec & 0.0002 sec \\   
				\hline
				$200,\ 10^{-8},\ 0.5,\ 0.9$ & 0.0000106 & 0.00000001 & 0.12 sec & 0.0004 sec \\
				\hline
				$400,\ 10^{-6},\ 0.4,\ 0.8$ & 0.0005573 & 0.0000018 & 0.13 sec & 0.002 sec \\
				\hline			
			\end{tabular}
			\caption{\footnotesize{Border rank tensor - final best performances.}} 
		\end{table}	
		
		\begin{table}
			\scriptsize
			\centering
  			\begin{tabular}{|c|c|c|c|c|} 
				\hline
				\textbf{Model} & \textbf{Average error} & \textbf{Variance error} & \textbf{Average time} & \textbf{Variance time} \\
				\hline
				$200,\ 10^{-6},\ 0.4,\ 0.9$ & 0.2384032 & 0.0007164 & 4.37 sec & 0.20 sec \\
				\hline
				$200,\ 10^{-6},\ 0.5,\ 0.9$ & 0.2396802 & 0.0005141 & 4.34 sec & 0.085 sec \\   
				\hline
				$200,\ 10^{-8},\ 0.5,\ 0.9$ & 0.2423333 & 0.0004590 & 4.33 sec & 0.067 sec \\
				\hline
				$400,\ 10^{-6},\ 0.4,\ 0.8$ & 0.2269175 & 0.0006688 & 9.35 sec & 5.06 sec \\
				\hline			
			\end{tabular}
			\caption{\footnotesize{Matrix multiplication tensor - final best performances.}} 
		\end{table}	
		
		\begin{table}
			\scriptsize
			\centering
  			\begin{tabular}{|c|c|c|c|c|} 
				\hline
				\textbf{Model} & \textbf{Average error} & \textbf{Variance error} & \textbf{Average time} & \textbf{Variance time} \\
				\hline
				$200,\ 10^{-6},\ 0.4,\ 0.9$ & 0.0288477 & 0.0000026 & 1.86 sec & 0.07 sec \\
				\hline
				$200,\ 10^{-6},\ 0.5,\ 0.9$ & 0.0288307 & 0.0000027 & 1.83 sec & 0.007 sec \\   
				\hline
				$200,\ 10^{-8},\ 0.5,\ 0.9$ & 0.0326175 & 0.0000003 & 2.63 sec & 0.07 sec \\
				\hline
				$400,\ 10^{-6},\ 0.4,\ 0.8$ & 0.0278611 & 0.0000026 & 1.79 sec & 0.003 sec \\
				\hline			
			\end{tabular}
			\caption{\footnotesize{Swamp 0.1 tensor - final best performances.}} 
		\end{table}	
		
		\begin{table}
			\scriptsize
			\centering
  			\begin{tabular}{|c|c|c|c|c|} 
				\hline
				\textbf{Model} & \textbf{Average error} & \textbf{Variance error} & \textbf{Average time} & \textbf{Variance time} \\
				\hline
				$200,\ 10^{-6},\ 0.4,\ 0.9$ & 0.1000538 & 0.0000001 & 1.78 sec & 0.003 sec \\
				\hline
				$200,\ 10^{-6},\ 0.5,\ 0.9$ & 0.1002456 & 0.0000001 & 1.86 sec & 0.013 sec \\   
				\hline
				$200,\ 10^{-8},\ 0.5,\ 0.9$ & 0.1010695 & 0.00000001 & 2.69 sec & 0.07 sec \\
				\hline
				$400,\ 10^{-6},\ 0.4,\ 0.8$ & 0.0999770 & 0.0000001 & 1.79 sec & 0.007 sec \\
				\hline			
			\end{tabular}
			\caption{\footnotesize{Swamp 0.5 tensor - final best performances.}} 
		\end{table}
		
		\begin{table}
			\scriptsize
			\centering
  			\begin{tabular}{|c|c|c|c|c|} 
				\hline
				\textbf{Model} & \textbf{Average error} & \textbf{Variance error} & \textbf{Average time} & \textbf{Variance time} \\
				\hline
				$200,\ 10^{-6},\ 0.4,\ 0.9$ & 0.1376466 & 0.0002145 & 1.96 sec & 0.30 sec \\
				\hline
				$200,\ 10^{-6},\ 0.5,\ 0.9$ & 0.1374688 & 0.0000559 & 1.83 sec & 0.010 sec \\   
				\hline
				$200,\ 10^{-8},\ 0.5,\ 0.9$ & 0.1337583 & 0.0000065 & 1.83 sec & 0.005 sec \\
				\hline
				$400,\ 10^{-6},\ 0.4,\ 0.8$ & 0.1366625 & 0.0000472 & 1.85 sec & 0.011 sec \\
				\hline			
			\end{tabular}
			\caption{\footnotesize{Swamp 0.9 tensor - final best performances.}} 
		\end{table}
		
		\begin{table}
			\scriptsize
			\centering
  			\begin{tabular}{|c|c|c|c|c|} 
				\hline
				\textbf{Model} & \textbf{Average error} & \textbf{Variance error} & \textbf{Average time} & \textbf{Variance time} \\
				\hline
				$200,\ 10^{-6},\ 0.4,\ 0.9$ & 0.2876770 & 0.0090001 & 1.86 sec & 0.016 sec \\
				\hline
				$200,\ 10^{-6},\ 0.5,\ 0.9$ & 0.2738036 & 0.0110769 & 1.85 sec & 0.009 sec \\   
				\hline
				$200,\ 10^{-8},\ 0.5,\ 0.9$ & 0.2562058 & 0.0129000 & 1.92 sec & 0.019 sec \\
				\hline
				$400,\ 10^{-6},\ 0.4,\ 0.8$ & 0.2624673 & 0.0085848 & 1.87 sec & 0.006 sec \\
				\hline			
			\end{tabular}
			\caption{\footnotesize{Bottleneck 0.1 tensor - final best performances.}} 
		\end{table}
		
		\begin{table}
			\scriptsize
			\centering
  			\begin{tabular}{|c|c|c|c|c|} 
				\hline
				\textbf{Model} & \textbf{Average error} & \textbf{Variance error} & \textbf{Average time} & \textbf{Variance time} \\
				\hline
				$200,\ 10^{-6},\ 0.4,\ 0.9$ & 0.2822554 & 0.0059065 & 1.87 sec & 0.011 sec \\
				\hline
				$200,\ 10^{-6},\ 0.5,\ 0.9$ & 0.2669849 & 0.0052620 & 1.84 sec & 0.007 sec \\   
				\hline
				$200,\ 10^{-8},\ 0.5,\ 0.9$ & 0.2447143 & 0.0048144 & 1.87 sec & 0.016 sec \\
				\hline
				$400,\ 10^{-6},\ 0.4,\ 0.8$ & 0.2853812 & 0.0038925 & 1.82 sec & 0.004 sec \\
				\hline			
			\end{tabular}
			\caption{\footnotesize{Bottleneck 0.5 tensor - final best performances.}} 
		\end{table}
		
		\begin{table}
			\scriptsize
			\centering
  			\begin{tabular}{|c|c|c|c|c|} 
				\hline
				\textbf{Model} & \textbf{Average error} & \textbf{Variance error} & \textbf{Average time} & \textbf{Variance time} \\
				\hline
				$200,\ 10^{-6},\ 0.4,\ 0.9$ & 0.2578984 & 0.0018668 & 1.80 sec & 0.003 sec \\
				\hline
				$200,\ 10^{-6},\ 0.5,\ 0.9$ & 0.2654803 & 0.0013860 & 1.94 sec & 0.005 sec \\   
				\hline
				$200,\ 10^{-8},\ 0.5,\ 0.9$ & 0.2446079 & 0.0022603 & 1.98 sec & 0.004 sec \\
				\hline
				$400,\ 10^{-6},\ 0.4,\ 0.8$ & 0.2568034 & 0.0029858 & 1.95 sec & 0.002 sec \\
				\hline			
			\end{tabular}
			\caption{\footnotesize{Bottleneck 0.9 tensor - final best performances.}} 
		\end{table}		
		
		We point out that the algorithm has a lot of room to improvements. The parameters choice may be changed, all updating strategies are wildly open to modifications, the preconditioner is surely not optimal (in fact, future work includes testing more preconditioners), and the compression strategies also can be changed.\vspace{28cm}

	\section{Tensor Fox vs. other implementations}	
		\subsection{Procedure}
			We have selected a set of very distinct tensors to test the known tensor implementations. Given a tensor $\mathcal{T}$ and a rank $R$, we compute the CPD of TFX (short for \emph{Tensor Fox}) with the default maximum number of iterations\footnote{Remember that the default for the maximum number of iterations of TFX is \texttt{maxiter}$=200$.} $100$ times and retain the best result, i.e., the CPD with the smallest relative error. Let $\varepsilon$ be this error. Now let ALG be any other algorithm implemented by some of the mentioned libraries. We set the maximum number of iterations to \verb|maxiter|, keep the other options with their defaults, run ALG with these options $100$ times. The only accepted solutions are the ones with relative error smaller that $\varepsilon + \varepsilon/100$. Between all accepted solutions we return the one with the smallest running time. If none solution is accepted, we increase it to \verb|maxiter| by a certain amount and repeat.
			
	We try the values \verb|maxiter| $ = 5, 10, 50, 100, 150, \ldots, 900, 950, 1000$, until there is an accepted solution. The running time associated with the accepted solution is the accepted time. Otherwise we consider that ALG failed. We also consider a fail if the computational time is exceedingly high, in which case we can stop increasing \verb|maxiter| earlier. These procedures favour all implementations against TFX since we are trying a solution close to the solution of TFX with the minimum number of iterations. This benchmark measures the computational effort that each program makes to achieve the precision of TFX's precision. We consider it to be a fair measure of performance, although other kind of benchmarks might be considered. We remark that the iteration process is always initiated with a random point. The option to generate a random initial point is offered by all libraries, and we use each one they offer (sometimes random initialization was already the default option). There is no much difference in their initializations, which basically amounts to draw random factor matrices from the standard Gaussian distribution. The time to perform the MLSVD, or any kind of preprocessing, is included in the time measurements. If one want to reproduce the tests presented here, they can be found at \url{https://github.com/felipebottega/Tensor-Fox/tree/master/tests}.
		
		\subsection{State of art implementations}
			Now we briefly describe what algorithms will be used in our tests together with their corresponding implementations source. Some general algorithms are repeated but with variations, which are particular of each implementation. The ALS algorithm, for instance, is the one with more variations because people are trying to improve it for decades. For all Tensorlab algorithm implementation we recommend reading \cite{tensorlab2}, for the Tensor Toolbox we recommend \cite{cp_opt}, for TensorLy we recommend \cite{tensorly}, and for TensorBox we recommend \cite{tensorbox, cichoki2013}. In these benchmarks we used Tensorlab version 3.0 and Tensor Toolbox version 3.1.	
			
			\subsubsection{TFX}
				The algorithm used in TFX's implementation is the nonlinear least squares scheme described in the previous section. There are more implementation details to be discussed, but the interested reader can check \cite{tensorfox} for more information.
			
			\subsubsection{ALS}
				This is the Tensorlab's implementation of ALS algorithm. Although ALS is remarkably fast and easy to implement, it is not very accurate specially in the presence of bottlenecks or swamps. It seems (see \cite{tensorlab2}) that this implementation is very robust while still fast.  
				
			\subsubsection{NLS} 
				This is the Tensorlab's implementation of NLS algorithm. This is the one we described in the previous section. We should remark that this implementation is similar to TFX's implementation at some points, but there are big differences when we look in more details. In particular the compression procedure, the preconditioner, the damping parameter update rule and the number of iterations of the conjugate gradient are very different.
				
			\subsubsection{MINF}
				This is the Tensorlab's implementation of the problem as an optimization problem. They use a quasi-Newton method, the limited-memory BFGS, and consider equation~\ref{min} just as a minimization of a function. 
				
			\subsubsection{OPT}
				Just as the MINF approach, the OPT algorithm is a implementation of Tensor Toolbox, which considers~\ref{min} as a minimization of a function. They claim that using the algorithm option ``lbfgs'' is the preferred option\footnote{ \url{https://www.tensortoolbox.org/cp_opt_doc.html}}, so we used this way.  
				
			\subsubsection{Tly-ALS}
				TensorLy has only one way to compute the CPD, which is a implementation of the ALS algorithm. We denote it by Tly-ALS, do not confuse with ALS, the latter is the Tensorlab's implementation.
				
			\subsubsection{fLMa}
				fLMA stands for \emph{fast Levenberg-Marquardt algorithm}, and it is a different version of the damped Gauss-Newton.\\			
				
			In Tensorlab it is possible to disable or not the option of refinement. The first action in Tensorlab is to compress the tensor and work with the compressed version. If we want to use refinement then the program uses the compressed solution to compute a refined solution in the original space. This can be more demanding but can improve the result considerably. In our experience working in the original space is not a good idea because the computational cost increases drastically and the gain in accuracy is generally very small. Still we tried all Tensorlab algorithms with and without refinement. Table 4.17 shows compares the shapes of the factor matrices for each implementation, and there we can see that most of them them doesn't compress the tensor. We will write ALSr, NLSr and MINFr for the algorithms ALS, NLS and MINF with refinement, respectively.
			
		\begin{table}
		\centering
		\begin{tabular}{|c|c|}
			\hline
			\textbf{Solver} & \textbf{Shape of factor matrices}\\
			\hline
			Tensor Toolbox & $m \times R$,\ \ $n \times R$,\ \ $p \times R$\\
			\hline
			Tensorbox & $m \times R$,\ \ $n \times R$,\ \ $p \times R$\\
			\hline
			Tensorly & $m \times R$,\ \ $n \times R$,\ \ $p \times R$\\
			\hline
			Tensorlab & $\min(m, R) \times R$,\ \ $\min(n, R) \times R$,\ \ $\min(p, R) \times R$\\
			\hline
			Tensor Fox & $R_1 \times R$,\ \ $R_2 \times R$,\ \ $R_3 \times R$\\
			\hline
		\end{tabular}
		\caption{\footnotesize{Shapes of the factor matrices of each implementation for rank-$R$ CPD computation of a tensor with shape $m \times n \times p$.}} 
		\end{table}			
			
		\subsection{Computational results}
			We used Linux Mint operational system in our tests. All tests were conducted using a processor Intel Core i7-4510U - 2.00GHz (2 physical cores and 4 threads) and 8GB of memory. The packages mentioned run in Python (Tensorly and TensorFox) or Matlab (Tensorlab, Tensor Toolbox and TensorBox). We use Python - version 3.6.5 and Matlab - version 2017a. In both platforms we used BLAS MKL-11.3.1. Finally, we want to mention that TFX is implemented in Python with Numba, a specialized library to make ``just-in-time compilation'' (JIT).\footnote{\url{http://numba.pydata.org/}} We used the version 0.41 of Numba in these tests. 
		
\begin{figure} 
			\includegraphics[scale=.36]{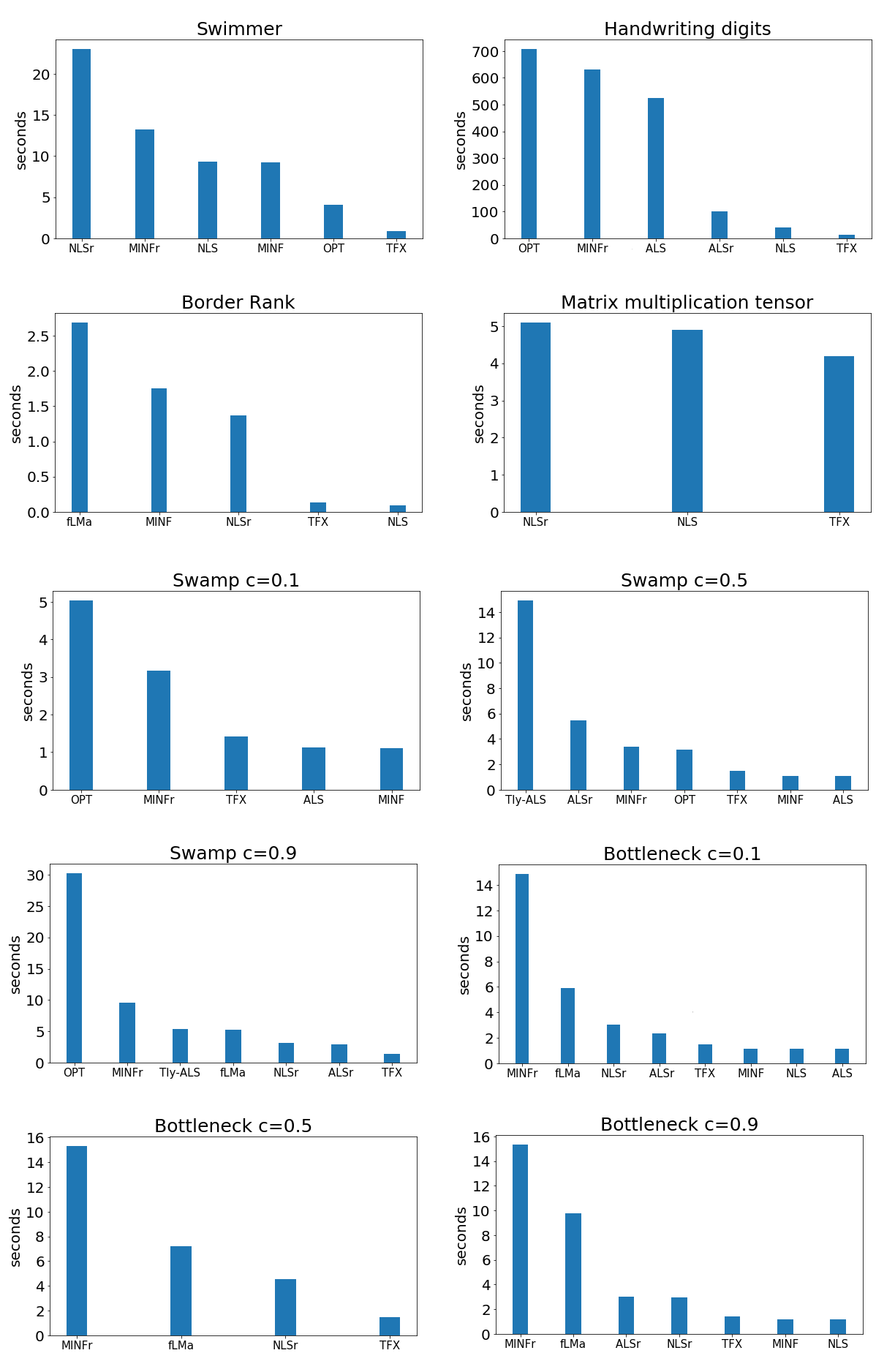}
			\caption{\footnotesize{Benchmarks of all tensors and all implementations.}}
			\label{benchmarks}
		\end{figure}

			In figure~\ref{benchmarks} there are some charts, each one shows the best running time of the algorithms with respect to each one of the tensors describe previously. If some algorithm is not included in a chart, it means that the algorithm was unable to achieve an acceptable error within the conditions described at the beginning if this section. 
		
		The first thing we should note is that there isn't a single algorithm which was able to match the speed of TFX in all tests. In the first four tests NLS algorithm could deliver an accepted solution in a reasonable time, but when collinearity was introduced this algorithm faced problems. When it produced a solution it was within a small time, but the algorithm can't be regarded as robust since it can fail to produce a solution in these cases. For this MINF and ALS looks more reliable while also being very fast. Still, these two algorithms couldn't produce an acceptable solution in some cases (for instance ALS only produced 4 solutions in the total of 12 tensors) and the existence of a solution can depend on the refinement option being active. Because of that their robustness is really only partial. About the timings, it is noticeable that TFX is the faster in most of the tests, and when it is not the faster it is close enough to the faster (in fact, every time an algorithm was faster than TFX it was by a difference of a fraction of a second). The algorithms NLS, ALS and MINF also are very fast (their version without refinement). NLS was faster than TFX two times, ALS was faster three times, and MINF was faster four times. Apart from that, none of the other algorithms seemed to stand out for some test.
		
		These results may lead the reader to the conclusion that Tensor Fox is the fastest implementation of all. This is partially true as we will see in the next experiment. This time, instead of making all algorithms to match Tensor Fox's accuracy, we run all with default options and show their box plots for accuracy and running time. With this we have a better understanding about the usual accuracy and usual running time of each algorithm. A box plot is used to summarize the main statistics of a random variable in a single image as~\ref{boxplot}.
		
		\begin{figure} 
			\hspace{3cm}\includegraphics[scale=.65]{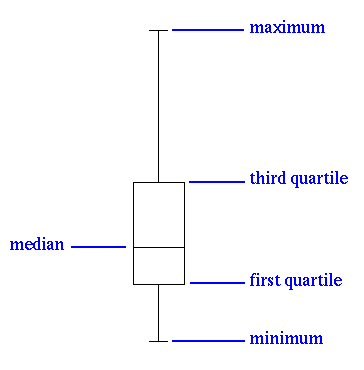}
			\caption{\footnotesize{The box plot is a standardized way of displaying the distribution of data based on the five number summary: minimum, first quartile, median, third quartile, and maximum.}}
			\label{boxplot}
		\end{figure}
		
		In the next figures we have the box plots of all algorithms, 20 runs for each tensor. An algorithm is omitted only when its running time is too big compared to the other or when the memory size to compute the CPD is too large (more than 7 GB in this context). Remember that in the previous benchmarks we manipulated the iteration of the algorithms in order to match Tensor Fox's accuracy with minimal time. Besides the Tensor Fox default, we considered the inverse procedure and manipulate Tensor Fox iterations\footnote{More precisely, we change the number of CG iterations and activate refinement sometimes, so the changes are similar to Tensorlab's.} to match the best solution (minimal error) of all other algorithms. This modified algorithm is labelled as TFXm.  
		
		\begin{figure} 
			\hspace{-3cm}\includegraphics[scale=.55]{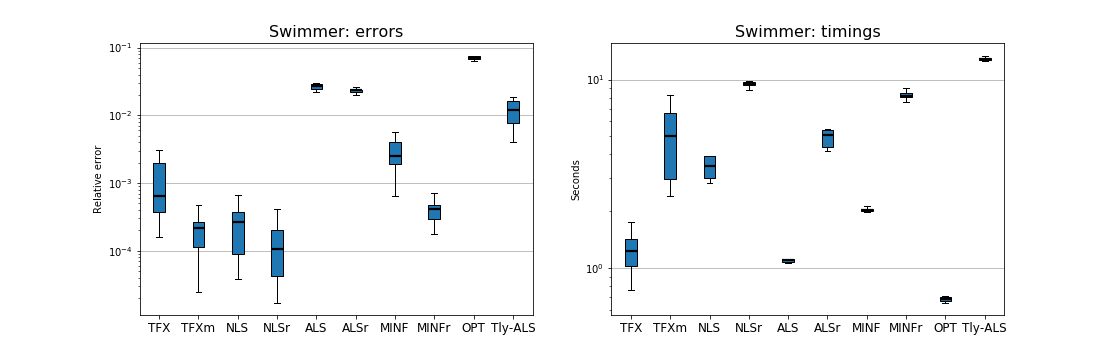}
			\caption{\footnotesize{Box plots of the swimmer tensor.}}
			\label{boxplot}
		\end{figure}
		
		\begin{figure} 
			\hspace{-3cm}\includegraphics[scale=.55]{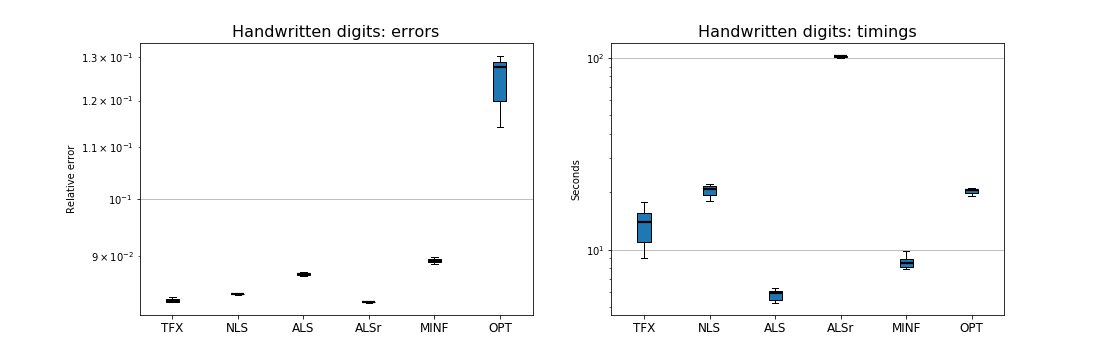}
			\caption{\footnotesize{Box plots of the handwritten digits tensor.}}
			\label{boxplot}
		\end{figure}
		
		\begin{figure} 
			\hspace{-3cm}\includegraphics[scale=.55]{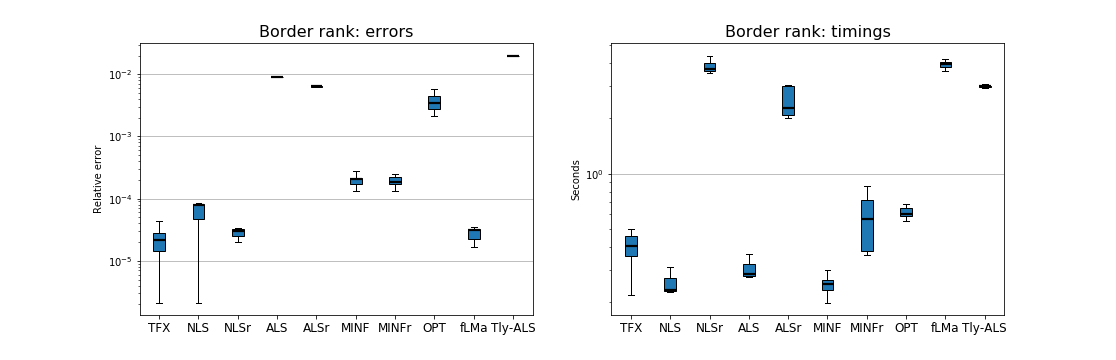}
			\caption{\footnotesize{Box plots of the border rank tensor.}}
			\label{boxplot}
		\end{figure}
		
		\begin{figure} 
			\hspace{-3cm}\includegraphics[scale=.55]{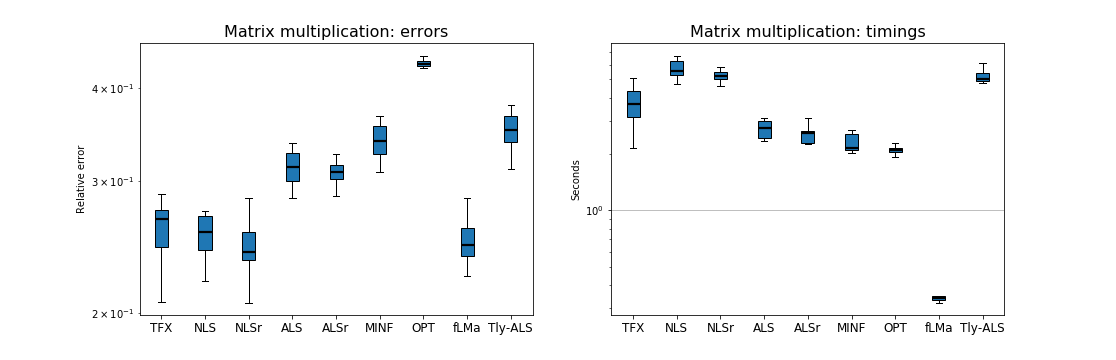}
			\caption{\footnotesize{Box plots of the $5 \times 5$ matrix multiplication tensor.}}
			\label{boxplot}
		\end{figure}
		
		\begin{figure} 
			\hspace{-3cm}\includegraphics[scale=.55]{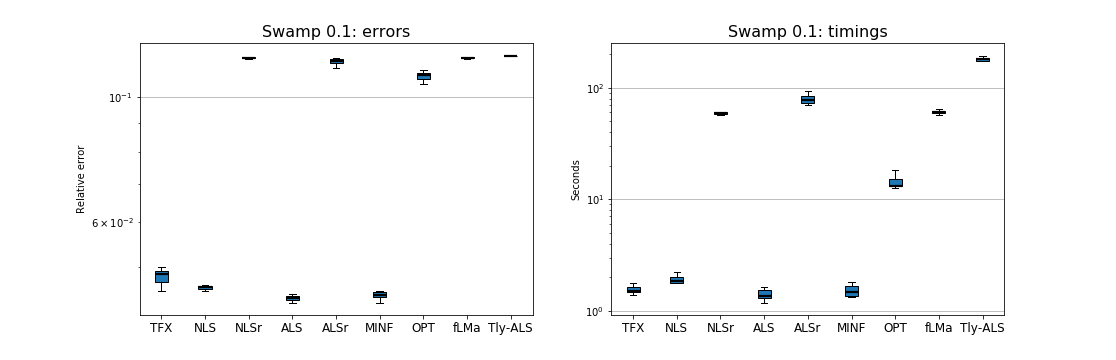}
			\caption{\footnotesize{Box plots of the swamp tensor with $c = 0.1$.}}
			\label{boxplot}
		\end{figure}
		
		\begin{figure} 
			\hspace{-3cm}\includegraphics[scale=.55]{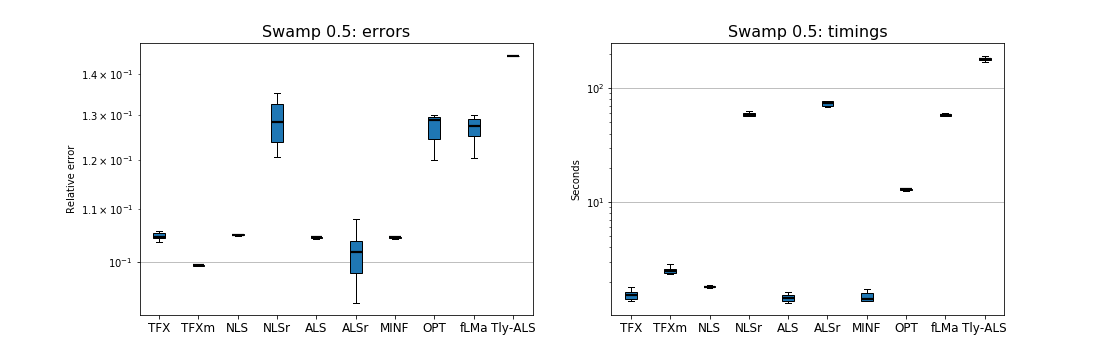}
			\caption{\footnotesize{Box plots of the swamp tensor with $c = 0.5$.}}
			\label{boxplot}
		\end{figure}
		
		\begin{figure} 
			\hspace{-3cm}\includegraphics[scale=.55]{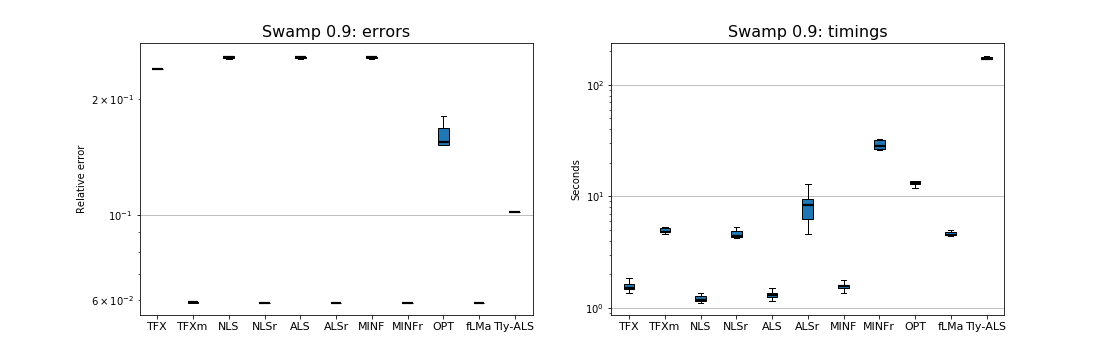}
			\caption{\footnotesize{Box plots of the swamp tensor with $c = 0.9$.}}
			\label{boxplot}
		\end{figure}
		
		\begin{figure} 
			\hspace{-3cm}\includegraphics[scale=.55]{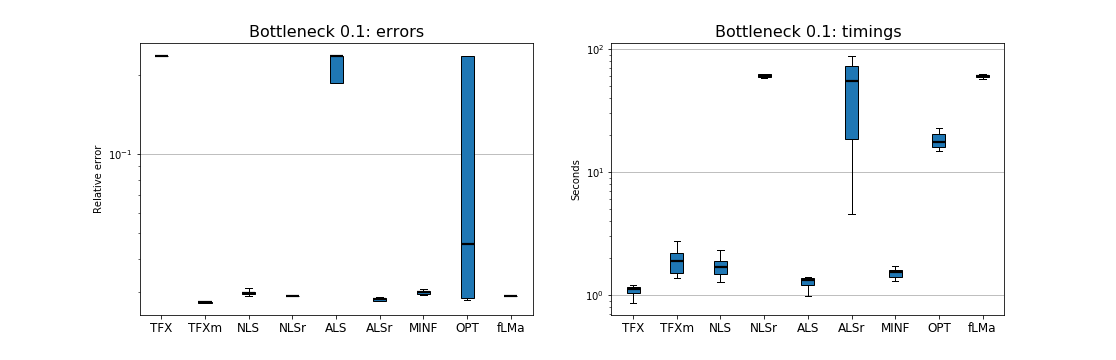}
			\caption{\footnotesize{Box plots of the double bottleneck tensor with $c = 0.1$.}}
			\label{boxplot}
		\end{figure}
		
		\begin{figure} 
			\hspace{-3cm}\includegraphics[scale=.55]{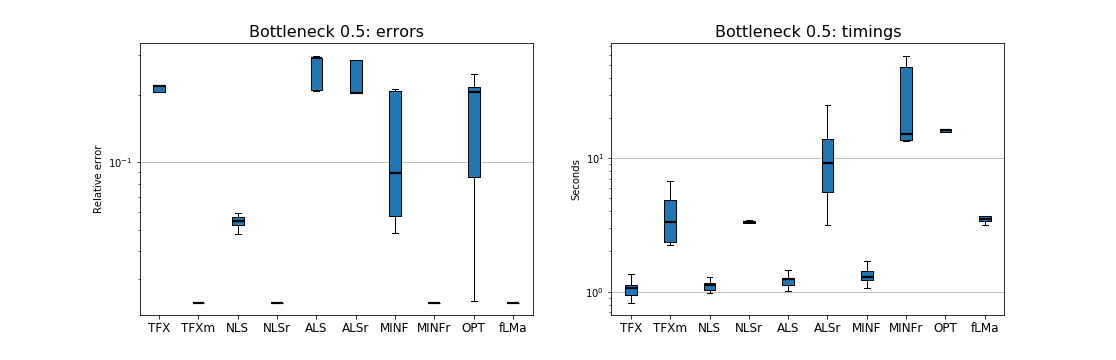}
			\caption{\footnotesize{Box plots of the double bottleneck tensor with $c = 0.5$.}}
			\label{boxplot}
		\end{figure}
		
		\begin{figure} 
			\hspace{-3cm}\includegraphics[scale=.55]{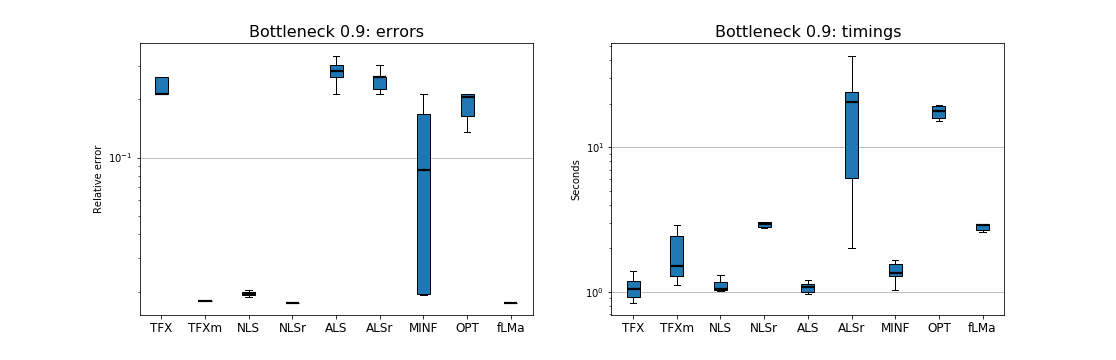}
			\caption{\footnotesize{Box plots of the double bottleneck tensor with $c = 0.9$.}}
			\label{boxplot}
		\end{figure}
		
		Tensor Fox default is competitive for all tensor except the double bottlenecks, where its accuracy is not optimal. In this case the modified versions are competitive, being as accurate and faster as the best algorithms. With these box plots we can see that no algorithm is optimal for all tensors, but the combinations TFX + TFXm or NLS + NLSr gives optimal solutions as fast as possible.
		
	\section{Tensor train and the CPD} 
		In table 4.4 we gave a summary of all costs necessary to compute a rank-$R$ CPD. The dominant cost comes from the dGN algorithm, which is dominated by a factor of $\mathcal{O}\left( \displaystyle L R \prod_{\ell=1}^L R_\ell \right)$ floats. With this cost we can see that the CPD computation suffers from the curse of dimensionality. More precisely, the cost increases exponentially with the tensor order. This limitation can be seen in every CPD implementation, and this is the main reason why this decomposition still is not widely used on large-scale problems. In the article \cite{cpd_tensortrain} they propose a way around this limitation, based on the \emph{tensor train decomposition} (TTD) \cite{tensor_train}. As a result, we get an algorithm capable of handling higher order tensors with low computational cost. First we will briefly introduce the TTD, following \cite{tensor_train}. After that we show how to connect the CPD with this new decomposition. This connection makes it possible to retrieve the CPD from the TTD. 
	
		\subsection{Tensor train decomposition}
			Let $\mathcal{T} \in \mathbb{R}^{I_1 \times \ldots I_L}$ be any order-$L$ tensor. The main idea of the TTD is that we can approximate $\mathcal{T}$ by a tensor $\tilde{\mathcal{T}}$ such that the entries of this tensor are given by
			\begin{equation}\label{ttd}
				\tilde{t}_{i_1 i_2 \ldots i_L} = \mathcal{G}^{(1)}(i_1) \mathcal{G}^{(2)}(i_2) \ldots \mathcal{G}^{(L)}(i_L)
			\end{equation}				
where each $\mathcal{G}^{(\ell)}(i_\ell)$ is a $r_{\ell-1} \times r_\ell$ full rank matrix for $\ell = 2 \ldots L$, with $r_1 = r_L = 1$. 

			Note that this definition generalizes the definition of rank one tensor. Instead of the coordinates of the tensor being given by a product of scalars, they are given by a product of matrices, where the first has one row and the last has one column so the result of the product still is a scalar. We can take one more step and consider each $\mathcal{G}^{(\ell)}(i_\ell)$ as the slice of a third order tensor with shape $r_{\ell-1} \times I_\ell \times r_\ell$. Denoting by $\mathcal{G}^{(\ell)}$ this tensor, we define $\mathcal{G}^{(\ell)}(i_\ell) = \mathcal{G}^{(\ell)}_{: \ i_\ell \ :}$ (the $i_\ell$-th vertical slice of $\mathcal{G}^{(\ell)}$). In coordinates we have that
			$$\tilde{t}_{i_1 i_2 \ldots i_L} = \sum_{j_0, j_1, \ldots, j_L} \mathcal{G}^{(1)}_{j_0 i_1 j_1} \ \mathcal{G}^{(2)}_{j_1 i_2 j_2} \ \ldots \ \mathcal{G}^{(L)}_{j_{L-1} i_L j_L}.$$ 
		
\begin{figure}[h] 
				\centering
				\includegraphics[scale=.35]{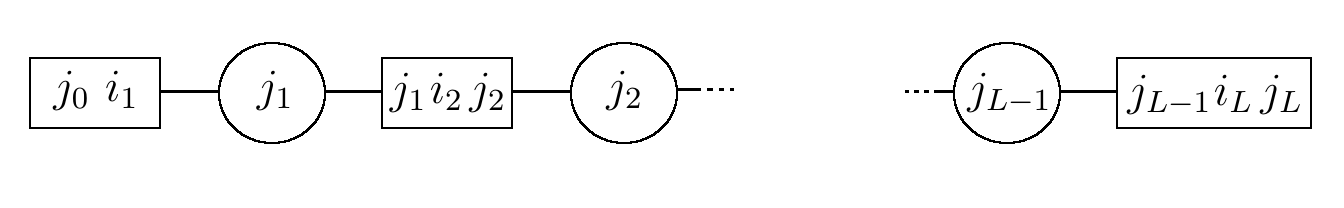}
				\caption{\footnotesize{Tensor-train network representation.}}
				\label{tensor-train}
			\end{figure}

			Since $\mathcal{G}^{(1)}$ and $\mathcal{G}^{(L)}$ are tensors of order $1 \times I_1 \times r_1$ and $r_{L-1} \times I_L \times 1$, respectively, they can be regarded as matrices (in particular, $j_0 = j_L = 1$). The ranks $r_\ell$ are called \emph{TT-ranks} and the three dimensional tensors $\mathcal{G}^{(\ell)}$ are called the \emph{cores} of the TTD. Note that the last index of each factor is the first of the next factor. This relation is illustrated in figure~\ref{tensor-train}. The rectangles contain spatial indexes (the indexes $j_{\ell-1}, i_\ell, j_\ell$), and the circles contain only the auxiliary nodes $j_\ell$, the ones representing the link between the cores. Quoting the original author of the TTD (Ivan Oseledets):\\
		
			``This picture looks like a train with carriages and links between them, and that justifies the name \emph{tensor train decomposition}, or simply \emph{TT-decomposition}.''\\ 
		
			The computation of the cores are made through successive SVDs and reshaping. The method used to reshape the matrices in the original article is the MATLAB reshape function, but it can be different as long as it is consistent with the other reshapes of the algorithm. In particular, note that in the first iteration of the TT-SVD loop we have $\textbf{M} = \mathcal{T}_{(1)}$, the first unfolding of $\mathcal{T}$. Thus the reshape function must be consistent with the unfolding function. The algorithm is presented without truncation in this context but the interested reader can refer to the mentioned article to see the proofs and the more general algorithm. The functions \verb|Unfolding| and \verb|SVD| were already introduce, while \verb|Reshape|$(\textbf{M}, dims)$ reshapes a certain array $\textbf{M}$ in order that the resulting array has the dimensions described by the tuple $dims$.\\
		
			\begin{algorithm}[TT-SVD]
				$ $\\
				\textbf{Input:} $\mathcal{T} \in \mathbb{R}^{I_1 \times \ldots \times I_L}$\vspace{.4cm}\\
				$ r_0 \leftarrow 1$\\
				$\textbf{M} \leftarrow \mathcal{T}$\\
				\verb|for | $\ell = 1 \ldots L-1$\\
					$\hspace{1cm} n_\ell \leftarrow \prod_{\ell' > \ell} I_{\ell'}$\\
					$\hspace{1cm} \textbf{M} \leftarrow \verb|Reshape|\left(\textbf{M}, (r_{\ell-1} I_\ell, n_\ell) \right)$\\
					$\hspace{1cm} \textbf{U}, \Sigma, \textbf{V} \leftarrow \verb|SVD|(\textbf{M})$\\
					$\hspace{1cm} r_\ell \leftarrow rank(\textbf{M})$\\
					$\hspace{1cm} \mathcal{G}^{(\ell)} \leftarrow \verb|Reshape|(\textbf{U}, (r_{\ell-1}, I_\ell, r_\ell))$\\
					$\hspace{1cm} \textbf{M} \leftarrow \Sigma \textbf{V}^T$\\
				$\mathcal{G}^{(L)} \leftarrow \textbf{M}$\vspace{.4cm}\\		
				\textbf{Output:} $\mathcal{G}^{(1)}, \ldots, \mathcal{G}^{(L)}$
			\end{algorithm}  
		
		\subsection{CPD tensor train}
			Given two tensors $\mathcal{T} \in \mathbb{R}^{I_1 \times \ldots \times I_\ell \times \ldots \times I_L}$ and $\mathcal{U} \in \mathbb{R}^{J_1 \times \ldots \times J_m \times \ldots \times J_M}$ such that $I_\ell = J_m$, where the modes $\ell,m$ may be different, the $\times_\ell^m$-contraction between $\mathcal{T}$ and $\mathcal{U}$ is the tensor $\mathcal{T} \times_m^\ell \mathcal{U} \in \mathbb{R}^{I_1 \times \ldots \times I_{\ell-1} \times I_{\ell+1} \times \ldots \times I_L \times J_1 \times \ldots \times J_{m-1} \times J_{m+1} \times \ldots \times J_M}$ given by
			$$\left( \mathcal{T} \times_\ell^m \mathcal{U} \right)_{i_1 \ldots i_{\ell-1} i_{\ell+1} \ldots i_L \ j_1 \ldots j_{m-1} j_{m+1} \ldots j_M} = \sum_{k=1}^{I_\ell} t_{i_1 \ldots i_{\ell-1} \ k \ i_{\ell+1} \ldots i_L} \cdot u_{j_1 \ldots j_{M-1} \ k \ j_{M+1} \ldots j_M}.$$			
Note that with the notion of contraction one can write the TTD in~\ref{ttd} more compactly as
			$$\mathcal{T} = \mathcal{G}^{(1)} \times_2^1 \mathcal{G}^{(2)} \times_3^1 \ldots \times_{L-1}^1\mathcal{G}^{(L-1)} \times_L^1 \mathcal{G}^{(L)},$$
where the operations are made from left to right but we suppressed the parenthesis to maintain a clean notation. Now we are able to prove our first theorem connecting the CPD and the TTD. Since the TTD assumes that all cores are of full rank, from this point we need to assume that all factor matrices of any CPD also are of full rank. This assumption is implicit in each of the next results. As we did in chapter 3, the objective of this part is to give the reader not only an understanding of the theoretical aspects, but also the computational aspects related of our problems, hence all proofs are given in coordinates.

			\begin{theorem}[I. Oseledets, 2010]\label{CPD_to_TT}
				Let $\mathcal{T} \in \mathbb{R}^{I_1 \times \ldots \times I_\ell \times \ldots \times I_L}$ be a rank-$R$ tensor with CPD given by
				$$\mathcal{T} = \Big( \textbf{W}^{(1)}, \ldots, \textbf{W}^{(L)} \Big) \cdot \mathcal{I}_{R \times \ldots \times R}.$$
			
				Then $\mathcal{T}$ admits a TTD of the form $\mathcal{T} = \mathcal{G}^{(1)} \times_2^1 \mathcal{G}^{(2)} \times_3^1 \ldots \times_{L-1}^1\mathcal{G}^{(L-1)} \times_L^1 \mathcal{G}^{(L)}$, where
				\begin{flalign*}			
					& \mathcal{G}^{(1)} = \textbf{W}^{(1)},\\
					& \mathcal{G}^{(\ell)} = \Big( \textbf{I}_R, \textbf{W}^{(\ell)}, \textbf{I}_R \Big) \cdot \mathcal{I}_{R \times R \times R}, \quad \ell = 2 \ldots L-1,\\
					& \mathcal{G}^{(L)} = \left( \textbf{W}^{(L)} \right)^T.\\
				\end{flalign*}
			\end{theorem} 
		
			\textbf{Proof:} Define the tensor $\mathcal{U} = \mathcal{G}^{(1)} \times_2^1 \mathcal{G}^{(2)} \times_3^1 \ldots \times_{L-1}^1\mathcal{G}^{(L-1)} \times_L^1 \mathcal{G}^{(L)}$ with cores as showed above. We will prove that $\mathcal{T} = \mathcal{U}$. First note that $\mathcal{G}_{ij}^{(1)} = \textbf{w}_{ij}^{(1)}$ and $\mathcal{G}_{ij}^{(L)} = \textbf{w}_{ji}^{(L)}$. For the other cores we have that
			$$\mathcal{G}^{(\ell)} = \Big( \textbf{W}^{(1)}, \ldots, \textbf{W}^{(L)} \Big) \cdot \mathcal{I}_{R \times \ldots \times R} = \sum_{r=1}^R \textbf{e}_r \otimes \textbf{w}_r^{(\ell)} \otimes \textbf{e}_r = $$
		
			$$ = \sum_{r=1}^R \textbf{e}_r \otimes \left( \sum_{i_\ell=1}^{I_\ell} w_{i_\ell r}^{(\ell)} \textbf{e}_{i_\ell} \right) \otimes \textbf{e}_r = \sum_{r=1}^R \sum_{i_\ell=1}^{I_\ell} w_{i_\ell r}^{(\ell)} \ \textbf{e}_r \otimes \textbf{e}_{i_\ell} \otimes \textbf{e}_r.$$
		
			In particular, this expansion of $\mathcal{G}^{(\ell)}$ shows that this tensor can be written as
			$$\mathcal{G}^{(\ell)} = \left\{ 
			\left[ \begin{array}{ccc}
				w_{1 1}^{(\ell)} & \ldots & w_{I_\ell 1}^{(\ell)}\\
				0 & \ldots & 0\\
				\vdots & & \vdots\\
				0 & \ldots & 0
			\end{array} \right], 
			\left[ \begin{array}{ccc}
				0 & \ldots & 0\\
				w_{1 2}^{(\ell)} & \ldots & w_{I_\ell 2}^{(\ell)}\\
				\vdots & & \vdots\\
				0 & \ldots & 0
			\end{array} \right],
			\ldots,
			\left[ \begin{array}{ccc}
				0 & \ldots & 0\\
				\vdots & & \vdots\\
				0 & \ldots & 0\\
				w_{1 r}^{(\ell)} & \ldots & w_{I_\ell r}^{(\ell)}\\
			\end{array} \right]
			\right\},$$
where each matrix is a frontal slice. Now we can just compute the coordinates of $\mathcal{U}$ and realize that they are equal to the $\mathcal{T}$ coordinates.

			$$u_{i_1 i_2 \ldots i_L} =  \sum_{j_0, j_1, \ldots, j_L} \mathcal{G}^{(1)}_{j_0 i_1 j_1} \ \mathcal{G}^{(2)}_{j_1 i_2 j_2} \ \ldots \ \mathcal{G}^{(L-1)}_{j_{L-2} i_{L-1} j_{L-1}} \ \mathcal{G}^{(L)}_{j_{L-1} i_L j_L} = $$
		
			\footnotesize
			$$\hspace{-0.4cm} = \sum_{j_1, \ldots, j_L} \left( \textbf{W}^{(1)} \right)_{i_1 j_1} \ \left( \sum_{r=1}^R \left( \textbf{e}_r \otimes \textbf{W}^{(2)} \otimes \textbf{e}_r \right)_{j_1 i_2 j_2} \right) \ \ldots \  \left( \sum_{r=1}^R \left( \textbf{e}_r \otimes \textbf{W}^{(L-1)} \otimes \textbf{e}_r \right)_{j_{L-2} i_{L-1} j_{L-1}} \right) \ \left( \textbf{w}_r^{(L)} \right)_{i_L j_{L-1}}.$$
			\normalsize
			Note that only the terms with $j_1 = j_2 = \ldots = j_{L_1} = r$, for $r = 1 \ldots R$, are not null. Therefore the expression above reduces to
		
			$$\sum_{r=1}^R w_{i_1 r}^{(1)} \ w_{i_2r}^{(2)} \ \ldots \ w_{i_{L-1} r}^{(L-1)} \ w_{i_L r}^{(L)},$$
from which we can conclude the equality, as desired. $\hspace{6cm} \square$\bigskip	

			The TTD based on the CPD given in this theorem will be called	\emph{CPD-train}. In particular, note that all $TT$-ranks of this decomposition are equal to $R$, the rank of the tensor. Now suppose we don't have the factors $\textbf{W}^{(\ell)}$ but the CPD-train holds. If we apply the TT-SVD algorithm to obtain the cores of the TTD we can expect to retrieve the factors of the CPD. It should be noted that the cores obtained are note necessarily in the form $\mathcal{G}^{(\ell)} = \Big( \textbf{I}_R, \textbf{W}^{(\ell)}, \textbf{I}_R \Big) \cdot \mathcal{I}_{R \times R \times R}$ since their representation are not unique (we can post multiply one core by an invertible matrix and pre multiply the next core for its inverse, therefore changing the cores but maintaining the TTD).
		
			The next lemma says that if we have the CPDs of two tensors $\mathcal{T}, \mathcal{U}$, and the last factor matrix of $\mathcal{T}$ is the transpose of the inverse of the first factor matrix of $\mathcal{U}$, then these factors are cancelled when contracted. This lemma is a preparation to prove next theorem, whose proof was only hinted in \cite{cpd_tensortrain}. Here we give a detailed proof with a computational flavor. To prove that lemma we rely on the fact that the $\times_L^1$-contraction of rank-1 terms (definition~\ref{contraction_def}) of the form $\textbf{a}^{(1)} \otimes \ldots \otimes \textbf{a}^{(L)}, \textbf{b}^{(1)} \otimes \ldots \otimes \textbf{b}^{(M)}$ (we are assuming that $\textbf{a}^{(L)}$ and $\textbf{b}^{(1)}$ have the same shape) is a tensor $\mathcal{T}$ given by
			$$\mathcal{T} = \langle \textbf{a}^{(L)}, \textbf{b}^{(1)} \rangle \ \textbf{a}^{(1)} \otimes \ldots \otimes \textbf{a}^{(L-1)} \otimes \textbf{b}^{(2)} \otimes \ldots \otimes \textbf{b}^{(M)}.$$ 
		
			\begin{lemma}
				Let two tensors $\mathcal{T} \in \mathbb{R}^{I_1 \times \ldots \times I_\ell \times \ldots \times I_L}$ and $\mathcal{U} \in \mathbb{R}^{J_1 \times \ldots \times J_m \times \ldots \times J_M}$ such that
				$$\mathcal{T} = \Big( \textbf{A}^{(1)}, \ldots, \textbf{A}^{(L)} \Big) \cdot \mathcal{I}_{R \times \ldots \times R}$$
for matrices $\textbf{A}^{(\ell)} \in \mathbb{R}^{I_\ell \times R}$, and 
				$$\mathcal{U} = \Big( \textbf{B}^{(1)}, \ldots, \textbf{B}^{(M)} \Big) \cdot \mathcal{I}_{R \times \ldots \times R}$$
for matrices $\textbf{B}^{(m)} \in \mathbb{R}^{I_m \times R}$. If $\textbf{A}^{(L)} = \big( \textbf{B}^{(1)} \big)^{-T}$ (which means $I_L = J_1 = R$), then
				$$\mathcal{T} \times_L^1 \mathcal{U} = \left( \Big( \textbf{A}^{(1)}, \ldots, \textbf{A}^{(L-1)}, \textbf{I}_R \Big) \cdot \mathcal{I}_{R \times \ldots \times R} \right) \times_L^1 \left( \Big( \textbf{I}_R, \textbf{B}^{(2)}, \ldots, \textbf{B}^{(M)} \Big) \cdot \mathcal{I}_{R \times \ldots \times R} \right).$$
			\end{lemma}
		
			\textbf{Proof:} First we expand the expression for the contraction, obtaining
			$$\mathcal{T} \times_L^1 \mathcal{U} = \left( \sum_{r=1}^R \textbf{a}_r^{(1)} \otimes \ldots \otimes \textbf{a}_r^{(L)} \right) \times_L^1 \left( \sum_{r'=1}^R \textbf{b}_{r'}^{(1)} \otimes \ldots \otimes \textbf{b}_{r'}^{(M)} \right) = $$
		
			$$ = \sum_{r, r'=1}^R \left( \textbf{a}_r^{(1)} \otimes \ldots \otimes \textbf{a}_r^{(L)} \right) \times_L^1 \left( \textbf{b}_{r'}^{(1)} \otimes \ldots \otimes \textbf{b}_{r'}^{(M)} \right) = $$
		
			$$ = \sum_{r, r' = 1}^R \langle \textbf{a}^{(L)}, \textbf{b}^{(1)} \rangle \ \textbf{a}_r^{(1)} \otimes \ldots \otimes \textbf{a}_r^{(L-1)} \otimes \textbf{b}_{r'}^{(2)} \otimes \ldots \otimes \textbf{b}_{r'}^{(M)}.$$
		
			For the moment, instead of writing $\textbf{a}_r^{(L)}$ for the $r$-th column of $\textbf{A}^{(L)}$, it will be convenient to use the Matlab convention and write $\textbf{a}_{:r}^{(L)}$. Furthermore, since $\textbf{a}_{:r}^{(L)}$ is the $r$-th row of $\left(\textbf{B}^{(L)}\right)^{-1}$ transposed, we have that
			$$\langle \textbf{a}_{:r}^{(L)}, \textbf{b}_{:r'}^{(1)} \rangle = \Big( \textbf{a}_{:r}^{(L)} \Big)^T \textbf{b}_{:r'}^{(1)} = $$
		
			$$ = \Big( \textbf{b}_{r:}^{(1)} \Big)^{-1} \textbf{b}_{:r'}^{(1)} = \delta_{r r'},$$
where $\textbf{b}_{r:}^{(1)}$ denotes the $r$-th row of $\textbf{B}^{(1)}$ and $\delta_{r r'}$ is the Kronecker delta of $r$ and $r'$. Thus the summation for $\mathcal{T} \times_L^1 \mathcal{U}$ reduces to
			$$\sum_{r=1}^R \langle \textbf{a}^{(L)}, \textbf{b}^{(1)} \rangle \ \textbf{a}_r^{(1)} \otimes \ldots \otimes \textbf{a}_r^{(L-1)} \otimes \textbf{b}_{r}^{(2)} \otimes \ldots \otimes \textbf{b}_{r}^{(M)}.$$
		
			It is not hard to see that 
			$$\left( \Big( \textbf{A}^{(1)}, \ldots, \textbf{A}^{(L-1)}, \textbf{I}_R \Big) \cdot \mathcal{I}_{R \times \ldots \times R} \right) \times_L^1 \left( \Big( \textbf{I}_R, \textbf{B}^{(2)}, \ldots, \textbf{B}^{(M)} \Big) \cdot \mathcal{I}_{R \times \ldots \times R} \right)$$ 
leads to the same expression, therefore they are the same tensor. $\hspace{4.2cm} \square$\bigskip

			\begin{theorem}[Y. Zniyed, R. Boyer, Andre L.F. de Almeida, G. Favier, 2018]
				Let $\mathcal{T} \in \mathbb{R}^{I_1 \times \ldots \times I_L}$ be a rank-$R$ tensor with TTD given by 
				$$\mathcal{T} = \mathcal{G}^{(1)} \times_2^1 \mathcal{G}^{(2)} \times_3^1 \ldots \times_{L-1}^1\mathcal{G}^{(L-1)} \times_L^1 \mathcal{G}^{(L)}$$ 
such that
				\begin{flalign*}
					& \mathcal{G}^{(1)} = \textbf{W}^{(1)} \big( \textbf{M}^{(1)} \big)^{-T}\\
					& \mathcal{G}^{(2)} = \Big( \textbf{M}^{(1)}, \textbf{W}^{(2)}, \textbf{M}^{(2)} \Big) \cdot \mathcal{I}_{R \times R \times R}\\
					& \mathcal{G}^{(3)} = \Big( \big( \textbf{M}^{(2)} \big)^{-T}, \textbf{W}^{(3)}, \textbf{M}^{(3)} \Big) \cdot \mathcal{I}_{R \times R \times R}\\
					& \vdots\\
					& \mathcal{G}^{(L-1)} = \Big( \big( \textbf{M}^{(L-2)} \big)^{-T}, \textbf{W}^{(L-1)}, \textbf{M}^{(L-1)} \Big) \cdot \mathcal{I}_{R \times R \times R}\\
					& \mathcal{G}^{(L)} = \big( \textbf{M}^{(L-1)} \big)^{-T} \big( \textbf{W}^{(L)} \big)^T.
				\end{flalign*}
Then $\mathcal{T} = \big( \textbf{W}^{(1)}, \ldots, \textbf{W}^{(L)} \big) \cdot \mathcal{I}_{R \times \ldots \times R}$ is a CPD for $\mathcal{T}$. All matrices $\textbf{M}^{(\ell)}$ are square $R \times R$ and invertible.
			\end{theorem}	
		
			\textbf{Proof:} The prove is based on successive applications of the last lemma. First note that
			$$\mathcal{G}^{(1)} \times_2^1 \mathcal{G}^{(2)} = \textbf{W}^{(1)} \big( \textbf{M}^{(1)} \big)^{-T} \times_2^1	 \left( \Big( \textbf{M}^{(1)}, \textbf{W}^{(2)}, \textbf{M}^{(2)} \Big) \cdot \mathcal{I}_{R \times R \times R} \right) = $$		
			$$ =  \textbf{W}^{(1)} \times_2^1 \left( \Big( \textbf{I}_R, \textbf{W}^{(2)}, \textbf{M}^{(2)} \Big) \cdot \mathcal{I}_{R \times R \times R} \right).$$
		
			Going to the next term gives
		
			$$\mathcal{G}^{(1)} \times_2^1 \mathcal{G}^{(2)} \times_3^1 \mathcal{G}^{(3)} = $$		
			$$ = \textbf{W}^{(1)} \times_2^1 \left( \Big( \textbf{I}_R, \textbf{W}^{(2)}, \textbf{M}^{(2)} \Big) \cdot \mathcal{I}_{R \times R \times R} \right) \times_3^1 \left( \Big( \big( \textbf{M}^{(2)} \big)^{-T}, \textbf{W}^{(3)}, \textbf{M}^{(3)} \Big) \cdot \mathcal{I}_{R \times R \times R} \right) = $$
			$$ = \textbf{W}^{(1)} \times_2^1 \left( \Big( \textbf{I}_R, \textbf{W}^{(2)}, \textbf{I}_R \Big) \cdot \mathcal{I}_{R \times R \times R} \right) \times_3^1 \left( \Big( \big( \textbf{I}_R, \textbf{W}^{(3)}, \textbf{M}^{(3)} \Big) \cdot \mathcal{I}_{R \times R \times R} \right).$$
		
			This pattern keeps going until we have
			$$\mathcal{G}^{(1)} \times_2^1 \mathcal{G}^{(2)} \times_3^1 \ldots \times_{L-1}^1 \mathcal{G}^{(L-1)} \times_L^1 \mathcal{G}^{(L)} = $$
			\footnotesize		
			$$\hspace{-0.5cm} = \textbf{W}^{(1)} \times_2^1 \left( \Big( \textbf{I}_R, \textbf{W}^{(2)}, \textbf{I}_R \Big) \cdot \mathcal{I}_{R \times R \times R} \right) \times_3^1 \ldots \times_{L-1}^1 \left( \Big( \textbf{I}_R, \textbf{W}^{(L-1)}, \textbf{M}^{(L-1)} \Big) \cdot \mathcal{I}_{R \times R \times R} \right) \times_L^1 \left( \big( \textbf{M}^{(L-1)} \big)^{-T} \big( \textbf{W}^{(L)} \big)^T \right) = $$
			\normalsize	
			$$ = \textbf{W}^{(1)} \times_2^1 \left( \Big( \textbf{I}_R, \textbf{W}^{(2)}, \textbf{I}_R \Big) \cdot \mathcal{I}_{R \times R \times R} \right) \times_3^1 \ldots \times_{L-1}^1 \left( \Big( \textbf{I}_R, \textbf{W}^{(L-1)}, \textbf{I}_R \Big) \cdot \mathcal{I}_{R \times R \times R} \right) \times_L^1 \big( \textbf{W}^{(L)} \big)^T.$$
		
			This is the decomposition showed in~\ref{CPD_to_TT}, which proves the theorem.$\hspace{3cm}\square$\bigskip
		
			We stated the theorem in a way it already gives insights about the computations to be made. More precisely, we should first compute a rank-$R$ CPD for the third order tensor $\mathcal{G}^{(2)}$. Since $\mathcal{G}^{(1)}$ is already know, we can obtain the first factor through the equality $\textbf{W}^{(1)} = \mathcal{G}^{(1)} \big( \textbf{M}^{(1)} \big)^T$. To compute $\textbf{W}^{(3)}$ we must compute a rank-$R$ CPD for $\mathcal{G}^{(3)}$ fixing the first factor to $\big( \textbf{M}^{(2)} \big)^{-T}$. We can keep going sequentially with this procedure, using the third factor of the previous CPD to construct the first factor of the next CPD. The last factor $\textbf{W}^{(L)}$ is then computed through the equality $\textbf{W}^{(L)} = \big( \mathcal{G}^{(L)} \big)^T \textbf{M}^{(L-1)}$ since $\textbf{M}^{(L-1)}$ will be already known at this point. This gives rise to the following algorithm.
		
			\begin{algorithm}[CPD-TTD]
				$ $\\
				\textbf{Input:} $\mathcal{T} \in \mathbb{R}^{I_1 \times \ldots \times I_L}, R$\vspace{.4cm}\\
				$\mathcal{G}^{(1)}, \ldots, \mathcal{G}^{(L)} \leftarrow \verb|TT-SVD|( \mathcal{T} )$\\
				$\textbf{M}^{(1)}, \textbf{W}^{(2)}, \textbf{M}^{(2)} \leftarrow \verb|CPD|( \mathcal{G}^{(2)}, R )$\\
				\verb|for | $\ell = 3 \ldots L-1$\\
					$\hspace{1cm} \textbf{W}^{(\ell)}, \textbf{M}^{(\ell)} \leftarrow \verb|Bi-CPD|\left( \big( \textbf{M}^{(\ell-1)} \big)^{-T}, \mathcal{G}^{(\ell)}, R \right)$\\
				$\textbf{W}^{(1)} \leftarrow \mathcal{G}^{(1)} \big( \textbf{M}^{(1)} \big)^T$\\
				$\textbf{W}^{(L)} \leftarrow \big( \mathcal{G}^{(L)} \big)^T \textbf{M}^{(L-1)}$\vspace{.4cm}\\		
				\textbf{Output:} $\textbf{W}^{(1)}, \ldots, \textbf{W}^{(L)}$
			\end{algorithm}
		
			\verb|CPD|$(\mathcal{T}, R)$ is any algorithm to compute a rank-$R$ CPD to $\mathcal{T}$ and \verb|Bi-CPD|$(\textbf{M}, \mathcal{T}, R)$ does the same job but with the first factor, $\textbf{M}$, fixed. Since we are talking about third order tensors, the latter algorithm just computes two factors, hence the name ``Bi-CPD''. Notice we can't use the conjugate descent gradient algorithm to compute each iteration of the dGN for the \verb|Bi-CPD| function. Suppose we are given the triple $\Big( \big( \textbf{M}^{(\ell-1)} \big)^{-T}, \textbf{W}^{(\ell)}, \textbf{M}^{(\ell)} \Big)$ respective to the CPD of $\mathcal{G}^{(\ell)}$ and now we have to fix $\big( \textbf{M}^{(\ell)} \big)^{-T}$ for the next CPD. At each iteration of the dGN we will be solving the system
			$$\left( \textbf{A}^T_{:,R^2 +1:} \textbf{A}_{:,R^2+1:} + \mu\textbf{D}_{:,R^2+1:} \right) \textbf{x}_{R^2+1:} = \textbf{A}^T \textbf{b} - \left( \textbf{A}^T_{:,:R^2} \textbf{A}_{:,:R^2} + \mu \textbf{D}_{:,:R^2} \right) \textbf{x}_{:R^2},$$
where $\textbf{x}_{: R^2} = vec\left( \big( \textbf{M}^{(\ell)} \big)^{-T} \right)$ is the fixed part. Notice we have to solve a system of $R^2 + R I_{\ell+1} + R^2$ equations and $R I_\ell + R^2$ variables, which can be solved by more than one algorithm.  
		
			The most relevant aspect of this algorithm is its cost, which is a great improvement compared to the costs showed in table 4.4. The higher costs comes from the TT-SVD algorithm, which amounts to computing $L-1$ SVDs of shapes $I_1 \times \displaystyle\prod_{\ell=2}^L$,  $R I_2 \times \displaystyle\prod_{\ell=3}^L I_\ell$,  $\ldots$, $R I_{L-1} \times I_L$, respectively. This has a total cost of
			$$\mathcal{O}\left( 2 I_1^2 \displaystyle\prod_{\ell=2}^L I_\ell + 2 I_1^3 \right) + \mathcal{O}\left( 2 R^2 I_2^2 \displaystyle\prod_{\ell=3}^L I_\ell + 2 R^3 I_2^3 \right) + \ldots + \mathcal{O}\left( 2 R^2 I_{L-1}^2 I_L + 2 R^3 I_{L-1}^3 \right) = $$
			$$ = \mathcal{O}\left( 2 \left( I_1^2 \prod_{\ell=2}^L I_\ell + I_1^3 + R^2 \left( \sum_{\ell=2}^{L-1} I_\ell^2 \prod_{\ell' = \ell+1}^L I_{\ell'} + R I_{\ell}^3 \right) \right) \right) = $$ 
			$$ = \mathcal{O}\left( 2 \left( I_1 \prod_{\ell=1}^L I_\ell + I_1^3 + R^2 \left( \sum_{\ell=2}^{L-1} I_\ell \prod_{\ell' = \ell}^L I_{\ell'} + R I_{\ell}^3 \right) \right) \right) \text{ flops}.$$
		
			Since $R$ is used as the rank for all cores of the TT-SVD, we must have $R \leq rank(T_{(\ell)}) \leq I_\ell$ for all $\ell$. We can put side by side the costs of computing these SVDs and the MLSVD computations of the previous algorithm, showed in table 4.18. At this point it should be clear the difference between both approaches. While the TT-SVD is always cutting a dimensions for each computed SVD, the MLSVD doesn't do this, it computes SVDs of unfolding with the size of the whole tensor $L$ times. The MLSVD suffers from the curse of dimensionality, whereas the TT-SVD doesn't (actually it still suffers, but very little and nothing compared to the MLSVD cost). 
		
			\begin{table}
				\centering
				\begin{tabular}{|c|c|c|}
					\hline
					& \textbf{TT-SVD} & \textbf{Compression}\\
					\hline
					\text{First SVD} & $2 I_1 \displaystyle\prod_{\ell=2}^L I_\ell + 2 I_1^3$ & $I_1 \displaystyle\prod_{\ell=1}^L I_\ell + I_1^3$\\
					\hline
					\text{Second SVD} & $2 I_1^2 I_2 \displaystyle\prod_{\ell=3}^L I_\ell + 2 I_1^3 I_2^3$ & $I_2 \displaystyle\prod_{\ell=1}^L I_\ell + I_2^3$\\
					\hline
					\vdots & \vdots & \vdots\\
					\hline
					$(L-1)$ \text{-th SVD} & $2 I_1^2 I_{L-1} I_L + 2 I_1^3 I_{L-1}^3$ & $I_{L-1} \displaystyle\prod_{\ell=1}^L I_\ell + I_{L-1}^3$\\
					\hline
					$L$ \text{-th SVD} & $0$ & $I_L \displaystyle\prod_{\ell=1}^L I_\ell + I_L^3$\\
					\hline
				\end{tabular}
				\caption{\footnotesize{Costs of TT-SVD vs. Compression}} \bigskip
			\end{table}		 
		
			The other costs consists of the computations of third order CPDs of shapes $R \times I_\ell \times R$, for $\ell = 2 \ldots L-1$.  With this we can see that the CPD-TTD approach is much faster than our previous algorithm, by orders of magnitude. To reinforce this claim we show some computational experiments with random tensors. We generate random tensors of shape $n \times n \times \ldots \times n$ and rank $R = 5$, where the entries of each factor matrix are drawn from the normal distribution (mean $0$ and variance $1$). First we consider fourth order tensors with shape $n \times n \times n \times n$, for $n = 10, 20, 30, 40, 50, 60, 70, 80$. Since the Tensorlab's NLS performed very well in the previous tests, we use only this one for Tensorlab and start with this algorithm, making $20$ computations for each dimension $n$ and averaging the errors and time.\footnote{In order to achieve the least possible difference between the errors, we accepted to discard 1 to 10 tests with bad approximations.} After that we run the other algorithms adjusting their tolerance in order to match the NLS results. 
  
\begin{figure}
				\centering
				\includegraphics[scale=.45]{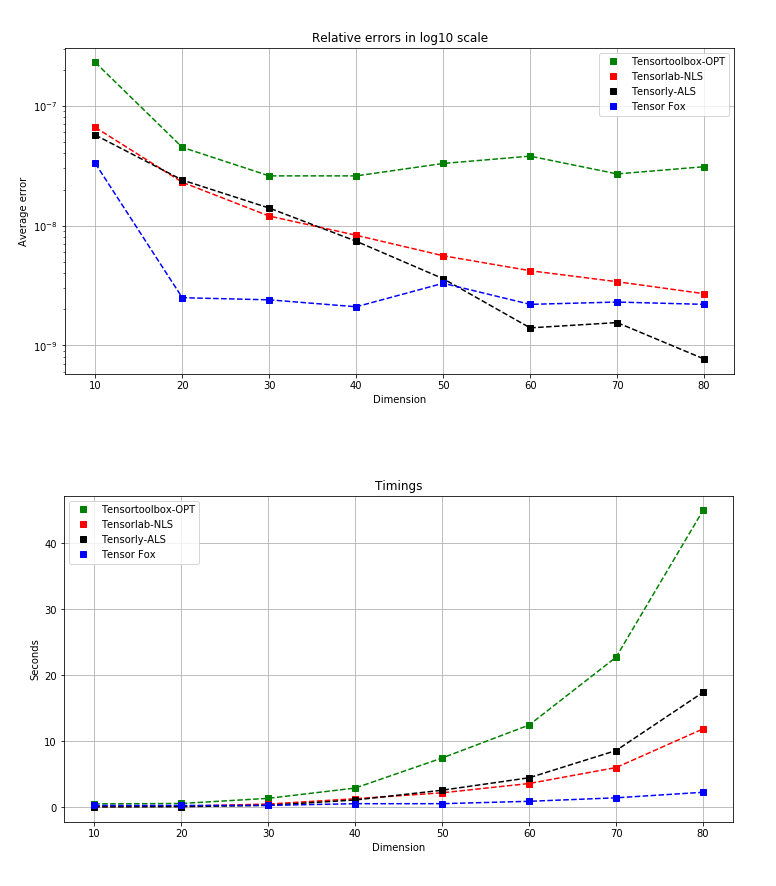}
				\caption{\footnotesize{Rank-$5$ CPD errors and timings of tensors with shape $n \times n \times n \times n$, for $n = 10, 20, \ldots, 70, 80$.}}
				\label{High order benchmarks varying the dimension}
			\end{figure}

 			In all tests we tried to choose the parameters in order to speed up the algorithms without losing accuracy. For example, we noticed that it was unnecessary to use compression, detection of structure and refinement for the NLS algorithm. These routines are very expensive and didn't bring much extra precision, so they were disabled in order to make the NLS computations faster. Similarly we used the initialization 'svd' for Tensorly because it proved to be faster than 'random', and we used the algorithm 'lbfgs' for Tensor Toolbox OPT. Finally, for Tensor Fox we just decreased its tolerance in order to match the precision given by the NLS algorithm. The results are showed in figure~\ref{High order benchmarks varying the dimension}.  
 		
			Next, we make the same procedure but this time we fixed $n$ to $n = 10$ and increased the order, from order $3$ to $8$. These last tests shows an important aspect of the CPD-TTD: it avoids the curse of dimensionality, whereas the other algorithms still suffers from that. We consider random rank-$5$ tensors of shape $10 \times 10 \times 10$, them $10 \times 10 \times 10 \times 10$, up to tensors of order $8$, i.e., with shape $\underbrace{10 \times 10 \times \ldots \times 10}_{8 \text{ times}}$, with the same distribution as before.\footnote{For anyone interested in reproducing these tests, the routine to generate these tensors can be found in \url{https://github.com/felipebottega/Tensor-Fox/blob/master/tests/gen_rand_tensor.py}.}
		
			\begin{figure}
				\centering
				\includegraphics[scale=.45]{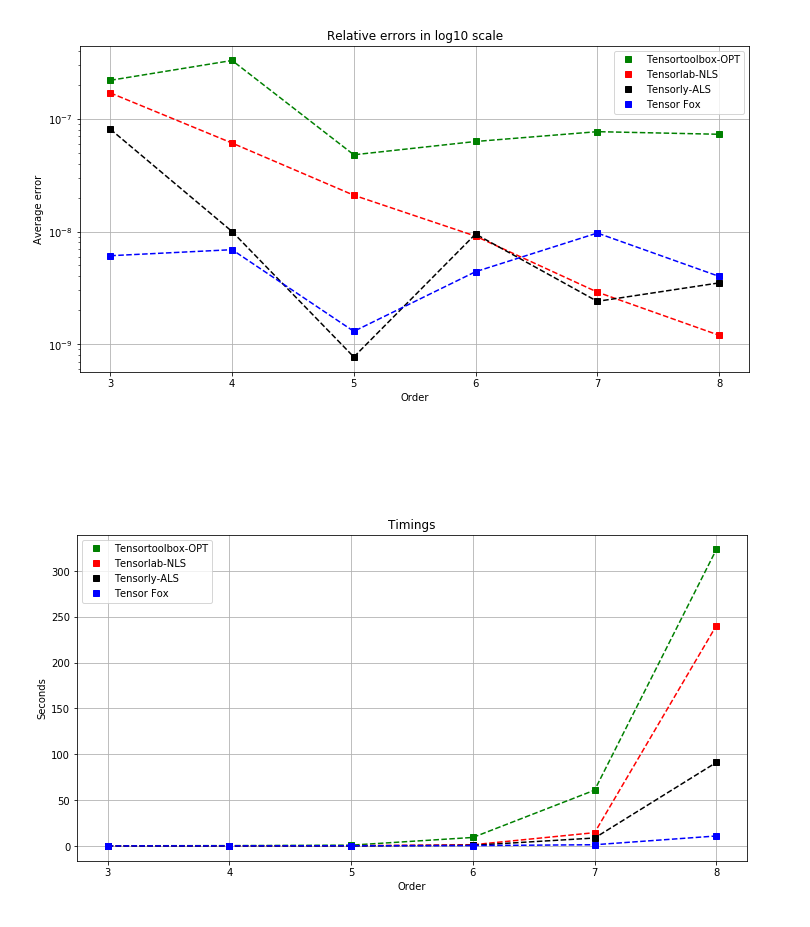}
				\caption{\footnotesize{Rank-$5$ CPD errors and timings of tensors with shape $10 \times 10 \times \ldots \times 10$ ($L$ times), for $L = 3, 4, \ldots, 8$.}}
				\label{High order benchmarks varying the order}
			\end{figure} 
		
			One limitation of the CPD-TTD is the rank itself. When constructing the cores of the tensor train we noted that we must have $R \leq \min_\ell I_\ell$, so this approach doesn't work for higher rank tensors. Since the CPD-TTD already is very fast compared to the other algorithms, we can afford to increase the cost a little when dealing with higher ranks. In the case $R > I_\ell$ for some $\ell$, we just increase this dimension to be of size $R$. The new values added to the tensor are random noise close to zero. For example, if we have a tensor $\mathcal{T}$ of shape $4 \times 5 \times 6 \times 7$ and $R = 6$, we increase the dimensions of $\mathcal{T}$ so we have a new tensor $\mathcal{T}'$ of shape $6 \times 6 \times 6 \times 7$ such that $\mathcal{T}'_{ijkl} = \mathcal{T}_{ijkl}$ for all $i=1 \ldots 4, j=1 \ldots 5, k=1 \ldots 6, l=1 \ldots 7$. The extra entries of $\mathcal{T}'$ won't affect too much the precision since they are very small, and since the CPD-TTD algorithm is already very fast, the CPD computation of this bigger tensor still is much faster than any other algorithm. We should point out that these new entries should be added after the MLSVD is computed, otherwise we are just introducing more complexity for nothing. After the CPD is computed we can truncate the CPD to its original dimensions. This approach works very well as long as $R$ is not too large or bigger than all dimensions (in practice we observed that we should have $R \leq I_\ell$ for at least one $\ell$). 	
			
	\section{Tensor Fox is not monotonic} \label{not_mon}
		Usually one expects that any minimization program produces a sequence of steps such that the corresponding sequence of errors decreases monotonically. In the ``warming up'' example, figure~\ref{warming_up-results} shows a sequence of errors which is not monotonic, sometimes the error increases, then it decreases back. This behavior was hinted at~\ref{unusual_steps} where we remarked that the way Tensor Fox handles the maximum number of CG iterations can introduce unusual steps sometimes.\footnote{These steps happen when \texttt{cg}\_\texttt{maxiter} is large, since the CG algorithm become unstable.} Now we make this statement clear. 
		
\begin{figure}[H]
				\centering
				\includegraphics[scale=.2]{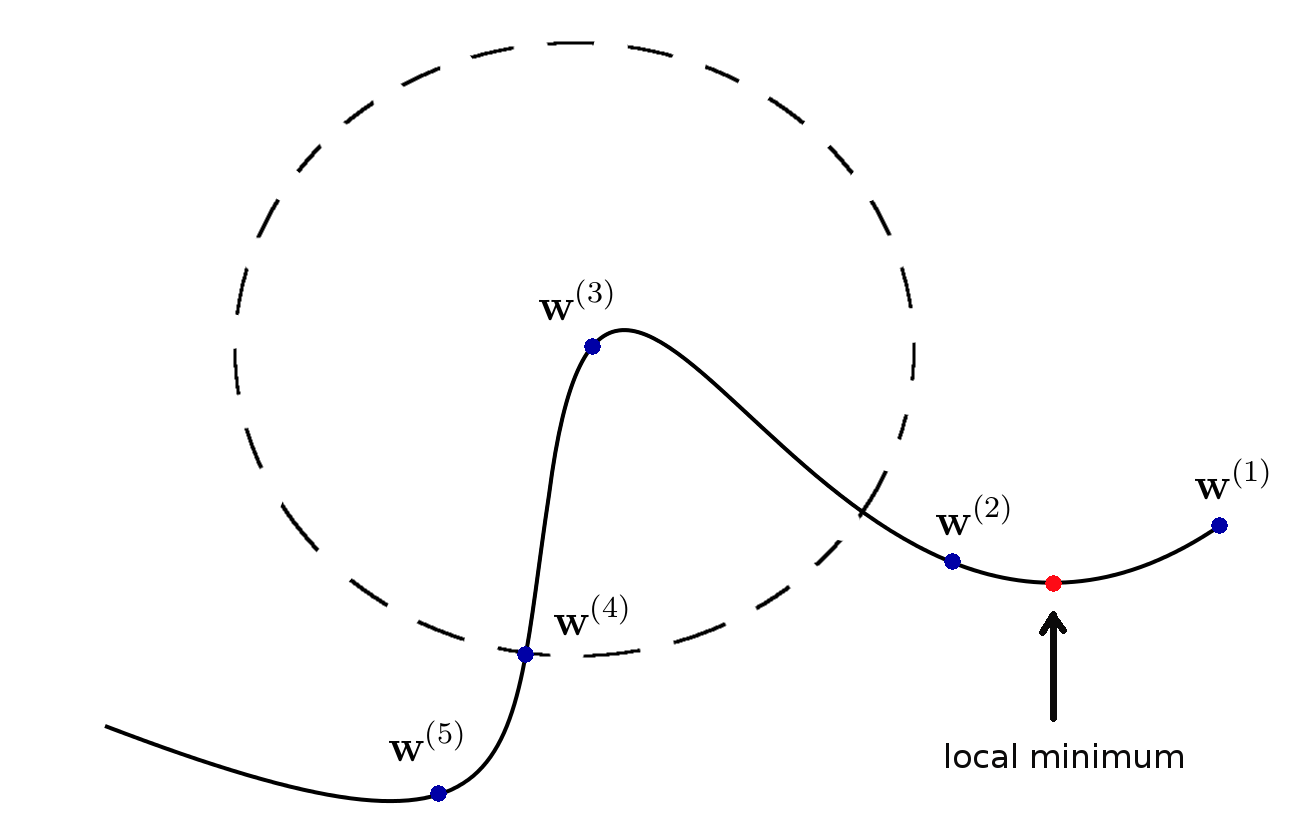}
				\caption{\footnotesize{After 2 iterations, the third one causes an increase in the error. However this allows the program to search in different regions (the circled region), which can lead the next iteration to be much better. In the worse case the next iteration is much likely to go back to the already located minimum.}}
				\label{Unusual steps}
			\end{figure}

		By allowing the CG method to perform a random number of steps, we observed that local minima are consistently avoided. Usually a poor step is followed by a very good step. The idea is that these steps go a little away from the local minimum, just enough so that the next iteration can search for better directions. If no new direction is found, the program goes back to the previous local minimum and probably will stay there. We say ``probably'' because there is always a little chance that an unusual step is taken, but most of the time the program are making usual steps. At the end of the day we conclude that this stochastic factor is successful when dealing with the high nonlinearity of the problem. This phenomenon is illustrated in figure~\ref{Unusual steps}. 	 
			
		  In~\ref{updates} we discussed the gain ratio $g$ and what it represents. Basically, large values of $g$ means that the error is being reduced substantially, $g = 0$ means that the error didn't change, and $g < 0$ means that the error increased. After the observations made here it is now clear that the unusual steps are responsible for making $g$ to be negative, and that the strategy to produce $\verb|cg_maxiter|$ (described in~\ref{dGN-cg}) is what causes this behavior. 
		  
\begin{figure}
				\centering
				\includegraphics[scale=.55]{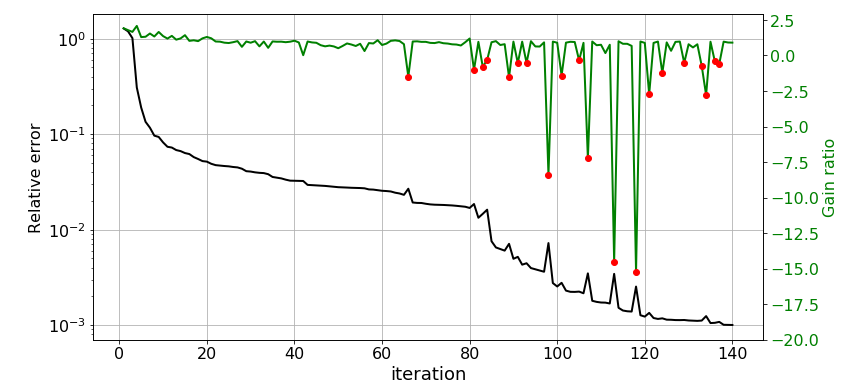}
				\caption{\footnotesize{The black curve is the represents the error, the green curve represents the gain ratio and the red dots are the points when the gain ratio became negative. Note that there is always a peak in the error when this happens. This is to be expected, but more remarkable is the fact that the error always decreases substantially after these points.}}
				\label{Gain ratio of the swimmer tensor}
			\end{figure}

		  Let's take a concrete example to observe this phenomenon happening in more details. In figure~\ref{Gain ratio of the swimmer tensor} we show the evolution of the error and gain ratio corresponding to a CPD computation of the swimmer tensor. We note that $g < 0$ precisely when the error increases, as expected. The interesting phenomenon here is the fact that the error always decreases substantially after these points. As explained, the program found a better direction to follow, which lead to better steps. Tensor Fox was tested with a fixed small number of CG iterations but this always lead to less accurate solutions, because the program got stuck at local minima.
		  
\begin{figure}
				\hspace{1cm}\includegraphics[scale=.8]{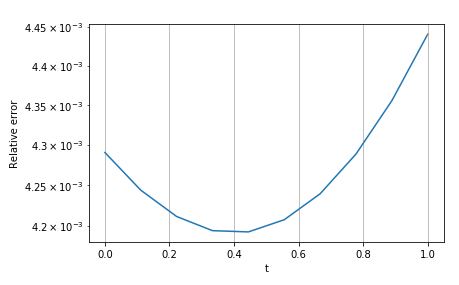}
				\caption{\footnotesize{Error curve $F(\textbf{w}_t)$. The error is minimal for $t \approx 0.4$ but we can see that the actual error (for $t = 1$) is bigger. At iteration 93 the gain ratio is $g = 0.974$, and at iteration 94 the gain ratio is $g = -0.518$. In both iterations the predicted error is $\mathcal{O}(10^{-8})$.}}
				\label{Error curve}
			\end{figure}

			Still with this same example, we can take iterations 93 and 94 to look in more detail. Iteration 93 is normal, with the error decreasing, but in iteration 94 the error increases. Let $F(\textbf{w}^{(93)})$ and $F(\textbf{w}^{(94)})$ be the errors at iterations 93 and 94, respectively. Remember that this notation was introduced in~\ref{min}. We can consider the line $\textbf{w}_t = \textbf{w}^{(93)} + t \left( \textbf{w}^{(94)} - \textbf{w}^{(93)} \right)$ between the approximated CPDs and analyze how is the curve $F(\textbf{w}_t)$, see figure~\ref{Error curve}
			
	\section{Regularization and preconditioning} \label{diag_reg}
		\subsection{Diagonal regularization}
		As suggested in~\ref{appenA}, our implementation of the preconditioner for Tensor Fox is the diagonal preconditioner. Remember we are using the following normal equations to compute the Gauss-Newton steps.
		$$(\textbf{J}_f^T \textbf{J}_f + \mu \textbf{D}) \textbf{x} = \textbf{J}_f^T \textbf{b}$$ 
where $\textbf{J}_f = \textbf{J}_f(\textbf{w}^{(k)})$, $\textbf{x} = \textbf{w} - \textbf{w}^{(k)}$, $\textbf{b} = - f(\textbf{w}^{(k)})$. For simplicity we start with the case of third order tensors, so $\textbf{w} = [vec(\textbf{X})^T, vec(\textbf{Y})^T, vec(\textbf{Z})^T ]^T$, where $\textbf{X} \in \mathbb{R}^{m \times R}, \textbf{Y} \in \mathbb{R}^{n \times R}, \textbf{Z}^{p \times R}$ are the corresponding factor matrices at the current step.

		This preconditioner is ideal to use when the approximated Hessian $\textbf{A}^T \textbf{A}$ is diagonally dominant. However this may not always be the case. We can overcome this problem by choosing a suitable regularization matrix $\textbf{D}$. The idea is to force $\textbf{J}_f^T \textbf{J}_f + \mu \textbf{D}$ to be diagonally dominant, then the preconditioner matrix is applied. The first row of $\textbf{J}_f^T \textbf{J}_f$ is the first row of the matrix
		\scriptsize
		$$\hspace{-.5cm} \left[ \langle \textbf{Y}_1, \textbf{Y}_1 \rangle \langle \textbf{Z}_1, \textbf{Z}_1 \rangle \textbf{I}_m \ \ \ldots \ \ \langle \textbf{Y}_1, \textbf{Y}_R \rangle \langle \textbf{Z}_1, \textbf{Z}_R \rangle \textbf{I}_m, \ \ \langle \textbf{Z}_1, \textbf{Z}_1 \rangle \textbf{X}_1 \textbf{Y}_1^T, \ \ \ldots \ \ \langle \textbf{Z}_1, \textbf{Z}_R \rangle \textbf{X}_R \textbf{Y}_1^T, \ \ \langle \textbf{Y}_1, \textbf{Y}_1 \rangle \textbf{X}_1 \textbf{Z}_1^T, \ \ \ldots \ \ \langle \textbf{Y}_1, \textbf{Y}_R \rangle \textbf{X}_R \textbf{Z}_1^T \right].$$
		\normalsize
where $\textbf{X}_r, \textbf{Y}_r, \textbf{Z}_r$ denotes the $r$-th column of $\textbf{X}, \textbf{Y}, \textbf{Z}$, respectively. We want to add a new term $\gamma$ to the first entry $\langle \textbf{Y}_1, \textbf{Y}_1 \rangle \langle \textbf{Z}_1, \textbf{Z}_1 \rangle$ such that 
		$$\langle \textbf{Y}_1, \textbf{Y}_1 \rangle \langle \textbf{Z}_1, \textbf{Z}_1 \rangle + \gamma \geq | \langle \textbf{Y}_1, \textbf{Y}_r \rangle \langle \textbf{Z}_1, \textbf{Z}_r \rangle |, \ \  | \langle \textbf{Z}_1, \textbf{Z}_r \rangle x_{1r} y_{j1} |, \ \  | \langle \textbf{Y}_1, \textbf{Y}_r \rangle x_{1r} z_{k1} |$$
for all $r, j, k$. For this it is sufficient to define
		$$\gamma = \| \textbf{Y}_1 \| \| \textbf{Z}_1 \| \max\left\{ \max_r \{ \| \textbf{Y}_r \| \| \textbf{Z}_r \| \}, \max_r \{ \| \textbf{X}_r \| \| \textbf{Z}_r \| \}, \max_r \{ \| \textbf{X}_r \| \| \textbf{Y}_r \| \} \right\}.$$
		This choice is cheap to compute and it also works for all first $m$ rows of $\textbf{J}_f^T \textbf{J}_f$. For the next $m$ rows we can apply a similar idea to obtain the term 
		$$\gamma = \| \textbf{Y}_2 \| \| \textbf{Z}_2 \| \max\left\{ \max_r \{ \| \textbf{Y}_r \| \| \textbf{Z}_r \| \}, \max_r \{ \| \textbf{X}_r \| \| \textbf{Z}_r \| \}, \max_r \{ \| \textbf{X}_r \| \| \textbf{Y}_r \| \} \right\}.$$
		
		Now we generalize this idea and introduce the notations
		$$\gamma_\textbf{X}^{(r)} = \|Y_r\| \|Z_r\| \max \Big\{ \max_{r'} \|Y_{r'}\| \|Z_{r'}\|, \ \max_{r'} \|X_{r'}\| \|Z_{r'}\|, \ \max_{r'} \|X_{r'}\| \|Y_{r'}\| \Big\},$$	
		$$\gamma_\textbf{Y}^{(r)} = \|X_r\| \|Z_r\| \max \Big\{ \max_{r'} \|Y_{r'}\| \|Z_{r'}\|, \ \max_{r'} \|X_{r'}\| \|Z_{r'}\|, \ \max_{r'} \|X_{r'}\| \|Y_{r'}\| \Big\},$$	
		$$\gamma_\textbf{Z}^{(r)} = \|X_r\| \|Y_r\| \max \Big\{ \max_{r'} \|Y_{r'}\| \|Z_{r'}\|, \ \max_{r'} \|X_{r'}\| \|Z_{r'}\|, \ \max_{r'} \|X_{r'}\| \|Y_{r'}\| \Big\}.$$
	
		From these values we finally define the regularization matrix		
		$$\textbf{D} = 
		\left[ \begin{array}{ccccccccc}
			\gamma_\textbf{X}^{(1)} \textbf{I}_m & & & & & & & &\\
			& \ddots & & & & & & &\\
			& & \gamma_\textbf{X}^{(R)} \textbf{I}_m & & & & & &\\
			& & & \gamma_\textbf{Y}^{(1)} \textbf{I}_n & & & & &\\
			& & & & \ddots & & & &\\
			& & & & & \gamma_\textbf{Y}^{(R)} \textbf{I}_n & & &\\
			& & & & & & \gamma_\textbf{Z}^{(1)} \textbf{I}_p & &\\
			& & & & & & & \ddots & \\
			& & & & & & & & \gamma_\textbf{Z}^{(R)} \textbf{I}_p
		\end{array} \right].$$
		
		Note that $\textbf{J}_f^T \textbf{J}_f + \textbf{D}$ is diagonally dominant, as desired, but now we are solving the regularized problem
		$$\min_\textbf{w} \| \textbf{J}_f^T \textbf{J}_f \textbf{w} - \textbf{J}_f^T \textbf{b} \|^2 + \| \textbf{D} \textbf{w} \|^2.$$
		This introduces a penalization over the solution, which can lead to local minimum, and this will have influence on the Gauss-Newton iterations. That is the reason for introducing the damping factor $\mu > 0$. Usually in the first iterations $\mu > 0$ has a moderate size, thus $\textbf{J}_f^T \textbf{J}_f + \mu \textbf{D}$ still is very regularized. On the other hand, the solution obtained with this system will be closer to the objective point than the initial guess anyway. Besides that, this will be done very fast since the preconditioner is being taken full advantage of. As the iteration continues, we have $\mu \to 0$, thus the regularization has less effect, and the solutions computed starts to get closer to the real objective solution. This also implies that at this stage the preconditioner will have less effect, which means we need to perform more CG iterations in order to obtain meaningful solutions. Note how this is in accordance with the number of CG iterations describe in~\ref{dGN-cg}, since the strategy adopted foresee the necessity to increase the number of CG iterations as we keep making Gauss-Newton iterations. 
		
		To finish, we conclude that the preconditioner $\textbf{M}$ is given by
		$$\textbf{M} = \left[ 
		\begin{array}{ccc}
			\textbf{D}_\textbf{X} & &\\
			& \textbf{D}_\textbf{Y} &\\
			& & \textbf{D}_\textbf{Z}
		\end{array} \right],$$ 
where
		$$\textbf{D}_\textbf{X} = \left[ 
		\begin{array}{ccc}
			\left( \|\textbf{Y}_1\|^2 \|\textbf{Z}_1\|^2 + \mu \gamma_\textbf{X}^{(1)} \right) \textbf{I}_m & &\\
			& \ddots & \\
			& & \left( \|\textbf{Y}_R\|^2 \|\textbf{Z}_R\|^2 + \mu \gamma_\textbf{X}^{(R)} \right) \textbf{I}_m
		\end{array} \right],$$
	
		$$\textbf{D}_\textbf{Y} = \left[ 
		\begin{array}{ccc}
			\left( \|\textbf{X}_1\|^2 \|\textbf{Z}_1\|^2 + \mu \gamma_\textbf{Y}^{(1)} \right) \textbf{I}_n & &\\
			& \ddots & \\
			& & \left( \|\textbf{X}_R\|^2 \|\textbf{Z}_R\|^2 + \mu \gamma_\textbf{Y}^{(R)} \right) \textbf{I}_n
		\end{array} \right],$$
	
		$$\textbf{D}_\textbf{Z} = \left[ 
		\begin{array}{ccc}
			\left( \|\textbf{X}_1\|^2 \|\textbf{Y}_1\|^2 + \mu \gamma_\textbf{Z}^{(1)} \right) \textbf{I}_p & &\\
			& \ddots & \\
			& & \left( \|\textbf{X}_R\|^2 \|\textbf{Y}_R\|^2 + \mu \gamma_\textbf{Z}^{(R)} \right) \textbf{I}_p
		\end{array} \right].$$

		For the reader interested in the general case, first we need to know what is the diagonal of $\textbf{J}_f^ T \textbf{J}_f$. From theorem~\ref{JfT_Jf} we know that
		
		$$\textbf{J}_f^T \textbf{J}_f = 
					\left[ \begin{array}{ccc}
						\textbf{H}_{11} & \ldots & \textbf{H}_{1L}\\
						\vdots & & \vdots\\
						\textbf{H}_{L1} & \ldots & \textbf{H}_{LL}
		\end{array} \right],$$
where
		$$\textbf{H}_{\ell' \ell'} = 
					\left[ \begin{array}{ccc}
					\displaystyle \prod_{\ell \neq \ell'} \omega_{11}^{(\ell)} \cdot \textbf{I}_{I_{\ell'}} & \ldots & \displaystyle \prod_{\ell \neq \ell'} \omega_{1R}^{(\ell)} \cdot \textbf{I}_{I_{\ell'}}\\
					\vdots & & \vdots\\
					\displaystyle \prod_{\ell \neq \ell'} \omega_{R1}^{(\ell)} \cdot \textbf{I}_{I_{\ell'}} & \ldots & \displaystyle \prod_{\ell \neq \ell'} \omega_{RR}^{(\ell)} \cdot \textbf{I}_{I_{\ell'}}
		\end{array} \right]$$
are the diagonal blocks of $\textbf{J}_f^T \textbf{J}_f$. The diagonal of $\textbf{H}_{\ell' \ell'}$ is 

		$$\textbf{D}_{\ell' \ell'} = 
					\left[ \begin{array}{ccc}
					\displaystyle \prod_{\ell \neq \ell'} \omega_{11}^{(\ell)} \cdot \textbf{I}_{I_{\ell'}} & & \\
					 & \ddots & \\
					 & & \displaystyle \prod_{\ell \neq \ell'} \omega_{RR}^{(\ell)} \cdot \textbf{I}_{I_{\ell'}}
		\end{array} \right] = \left[ \begin{array}{ccc}
					\displaystyle \prod_{\ell \neq \ell'} \| \textbf{w}_1^{(\ell)} \|^2 \cdot \textbf{I}_{I_{\ell'}} & & \\
					 & \ddots & \\
					 & & \displaystyle \prod_{\ell \neq \ell'} \| \textbf{w}_R^{(\ell)} \|^2 \cdot \textbf{I}_{I_{\ell'}}
		\end{array} \right],$$		
hence the diagonal of $\textbf{J}_f^T \textbf{J}_f$ is the matrix
		
		$$\textbf{D} = \left[ 
		\begin{array}{ccc}
			\textbf{D}_{11} & &\\
			& \ddots &\\
			& & \textbf{D}_{LL}
		\end{array} \right].$$
		
		In the case the factors are norm-balanced, we know that $\| \textbf{w}_r^{(1)} \| = \ldots = \| \textbf{w}_r^{(L)} \|$ for all $r=1 \ldots R$. Denote this norm by $\alpha_r$. Then we have that $\displaystyle\prod_{\ell \neq \ell'} \| \textbf{w}_r^{(\ell)} \|^2 = \alpha_r^{2(L-1)}$. Therefore all matrices $\textbf{D}_{\ell' \ell'}$ are equal and
		
		$$\textbf{D} = \left[
		\begin{array}{ccc}
			\left[
			\begin{array}{ccc}
				\alpha_1^{2(L-1)} \cdot \textbf{I}_{I_1} & &\\
				& \ddots &\\
				& & \alpha_R^{2(L-1)} \cdot \textbf{I}_{I_1}
			\end{array} \right] & &\\
			& \ddots &\\
			& & \left[
			\begin{array}{ccc}
				\alpha_1^{2(L-1)} \cdot \textbf{I}_{I_L} & &\\
				& \ddots &\\
				& & \alpha_R^{2(L-1)} \cdot \textbf{I}_{I_L}
			\end{array} \right]
		\end{array} \right].$$
		
		This may be useful if one want to use a diagonal preconditioner $\textbf{M}$ such that the matrix $\textbf{M}^{-\frac{1}{2}} ( \textbf{J}_f^T \textbf{J}_f + \mu \textbf{D} ) \textbf{M}^{-\frac{1}{2}}$ is unit diagonal (but not necessarily diagonally dominant). We remark that this kind of preconditioner may be enough since the balancing of the factor matrices (see section~\ref{updates}) already helps $\textbf{J}_f^T \textbf{J}_f$ to be more diagonally dominant.	
		
		\subsection{Computational experiments}
			To simplify notation let $\textbf{A} = \textbf{J}_f, \textbf{H} = \textbf{A}^T \textbf{A}, \textbf{b} = - f$, where we suppress the argument $\textbf{w}^{(k)}$. Tensor Fox use equation~\ref{dGN-eq} for the CG algorithm, however it is more usual to use the Levenberg-Marquardt equation, 
			\begin{equation} \label{LM}
				\text{diag}(\textbf{H} + \mu \textbf{I})^{-1/2} \cdot (\textbf{H} + \mu \textbf{I}) \cdot \text{diag}(\textbf{H} + \mu \textbf{I})^{-1/2} \textbf{x} = \text{diag}(\textbf{H} + \mu \textbf{I})^{-1/2} \textbf{A}^T \textbf{b}. 
			\end{equation}
			
			The only difference is that Tensor Fox uses $\textbf{D}$ instead of the identity matrix. We know that the condition number of $\text{diag}(\textbf{H} + \mu \textbf{I})^{-1/2} \cdot (\textbf{H} + \mu \textbf{I}) \cdot \text{diag}(\textbf{H} + \mu \textbf{I})^{-1/2}$ must be small to solve equation~\ref{LM} efficiently. Although the Levenberg-Marquardt formulation is used in many implementations of the dGN, we will see that our diagonal matrix produces better results.
			
\begin{figure}[H]
				\hspace{-1cm}\includegraphics[scale=.45]{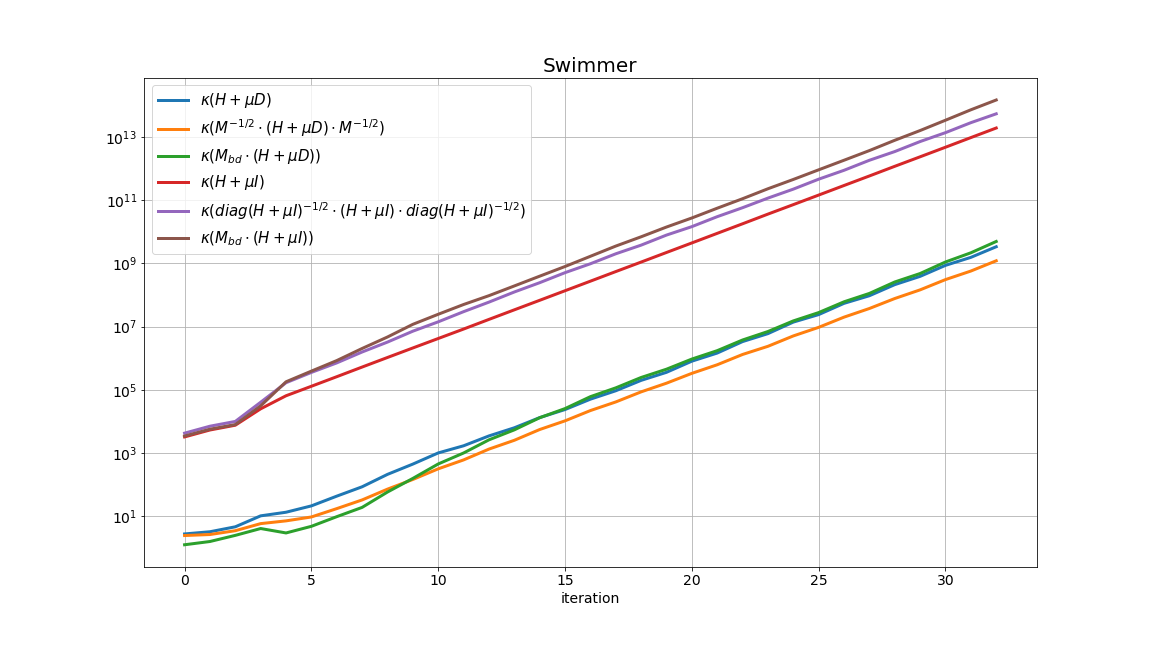}
				\caption{\footnotesize{Condition number for each iteration of several approaches to compute a CPD for the swimmer tensor. Note that it is showed the condition number of the approximated Hessian (always regularized) with and without regularization. The condition number of Tensor Fox is the orange one.}}
				\label{Swimmer conditioning}
			\end{figure}

\begin{figure}[H]
				\hspace{-1cm}\includegraphics[scale=.45]{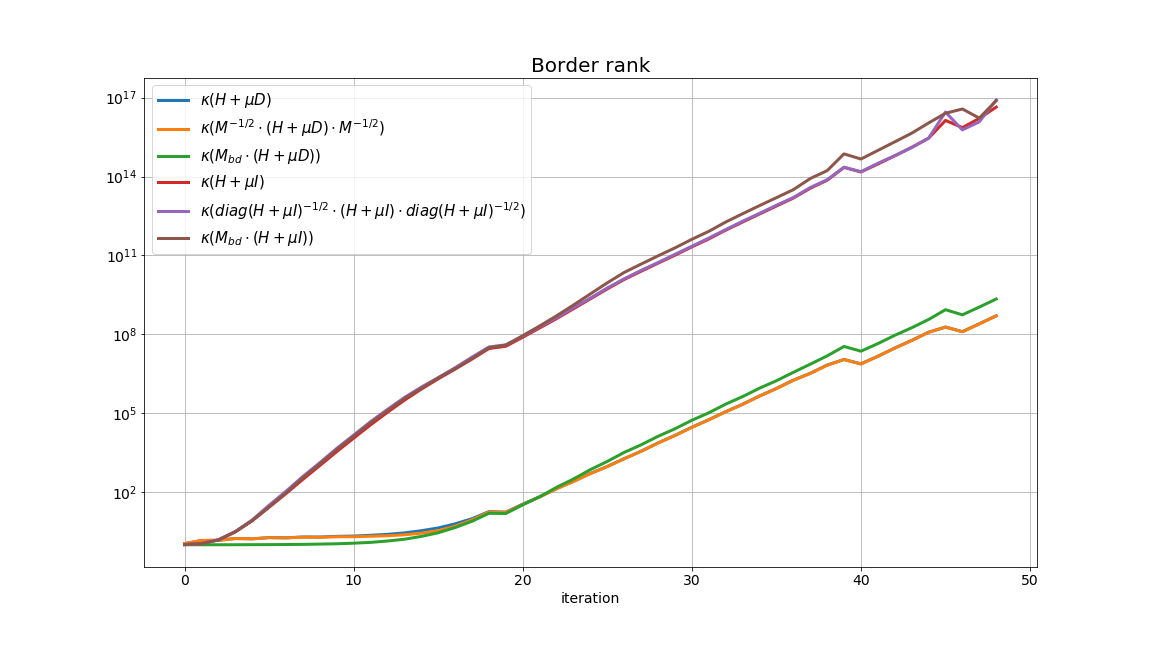}
				\caption{\footnotesize{Condition number for each iteration of several approaches to compute a CPD for the border rank tensor.}}
				\label{Border rank conditioning}
			\end{figure}

\begin{figure}[H]
				\hspace{-1cm}\includegraphics[scale=.45]{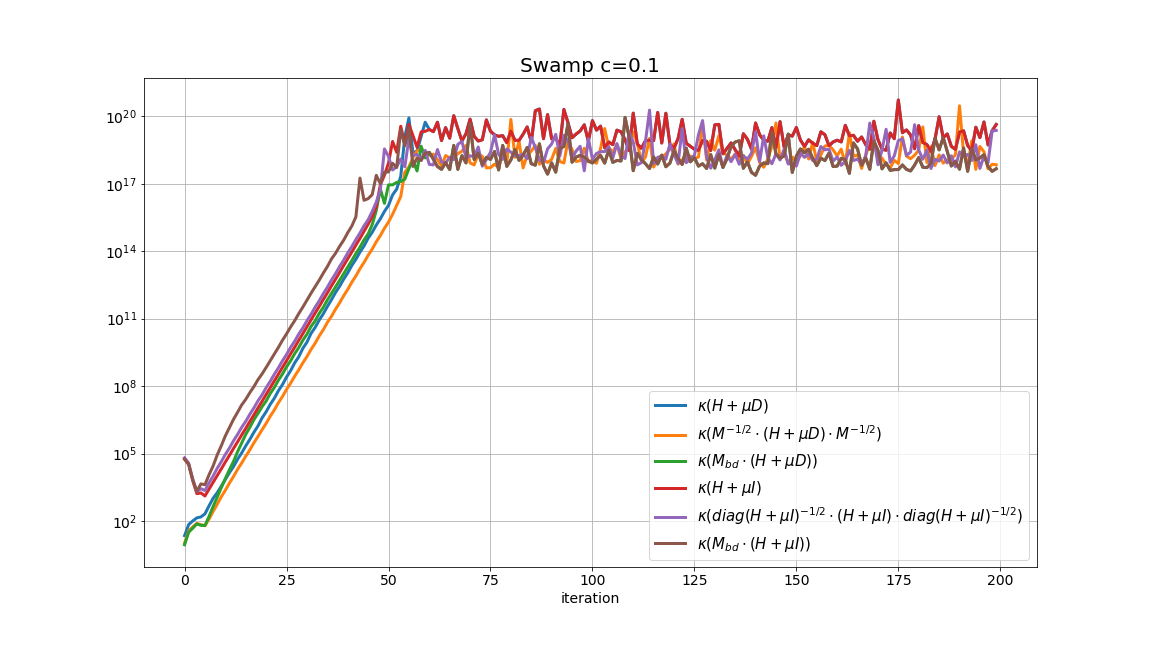}
				\caption{\footnotesize{Condition number for each iteration of several approaches to compute a CPD for the swamp tensor with $c = 0.1$.}}
				\label{Swamp conditioning}
			\end{figure}

\begin{figure}[H]
				\hspace{-1cm}\includegraphics[scale=.45]{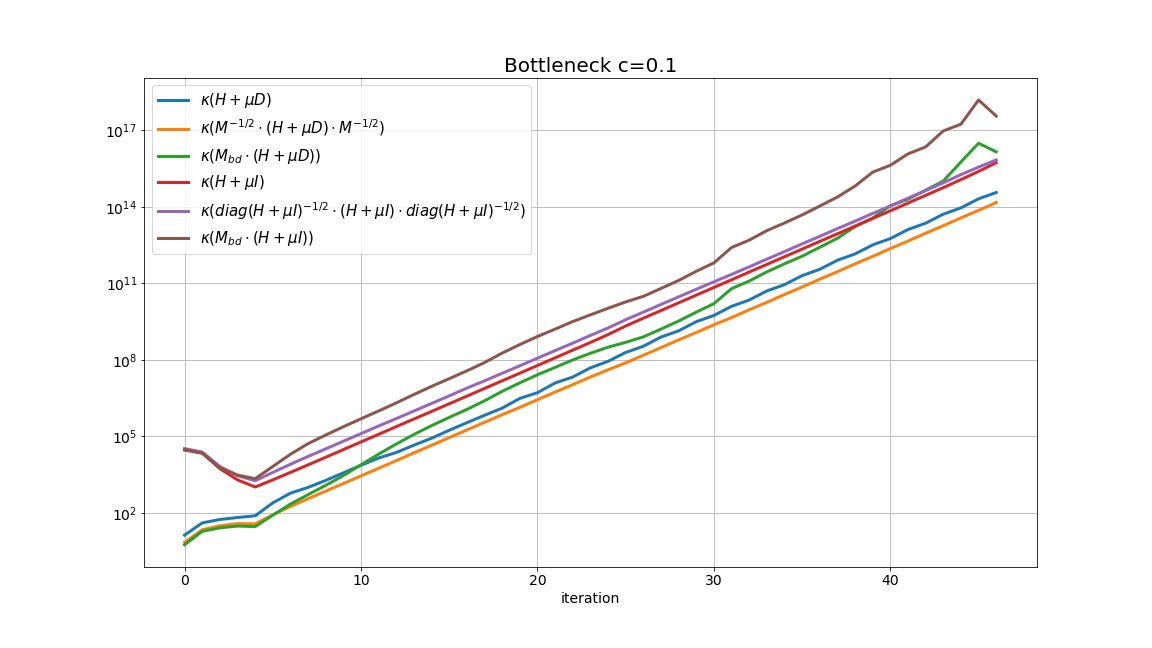}
				\caption{\footnotesize{Condition number for each iteration of several approaches to compute a CPD for the double bottleneck tensor with $c = 0.1$.}}
				\label{Bottleneck conditioning}
			\end{figure}

			Besides that formulation, we also have the Tensorlab's formulation  
			\begin{equation} \label{bd}
				\textbf{M}_{bd} \cdot (\textbf{H} + \mu \textbf{I}) \textbf{x} = \textbf{M}_{bd} \cdot \textbf{A}^T \textbf{b}, 
			\end{equation}
where $\textbf{M}_{bd}$ is the block diagonal preconditioner as explained in~\ref{comparison}. In figures~\ref{Swimmer conditioning}, ~\ref{Border rank conditioning}, ~\ref{Swamp conditioning}, ~\ref{Bottleneck conditioning} we show the evolution of the condition number of the approximated (with and without preconditioning) Hessian for some of the test tensors used before.

			Note that we tried the block diagonal preconditioner $\textbf{M}_{bd}$ both on the Levenberg-Marquardt and Tensor Fox formulations. It is clear that the diagonal preconditioner we are using here brings the best results in terms of lowering the condition number. Maybe the only exception are the swamp tensors, where after some iterations all approaches are as large as possible. 
			
\begin{figure}[H]
				\hspace{-1cm}\includegraphics[scale=.45]{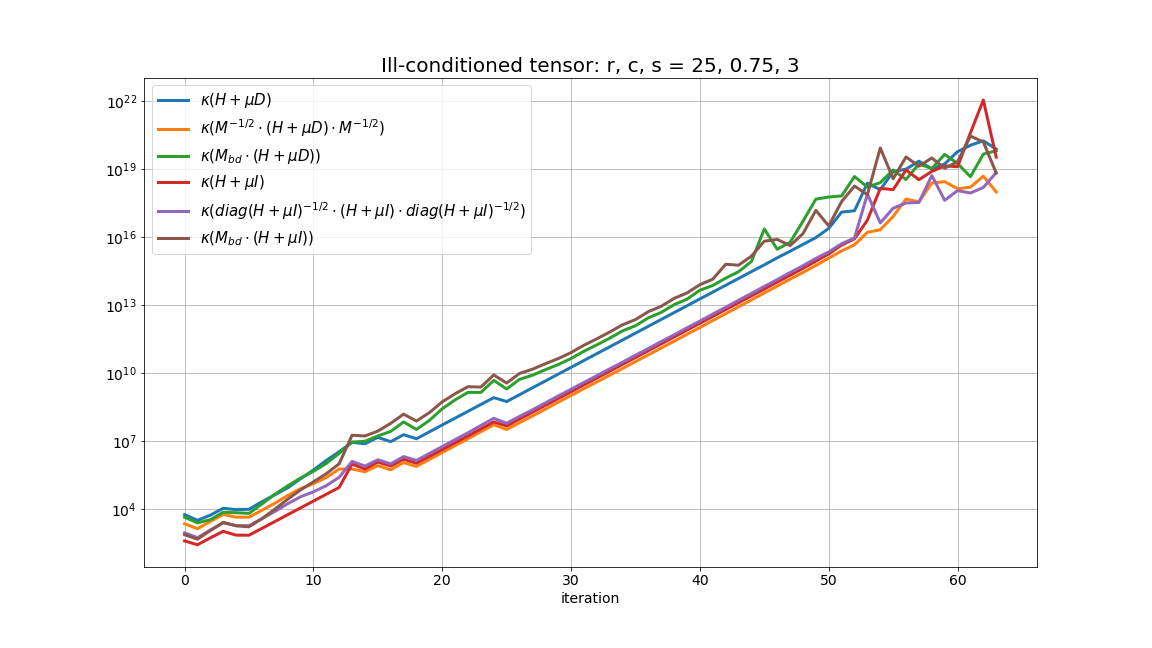}
				\caption{\footnotesize{Condition number for each iteration of several approaches to compute a CPD for an ill-conditioned tensor with $r = 25, c = 0.75, s = 3$.}}
				\label{Ill-conditioned tensor}
			\end{figure}

			In the next section we talk about the condition number of tensors and introduce a family of difficult tensors to handle. These tensors are such that most close CPDs are ill-conditioned even when the original CPD is well-conditioned. Figure~\ref{Ill-conditioned tensor} shows the evolution of the condition number of the approximated Hessian for the parameters $r = 25, c = 0.75, s = 3$ (they are explained in the next section). In this case, again, we observe that the approach of Tensor Fox is better than the others.
			
		\section{Conditioning} \label{cond_section}
			In this section we briefly mention the main concepts and results of conditioning of tensors. For more about the subject the reader is invited to check the references \cite{nick, join_cond, rgh, breiding}. 
			
			\subsection{Definitions and results}			
			\begin{definition}
				The \emph{Segre map} is the map given by
			$$Seg:\underset{\hspace{-.1cm}\displaystyle (\textbf{x}^{(1)}, \ldots, \textbf{x}^{(L)}) \mapsto \textbf{x}^{(1)} \otimes \ldots \otimes \textbf{x}^{(L)}}{\hspace{-.5cm} \mathbb{R}^{I_1} \times \ldots \times \mathbb{R}^{I_L} \to \mathbb{R}^{I_1 \times \ldots \times I_L} }.$$
			\end{definition}
			
			The image of this map without its origin is called the \emph{Segre variety}. It will be denoted by $S= \{\textbf{x}^{(1)} \otimes \ldots \otimes \textbf{x}^{(L)}: \ \textbf{x}^{(\ell)} \neq 0 \text{ for all } \ell=1 \ldots L \}$. Now consider the additive map given by
	$$\Phi_R: \underset{\displaystyle(\mathcal{T}_1, \ldots, \mathcal{T}_R) \mapsto \mathcal{T}_1 + \ldots + \mathcal{T}_R}{S \times \ldots \times S \to \mathbb{R}^{I_1 \times \ldots \times I_L} }.$$
	
			Let $\sigma_R$ be the set of tensors with rank $\leq R$. From the definition we have that $\text{Im}(\Phi_R) = \sigma_R$ and $rank(T) = \min \{ r: \ \mathcal{T} \in \sigma_r \}$ for all $T \in \mathbb{R}^{I_1 \times \ldots \times I_L}$. The derivative of $\Phi_R$ at $(\mathcal{T}_1, \ldots, \mathcal{T}_R)$ is the map
			$$\text{d}_{(\mathcal{T}_1, \ldots, \mathcal{T}_R)} \Phi_R: \text{T}_{\mathcal{T}_1}S \times \ldots \times \text{T}_{\mathcal{T}_R}S \to \text{T}_{\Phi_R(\mathcal{T}_1, \ldots, \mathcal{T}_R)} \mathbb{R}^{I_1 \times \ldots \times I_L}$$ 
given by
		$$\text{d}_{(\mathcal{T}_1, \ldots, \mathcal{T}_R)} \Phi_R (\dot{\mathcal{T}}_1, \ldots, \dot{\mathcal{T}}_R) = \dot{\mathcal{T}}_1 + \ldots + \dot{\mathcal{T}}_R,$$
where $\text{T}_{\mathcal{T}_r}S$ is the tangent space of S at $\mathcal{T}_r$. Then the condition number \cite{cucker} of the CPD at $(\mathcal{T}_1, \ldots, \mathcal{T}_R)$ is
			$$\kappa\left( (\mathcal{T}_1, \ldots, \mathcal{T}_R), \mathcal{T} \right) = \lim_{\varepsilon \to 0} \  \max_{\mathcal{T}' \in B_\varepsilon(\mathcal{T}) \cap \sigma_R} \frac{ \| \tilde{\Phi}^{-1}(\mathcal{T}) - \tilde{\Phi}^{-1}(\mathcal{T}') \| }{ \| \mathcal{T} - \mathcal{T}' \| },$$
where $B_\varepsilon(\mathcal{T}')$ is an $\varepsilon$-ball centered at $\mathcal{T}$ and $\tilde{\Phi}^{-1}$ is a local inverse of $\Phi$ at $(\mathcal{T}_1, \ldots, \mathcal{T}_R)$. If this local inverse doesn't exist then we define $\kappa\left( (\mathcal{T}_1, \ldots, \mathcal{T}_R), \mathcal{T} \right) = \infty$. Both norm are the Frobenius norm as it is the case in all this work. Note that $\kappa\left( (\mathcal{T}_1, \ldots, \mathcal{T}_R), \mathcal{T} \right)$ is completely determined by the choice of $\mathcal{T}_1, \ldots, \mathcal{T}_R$ since $\mathcal{T} = \Phi_R(\mathcal{T}_1, \ldots, \mathcal{T}_R)$, therefore we set $\kappa\left( \mathcal{T}_1, \ldots, \mathcal{T}_R \right) = \kappa\left( (\mathcal{T}_1, \ldots, \mathcal{T}_R), \mathcal{T} \right)$. In \cite{join_cond} it is showed that the condition number is the inverse of the smallest singular value of $\text{d}_{(\mathcal{T}_1, \ldots, \mathcal{T}_R)} \Phi_R$.
				
			Remember that the condition number measures the sensitivity of $(\mathcal{T}_1, \ldots, \mathcal{T}_R)$ to perturbations of $\Phi_R(\mathcal{T}_1, \ldots, \mathcal{T}_R)$. A rule of thumb in numerical analysis is the inequality
			$$\text{forward error } \lessapprox \text{ condition number } \cdot \text{ backward error},$$
which in this context translates to
			\begin{equation} \label{cond}
				\| (\mathcal{T}_1', \ldots, \mathcal{T}_R') - (\mathcal{T}_1, \ldots, \mathcal{T}_R) \| \lessapprox \kappa(\mathcal{T}_1, \ldots, \mathcal{T}_R) \cdot \| \mathcal{T}' - \mathcal{T} \|,
			\end{equation}
where 
			$$\| (\mathcal{T}_1', \ldots, \mathcal{T}_R') - (\mathcal{T}_1, \ldots, \mathcal{T}_R) \| = \min_{\sigma \in S_R} \displaystyle\sqrt{\sum_{r=1}^R \| \mathcal{T}_r' - \mathcal{T}_{\sigma(r)} \|^2}$$
and $S_R$ is the group of permutations of $R$ elements. From this inequality we can see that all algorithms showed so far are minimizing the backward error instead of the forward error. Relying on the backward error alone can be a dangerous practice in the presence of ill-conditioned CPDs.

			\subsection{A special family of tensors}
			In \cite{rgh} this aspect of conditioning is explored and a family of $15 \times 15 \times 15$ tensors is proposed for testing. This family is such that the best algorithm of Tensorlab (NLS) failed to produce well conditioned approximations even when the original tensor is well conditioned. The parameters $(r, c, s) \in \{15, 20, 25, 30\} \times \{0, 0.25, 0.5, 0.75\} \times \{1, 2, 3, 4\}$ are used to generate each tensor of the family. Let $\textbf{R}_c$ be the upper triangular factor in the Cholesky decomposition $\textbf{R}_c^T \textbf{R}_c = c \textbf{1} \textbf{1}^T + (1-c) \textbf{I}_r$, where $\textbf{1} \in \mathbb{R}^r$ is the vector of ones and $\textbf{I}_r$ is the $15 \times 15$ identity matrix. Then each CPD of this family with parameters $r, c, s$ is given by $\textbf{A}_1, \textbf{A}_2, \textbf{A}_3 \in \mathbb{R}^{15 \times r}$, where
			$$\textbf{A}_i = \textbf{N}_i \cdot \textbf{R}_c \cdot \text{diag}\left( 10^{\frac{s}{3r}}, 10^{\frac{2s}{3r}}, \ldots, 10^{\frac{rs}{3r}} \right),$$
where $\textbf{N}_i$ is a $15 \times r$ random matrix with its entries draw from the normal distribution $\mathcal{N}(0,1)$ with mean 0 and variance 1.	

			Their experiment proceed as follow:	
			\begin{enumerate}
				\item Randomly sample $\textbf{A}_1, \textbf{A}_2, \textbf{A}_3$ as described above, and let $\mathcal{T}' = \Phi_r(\mathcal{T}_1', \ldots, \mathcal{T}_r')$, where the rank one tensor are obtained from the factors $\textbf{A}_1, \textbf{A}_2, \textbf{A}_3$;
				\item create a perturbed tensor $\mathcal{T} = \frac{\mathcal{T}'}{\| \mathcal{T}' \|} + 10^{-3} \frac{\mathcal{E}}{\| \mathcal{E} \|}$, where $\mathcal{E}$ is a $15 \times 15 \times 15$ random tensor with its entries draw from the normal distribution $\mathcal{N}(0,1)$;
				\item randomly samples factor matrices $\textbf{M}_1, \textbf{M}_2, \textbf{M}_3 \in \mathbb{R}^{15 \times r}$ to be used as initialization;
				\item compute an approximated CPD for $\mathcal{T}$ from the initialization $\textbf{M}_1, \textbf{M}_2, \textbf{M}_3$.
			\end{enumerate}		
			
			In the original experiment they used some specific parameter configuration for Tensorlab NLS, but we won't go into such details here. The reader may read section 7 of \cite{rgh} for more information. For each triple $(r, c, s)$ one tensor is generated by the procedure above, then 25 random initializations are generated and 25 CPDs are computed, on for each different initialization. For a very significant fraction of initializations, the state-of-the-art method halt at extremely ill-conditioned CPDs. 
			
			Table 4.19 shows a part of their results, where Tensorlab had more trouble to produce well-conditioned solutions. The percentage showed is the fraction of cases where the condition number of the approximated CPD is bigger than $10^3$ among those CPDs whose backward error is very small (less than $1.1 \cdot 10^{-3}$). From~\ref{cond} we can conclude that these CPDs are completely uninterpretable, that is, no correct significant digits are present in the individual rank one terms. RGN-HR stands for \emph{Riemannian Gauss-Newton with hot restarts} and RGN-Reg stands for \emph{Riemannian Gauss-Newton with regularization}. The former is their main method, while the latter was included only to illustrate that the proposed hot restarts mechanism provides a superior way of handling ill-conditioning. RGN-HR is not just more well-conditioned than Tensorlab's NLS, but it is also faster in most cases as showed in the paper. 
			
			\begin{table}[H]
			\centering
				\begin{tabular}{|c|c|c|}
					\hline
					\bf{RGN-HR} & \bf{RGN-Reg} & \bf{Tensorlab}\\
					\hline
					$23.6 \%$ & $42.3 \%$ & $51.0 \%$\\
					\hline
				\end{tabular}
				\caption{\footnotesize{Fraction of ill-conditioned CPDs.}} \bigskip
			\end{table}
			
			\subsection{Results}
			Here we tried to make the parameter setting for Tensor Fox as close as possible to the choices made for Tensorlab. First it was observed that better results were obtained with the damping parameter fixed at $\mu = 10^{-8}$. The most relevant parameter is the number of CG iterations. We noted that increasing this value lead to more well-conditioned solutions consistently, which makes sense because the ill-conditioning comes from the approximated Hessian formulation (which is singular). In the paper they used a maximum of $\verb|cg| \_ \verb|maxiter| = 75$ number of iterations for the CG. Also, although Tensor Fox normally does not use the dogleg method\footnote{See \cite{madsen} for more about this method.}, for this problem in particular this method showed to be handful. 
						
\begin{figure}[H]
				\centering
				\includegraphics[scale=.3]{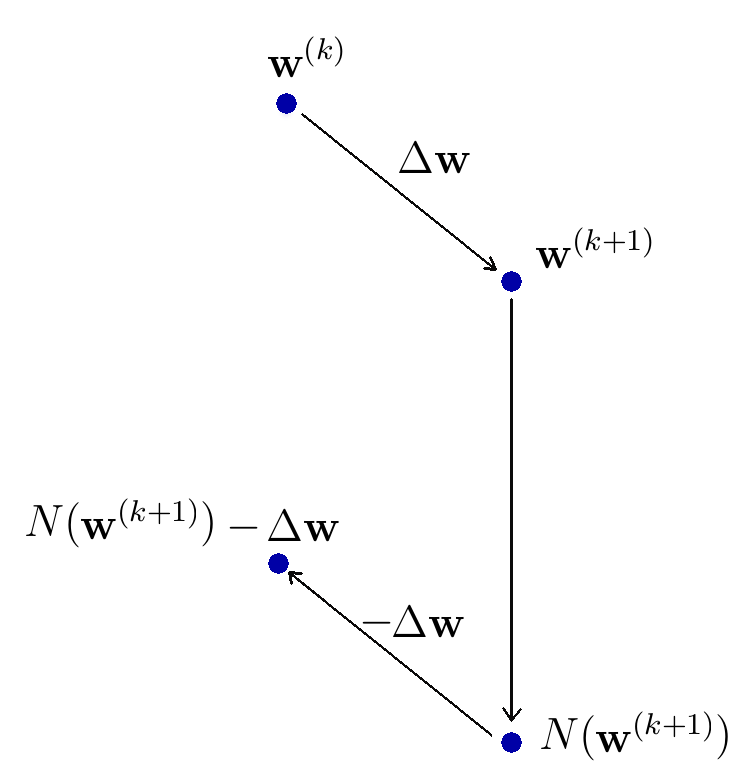}
				\caption{\footnotesize{After making a bad step, the program draws back in a way that the new point is close to norm-balanced, then it applies the dogleg method.}}
				\label{Draw back}
			\end{figure}

			The dogleg method is used when the error of some iteration is too big compared to the previous error. This is a case when the unusual step is not considered unusual but some kind of divergence. This rarely is used in Tensor Fox because the program only recognizes divergence when the error is 100 times bigger than the previous error, and this almost never happens. In any case, when this does happen, the program draws back to a ``version'' of the previous step and perform the dogleg method to produce the new step. We explain what is this version of the previous step. Let $\textbf{w}^{(k)}$ be the $k$-th dGN step. Remember that before computing the next iteration, first we make the factor norm balanced, let $N(\textbf{w}^{(k)})$ be this norm-balanced representation. Then the next dGN iteration is computed with $N(\textbf{w}^{(k)})$ instead of $\textbf{w}^{(k)}$. We already observed in~\ref{updates} that this procedure always improves the conditioning of the iteration. In the case the next step, $\textbf{w}^{(k+1)} = \textbf{w}^{(k)} + \Delta \textbf{w}$, has a large error, we take its norm-balanced representation and shift it back with the dGN step $\Delta \textbf{w}$. The new point obtained, $N(\textbf{w}^{(k+1)}) - \Delta \textbf{w}$ should be close to the norm-balanced representation of $\textbf{w}^{(k)}$, but not quite the same. Then the dogleg method takes action and improves the step. Figure~\ref{Draw back} illustrates the process of drawing back a step. We remark that $N(\textbf{w}^{(k+1)}) - \Delta \textbf{w}$ usually is different than $N(\textbf{w}^{(k)})$.
			
\begin{figure}[H]
				\centering
				\includegraphics[scale=.8]{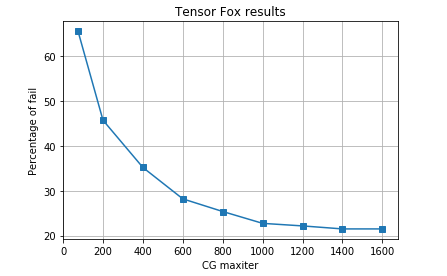}
				\caption{\footnotesize{Percentage of fails as we increase the maximum number of CG iterations.}}
				\label{Condition number fail percentage}
			\end{figure}

			For this particular family of tensors we could verify that it is better to perform draw back and apply the dogleg method in every iteration. We repeat the previous table plus the Tensor Fox results, shown in table 4.20. As we can see, Tensor Fox is no better than Tensorlab for $\verb|cg| \_ \verb|maxiter| = 75$ and it is better than RGN-HR for $\verb|cg| \_ \verb|maxiter| = 1600$. The situation changes as we increase $\verb|cg| \_ \verb|maxiter|$. We tested the values $\verb|cg| \_ \verb|maxiter| = 75, 200, 400, 600, 800, 1000, 1200, 1400, 1600$, see figure~\ref{Condition number fail percentage}.\bigskip
			
			\begin{table}[H]
			\centering
				\begin{tabular}{|c|c|c|c|c|}
					\hline
					\bf{RGN-HR} & \bf{RGN-Reg} & \bf{Tensorlab} & \bf{Tensor Fox} $(75)$ & \bf{Tensor Fox} $(1600)$\\
					\hline
					$23.6 \%$ & $42.3 \%$ & $51.0 \%$ & $65.5 \%$ & $21.4 \%$\\
					\hline
				\end{tabular}
				\caption{\footnotesize{Fraction of ill-conditioned CPDs, including Tensor Fox.}} \bigskip
			\end{table}
			
			Of course the running time increases as we increase $\verb|cg| \_ \verb|maxiter|$, however the expected time to get a well-conditioned CPD decreases. The same does not occurs with Tensorlab. We increased $\verb|cg| \_ \verb|maxiter|$ for Tensorlab and observed an improvement with regard the conditioning but the running time in this case increased substantially. For instance, when $\verb|cg| \_ \verb|maxiter| = 400$, the slowest expected time for Tensor Fox obtained was $33.3$ seconds, with $r=30, c=0.75, s=1$, and for Tensorlab it was $101.2$ seconds, with $r=30, c=0.25, s=1$. Even more, in all tests with Tensor Fox the expected time is always less than $40$ seconds. This happens because the CG iterations of Tensor Fox are lighter than Tensorlab's, which is something already observed in~\ref{comparison}.   
			
			We finish this section observing that the algorithm RGN-HR is able to obtain many well-conditioned solutions in less time than Tensor Fox and Tensorlab. Their approach is very different in some aspects and it is worth reading about it. The speed of RGN-HR tends to be slower for larger tensors, which is a major drawback. However if the rank or the multilinear rank is small enough, using tolerance-based compression may reduced the problem substantially, and then RGN-HR can be efficiently used.  
			
	\section{Parallelism} \label{par}
		Several parts of Tensor Fox are open to parallelism: unfoldings, SVDs, matrix multiplications, Khatri-Rao products, and so on. With regard to the dGN algorithm, the main costs comes from the Khatri-Rao products to obtain the gradient, and the matrix-matrix multiplications. Both routines run in parallel in Tensor Fox. In particular, the matrix-matrix multiplication uses the BLAS parallelism, which is as good as one can get. 
		
\begin{figure}[H]
			\centering
			\includegraphics[scale=.75]{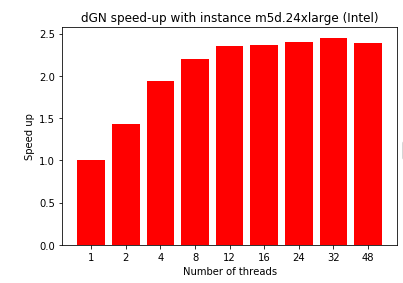}
			\caption{\footnotesize{Speed-up with the number of threads when computing a rank-$15$ CPD of a $2000 \times 2000 \times 2000$ random tensor with rank $15$ without compression.}}
			\label{Speed-up}
		\end{figure}

		Still, the dGN algorithm, by design is not highly parallel, that is, each iteration must be computed sequentially. Besides that, within each iteration we must run the CG algorithm, which run each iteration sequentially again. Figure~\ref{Speed-up} shows the speed-up obtained by the addition of threads to compute a rank-$15$ CPD of a $2000 \times 2000 \times 2000$ random tensor with rank $15$. We are not compressing in this example. The machine used in this experiment was the AWS instance m5d.24xlarge, consisting of a Intel Xeon Platinum 8000 series (Skylake-SP) processor with a sustained all core Turbo CPU clock speed of up to 3.1 GHz, and 384 GB of memory. This computational and memory power was necessary due to the size of the tensor. The reason for using such a big tensor is because the parallelism increases as we increase the dimensions and the rank. 
		
\begin{figure}[H]
			\centering
			\includegraphics[scale=.75]{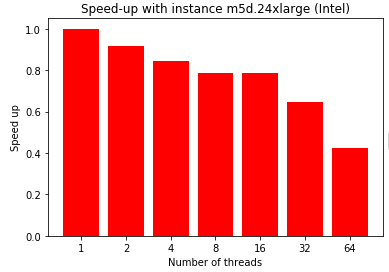}
			\caption{\footnotesize{Speed-up with the number of threads when computing a rank-$50$ CPD of the swimmer tensor.}}
			\label{Speed-up of the swimmer tensor}
		\end{figure}

		Even with a big tensor as this one, the speed-up curve gets flat already for 16 threads. This shows how the dGN algorithm, the way it is presented in this thesis, is not highly parallel. For smaller tensors as the swimmer tensor (see figure~\ref{Speed-up of the swimmer tensor}) one can expect the speed-up to drop as the number of threads increases due to the overhead communication. There is just too much data being passed for a small number of operations.  
		
		\begin{figure}[H]
			\hspace{-1cm}\includegraphics[scale=.5]{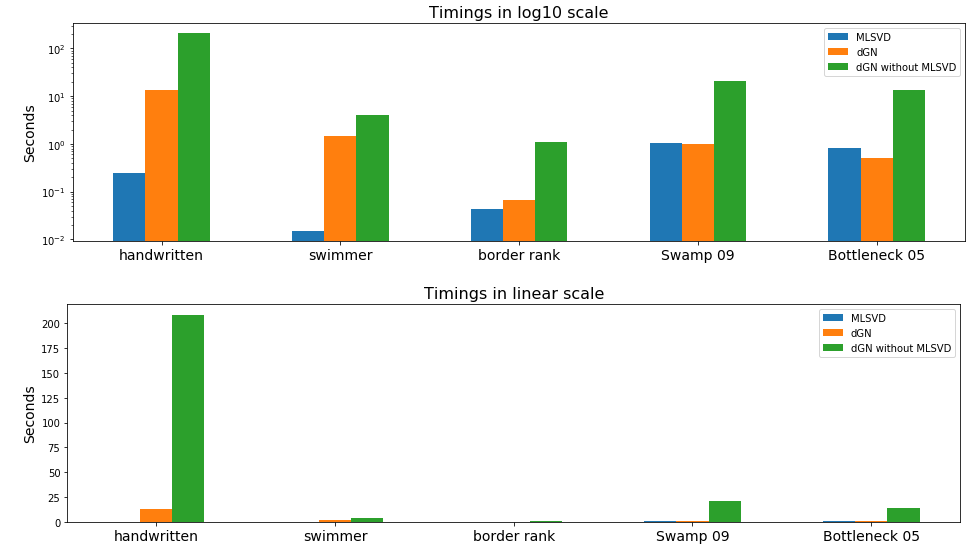}
			\caption{\footnotesize{Timings to compute parts of the CPDs of known tensors. Note that summing the bars of the MLSVD and dGN timings is not the correct comparison in log scale. That is why we also showed the linear scale, so the reader have a notion of the real difference in time when there is no compression.}}
			\label{Compression vs. no compression}
		\end{figure} 
		
		The main source of parallelism comes from the MLSVD, which relies on several randomized SVDs \cite{rand_svd}, and this method is highly parallel. \footnote{In fact, the randomized method for SVDs are also very efficient if running in a GPU.} As we stressed before, compressing before running the dGN is a step we should always take. Figure~\ref{Compression vs. no compression} reinforce this claim by showing that the compression time plus the dGN time is less than the dGN time without compression.		
		
	\section{What are the main features of Tensor Fox?} \label{main_feat}
		In the previous sections we introduced Tensor Fox in detail and conducted a lot of experiments proving that Tensor Fox is competitive in a wide range of examples. After that the reader is left with the impression that all new features are relevant and indispensable. In this section we investigate the main features, turning them off and running the test tensors. The features to be investigated are:
		
		\begin{itemize}
			\item Stopping conditions
			\item Regularization
			\item Preconditioning
			\item Number of CG iterations
		\end{itemize}	
			
		Before anything, we want to remark that one of the most relevant features of Tensor Fox is the tolerance based compression, but this feature is already investigated enough. 
		
\begin{figure}[H]
			\centering
			\includegraphics[scale=.31]{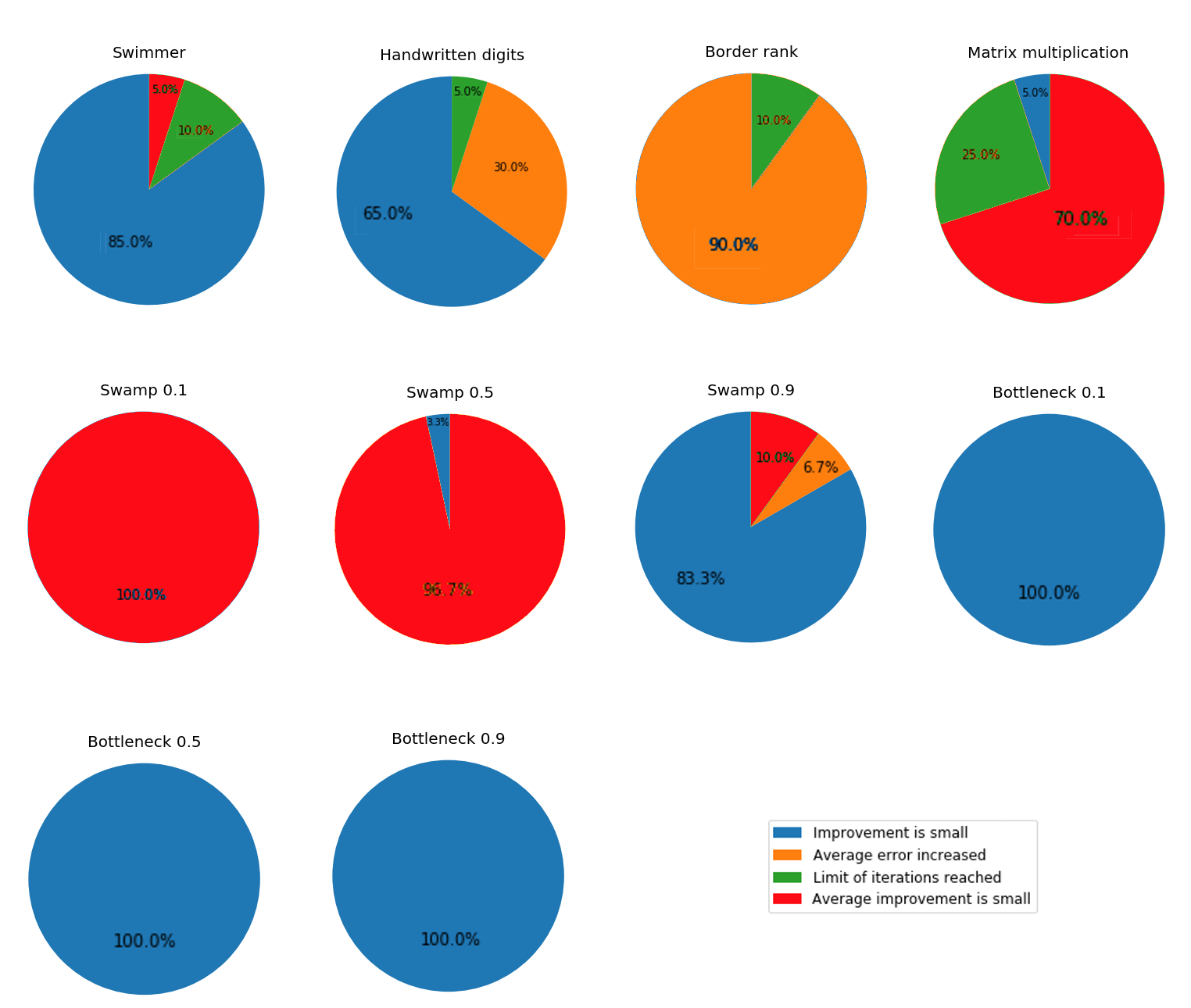}
			\caption{\footnotesize{Pie charts with the stopping condition frequencies.}}
			\label{Pie charts with the stopping condition frequencies}
		\end{figure}

		One thing the reader may notice is the large number of stopping conditions introduced in~\ref{Pie charts with the stopping condition frequencies}. Below we show their frequency when running 20 CPDs for each tensor. As we can see, most of the time the stopping condition triggered is the error improvement, which actually is the most common in all solvers. Sometimes none condition is triggered and the program just reach the maximum number of iterations. The novelty here are the conditions with averages, which really helps to stop the iterations earlier. 
		
\begin{figure}[H]
			\centering
			\includegraphics[scale=.5]{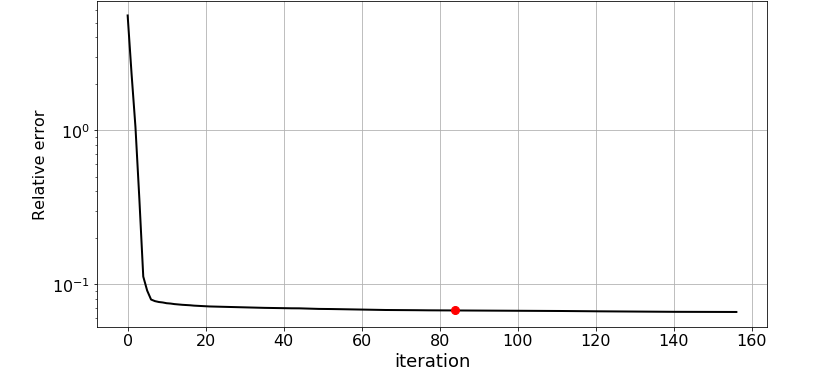}
			\caption{\footnotesize{Disabling stopping condition for the swamp tensor with $c = 0.1$.}}
			\label{Disabling stopping condition}
		\end{figure}

		For instance, we can note that the averages of the improvements are the unique stopping condition triggered for the swamp tensor with $c = 0.1$. In figure~\ref{Disabling stopping condition} we show the error evolution to compute a CPD for this tensor with this stopping condition disabled. The red dot is the iteration where the program would stop if the condition were enabled. We can see that there is almost no reason to continue iterating, therefore the stopping condition is appropriate. 
		
\begin{figure}[H]
			\hspace{-1.5cm}\includegraphics[scale=.53]{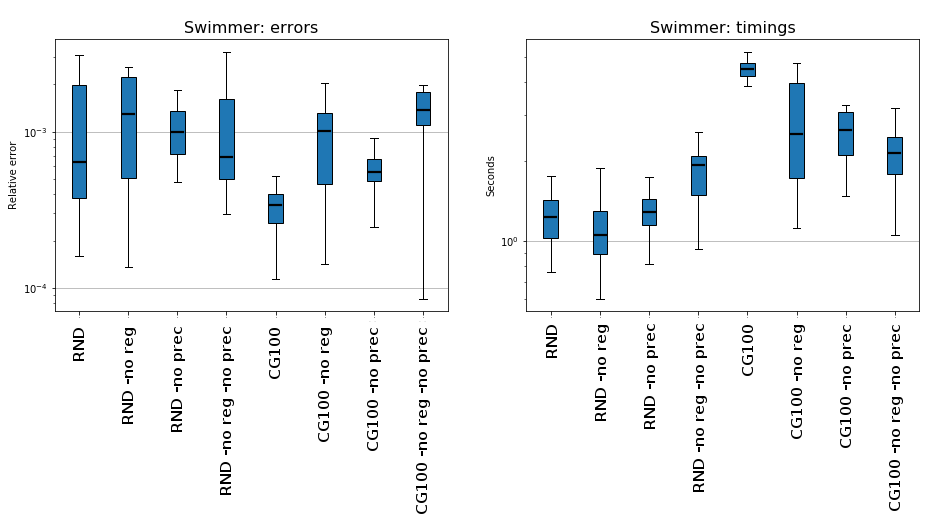}
			\caption{\footnotesize{Box plots with errors and timings for the swimmer tensor.}}
			\label{Box plot: swimmer tensor}
		\end{figure}

\begin{figure}[H]
			\hspace{-1.5cm}\includegraphics[scale=.53]{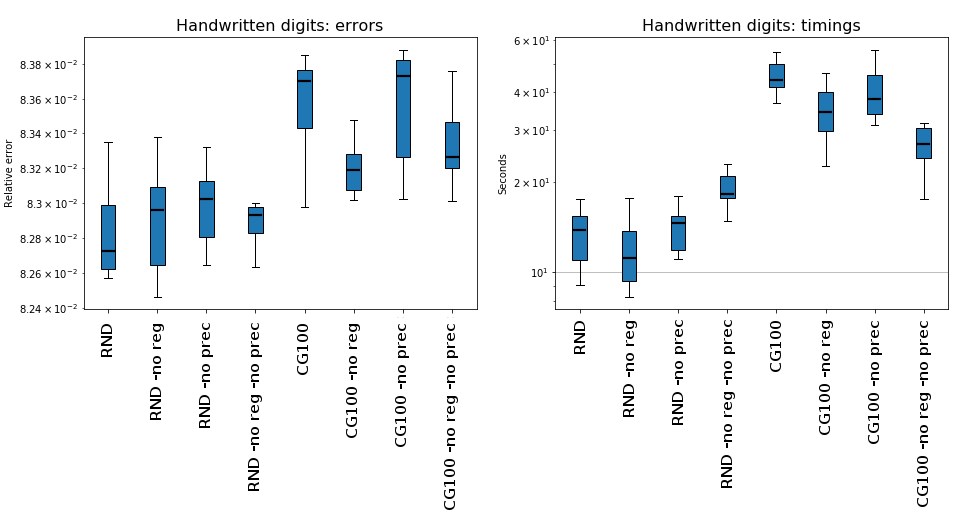}
			\caption{\footnotesize{Box plots with errors and timings for the handwritten digits tensor.}}
			\label{Box plot: handwritten digits tensor}
		\end{figure}

\begin{figure}[H]
			\hspace{-1.5cm}\includegraphics[scale=.53]{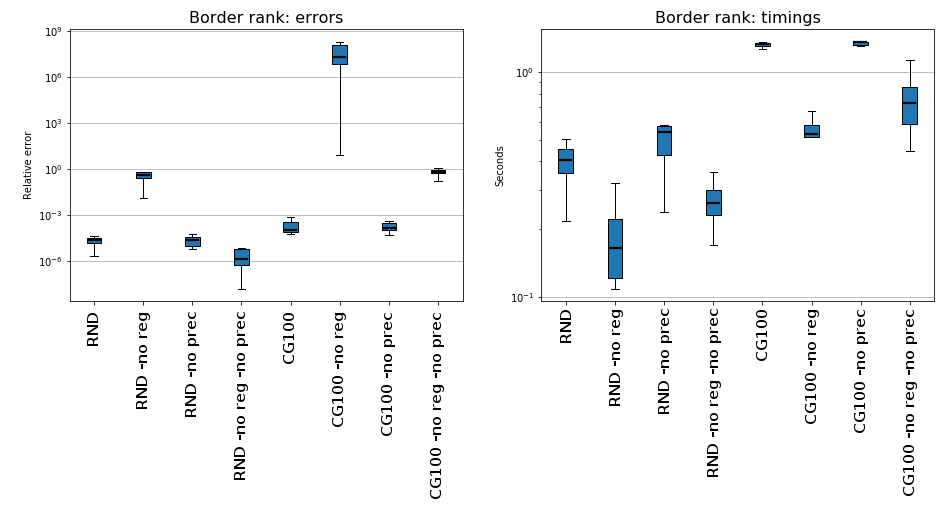}
			\caption{\footnotesize{Box plots with errors and timings for the border rank tensor.}}
			\label{Box plot:border rank tensor}
		\end{figure}

\begin{figure}[H]
			\hspace{-1.5cm}\includegraphics[scale=.53]{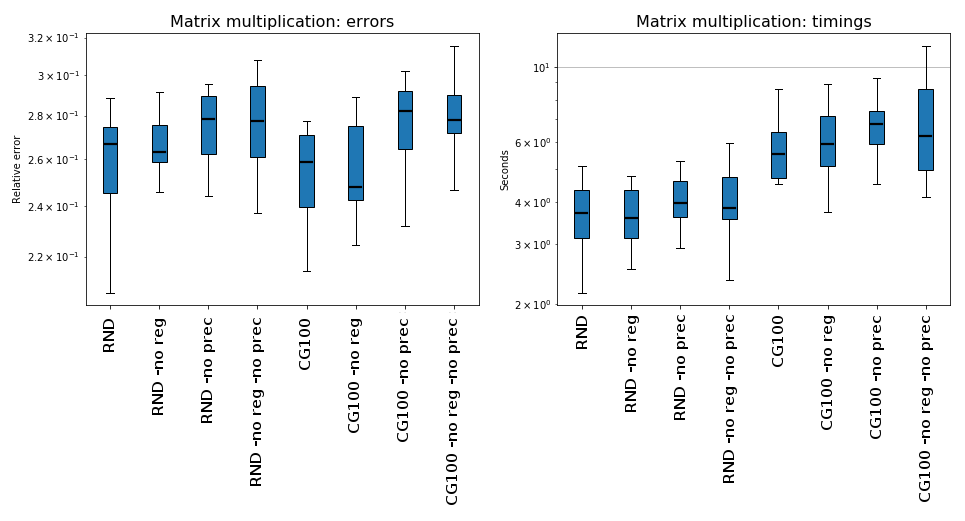}
			\caption{\footnotesize{Box plots with errors and timings for the matrix multiplication tensor.}}
			\label{Box plot: matrix multiplication tensor}
		\end{figure}

\begin{figure}[H]
			\hspace{-1.5cm}\includegraphics[scale=.53]{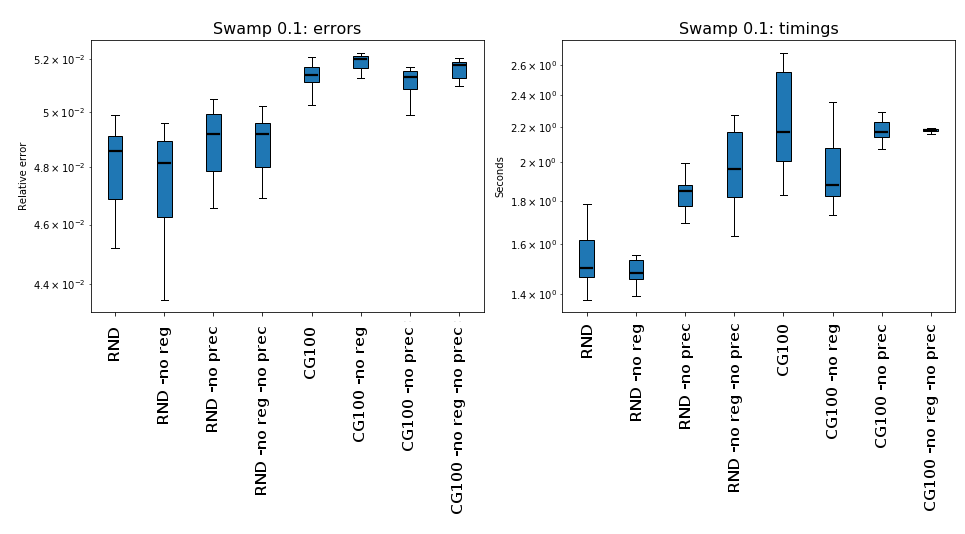}
			\caption{\footnotesize{Box plots with errors and timings for the swamp $c = 0.1$ tensor.}}
			\label{Box plot: swamp $c = 0.1$ tensor}
		\end{figure}

\begin{figure}[H]
			\hspace{-1.5cm}\includegraphics[scale=.53]{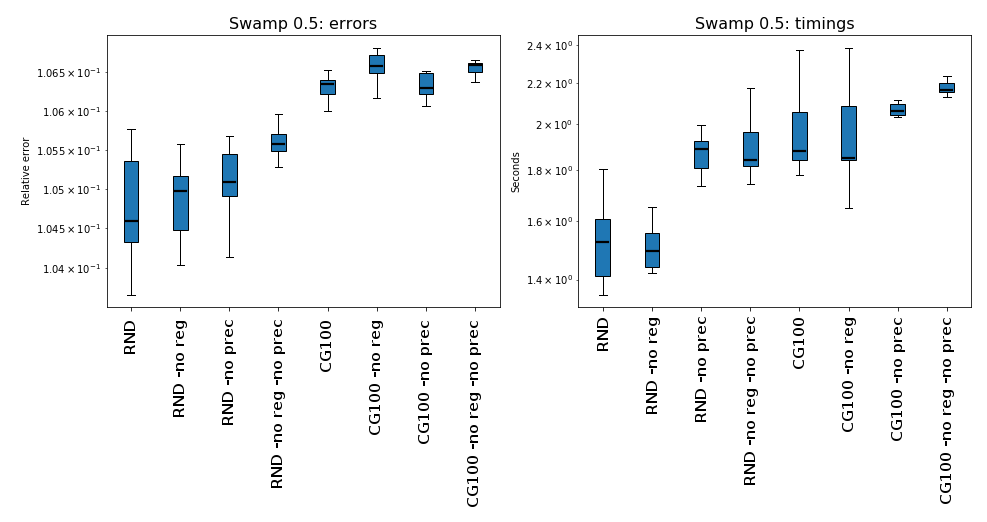}
			\caption{\footnotesize{Box plots with errors and timings for the swamp $c = 0.5$ tensor.}}
			\label{Box plot: swamp $c = 0.5$ tensor}
		\end{figure}

\begin{figure}
			\hspace{-1.5cm}\includegraphics[scale=.53]{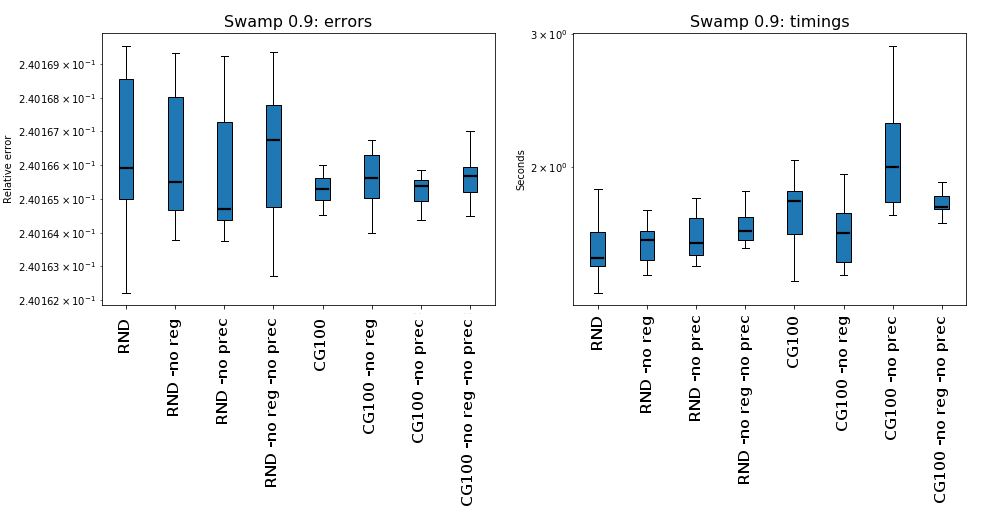}
			\caption{\footnotesize{Box plots with errors and timings for the swamp $c = 0.9$ tensor.}}
			\label{Box plot: swamp $c = 0.9$ tensor}
		\end{figure}

\begin{figure}[H]
			\hspace{-1.5cm}\includegraphics[scale=.53]{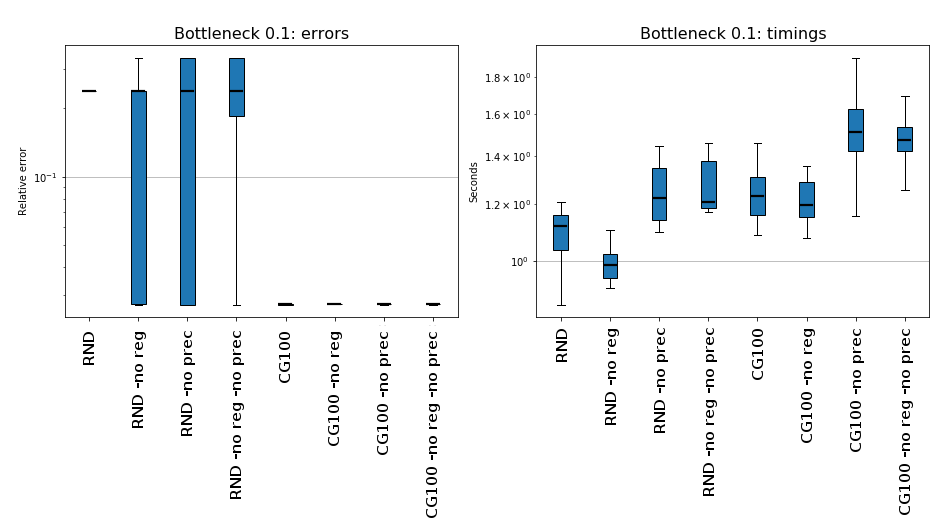}
			\caption{\footnotesize{Box plots with errors and timings for the bottleneck $c = 0.1$ tensor.}}
			\label{Box plot: bottleneck $c = 0.1$ tensor}
		\end{figure}

\begin{figure}[H]
			\hspace{-1.5cm}\includegraphics[scale=.53]{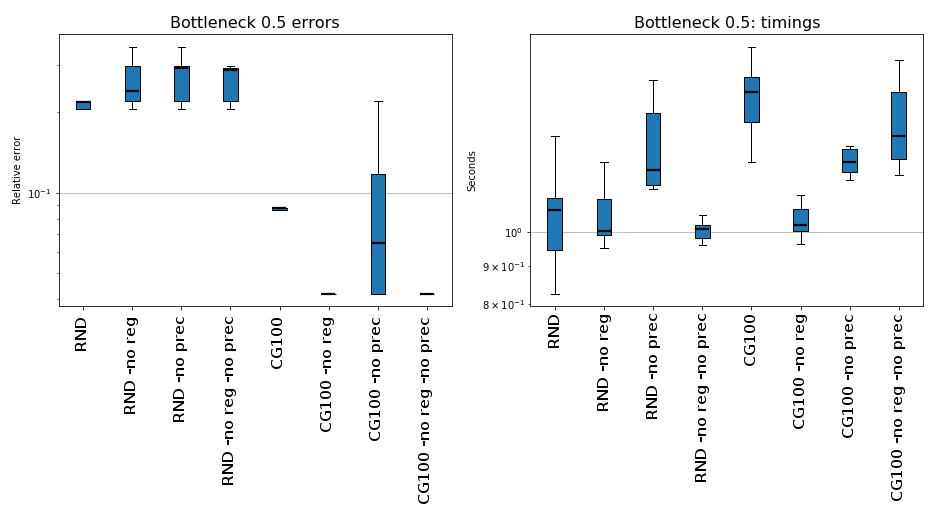}
			\caption{\footnotesize{Box plots with errors and timings for the bottleneck $c = 0.5$ tensor.}}
			\label{Box plot: bottleneck $c = 0.5$ tensor}
		\end{figure}

\begin{figure}[H]
			\hspace{-1.5cm}\includegraphics[scale=.53]{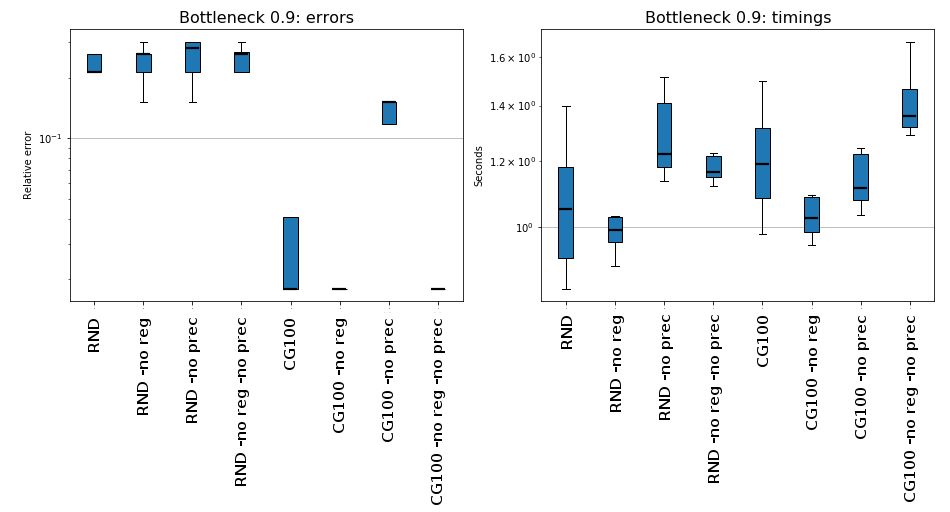}
			\caption{\footnotesize{Box plots with errors and timings for the bottleneck $c = 0.9$ tensor.}}
			\label{Box plot: bottleneck $c = 0.9$ tensor}
		\end{figure}

		Now we will analyze the other features of Tensor Fox at the same time. Basically we will repeat the experiments as we did before to produce the box plots, but the models compared are variations of Tensor Fox, with the relevant features turned on or off. We denote by \textbf{RND} the Tensor Fox default model, using random \verb|cg_maxiter| as described in~\ref{dGN-cg}, and \textbf{CG100} the model where \verb|cg_maxiter|$= 100$. Additionally, we write ``no reg'' when regularization is removed from the model and ``no prec'' when preconditioning is removed from the model. All results are summarized in figures~\ref{Box plot: swimmer tensor} to~\ref{Box plot: bottleneck $c = 0.9$ tensor}.
		
		The first thing to note is that \textbf{CG100} is not always more accurate than \textbf{RND}, that is, making more CG iterations does not necessarily improve the results. Also it is evident that all \textbf{CG100} based models are slower than their corresponding \textbf{RND} models. We observe that the use of regularization is not always necessary, and in most cases the performance is the same. The exception are the border rank tensor, where removing regularization made the accuracy much worse, and in all bottleneck tensors, where removing regularization increased the variance of the solutions. In fact, there are several other cases where the removal of regularization increases the variance of the errors. This is true for the \textbf{RND} and \textbf{CG100} models. Removing the preconditioning has the expected result, that is, the increase of the running time. Because the main role of the preconditioning is to decrease the number of CG iterations. In some cases removing the preconditioning even increased the error. 
		
		We conclude that removing regularization is feasible sometimes depending on the tensor. For this reason Tensor Fox has the option \verb|init_damp|, which can be initialized to zero, thus removing the effect of regularization. Removing the preconditioning is never good, therefore the preconditioning is always performed.

%% file: CHAPTER_5.tex
\chapter{Tensor learning}\label{cap-5} 
	We call \emph{tensor learning} any tensor technique applied to a machine learning problem. Here we will some of these techniques and their performance will be evaluated. The main tools are the MLSVD and the CPD. Usually the MLSVD is used to reduce dimensionality and the CPD is used to find latent parameters. In order to solve machine learning problems with tensor techniques, the biggest problem challenge is the modelling part, that is, transforming a machine learning problem into a tensor problem. 
	
	The main contribution is in~\ref{cpd-nn}, where a model of neural network based on the CPD is introduced. Although this neural network architecture still acts like a black box, the outputs consists of multilinear maps, which are more understandable than classical neural networks. The work here only takes in account dense neural networks, so the comparison is limited. There is still a lot of research to be done before discarding or not this new model, but for the moment this model seems to be reasonable.

	\section{Classification with the MLSVD} 
\begin{figure}[h] 
			\includegraphics[scale=.3]{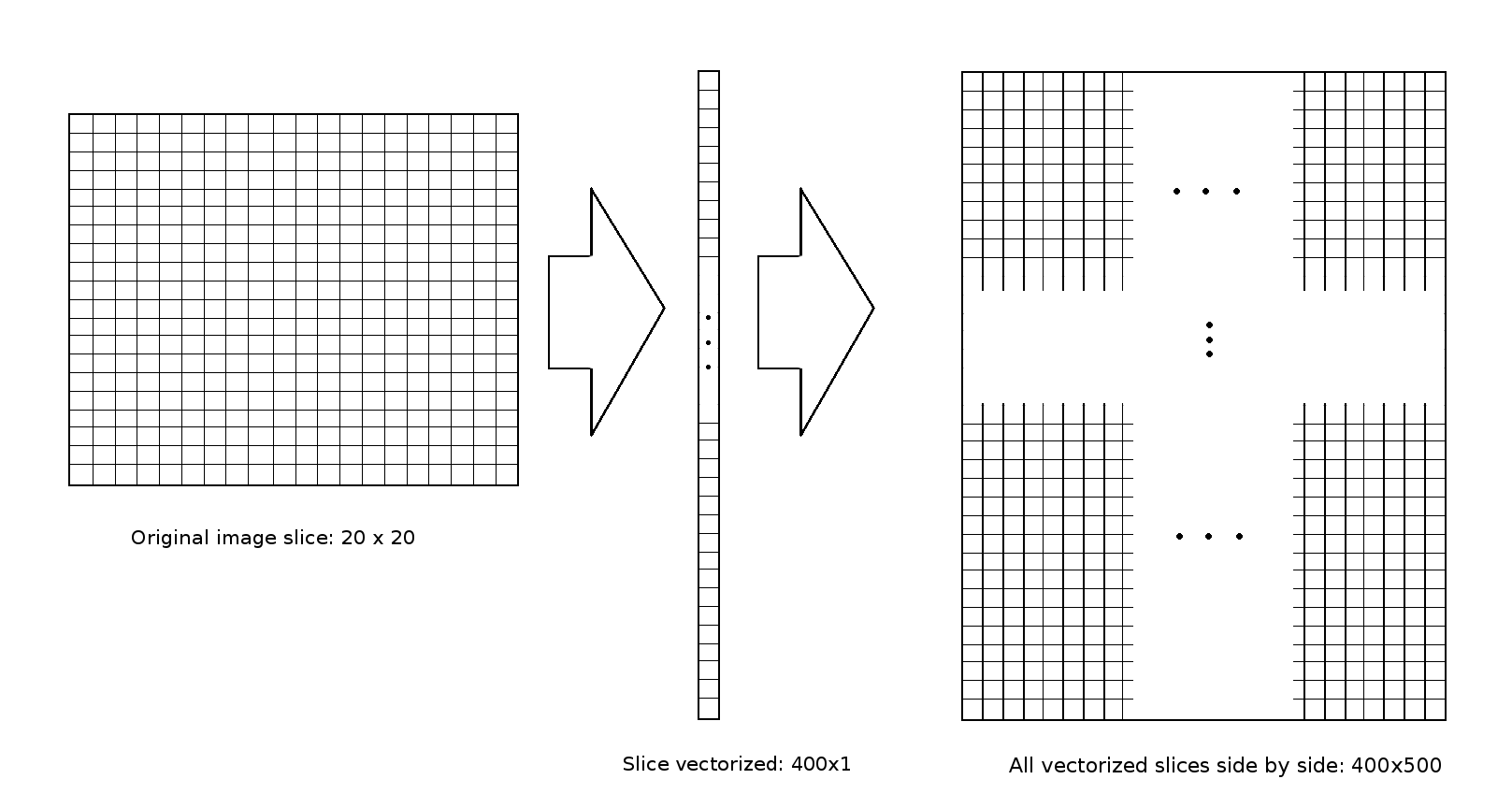}
			\caption{\footnotesize{Transformation of a frontal slice into a vector and attaching dimension responsible to the class of the vector (the image represents the class $\textbf{e}_2$, which corresponds to the digit 2).}}
			\label{hw-transform}
		\end{figure}

		The MLSVD can be considered as a high order PCA. In this section we present a method to solve classification problems using the MLSVD, following \cite{savas} and \cite{hw}. The problem of handwritten digit classification will be used to test the performance of the method. We start with the same tensor described in~\ref{hw} and apply some transformations to it in order to take the classes (the digits) in account. Each frontal slice (a $20 \times 20$ matrix) is vectorized to to form a column vector of size $400$. Then we gather all these vectors side by side to form a $400 \times 500$ matrix corresponding to a certain class. In total we have 10 classes, each one corresponds to a frontal slices of the tensor $\mathcal{T} \in \mathbb{R}^{400 \times 500 \times 10}$. The first slice correspond to the digit 0, the second correspond to the digit 1, and so on. Figure~\ref{hw-transform} describes the procedure to form each frontal slice of $\mathcal{T}$.
		
		Let $\mathcal{T} = (\textbf{U}^{(1)}, \textbf{U}^{(2)}, \textbf{U}^{(3)}) \cdot \mathcal{S}$ be a truncated MLSVD of $\mathcal{T}$, where $\textbf{U}^{(1)} \in \mathbb{R}^{400 \times R_1}$, $\textbf{U}^{(2)} \in \mathbb{R}^{500 \times R_2}$, $\textbf{U}^{(3)} \in \mathbb{R}^{10 \times R_3}$ and $\mathcal{S} \in \mathbb{R}^{R_1 \times R_2 \times R_3}$. Since we don't want to truncate the number of classes, it is necessary that $R_3 = 10$. Remember that $\textbf{U}^{(1)}$ acts on the mode-$1$ fibers of $\mathcal{T}$, which are the vectors $\mathcal{T}_{:jk}$. Each one of these vectors correspond to a vectorized image, so we can regard $\textbf{U}^{(1)}$ as a compression matrix, which compress images from $400$ pixels to $R_1$ pixels. The matrix $\textbf{U}^{(2)}$ acts of the fibers $\mathcal{T}_{i:k}$, which correspond to the pixel $i$ of all images in class $k$. This matrix ``decides'' if all images are really necessary to describe this pixel, so we can regard $\textbf{U}^{(2)}$ as a matrix of redundancy detection, which reduces the number of images from $500$ to $R_2$, for all classes. The matrix $\textbf{U}^{(3)}$ is a square nonsingular so it doesn't change the dimension. It just makes some rotations of minor importance. We want to work with $\mathcal{S}$ but without the rotated classes. For this we just disconsider $\textbf{U}^{(3)}$ from the definition of $\mathcal{S} = \big( \textbf{U}^{(1)^T}, \textbf{U}^{(2)^T}, \textbf{U}^{(3)^T} \big) \cdot \mathcal{T}$. Thus we define the tensor $\mathcal{F} = \big( \textbf{U}^{(1)^T}, \textbf{U}^{(2)^T}, \textbf{I}_{10} \big) \cdot \mathcal{T} \in \mathbb{R}^{R_1 \times R_2 \times 10}$. The training stage basically consists in computing the MLSVD and constructing the tensor $\mathcal{F}$.
		
		The computation of a truncated MLSVD in Tensor Fox needs some value $R$ for the rank as parameter, even if we are not interested in computing the CPD (as is the case here). The reason for this was mentioned in~\ref{compression}, we can't have $R_\ell > R$ for $\ell=1,2,3$, and in order to assure this it is necessary to give a rank estimate. Furthermore, if the tensor has enough noise, it will be hard to truncate and the program will choose the truncation $R \times R \times 10$. After some tests we choose to use $R = 40$.
 
		Now let's talk about the test stage. For each new image $\textbf{z} \in \mathbb{R}^{400}$, convert it to the $R_1$ dimensional space using the transformation $\textbf{z}_{new} = \textbf{U}^{(1)^T} \cdot \textbf{z}$. In this space we compare $\textbf{z}_{new}$ to the space generated by the columns of each frontal slice of $\mathcal{F}$ (remember that each frontal slice correspond to a class). The best match will be the class of $\textbf{z}_{new}$. This stage of comparison is not made using the slice directly, but instead we compute the SVD of the slice and truncate it to another with lower rank. The overall procedure (training + test) is summarized in the steps below.

		\begin{enumerate}	
			\item Compute truncated MLSVD $\mathcal{T} = (\textbf{U}^{(1)}, \textbf{U}^{(2)}, \textbf{U}^{(3)}) \cdot \mathcal{S}$
			\item Define $\mathcal{F} = \big( \textbf{U}^{(1)^T}, \textbf{U}^{(2)^T}, \textbf{I}_{10} \big) \cdot \mathcal{T}$
			\item Input: $\textbf{z} \in \mathbb{R}^{400}$   
			\item Compress input: $\textbf{z}_{new} = \textbf{U}^{(1)^T} \cdot \textbf{z} \in \mathbb{R}^{R_1}$
			\item Compute SVD of slices: $\mathcal{F}_{::k} = \textbf{W}^{(k)} \Sigma^{(k)} \textbf{V}^{(k)^T}$ for $k = 1, \ldots, 10$
			\item Truncate $\textbf{W}^{(k)}$ to have $R_1' < R_1$ columns: $\tilde{\textbf{W}}^{(k)} = \left[ \textbf{W}^{(k)}_1, \ldots, \textbf{W}^{(k)}_{R_1'} \right]$
			\item Solve least squares problems: $\displaystyle \min_{\textbf{x}} \| \tilde{\textbf{W}}^{(k)} \cdot \textbf{x} - \textbf{z}_{new} \|$ for $k = 1, \ldots, 10$
			\item The class $k$ associated with the smallest error (of all least squares problems solved) is chosen to be the class of $\textbf{z}$
		\end{enumerate}

		Note that, since the columns of $\tilde{\textbf{W}}^{(k)}$ are orthonormal, the solution of the least squares problem is given by
		$$\textbf{x}_\ast = \tilde{\textbf{W}}^{(k)^T} \cdot \textbf{z}_{new},$$
which is easy to compute at this point. It should be noted that the matrices $\tilde{\textbf{W}}^{(k)}$ doesn't have to be computed for each new input. Compute them just once and use them for all new incoming inputs. Finally, the choice of the SVD truncations $R_1', R_2'$ depends on each case. There are some known procedures to guess good values but we won't address this issue here.

		We used 2000 images for training and 3000 for testing. The dimensions of $\mathcal{S}$ are not bigger than $40 \times 40 \times 10$. Since the training data is $400 \times 200 \times 10$, we see that $100 \left(1 - \frac{40 \cdot 40 \cdot 10}{400 \cdot 200 \cdot 10} \right) = 98 \%$ of the data is compressed. Working with only $2 \%$ of the information may sound too few, but only with this much of data we could obtain more than $95 \%$ of accuracy in just a few seconds. Furthermore, this technique clearly can be generalized for practically any classification problem with minor changes. This example can be reproduced with the code in the following link: \url{https://github.com/felipebottega/Tensor-Fox/blob/master/examples/handwritten_digit_mlsvd.ipynb}.
  
	\section{Tensor learning vs. neural network} \label{cpd-nn}
		Given an input data $\textbf{x} = [x_1, \ldots, x_n]^T \in \mathbb{R}^n$, the usual neural network uses a certain number of layers with weights and activation functions in order to compute the respective output $\textbf{y} \in \mathbb{R}^m$. Each layer of weights forms is a matrix, usually bigger than $n \times n$ (except in the last layer). As an alternative to the neural network model, we propose to use $m$ tensors $\mathcal{T}_1, \ldots, \mathcal{T}_m \in \mathbb{R}^{I_1 \times \ldots \times I_L}$, given by

		$$\mathcal{T}_k = \sum_{r=1}^R \textbf{w}_r^{(k,1)} \otimes \ldots \otimes \textbf{w}_r^{(k,L)},$$
for $I_1 = I_2 = \ldots = I_L = n$ and $k = 1 \ldots m$.

		Each $\textbf{w}_r^{(k,\ell)} \in \mathbb{R}^n$ is a vector of weights. These weights are computed in the same way we compute the weights of a neural network, that is, by minimizing a cost function. In this model, the input-output function, also called the \emph{hypothesis}, is the function $h: \mathbb{R}^n \to \mathbb{R}^m$ defined as

		\begin{equation} \label{hypothesis}			
			h(\textbf{x}) = \textbf{f}\left[
			\begin{array}{c}
    			\mathcal{T}_1(\underbrace{\textbf{x}, \ldots, \textbf{x}}_{L \text{ times}})\\
    			\vdots\\
    			\mathcal{T}_m(\underbrace{\textbf{x}, \ldots, \textbf{x}}_{L \text{ times}})\\
			\end{array}
			\right]
		\end{equation}
where $\textbf{f}:\mathbb{R}^m \to \mathbb{R}$ is defined as $\textbf{f} = (\underbrace{f, \ldots, f}_{m \text{ times}})$ for a nonlinear function $f:\mathbb{R} \to \mathbb{R}$. If one wants to write~\ref{hypothesis} in explicit form, then we have that

		$$h(\textbf{x}) = \textbf{f} \left[
		\begin{array}{c}
   			\sum_{r=1}^R \prod_{\ell=1}^L \langle \textbf{w}_r^{(1,\ell)}, \textbf{x} \rangle  \\
   			\vdots\\
   			\sum_{r=1}^R \prod_{\ell=1}^L \langle \textbf{w}_r^{(m,\ell)}, \textbf{x} \rangle  \\
		\end{array}
		\right] = $$

		$$ = \left[
		\begin{array}{c}
   			f\left( \sum_{r=1}^R \sum_{i_1, \ldots, i_L=1}^n w_{i_1 r}^{(1,1)} \ldots w_{i_\ell r}^{(1,L)} x_{i_1} \ldots x_{i_L} \right) \\
   			\vdots\\
   			f\left( \sum_{r=1}^R \sum_{i_1, \ldots, i_L=1}^n w_{i_1 r}^{(m,1)} \ldots w_{i_\ell r}^{(m,L)} x_{i_1} \ldots x_{i_L} \right) \\
		\end{array}
		\right].$$

		Notice that $rank(\mathcal{T}_k) \leq R$ by construction. As is usual in machine learning algorithms, we use the convention $\textbf{x} = [1, x_2, \ldots, x_n]^T$, that is, $x_1 = 1$ is the bias of the model. Now let $S = \{ (i_1, \ldots, i_L): 1 \leq i_1, \ldots, i_L \leq n, \ (i_1, \ldots, i_L) \neq (1, \ldots, 1) \}$, then we can write the above expression as 

		$$ \left[
		\begin{array}{c}
   			f\left( b_1 + \sum_{r=1}^R \sum_{(i_1, \ldots, i_L) \in S} w_{i_1 r}^{(1,1)} \ldots w_{i_L r}^{(1,L)} x_{i_1} \ldots x_{i_L} \right) \\
   			\vdots\\
   			f\left( b_m + \sum_{r=1}^R \sum_{(i_1, \ldots, i_L) \in S} w_{i_1 r}^{(m,1)} \ldots w_{i_L r}^{(m,L)} x_{i_1} \ldots x_{i_L} \right) \\
		\end{array}
		\right]$$
where each $b_k = \sum_{r=1}^R (w_{1 r}^{(k,1)})^{L}$ is the bias parameter of the model.

		In a neural network, each activation function receives a weighted sum of the form 
		$$\sum_{j=1}^n w_j x_j,$$
whereas in the tensor learning we have multilinear weighted sums of the form 

		$$\sum_{i_1, \ldots, i_L=1}^n w_{i_1 r}^{(k,1)} \ldots w_{i_L r}^{(k,L)} x_{i_1} \ldots x_{i_L}$$ 
and just one activation per output. Contrary to the neural network model where we have lots of activations and few direct interactions between the weights\footnote{By \emph{direct interactions} we mean multiplications between the weights.}, in the tensor learning approach we have few activations and lots of direct interactions between the weights, through the multilinearity. Our hope is that this multilinearity compensates the lack of activation functions.

		\subsection{Tensor learning as a special neural network}
			For a fixed input, note that each multilinear weighted sum shown above is just a polynomial of degree $L$ with some weights as variables. This polynomial is such that all monomials have degree $L$. From this point of view one can interpret this model as a neural network (with no hidden layers) where we just changed from dot products to polynomials, a specific class of polynomials. Although this is a totally valid interpretation, it misses the importance of what these polynomials represent. They are tensors evaluated at the input, and their structure is such that all weights are directly mixed during the learning, whereas in the classic neural network\footnote{By \emph{classic neural network} I mean the neural network model in which only dot products are passed to the activation functions.} the weights are not so much mixed. This lack of interaction between the weights in the classic neural network is addressed by the tensor learning. 
				
			It should be noted that this model probably wouldn't be noticed by the usual machine learning practitioner even if it is desired to use polynomial weights as input to the activation functions. The usual way to model weights with polynomials is by introducing a weight per monomial. Therefore, one used to work with machine learning models would probably come with a polynomial of the form
				
			$$\sum_{i_1, \ldots, i_L=1}^n w_{i_1 \ldots i_L} x_{i_1} \ldots x_{i_L},$$
for $w_{i_1 \ldots i_L} \in \mathbb{R}$. This would cause a lost of interactions between the weights. The tensor approach is not one usually would think, unless one is working with tensors, as it is the case of this thesis. 
				
\begin{figure}[h] 
				\includegraphics[scale=.45]{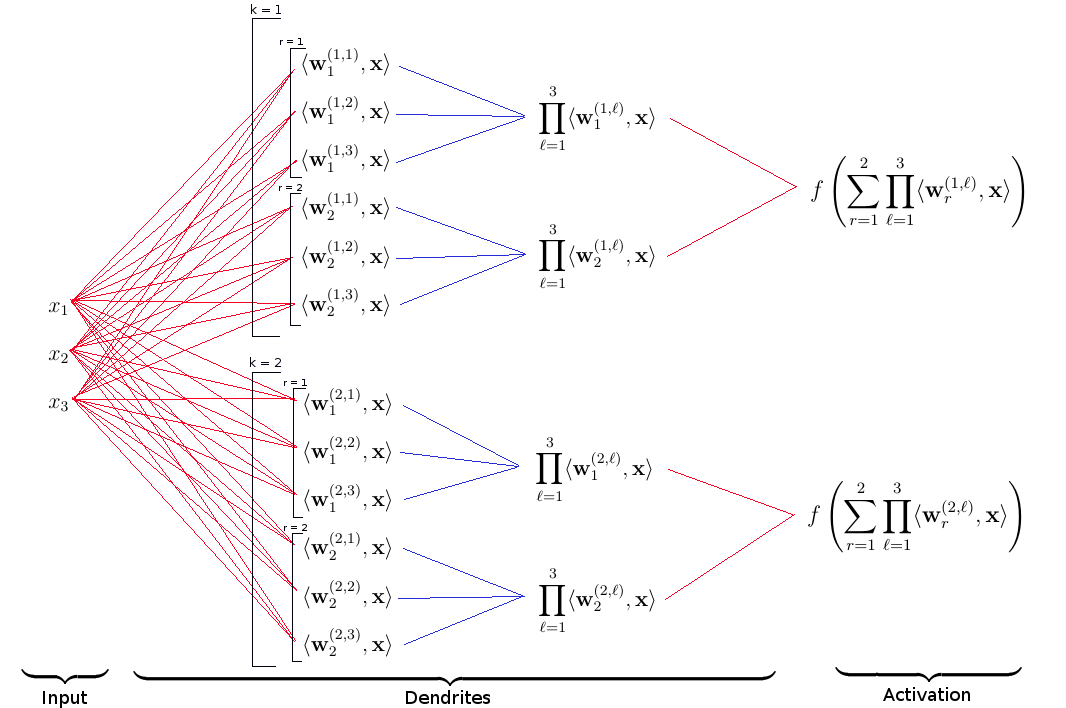}
				\caption{\footnotesize{Tensor learning as a special neural network.}}
				\label{neural}
			\end{figure}

			It is possible to interpret tensor learning as a neural network (with 2 hidden layers) which relies almost entirely on the classic dot product plus activation function recipe. This interpretation shades more light on what is really happening and why tensor learning may be a valid alternative to the classic neural network. Figure~\ref{neural} shows the architecture of such interpretation for $n=3, m=2, L=2, R=2$. The red lines represents classic neural network connections, that is, each line has a weight which is multiplied by the input and summed with all other lines reaching the same neuron. Blue lines means to multiply the inputs by the inputs of all other lines reaching the same neuron. The first layer of lines contains all the weights in our model. The second layer of lines makes the multiplications, and is in this part where we the model is no more the classic one. The third layer of lines are classic neural network connections with all weights equal to 1. Except for the last layer of neurons (the output layer), all other neurons have activation function equals the identity function. In the last layer the activation function is the nonlinear function $f$. We can note that the blue lines are the great novelty in this approach compared to the classic neural network. They represent the part where the weights are mixed together. Also note that this special neural network has $m R L + m R + m$ neurons. This usually will be much smaller than any neural network applied to the same problem, because a neural network uses $\mathcal{O}(n)$ neurons at the first layers, and we almost always have that $m, R, L \ll n$.
				
			One could argue that these activations are too artificial and doesn't match any notion of a neuron. A new, and very different interpretation, is given in figure~\ref{neural}. The idea is that the dendrites intercept sometimes before the layer of activations, which is the last one. Each node of the dendrite layers correspond to a point of intersection between the dendrites. The information travels through a dendrite while also touching other dendrites, and this causes the interaction between the weights and the data. When it finally gets to the activation, we will have a multilinear combination to process. This interpretation correspond to a scenario where the dendrites path are not isolated, they touch each other several times and the information between one neuron and the other is affected by these interactions.
				
		\subsection{Cost function}
			Denote $\textbf{W} = \left( \textbf{W}^{(1,1)}, \ldots, \textbf{W}^{(1,L)}, \ldots, \textbf{W}^{(m,1)}, \ldots, \textbf{W}^{(m,L)} \right)$, where we have that $\mathcal{T}_k = \Big( \textbf{W}^{(k,1)}, \ldots, \textbf{W}^{(k,L)} \Big) \cdot \mathcal{I}_{R \times \ldots \times R}$, using the multilinear multiplication notation. Let $(\textbf{x}^{(j)}, \textbf{y}^{(j)})$, $j=1 \ldots N$, be the training dataset, where we expect to have $n, m, L, R \ll N$. We will use the classic cost function $J: \mathbb{R}^{L m n R} \to \mathbb{R}$ defined as

			$$J(\textbf{W}) = \frac{1}{2N} \sum_{j=1}^N \Big\| h(\textbf{x}^{(j)}) - \textbf{y}^{(j)} \Big\|^2 = $$

			$$ = \frac{1}{2N} \sum_{j=1}^N \sum_{k=1}^m \left( h_k(\textbf{x}^{(j)}) - y_k^{(j)} \right)^2 = $$

			$$ = \frac{1}{2N} \sum_{j=1}^N \sum_{k=1}^m \left( f \left( \sum_{r=1}^R \prod_{\ell=1}^L \langle \textbf{w}_r^{(k,\ell)}, \textbf{x}^{(j)} \rangle \right) - y_k^{(j)} \right)^2.$$
				
			\begin{theorem}
				For $k' = 1 \ldots m$, $\ell' = 1 \ldots L$, $i' = 1 \ldots n$ and $r' = 1 \ldots R$, we have that
					
				$$ \frac{\partial J}{\partial w_{i' r'}^{(k',\ell')}} = \frac{1}{N} \sum_{j=1}^N \nabla J_{k',\ell',i',r'}^{(j)},$$
where	
				$$\hspace{-0.5cm} \nabla J_{k',\ell',i',r'}^{(j)} = \left( f \left( \sum_{r=1}^R \prod_{\ell=1}^L \langle \textbf{w}_r^{(k',\ell)}, \textbf{x}^{(j)} \rangle \right) - y_{k'}^{(j)} \right) \cdot f'\left( \sum_{r=1}^R \prod_{\ell=1}^L \langle \textbf{w}_r^{(k',\ell)}, \textbf{x}^{(j)} \rangle \right) \cdot \left( x_{i'}^{(j)} \prod_{\ell=1, \ell \neq \ell'}^L \langle \textbf{w}_{r'}^{(k',\ell)}, \textbf{x}^{(j)} \rangle \right).$$					
			\end{theorem}

			\textbf{Proof:}

			$$ \frac{\partial J}{\partial w_{i' r'}^{(k',\ell')}} = \frac{1}{2N} \sum_{j=1}^N \sum_{k=1}^m \frac{\partial}{\partial w_{i' r'}^{(k',\ell')}} \left( h_k(\textbf{x}^{(j)}) - y_k^{(j)} \right)^2 = $$

			$$ = \frac{1}{2N} \sum_{j=1}^N \sum_{k=1}^m \left[ 2 \left( h_k(\textbf{x}^{(j)}) - y_k^{(j)} \right) \cdot \frac{\partial}{\partial w_{i' r'}^{(k',\ell')}} h_k(\textbf{x}^{(j)}) \right] = $$

			$$ = \frac{1}{2N} \sum_{j=1}^N \sum_{k=1}^m \left[ \left( h_k(\textbf{x}^{(j)}) - y_k^{(j)} \right) \cdot f'\left( \sum_{r=1}^R \prod_{l=1}^\ell \langle \textbf{w}_r^{(k,l)}, \textbf{x}^{(j)} \rangle \right) \cdot \frac{\partial}{\partial w_{i' r'}^{(k',\ell')}} \left( \sum_{r=1}^R \prod_{l=1}^\ell \langle \textbf{w}_r^{(k,l)}, \textbf{x}^{(j)} \rangle \right) \right] = $$

			$$ = \frac{1}{N} \sum_{j=1}^N \left[ \left( h_{k'}(\textbf{x}^{(j)}) - y_{k'}^{(j)} \right) \cdot f'\left( \sum_{r=1}^R \prod_{l=1}^\ell \langle \textbf{w}_r^{({k'},l)}, \textbf{x}^{(j)} \rangle \right) \cdot \left( \frac{\partial}{\partial w_{i' r'}^{(k',\ell')}} \prod_{l=1}^\ell \langle \textbf{w}_{r'}^{(k',l)}, \textbf{x}^{(j)} \rangle \right) \right] = $$

			$$\hspace{-1cm} = \frac{1}{N} \sum_{j=1}^N \left[ \left( h_{k'}(\textbf{x}^{(j)}) - y_{k'}^{(j)} \right) \cdot f'\left( \sum_{r=1}^R \prod_{l=1}^\ell \langle \textbf{w}_r^{(k',l)}, \textbf{x}^{(j)} \rangle \right) \cdot \left( \frac{\partial}{\partial w_{i' r'}^{(k',\ell')}} \langle \textbf{w}_{r'}^{(k',\ell')}, \textbf{x}^{(j)} \rangle \prod_{l=1, l \neq \ell'}^\ell \langle \textbf{w}_{r'}^{(k',l)}, \textbf{x}^{(j)} \rangle \right) \right] = $$

			$$ = \frac{1}{N} \sum_{j=1}^N \left[ \left( h_{k'}(\textbf{x}^{(j)}) - y_{k'}^{(j)} \right) \cdot f'\left( \sum_{r=1}^R \prod_{l=1}^\ell \langle \textbf{w}_r^{(k',l)}, \textbf{x}^{(j)} \rangle \right) \cdot \left( x_{i'}^{(j)} \prod_{l=1, l \neq \ell'}^\ell \langle \textbf{w}_{r'}^{(k',l)}, \textbf{x}^{(j)} \rangle \right) \right] = $$

			$$\hspace{-0.4cm} = \frac{1}{N} \sum_{j=1}^N \left( f \left( \sum_{r=1}^R \prod_{l=1}^\ell \langle \textbf{w}_r^{(k',l)}, \textbf{x}^{(j)} \rangle \right) - y_{k'}^{(j)} \right) \cdot f'\left( \sum_{r=1}^R \prod_{l=1}^\ell \langle \textbf{w}_r^{(k',l)}, \textbf{x}^{(j)} \rangle \right) \cdot \left( x_{i'}^{(j)} \prod_{l=1, l \neq \ell'}^\ell \langle \textbf{w}_{r'}^{(k',l)}, \textbf{x}^{(j)} \rangle \right) = $$

			$$\hspace{6cm} = \frac{1}{N} \sum_{j=1}^N \nabla J_{k',\ell',i',r'}^{(j)}. \hspace{6cm}\square$$
				
		\subsection{Regularization}
			If one want to use regularization, just introduce a regularization parameter $\lambda > 0$ and work with the cost function

			$$J(\textbf{W}) = \frac{\lambda}{2N} \sum_{j=1}^N \| \textbf{W} \|^2 + \frac{1}{2N} \sum_{j=1}^N \Big\| h(\textbf{x}^{(j)}) - \textbf{y}^{(j)} \Big\|^2.$$

			Now we have that

			$$\nabla J_{k',\ell',i',r'}^{(j)} = $$
				
			$$\hspace{-0.5cm} = \lambda w_{i' r'}^{(k', \ell')} + \left( f \left( \sum_{r=1}^R \prod_{\ell=1}^L \langle \textbf{w}_r^{(k',\ell)}, \textbf{x}^{(j)} \rangle \right) - y_{k'}^{(j)} \right) \cdot f'\left( \sum_{r=1}^R \prod_{\ell=1}^L \langle \textbf{w}_r^{(k',\ell)}, \textbf{x}^{(j)} \rangle \right) \cdot \left( x_{i'}^{(j)} \prod_{\ell=1, \ell \neq \ell'}^L \langle \textbf{w}_{r'}^{(k',\ell)}, \textbf{x}^{(j)} \rangle \right).$$
				
		\subsection{Stochastic gradient}
			After some initial guess $\textbf{W} = \left( \textbf{W}^{(1,1)}, \ldots, \textbf{W}^{(1,L)}, \ldots, \textbf{W}^{(m,1)}, \ldots, \textbf{W}^{(m,L)} \right)$, we could use the gradient descent iteration

			$$\textbf{W} \leftarrow \textbf{W} - \alpha \nabla J(\textbf{W})$$
to obtain sequential improvements for the weights. This procedure is costly and we will prefer to use the stochastic gradient iteration

			$$\textbf{W} \leftarrow \textbf{W} - \alpha \nabla J^{(j)}(\textbf{W}),$$
where $\nabla J^{(j)} = \left[ \nabla J_{k,\ell,i,r}^{(j)} \right]_{k,\ell,i,r}$. After making this iteration for each $j = 1 \ldots N$, we probably will have a good set of weights. Otherwise we may repeat the procedure above sometimes. It should be that this model of learning, by construction, avoids the problem of vanishing gradients, which is a common problem in neural networks.

		\subsection{Computational cost}
			Now we analyze the cost of the stochastic gradient descent algorithm. Fix a sample $\textbf{x}^{(j)}$. First we compute and save each $\langle \textbf{w}_r^{(k,\ell)}, \textbf{x}^{(j)} \rangle$. Since each one costs $n$ flops and we have to compute it $LmR$ times, we have a total cost of $LmnR$ flops. The computational cost of $\nabla J_{k,\ell,i,r}^{(j)}$ is summarized below.

			\begin{flalign*}
   				& \lambda w_{i' r'}^{(k', \ell')}: \quad 1 \text{ flop} \\
    			& f \left( \sum_{r=1}^R \prod_{\ell=1}^L \langle \textbf{w}_r^{(k',\ell)}, \textbf{x}^{(j)} \rangle \right) - y_{k'}^{(j)}: \quad \ell R + 1 \text{ flops} \\
    			& f'\left( \sum_{r=1}^R \prod_{\ell=1}^L \langle \textbf{w}_r^{(k',\ell)}, \textbf{x}^{(j)} \rangle \right): \quad \ell R +1 \text{ flops} \\
    			& \left( x_{i'}^{(j)} \prod_{\ell=1, \ell \neq \ell'}^L \langle \textbf{w}_{r'}^{(k',\ell)}, \textbf{x}^{(j)} \rangle \right): \quad L \text{ flops} \\
    			& ---------------------\\
    			& \textbf{Total: } 1 + 2(L R + 1) + L \text{ flops} 
			\end{flalign*}
			
			To compute all $\nabla J_{k,\ell,i,r}^{(j)}$ we have make these computations $LmnR$ times, so we have a total of $L m n R + 2(L m n R)(L R +1) + L^2 m n R$ flops. Taking in account also the computations of the dot products, we have a total of $3L m n R + 2L^2 m n R^2 + L^2 m n R$ flops. This is for one iteration of the stochastic gradient. Running the program over all the samples has a cost of $(3L m n R + 2L^2 m n R^2 + L^2 m n R)N$ flops. Assuming we have to repeat the whole process $E$ times (number of epochs), we have a final cost of 

			$$(3L m n R + 2L^2 m n R^2 + L^2 m n R)NE \text{ flops}.$$

			Now assume a neural network with $\tilde{L}$ layers and $Cn$ neurons in each layer (except the last one), where $C > 1$. For each layer we have to compute $Cn$ dot products (one per neuron), except the last layer which we have to compute only $m$ dot products. This gives us a total of $Cn(\tilde{L}-1) + m$ dot products. Since each dot product costs $n$ flops, this amounts to a total of $Cn^2(\tilde{L}-1) + mn\tilde{L}$ flops. We also have to evaluate as many function evaluations as neurons, which gives us a total of $Cn(\tilde{L}-1) + m$ flops. Therefore we need to perform 

			$$Cn^2(\tilde{L}-1) + mn\tilde{L} + Cn(\tilde{L}-1) + m = (n^2 + n)C(\tilde{L}-1) + mn\tilde{L} + m \text{ flops}$$
in order to make the forward propagation of one sample. Running this over all the samples has a cost of $((n^2 + n)C(\tilde{L}-1) + mn\tilde{L} + m)N$ flops. Assuming we have to repeat the whole process $E$ times, we have a final cost of 

			$$((n^2 + n)C(\tilde{L}-1) + mn\tilde{L} + m)NE \text{ flops}.$$ 

			Below we put both costs together for comparison.

			$$\textbf{Tensor learning: } (3L m n R + 2L^2 m n R^2 + L^2 m n R)NE \text{ flops}$$
			$$\textbf{Neural network: } ((n^2 + n)C(\tilde{L}-1) + mn\tilde{L} + m)NE \text{ flops}$$ 

			We expect to have $m, R \ll n$ and $L < \tilde{L}$ so the cost of the tensor learning is, in general, much lower than the neural network. There is no rule of thumb about the number of layers and the number of neurons at each layer, but some reasonable values, for example, would be $R=3, C=2$ and $\tilde{L} = L = 3$. In this case we have the following costs.

			$$\hspace{-1.9cm}\textbf{Tensor learning: } 3^3mn + 2 \cdot 3^4mn + 3^3mn = 216mn  \text{ flops}$$
			$$\hspace{0.3cm}\textbf{Neural network: } (n^2 + n)2^2 + 3mn + m = 4n^2 + 4n + 3mn + m  \text{ flops}$$
			
			After the weights are trained we have a CPD format $(\textbf{W}^{(1)}, \ldots, \textbf{W}^{(L)}) \cdot \mathcal{I}_{R \times \ldots \times R}$ to makes predictions. As already observed, $rank(\mathcal{T}_k) \leq R$ by construction. Therefore it is interesting to compute low rank CPD approximations for all the tensors $\mathcal{T}_k$, hence simplifying the model. We tested this formulation over the MNIST handwritten digits, the same tensor described in the previous section. Making $L=3$ (third order tensors) and $R=20$ we could find a set of weights which predicts the entire dataset with $97.8 \%$ of accuracy. After that we start computing low rank approximations and the new rank $R' = 17$ worked very well, producing a smaller set of weights which models the entire dataset with the practically the same accuracy (it achieved $97.628 \%$).

%% file: Conclusion.tex
\chapter{Conclusions}\label{conclusion}
	Chapters 1 and 2 were devoted to the study of tensor decompositions and some geometrical aspects of the tensor rank. In chapter 3 we introduced the Gauss-Newton algorithm and showed how it can be used to the problem of finding an approximated CPD. The block structure of the approximated Hessian is of big importance. The study of this structure allowed us to make fast matrix-vector products for the CG algorithm. Chapter 4 was devoted to the development, in details, of the algorithms used within Tensor Fox. These algorithms now are part of a tensor package which includes not only a CPD solver, but also many routines of multilinear algebra and machine learning tools. The CPD algorithm developed in this work was compared to the state of art and it proved to be competitive, both in terms of speed and accuracy. Finally, in chapter 5 we introduced the concept of tensor learning and showed how tensor techniques can be used to machine learning problems.
	
	This thesis now is part of an undergoing project. The solver we develop is fast, accurate and scalable, but it still has room to improvements. Below we list the points to be addressed in the future.
	
	\begin{itemize}
		\item The diagonal preconditioner is certainly not the best preconditioner one could use to solve the subproblems of the dGN algorithm. For instance Tensorlab uses a block diagonal preconditioner, thus reducing the number of CG iterations much more times than Tensor Fox. The drawback of this approach is the necessity of solving more linear systems in order to obtain these preconditioners. It would be ideal to exploit the block structure of the approximated Hessian somehow to obtain a better preconditioner with few computations.
		
		\item The rule of the current damping parameter update doesn't seem to be optimal. Designing and testing new rule updates could bring more improvements.
		
		\item In recent work \cite{rgh} it was verified that Tensorlab can produce very accurate solutions that are extremely ill-conditioned, even when the objective solution is well-conditioned. In~\ref{cond_section} we could verify that Tensor Fox suffers from the same problem. It is not a secret that problems in high-dimensional spaces usually are full of \emph{degenerate saddle points} \cite{anandkumar3}. There is the possibility that both solvers are converging to these points, which happen to be close to the objective but are ill-conditioned. Since they are ``looking'' only for the forward error of the approximation, nothing are being made to avoid ill-conditioned solutions. Increasing the number of CG iterations could address this problem but this is costly.
		
		\item Currently Tensor Fox is able to handle large sparse tensors with low memory cost. It is of interest to make use of GPUs to perform some computations for these kind of tensors. 
		
		\item Expanding the technique of tensor learning present in~\ref{cpd-nn} can be fruitful for machine learning problems. Although the architecture of the tensor approach still acts like a black box, the outputs consists of multilinear maps, which are understandable. Additionally, comparing the tensor model with other kind of neural networks is something to be addressed. It is also of interest to test the tensor model with different orders, different datasets, a wider range of ranks and so on. In short, there are a lot to be done in this subject, and it is the hope of the author that this first step is a step in a rich direction 
	\end{itemize}

%% file: appenA.tex
\chapter{Numerical linear algebra} \label{appenA}

	\section{Singular value decomposition}

		In this section $\mathbb{K}$ can be $\mathbb{R}$ or $\mathbb{C}$. We start with the singular value decomposition (SVD). The condition $m \geq n$ won't be assumed.
	
		\begin{theorem}[Reduced SVD] 
			Let $\textbf{M} \in \mathbb{K}^{m \times n}$. Then there exists unitary (orthogonal) matrices\footnote{Here we are talking about rectangular unitary (orthogonal) matrices. There are two definitions for a rectangular matrix $\textbf{U} \in \mathbb{K}^{m \times n}$. In the case $m \geq n$ we say $\textbf{U}$ is unitary (orthogonal) if its columns are orthonormal (which is equivalent to $\textbf{U}^\ast \textbf{U} = \textbf{I}_n$). In the case $m \leq n$ we say $\textbf{U}$ is unitary (orthogonal) if its rows are orthonormal (which is equivalent to $\textbf{U} \textbf{U}^\ast = \textbf{I}_m$).} $\textbf{U} \in \mathbb{K}^{m \times n}, \textbf{V} \in \mathbb{K}^{n \times n}$ and a diagonal matrix $\Sigma = \text{diag}(\sigma_1, \ldots, \sigma_n) \in \mathbb{R}^{n \times n}$, with $\sigma_1 \geq \ldots \geq \sigma_n$, such that $\textbf{M} = \textbf{U} \Sigma \textbf{V}^\ast$.
		\end{theorem}

		The columns of $\textbf{U}$ are called \emph{left singular vector}, the columns of $\textbf{V}$ are called \emph{right singular vectors}, and the values $\sigma_i$ are called \emph{singular values}. The decomposition given in the theorem is the \emph{reduced SVD of $\textbf{M}$}. In the \emph{full SVD of $\textbf{M}$} we can write $\textbf{M} = \textbf{U} \Sigma \textbf{V}^\ast$, where $\textbf{U} \in \mathbb{K}^{m \times m}$ and $\textbf{V} \in \mathbb{K}^{n \times n}$ are unitary (orthogonal), and $\Sigma = \text{diag}(\sigma_1, \ldots, \sigma_p) \in \mathbb{R}^{m \times n}$, where $p = \min\{ m,n\}$. Note that $\Sigma$ has the same dimensions of $\textbf{M}$ so we are calling it ``diagonal'' in accord with the structures below. There are three possible cases.\\
	
		\underline{\textbf{$m = n$}}
		$$ \textbf{M} =\textbf{U}
		\left[ \begin{array}{ccc}
			\sigma_1 & & \\
			& \ddots & \\
			& & \sigma_ n\\
		\end{array}\right]
		\textbf{V}^\ast$$	

		\underline{\textbf{$m > n$}}
		$$ \textbf{M} = \textbf{U}
		\left[ \begin{array}{ccc}
			\sigma_1 & & \\
			& \ddots & \\
			& & \sigma_n \\
			& & \\
			& & \\
		\end{array}\right]
		\textbf{V}^\ast$$	

		\underline{\textbf{$m < n$}}
		$$ \textbf{M} = \textbf{U}
		\left[ \begin{array}{cccccc}
			\sigma_1 & &\ &\ &\ & \\
			& \ddots & &\ &\ &\ \\
			& & \sigma_m &\ &\ &\ \\
		\end{array}\right]
		\textbf{V}^\ast$$\\	
	
		The following results will be useful.
	
		\begin{lemma} \label{svd-transpose-lemma}
			Let $\textbf{M} \in \mathbb{K}^{m \times m}$ be symmetric with eigendecomposition $\textbf{M} = \textbf{U} \Lambda \textbf{U}^\ast$ where $\textbf{U} \in \mathbb{K}^{m \times m}$ is unitary (orthogonal) and  $\Lambda = \text{diag}(\lambda_1, \ldots, \lambda_m) \in \mathbb{K}^{m \times m}$. Furthermore, assume the eigenvalues $\lambda_i$ are such that $|\lambda_1| \geq \ldots \geq |\lambda_m|$. Then an SVD of $\textbf{M}$ is $\textbf{M} = \textbf{U} \Sigma \textbf{V}^T$, where $\sigma_i = |\lambda_i|$ and $\textbf{V}_{:i} = \text{sign}(\lambda_i) \textbf{U}_{:i}$.
		\end{lemma}
	
		\begin{lemma} 
			Let $\textbf{M} \in \mathbb{K}^{m \times n}$ such that $m \leq n$ and let $\sigma_1 \geq \ldots \geq \sigma_m$ be the singular values of $\textbf{M}$. Then the singular values of $\textbf{M} \textbf{M}^T$ are $\sigma_1^2, \ldots, \sigma_m^2$.
		\end{lemma}
	
		\begin{theorem} \label{svd-transpose}
			Let $\textbf{M} \in \mathbb{K}^{m \times n}$ such that $m \leq n$ and let $\lambda_1, \ldots, \lambda_m$ be the eigenvalues of $\textbf{M} \textbf{M}^T$, where $|\lambda_1| \geq \ldots \geq |\lambda_m|$. Then the singular values of $\textbf{M}$ are $\sqrt{|\lambda_1|} \geq \ldots \geq \sqrt{|\lambda_m|}$.
		\end{theorem}
		
		\begin{remark}\label{MLSVD-SVD}
			This theorem is particularly useful when $m < n$. Suppose we want to compute a SVD of $\textbf{M}$ in this situation. Then the usual cost is at least of $2m^2n + 2m^3$ flops (see lecture 31 of \cite{trefethen}). Instead of this we can compute $\textbf{M} \textbf{M}^T$ first, at the cost of $m^2n$ flops. Then the cost of computing the eigenvalues of a $m \times m$ symmetric matrix is of $\mathcal{O}( m^3 )$ flops with a low constant, giving a total of $\mathcal{O}\left( m^2n + m^3 \right)$ flops. 
		\end{remark}
		
	\section{Conjugate gradient}
		 Let $\textbf{A} \in \mathbb{R}^{n \times n}$ be symmetric positive definite and $\textbf{b} \in \mathbb{R}^n$. Suppose we want to solve a linear system
		\begin{equation} \label{CG}		 
			\textbf{A} \textbf{x} = \textbf{b}
		\end{equation}
for $\textbf{x}$. The solution of this system is $\textbf{x}_\ast = \textbf{A}^{-1} \textbf{b}$. The \emph{$i$-th Krylov space generated by $\textbf{A}$ and $\textbf{b}$} is defined as 

		$$ \mathcal{K}_i = \text{span}(\textbf{b}, \textbf{A} \textbf{b}, \ldots, \textbf{A}^{i-1} \textbf{b}).$$
		
		Note that $\mathcal{K}_i$ is a subspace of $\mathbb{R}^{n \times n}$. Moreover. the vectors generating the space are orthonormal so they form a orthonormal basis for $\mathcal{K}_i$. Any algorithm that constructs these vectors sequentially is called a \emph{Krylov method}. One of the advantages of Krylov methods is that sometimes they can be carried on without actually storing the matrix $\textbf{A}$. Depending on the structure of $\textbf{A}$ it is possible to accomplish this. When this is possible we way the method is \emph{matrix-free}. For example, one can see that theorem~\ref{Hv} allow us to compute the product $\textbf{J}_f^T \textbf{J}_f \cdot \textbf{x}$ without needing to construct $\textbf{J}_f^T \textbf{J}_f$ explicitly.
		
		The conjugate gradient method is a Krylov method made specifically to deal with linear systems with symmetric positive definite matrices. Given the system~\ref{CG}, this method iteratively constructs a sequence of approximations $\textbf{x}^{(i)}$ which minimizes $\varepsilon^{(i)} = \| \textbf{x}_\ast - \textbf{x}^{(i)} \|_{\textbf{A}}$, the \emph{$\textbf{A}$-residual of $\textbf{x}^{(i)}$}, where $\| \textbf{x} \|_{\textbf{A}} = \sqrt{\langle \textbf{A} \textbf{x}, \textbf{x} \rangle}$. We won't go into the details of this method. This is just a brief introduction of the algorithm and some of its properties.\\ 
		
		\begin{algorithm}[Conjugate gradient] \label{cg-alg}
			$ $\\
			\textbf{Input:} $\textbf{A} \in \mathbb{R}^{n \times n}$ symmetric positive definite, $\textbf{b} \in \mathbb{R}^n$ \vspace{.4cm}\\
			$\textbf{x}^{(0)} \leftarrow 0 \in \mathbb{R}^n$\\
			$\textbf{r}^{(0)} \leftarrow \textbf{b}$\\
			$\textbf{p}^{(0)} \leftarrow \textbf{r}^{(0)}$\\
			$i = 1$\\
			\verb|repeat | \\
				$\hspace{1cm} \alpha^{(i)} \leftarrow \frac{ \| \textbf{r}^{(i-1)} \|^2 }{ \| \textbf{p}^{(i-1)} \|_{\textbf{A}}^2 }$\\
				$\hspace{1cm} \textbf{x}^{(i)} \leftarrow \textbf{x}^{(i-1)} + \alpha^{(i)} \textbf{p}^{(i-1)}$\\
				$\hspace{1cm} \textbf{r}^{(i)} \leftarrow \textbf{r}^{(i-1)} - \alpha^{(i)} \textbf{A} \textbf{p}^{(i-1)}$\\
				$\hspace{1cm} \beta^{(i)} \leftarrow \frac{ \| \textbf{r}^{(i)} \|^2 }{ \| \textbf{r}^{(i-1)} \|^2 }$\\
				$\hspace{1cm} \textbf{p}^{(i)} \leftarrow \textbf{r}^{(i)} + \beta^{(i)} \textbf{p}^{(i-1)}$\\
				$\hspace{1cm} i \leftarrow i+1$
			\verb|until stopping criteria is met |\vspace{.4cm}\\
			\textbf{Output:} $\textbf{x}^{(i')}$, where $i'$ is the last index of the iterations
		\end{algorithm}
		
		Each iteration involves several vector manipulations and one matrix-vector multiplication, which appears twice in the algorithm but need to be computed only once. In general the matrix-vector multiplication dominates the computational cost so we can consider that each iteration costs $\mathcal{O}(n^2)$ flops. If $\textbf{A}$ has some structure to be exploited this cost may be lower. The term ``conjugate'' in the name of the algorithm comes from the fact that the ``search direction'' vectors $\textbf{p}^{(i)}$ are $\textbf{A}$-conjugate with respect to each other, that is, $\langle \textbf{p}^{(i)}, \textbf{p}^{(j)} \rangle_{\textbf{A}} = 0$ whenever $i \neq j$.
		
		Usually one can use as stopping condition the distance between the current approximation, $\textbf{x}^{(i)}$, and the objective solution, $\textbf{b}$. The error is measured by $\varepsilon = \| \textbf{b} - \textbf{A} \textbf{x}^{(i)}\|^2$, but a faster way to compute this error is through the identity $\varepsilon = \| \textbf{r}^{(i)} \|^2$. To see this identity holds just note that 
				
				$$\varepsilon = \| \textbf{r}^{(i)} \|^2 = \| \textbf{r}^{(i-1)} - \alpha^{(i)} \textbf{A} \textbf{p}^{(i-1)} \|^2 = $$
				
				$$ = \| \textbf{r}^{(i-2)} - \alpha^{(i-1)} \textbf{A} \textbf{p}^{(i-2)} - \alpha^{(i)} \textbf{A} \textbf{p}^{(i-1)} \|^2 = \ldots $$
				
				$$\ldots = \| \textbf{r}^{(0)} - \alpha^{(1)} \textbf{A} \textbf{p}^{(0)} - \alpha^{(2)} \textbf{A} \textbf{p}^{(1)} - \ldots -  \alpha^{(i)} \textbf{A} \textbf{p}^{(i-1)} \|^2 =$$
				
				$$ = \left\| \textbf{b} - \textbf{A} \sum_{j=1}^i \alpha^{(j)} \textbf{p}^{(j-1)} \right\|^2 = \left\| \textbf{b} - \textbf{A} \textbf{x}^{(i)} \right\|^2.$$
		
		\begin{theorem}
			Consider conjugate gradient algorihtm applied to the system~\ref{CG}. If the iteration $i$ has not already converged ($\textbf{r}^{(i-1)} \neq 0$), then $\textbf{x}^{(i)}$ is the unique point in $\mathcal{K}_i$ that minimizes $\| \varepsilon^{(i)} \|_{\textbf{A}}$. The convergence is monotonic (that is, $\| \varepsilon^{(i)} \|_{\textbf{A}} \leq \| \varepsilon^{(i-1)} \|_{\textbf{A}}$) and $\varepsilon^{(i)} = 0$ is achieved for some $i \leq n$.
		\end{theorem}		    
		
		This theorem highlights some of the most known properties of the conjugate gradient. The last property, in particular, assures the algorithm will find the exact solution so it can be regarded as a direct algorithm, not an iterative one. We say it is iterative because in practice (working with finite precision) this exact solution is not necessarily found in less than $n$ iterations. Usually one just compute some iterations of the algorithm, stopping before if the residual increases, because the monotonic property also can be lost when working in finite precision. 
		
		As we will see now, the rate of convergence depends a lot on the distribution of the eigenvalues of $\textbf{A}$. For the next theorem denote $P_i$ as the set of univariate polynomials $p$ with degree $\leq i$ and $p(0) = 1$. More precisely $P_i$ is the set of polynomials of the form $a_i x^i + a_{i-1} x^{i-1} + \ldots + a_1 x + 1$. We also use the notation $\Lambda( \textbf{A} )$ to be the spectrum of $\textbf{A}$.
		
		\begin{theorem}
			At iteration $i$ of the conjugate gradient algorithm, we have 
			$$\frac{ \| \varepsilon^{(i)} \|_{\textbf{A}} }{ \| \varepsilon^{(0)} \|_{\textbf{A}} } \leq \inf_{p \in P_i} \max_{\lambda \in \Lambda(\textbf{A})} |p(\lambda)|.$$ 
		\end{theorem}
		
		\begin{corollary}
			If $\textbf{A}$ has at most $i$ distinct eigenvalues, then the conjugate gradient converges in at most $i$ iterations.
		\end{corollary}
		
		This corollary basically says that matrices with repeated eigenvalues may perform better. the general idea is that if the eigenvalues of $\textbf{A}$ are group in $i$ small clusters, then we can expect to have convergence in at most $i$ iterations. The corollary is the extreme case where these clusters collapse in just points. 
		
		Let $\kappa(\textbf{A}) = \| \textbf{A} \|_2 \| \textbf{A}^{-1} \|_2$ be the condition number of $\textbf{A}$, where $\| \cdot \|_2$ is the spectral norm. The next result shows that the error is also bounded by a function of the condition number.
		
		\begin{corollary}
			At iteration $i$ of the conjugate gradient algorithm, we have 
			$$\frac{ \| \varepsilon^{(i)} \|_{\textbf{A}} }{ \| \varepsilon^{(0)} \|_{\textbf{A}} } \leq 2 \left( \frac{\sqrt{ \kappa(\textbf{A}) } -1 }{ \sqrt{ \kappa(\textbf{A}) } +1 } \right)^i.$$ 
		\end{corollary}
		
		Since $\kappa(\textbf{A}) = \frac{\lambda_{\max}}{\lambda_{\min}}$, this theorem implies that the largest eigenvalue shouldn't be too far from the smallest. Again, this reinforce the idea that the eigenvalues should be grouped together. Finally, since $\frac{\sqrt{ \kappa(\textbf{A}) } -1 }{ \sqrt{ \kappa(\textbf{A}) } +1 } \approx 1 - \frac{2}{ \sqrt{ \kappa(\textbf{A}) } }$, we have that the convergence to a specified tolerance is expected to be done in $\mathcal{O}(\sqrt{ \kappa(\textbf{A}) })$ iterations. This is just an upper bound, convergence may be faster if the spectrum is clustered. For more information about this subject check lecture 38 of \cite{trefethen} and 6.6.3 of \cite{demmel}.  
		
	\section{Preconditioning}
		We assume the reader is familiar with the subject of condition numbers in linear algebra. Let $\textbf{A} \in \mathbb{R}^{n \times n}$ and $\textbf{b} \in \mathbb{R}^n$. Given a linear system
		\begin{equation} \label{precond}		 
			\textbf{A} \textbf{x} = \textbf{b}
		\end{equation}
to be solved for $\textbf{x}$, one may encounter numerical issues depending on the conditioning of $\textbf{A}$. If $\textbf{A}$ is ill-conditioned there will be problems to compute an accurate solution regardless of the algorithm used. When facing this situation we can use the technique of \emph{preconditioning} to transform the ill-conditioned problem~\ref{precond} into a well-conditioned one. First, suppose we have at our disposal an invertible matrix $\textbf{M} \in \mathbb{R}^{n \times n}$ such that $\textbf{M}$ is well-conditioned and $\textbf{M}^{-1} \textbf{A}$ is easily computable and well-conditioned. $\textbf{M}$ is called the \emph{preconditioner}. Note that the system

		\begin{equation} \label{precond2}
			\textbf{M}^{-1} \textbf{A} \textbf{x} = \textbf{M}^{-1} \textbf{b}
		\end{equation}		  	 
has the same solution of~\ref{precond}. Therefore, we can solve this well-conditioned system instead of an ill-conditioned and obtain the desired solution. 

		We don't actually need to construct the inverse $\textbf{M}^{-1}$. First we solve the system $\textbf{M} \textbf{y} = \textbf{b}$ which should be easy since $\textbf{M}$ is well-conditioned. The solution is $\textbf{y}_\ast = \textbf{M}^{-1} \textbf{b}$. With this solution now we solve the system
		
		\begin{equation} \label{precond3}
			(\textbf{M}^{-1} \textbf{A}) \textbf{x} = \textbf{y}_\ast.
		\end{equation}	
		
		We assume that $\textbf{M}^{-1} \textbf{A}$ is already computed and that it is a well behaved matrix so we can use it to solve the system without any problems. The respective solution is $\textbf{x}_\ast = \textbf{A}^{-1} \textbf{M} \textbf{M}^{-1} \textbf{A}^{-1} \textbf{y}_\ast = \textbf{A}^{-1} \textbf{M} \textbf{M}^{-1} \textbf{b} = \textbf{A}^{-1} \textbf{b}$, which is the desired solution. This is the procedure of transforming and solving an ill-conditioned problem into a well-conditioned one. In the case $\textbf{M}$ is diagonal we don't need to worry since $\textbf{M}^{-1}$ is very easy to obtain.
		
		This approach has to be changed when $\textbf{A}$ is symmetric positive definite and we want to use the conjugate gradient algorithm. The reason for this is that $\textbf{M}^{-1} \textbf{A}$ is generally not symmetric even if $\textbf{M}$ is symmetric. To derive the correct preconditioning assume $\textbf{M}^{-1}$ is symmetric positive definite. We can rewrite the preconditioned system $\textbf{M}^{-1} \textbf{A} \textbf{x} = \textbf{M}^{-1} \textbf{b}$ as 
		
		\begin{equation} \label{precond4}		
			\underbrace{(\textbf{M}^{-1/2} \textbf{A} \textbf{M}^{-1/2})}_{\hat{\textbf{A}}} \underbrace{(\textbf{M}^{1/2} \textbf{x})}_{\hat{\textbf{x}}} = \underbrace{\textbf{M}^{-1/2} \textbf{b}}_{\hat{\textbf{b}}}.
		\end{equation}
		
		Note that $\hat{\textbf{A}}$ is symmetric positive definite so we can apply the conjugate gradient algorithm. Moreover $\hat{\textbf{A}}$ and $\textbf{M}^{-1} \textbf{A}$ have the same eigenvalues since they are similar, therefore the performance of the gradient conjugate algorithm is not altered.  	
		
		One of the simplest preconditioners is the \emph{Jacobi preconditioner}, also called \emph{diagonal preconditioner}. Its defined as $\textbf{M} = \text{diag}(a_{ii})$, where we suppose $\textbf{A}$ has positive diagonal. The result of this preconditioning is the matrix
		
		$$\textbf{M}^{-1} \textbf{A} = 
		\left[ \begin{array}{cccc}
			1 & \frac{a_{12}}{a_{11}} & \ldots & \frac{a_{1n}}{a_{11}}\\
			\frac{a_{21}}{a_{22}} & 1 & \ldots & \frac{a_{2n}}{a_{22}}\\
			\vdots & \vdots & \ddots & \vdots\\
			\frac{a_{n1}}{a_{nn}} & \frac{a_{n2}}{a_{nn}} & \ldots & 1
		\end{array}	\right].$$	
		
		This is particularly useful if $\textbf{A}$ is diagonally dominant. Now suppose $\textbf{A}$ is symmetric positive definite. Then we have
		
		$$\textbf{M}^{-1/2} \textbf{A} \textbf{M}^{-1/2} = 
		\left[ \begin{array}{cccc}
			1 & \frac{a_{12}}{\sqrt{a_{11}} \sqrt{a_{22}}} & \ldots & \frac{a_{1n}}{\sqrt{a_{11}} \sqrt{a_{nn}}}\\
			\frac{a_{21}}{\sqrt{a_{22}} \sqrt{a_{11}}} & 1 & \ldots & \frac{a_{2n}}{\sqrt{a_{22}} \sqrt{a_{nn}}}\\
			\vdots & \vdots & \ddots & \vdots\\
			\frac{a_{n1}}{\sqrt{a_{nn}} \sqrt{a_{11}}} & \frac{a_{n2}}{\sqrt{a_{nn}} \sqrt{a_{22}}} & \ldots & 1
		\end{array}	\right].$$		
		
		Even if $\textbf{A}$ is not diagonally dominant this preconditioner may improve the conditioning by lowering the condition number. It can be shown \cite{sluis, bradley} that 
		$$\kappa(\textbf{A}) \leq n \cdot \kappa(\textbf{M}^{-1/2} \textbf{A} \textbf{M}^{-1/2})$$ 
and
		$$\kappa(\textbf{M}^{-1/2} \textbf{A} \textbf{M}^{-1/2}) \leq n \cdot \min_{\textbf{D} \in \textbf{D}_n} \kappa(\textbf{D} \textbf{A} \textbf{D}),$$
where $\kappa(\textbf{A}) = \| \textbf{A} \|_2 \| \textbf{A}^{-1} \|_2$ and $\textbf{D}_n$ is the set of nonsingular diagonal matrices. With this result we have a lower and an upper bound for $\kappa(\textbf{M}^{-1/2} \textbf{A} \textbf{M}^{-1/2})$, namely,
		
		$$ \frac{1}{n} \kappa(\textbf{A}) \leq \kappa(\textbf{M}^{-1/2} \textbf{A} \textbf{M}^{-1/2}) \leq \kappa(\textbf{M}^{-1/2})^2 \cdot \kappa(\textbf{A}) = \frac{\max_i a_{ii}}{\min_i a_{ii}} \cdot \kappa(\textbf{A}).$$		
		
		 A good survey on preconditioning I would recommend is \cite{precond}, also part 6.6.5 of \cite{demmel} is worth reading. 
			

%% file: appenB.tex
\chapter{Tensor algebra} \label{appenB}
	This appendix section can be viewed as a complement of chapter 1.

	\section{Tensor product properties}
		
		The main algebraic properties about tensor product are summarized below.
		
		\begin{theorem}
		Consider the vector spaces $\mathbb{V}^{(1)}, \mathbb{V}^{(2)}, \mathbb{V}^{(1)}$, the vectors $\textbf{v}^{(1)}, \textbf{u}^{(1)} \in \mathbb{V}^{(1)}$, $\textbf{v}^{(2)}, \textbf{u}^{(2)} \in \mathbb{V}^{(2)}$, $\textbf{v}^{(3)} \in \mathbb{V}^{(3)}$ and a scalar $\alpha \in \mathbb{K}$. Then
			\begin{enumerate}
				\item $\textbf{v}^{(1)} \otimes \textbf{v}^{(2)} = 0 \implies \textbf{v}^{(1)} = 0$ or $\textbf{v}^{(2)} = 0$.
				\item $(\alpha \textbf{v}^{(1)}) \otimes \textbf{v}^{(2)} = \textbf{v}^{(1)} \otimes (\alpha \textbf{v}^{(2)}) = \alpha (\textbf{v}^{(1)} \otimes \textbf{v}^{(2)})$
				\item $\textbf{v}^{(1)} \otimes (\textbf{v}^{(2)} + \textbf{u}^{(2)}) = \textbf{v}^{(1)} \otimes \textbf{v}^{(2)} + \textbf{v}^{(1)} \otimes \textbf{u}^{(2)}$
				\item $(\textbf{v}^{(1)} + \textbf{u}^{(1)}) \otimes \textbf{v}^{(2)} = \textbf{v}^{(1)} \otimes \textbf{v}^{(2)} + \textbf{u}^{(1)} \otimes \textbf{v}^{(2)}$
				\item $\mathbb{V}^{(1)} \otimes \mathbb{K} \cong \mathbb{V}^{(1)}$
				\item $(\mathbb{V}^{(1)} \otimes \mathbb{V}^{(2)})^\ast \cong (\mathbb{V}^{(1)})^\ast \otimes (\mathbb{V}^{(1)})^\ast$
				\item $\mathbb{V}^{(1)} \otimes \mathbb{V}^{(2)} \cong \mathbb{V}^{(2)} \otimes \mathbb{V}^{(1)}$
				\item $(\mathbb{V}^{(1)} \otimes \mathbb{V}^{(2)}) \otimes \mathbb{V}^{(3)} \cong \mathbb{V}^{(1)} \otimes (\mathbb{V}^{(2)} \otimes \mathbb{V}^{(3)})$
				\item $\textbf{v}^{(1)} \otimes (\textbf{v}^{(2)} \otimes \textbf{v}^{(3)}) = (\textbf{v}^{(1)} \otimes \textbf{v}^{(2)}) \otimes \textbf{v}^{(3)}$\bigskip
			\end{enumerate}
		\end{theorem}
	
		\begin{remark}
			The generalization of these results to more products instead 3 is immediate. With respect to property 7, it is important to note that it does not imply commutativity. In general tensor spaces are associative but not commutative. Finally, these properties implies that tensor spaces are also an algebra\footnote{An \emph{algebra} is a vector spaces $V$ with a multiplication $\cdot$ satisfying the following properties: 
			$$(\textbf{x} + \textbf{y}) \cdot \textbf{z} = \textbf{xz} + \textbf{yz} \quad \text{(right distributivity)}$$			
			$$\textbf{x} \cdot (\textbf{y} + \textbf{z}) = \textbf{xy} + \textbf{xz} \quad \text{(left distributivity)}$$			
			$$(\alpha \textbf{x}) \cdot (\beta \textbf{y}) = (\alpha \beta) (\textbf{x} \cdot \textbf{y}) \quad \text{(compatibility with scalars)}$$}, 
and for this reason it is common to call them a \emph{tensor algebra}.  
		\end{remark}
		
		For $L = 1$ one can derive formula~\ref{formula_coordinates} directly from the computation of $\langle \textbf{x}, \mathcal{T} \rangle_\mathbb{R}$ (identifying $\mathcal{T}$ with a vector), and in the case $L = 2$ the formula is derived from the computation of $\langle \mathcal{T} \textbf{x}, \textbf{y} \rangle_\mathbb{R}$ (identifying $\mathcal{T}$ with a matrix). For $L = 3$ we want to use the coordinate representation of $\mathcal{T}$ (that is, identify $\mathcal{T}$ with a 3-D matrix) to obtain some kind of generalization. To accomplish this, let's take a step back and look at the case $L = 2$ again. Observe that it is possible to reduce the expression to a particular case of $L = 1$. More precisely, we can write 	
		$$\langle \mathcal{T} \textbf{x}, \textbf{y} \rangle_\mathbb{R} = 
		\textbf{y}^T \left( \sum_{j=1}^m x_j \mathcal{T}_{:j} \right) = 
		\sum_{j=1}^m x_j ( \textbf{y}^T \mathcal{T}_{:j} ) = 
		\langle \textbf{x}, \textbf{w} \rangle_\mathbb{R},$$
where $\textbf{w} = \left[ \textbf{y}^T \mathcal{T}_{:1}, \ldots, \textbf{y}^T \mathcal{T}_{:m} \right]^T$. The idea is to compute the Euclidean inner product between $\textbf{y}$ and the columns of $\mathcal{T}$ and use the results to construct a vector, which have the same size of $\textbf{x}$. Then we compute the Euclidean inner product between $\textbf{x}$ and this vector. We can obtain the same result when starting with $\textbf{x}$ instead of $\textbf{y}$, but in this case we compute the inner product between the rows of $\mathcal{T}$. In the case $L = 3$ we proceed similarly, first we fix all dimensions of $\mathcal{T}$, except the one associated to $\textbf{x}$ (the first dimension), which will vary to create a vector $\mathcal{T}_{:jk}$. Then we compute the inner product between $\textbf{x}$ and this vector. After doing this for all possible vectors, we use the results of these computations to produce a matrix such that the entry $(j,k)$ is $\langle \textbf{x}, \mathcal{T}_{:jk} \rangle_\mathbb{R}$, in which we repeat the process using the previous procedure ($L = 2$) to obtain a scalar. The ordering in which we use the vectors doesn't affect the final result, so we will work first with $\textbf{x}$, then $\textbf{y}$ and then $\textbf{z}$. After computing all inner products $\langle \textbf{x}, \mathcal{T}_{:jk} \rangle_\mathbb{R}$, we construct the matrix		
		$$\textbf{W} = \left[
		\begin{array}{ccc}
			\langle \textbf{x},\mathcal{T}_{:11} \rangle_\mathbb{R} & \ldots & \langle \textbf{x}, \mathcal{T}_{:1p} \rangle_\mathbb{R}\\
			\vdots & & \vdots\\
			\langle \textbf{x}, \mathcal{T}_{:n1} \rangle_\mathbb{R} & \ldots & \langle \textbf{x}, \mathcal{T}_{:np} \rangle_\mathbb{R}
		\end{array}
		\right]$$
respecting the ordering of the indexes. Now we repeat the process for $\textbf{y}$, that is, compute all inner products $\langle \textbf{y}, \textbf{W}_{:k} \rangle_\mathbb{R}$ and form the vector
		$$\textbf{w} = \left[
		\begin{array}{c}
			\langle \textbf{y}, \textbf{W}_{:1} \rangle_\mathbb{R}\\
			\vdots\\
			\langle \textbf{y}, \textbf{W}_{:p} \rangle_\mathbb{R}
		\end{array}
		\right].$$
Finally, compute the inner product $\langle \textbf{z}, \textbf{w} \rangle_\mathbb{R}$. It is not difficult to verify that this value is equal to $\mathcal{T}(\textbf{x}, \textbf{y}, \textbf{z})$. This procedure can be generalized for any $L$, the idea is always to compute the inner products with respect to one dimension, thus cutting the order by one. This is repeated until order one, when we compute a single inner product. We can formalize this procedure with the following definition.

		\begin{definition}\label{contraction_def}
			Let two tensors $\mathcal{T} \in \mathbb{K}^{I_1} \otimes \ldots \otimes \mathbb{K}^{I_L}$ and $\mathcal{W} \in \mathbb{K}^{J_1} \otimes \ldots \otimes \mathbb{K}^{J_M}$ such that $\mathcal{T} = \textbf{v}^{(1)} \otimes \ldots \otimes \textbf{v}^{(L)}$, $\mathcal{W} = \textbf{w}^{(1)} \otimes \ldots \otimes \textbf{w}^{(M)}$ and $\dim(\mathbb{K}^{I_\ell}) = \dim(\mathbb{K}^{I_m})$ for some $\ell,m$. Then the \emph{$\times_\ell^m$-contraction} between $\mathcal{T}$ and $\mathcal{W}$ is the tensor $\mathcal{T} \times_m^\ell \mathcal{W} \in \mathbb{K}^{I_1} \otimes \ldots \otimes \mathbb{K}^{I_{\ell-1}} \otimes \mathbb{K}^{I_{\ell+1}} \otimes \ldots \otimes \mathbb{K}^{I_L} \otimes \mathbb{K}^{J_1} \otimes \ldots \otimes \mathbb{K}^{J_{m-1}} \otimes \mathbb{K}^{J_{m+1}} \otimes \ldots \otimes \mathbb{K}^{J_M}$ defined as
			$$\hspace{-.3cm}\mathcal{T} \times_\ell^m \mathcal{W} = \langle \textbf{v}^{(\ell)}, \textbf{w}^{(m)} \rangle_\mathbb{R} \  \textbf{v}^{(1)} \otimes \ldots \otimes \textbf{v}^{(\ell-1)} \otimes \textbf{v}^{(\ell+1)} \otimes \ldots \otimes \textbf{v}^{(L)} \otimes \textbf{w}^{(1)} \otimes \ldots \otimes \textbf{w}^{(m-1)} \otimes \textbf{w}^{(m+1)} \otimes \ldots \otimes \textbf{w}^{(L)}.$$
		\end{definition}
		
		Consider the tensor $\mathcal{T}$ of the definition above. Given $L$ vectors $\textbf{x}^{(1)}, \ldots, \textbf{x}^{(L)}$ we can use equation~\ref{rank1_map} and the definition above to write
		$$\mathcal{T}(\textbf{x}^{(1)}, \ldots, \textbf{x}^{(L)}) = \textbf{v}^{(1)} \otimes \ldots \otimes \textbf{v}^{(L)} (\textbf{x}^{(1)}, \ldots, \textbf{x}^{(L)}) = $$
		$$ = \langle \textbf{x}^{(1)}, \textbf{v}^{(1)} \rangle_\mathbb{R} \cdot \ldots \cdot \langle \textbf{x}^{(L)}, \textbf{v}^{(L)} \rangle_\mathbb{R} = $$
		\begin{equation}\label{contraction}		
			= \mathcal{T} \times _1^1 \textbf{x}^{(1)} \times_2^1 \textbf{x}^{(2)} \times_3^1 \ldots \times_L^1 \textbf{x}^{(L)}.
		\end{equation}
		
		For a generic tensor this gives an immediate expression for $\mathcal{T}(\textbf{x}^{(1)}, \ldots, \textbf{x}^{(L)})$ since it is a linear combination of tensors of the form $\textbf{e}_{i_1}^{(1)} \otimes \ldots \otimes \textbf{e}_{i_L}^{(L)}$, whose contraction is described by the formula~\ref{contraction}.
		
		Another useful definition is the \emph{mode $\ell$ tensor-vector} product, which is a contraction along one dimension. Given a vector $\textbf{x} \in \mathbb{R}^{I_\ell}$, we define $\mathcal{T} \times_\ell \textbf{x} \in \mathbb{R}^{I_1 \times \ldots \times I_{\ell-1} \times I_{\ell+1} \times \ldots \times I_L}$ by the relation 
		$$\mathcal{T} \times_\ell \textbf{x} = \left( \textbf{I}_1, \ldots, \textbf{I}_{\ell-1}, \textbf{x}^T, \textbf{I}_{\ell+1}, \ldots, \textbf{I}_L \right) \cdot \mathcal{T}.$$
		
	\section{Rank properties}
		We begin with some basic result regard to the tensor rank.
	
		\begin{theorem}
			Let $\mathcal{T} \in \mathbb{V}^{(1)} \otimes \ldots \otimes \mathbb{V}^{(L)}$ and $\mathcal{S} \in \mathbb{W}^{(1)} \otimes \ldots \otimes \mathbb{W}^{(M)}$. Then the following holds.
			\begin{enumerate}
				\item If $rank(\mathcal{S}) = 1$, then $rank(\mathcal{T} \otimes \mathcal{S}) = rank(\mathcal{T})$.
				\item If $rank(\mathcal{T}), rank(\mathcal{S}) > 1$, then $rank(\mathcal{T} \otimes \mathcal{S}) \leq rank(\mathcal{T}) \cdot rank(\mathcal{S})$.
				\item Let $\mathcal{T} = \displaystyle \sum_{r=1}^R \textbf{v}_r^{(1)} \otimes \ldots \otimes \textbf{v}_r^{(L)}$ be such that $\textbf{v}_1^{(\ell)}, \ldots, \textbf{v}_R^{(\ell)} \in \mathbb{V}^{(\ell)}$ are linearly independent vectors for each $\ell = 1 \ldots L$. Then $rank(\mathcal{T}) = R$.
			\end{enumerate}			 
		\end{theorem}	
	
		Note that in the last item of the theorem the hypothesis about the vectors is equivalent to saying that each factor matrix $\textbf{V}^{(\ell)} = [ \textbf{v}_1^{(\ell)}, \ldots, \textbf{v}_R^{(\ell)} ]$ is of full rank. 
		
		It is of interest to understand how often the CPD is unique, and for this we will need to introduce other types of rank.
		
		\begin{definition}
			We say $R$ is a \emph{typical rank of $\mathbb{V}^{(1)} \otimes \ldots \otimes \mathbb{V}^{(L)}$} if the set of rank-$R$ tensors in this space has positive probability (positive Lebesgue measure).  
		\end{definition}
		
		\begin{definition}
			We say $R$ is a \emph{generic rank of $\mathbb{V}^{(1)} \otimes \ldots \otimes \mathbb{V}^{(L)}$} if the set of rank-$R$ tensors in this space has probability 1.  
		\end{definition}
		
		These definitions assume that the tensors are being drawn according to a continuous probability distribution. Note that $R$ is the generic rank if, and only if there is only one typical rank in $\mathbb{V}^{(1)} \otimes \ldots \otimes \mathbb{V}^{(L)}$. In this case all other ranks occur with zero probability. It is important to call attention to the fact that either type of rank depends on the field $\mathbb{K}$. In \cite{kruskal2} there is a example of a real third order tensor which has rank 3 over $\mathbb{R}$ and rank 2 over $\mathbb{C}$. Another big difference between the real case and the complex case is given by the following theorem.
		
		\begin{theorem}
			If $\mathbb{V}^{(1)} \otimes \ldots \otimes \mathbb{V}^{(L)}$ is a $\mathbb{C}$-vector space, then it has a single typical rank.
		\end{theorem} 
		
		As this result shows, complex tensor spaces only have one typical rank, which is the generic rank, whereas real tensor spaces usually have more than one typical rank, so they does not have a generic rank. In the real case, some authors calls the least typical rank as the generic rank, but we won't use this terminology here. There is only one more type of rank we need to introduce. 
		
		\begin{example}
			Consider the space $\mathbb{K}^m \otimes \mathbb{K}^n \otimes \mathbb{K}^p$. Let $\textbf{x} \otimes \textbf{y} \otimes \textbf{z}$ be a rank one tensor such that $x_1, y_1,z_1$ are not zero. Then we can rewrite this tensor as 
			
			$$\textbf{x} \otimes \textbf{y} \otimes \textbf{z} = [x_1, x_2, \ldots, x_m]^T \otimes [y_1, y_2, \ldots, y_n]^T \otimes [z_1, z_2, \ldots, z_p]^T = $$
			
			$$ = \left( x_1\left[1, \frac{x_2}{x_1}, \ldots, \frac{x_m}{x_1}\right]^T \right) \otimes \left( y_1\left[1, \frac{y_2}{y_1}, \ldots, \frac{y_n}{y_1}\right]^T \right) \otimes \left( z_1\left[1, \frac{z_2}{z_1}, \ldots, \frac{z_p}{z_1}\right]^T \right) = $$
			
			$$ = (x_1 y_1 z_1)\cdot \left[1, \tilde{x}_1, \ldots, \tilde{x}_{m-1}\right]^T \otimes \left[1, \tilde{y}_1, \ldots, \tilde{y}_{n-1}\right]^T \otimes \left[1, \tilde{z}_1, \ldots, \tilde{z}_{p-1}\right]^T = $$
			
			$$ = \lambda \cdot \left[1, \tilde{x}_1, \ldots, \tilde{x}_{m-1}\right]^T \otimes \left[1, \tilde{y}_1, \ldots, \tilde{y}_{n-1}\right]^T \otimes \left[1, \tilde{z}_1, \ldots, \tilde{z}_{p-1}\right]^T.$$
			
			With this we can conclude that any tensor of this form can be described using only $1 + (m-1) + (n-1) + (p-1)$ parameters instead of $m+n+p$ if we used the original form. Now let $M_{ijk}$ be the set of tensors $\textbf{x} \otimes \textbf{y} \otimes \textbf{z}$ such that $x_i \neq 0, y_j \neq 0, z_k \neq 0$. We have that $M_{ijk}$ has dimension $1 + (m-1) + (n-1) + (p-1)$ since we can parametrize it with $1 + (m-1) + (n-1) + (p-1)$ parameters and no less than that. Note that the set of all rank one tensors also has dimension $1 + (m-1) + (n-1) + (p-1)$ since this set is the (finite) union 
			
			$$\bigcup_{\substack{1 \leq i \leq m \\ 1 \leq j \leq n \\ 1 \leq k \leq p}} M_{ijk}$$ 
of sets with dimension $1 + (m-1) + (n-1) + (p-1)$. 			
		\end{example} 
		
		It is not hard to generalize this example to the space $\mathbb{V}^{(1)} \otimes \ldots \otimes \mathbb{V}^{(L)}$ and then conclude that the set of rank one tensors has dimension $1 + (I_1 - 1) + \ldots + (I_L - 1) = 1 + \displaystyle \sum_{\ell=1}^L (I_\ell-1)$. Given a value $1 \leq R \leq \displaystyle \prod_{\ell=1}^L I_\ell$ for the rank, the set of rank-$R$ tensors can be written as the sum of $R$ copies of the set of rank one tensors. Therefore we could conclude that the dimension of this set is $R \left( 1 + \displaystyle \sum_{\ell=1}^L (I_\ell-1) \right)$. The parametrization given is not necessarily optimal, it may be possible to obtain other with less parameters. The best we can say is that this parametrization may be close to optimal, so  we concluded that the dimension is close to the actual one. Furthermore, note that this set actually contains not only the rank-$R$ tensors, it also contains all tensors of rank $\leq R$. 
		
		Denote $\sigma_R \left( \mathbb{V}^{(1)} \otimes \ldots \otimes \mathbb{V}^{(L)} \right) = \left\{ \mathcal{T} \in \mathbb{V}^{(1)} \otimes \ldots \otimes \mathbb{V}^{(L)}: \ rank(\mathcal{T}) \leq R \right\}$ for the set of tensors with rank $\leq R$. If the space is understood from the context we may just denote $\sigma_R$. With this notation and the previous observation about the parametrization one can write 
		
		$$\dim ( \sigma_R ) \lessapprox R \left( 1 + \displaystyle \sum_{\ell=1}^L (I_\ell-1) \right).$$
		
		We are talking about dimension of sets in a informal manner to refer to the number of parameters necessary to parameterize the set but in general it is not true that $\sigma_R$ is a manifold. Now suppose $R$ is such that our parametrization is optimal and it generates the entire space. Then we have 
		
		$$\dim ( \sigma_R ) = \dim\left( \mathbb{V}^{(1)} \otimes \ldots \otimes \mathbb{V}^{(L)} \right) \implies$$
		
		$$\implies R \left( 1 + \sum_{\ell=1}^L (I_\ell-1) \right) = \prod_{\ell=1}^L I_\ell \implies R = \left\lceil \frac{ \prod_{\ell=1}^L I_\ell }{ 1 + \sum_{\ell=1}^L (I_\ell-1) } \right\rceil.$$ 
		
		By definition, this $R$ is a good candidate to be the generic rank of $\mathbb{V}^{(1)} \otimes \ldots \otimes \mathbb{V}^{(L)}$. This motivates the following definition.
		
		\begin{definition}
			The \emph{expected generic rank of $\mathbb{V}^{(1)} \otimes \ldots \otimes \mathbb{V}^{(L)}$} is defined as being the value
			
			$$R_E = \left\lceil \frac{ \prod_{\ell=1}^L I_\ell }{ 1 + \sum_{\ell=1}^L (I_\ell-1) } \right\rceil.$$
		\end{definition}
		
		Now we are ready to use these definitions to state some useful theorems about tensor rank.
		
		\begin{theorem}
			Let $\mathbb{V}^{(1)} \otimes \ldots \otimes \mathbb{V}^{(L)}$ be such that $I_1 \geq I_2 \geq \ldots \geq I_L$ and let $R$ be its least typical rank. Then $rank(\mathcal{T}) \leq \min\left\{ \prod_{\ell=2}^L I_\ell, 2R \right\}$ for all $\mathcal{T} \in \mathbb{V}^{(1)} \otimes \ldots \otimes \mathbb{V}^{(L)}$.
		\end{theorem}
		
		\begin{theorem}
			Let $\mathbb{C}^n \otimes \mathbb{C}^n \otimes \mathbb{C}^n$ be such that $n \neq 3$, let $R$ be its generic rank and $R_E$ its expected generic rank. Then $R = R_E = \left\lceil \frac{n^3 - 1}{3n - 2} \right\rceil$.
		\end{theorem}
		
		\begin{theorem}
			Let $\mathbb{C}^n \otimes \mathbb{C}^n \otimes \mathbb{C}^n$ be such that $n \neq 3$ and let $R_E$ be its expected generic rank. If $\mathcal{T} \in \mathbb{C}^n \otimes \mathbb{C}^n \otimes \mathbb{C}^n$ is such that $rank(\mathcal{T}) < R_E$, then, with probability 1, $\mathcal{T}$ has finite CPD's.
		\end{theorem}
		
		\begin{definition} \label{unbalanced}
			Let $\mathbb{V}^{(1)} \otimes \ldots \otimes \mathbb{V}^{(L)}$ such that $I_1 \geq I_2 \geq \ldots \geq I_L$. We say this space is \emph{unbalanced} if 			
			$$\sum_{\ell=1}^L (I_\ell - 1) > \prod_{\ell=2}^L I_\ell.$$ 
Otherwise we say it is \emph{balanced}.
		\end{definition}
		
		Intuitively, a tensor space is unbalanced when the largest dimension is much larger than the other dimensions of the tensor product. In this case the tensor space essentially behaves as a matrix space. As the next result shows, the notion of expected generic rank is not necessary for unbalanced spaces because now we have an explicit formula for the generic rank.
		
		\begin{theorem}
			Let $\mathbb{V}^{(1)} \otimes \ldots \otimes \mathbb{V}^{(L)}$ be unbalanced and such that $I_1 \geq I_2 \geq \ldots \geq I_L$. If $R$ is its generic rank, then $R = \min \left\{ I_1, \prod_{\ell=2}^L I_\ell \right\}$.
		\end{theorem}	
		
	\section{Special products}
			In addition to the tensor product and multilinear multiplication, there are a few other products that will be relevant to us. 
	
			\begin{definition}
				Let two matrices $\textbf{A} \in \mathbb{K}^{k \times \ell}, \textbf{B} \in \mathbb{K}^{m \times n}$. The \emph{Kronecker product} between $\textbf{A}$ and $\textbf{B}$ is defined by
				$$\textbf{A} \tilde{\otimes} \textbf{B} = \left[
				\begin{array}{cccc}
					a_{11} \textbf{B} & a_{12} \textbf{B} & \ldots & a_{1\ell} \textbf{B}\\
					a_{21} \textbf{B} & a_{22} \textbf{B} & \ldots & a_{2\ell} \textbf{B}\\
					\vdots & \vdots & \ddots & \vdots\\
					a_{k1} \textbf{B} & a_{k2} \textbf{B} & \ldots & a_{k\ell} \textbf{B}
				\end{array}
				\right].$$
			\end{definition}
	
			The matrix given in the definition is a block matrix such that each block is a $m \times n$ matrix, so $\textbf{A} \tilde{\otimes} \textbf{B}$ is a $km \times \ell n$ matrix. We would like to point out that some texts uses $\otimes$ for the Kronecker product and $\circ$ for the tensor product. 
			
			\begin{lemma} \label{kronecker-cpd}
				For any vectors $\textbf{x}^{(1)} \in \mathbb{R}^{I_1}, \ldots, \textbf{x}^{(L)} \in \mathbb{R}^{I_L}$ we have that
				
				$$\textbf{x}^{(1)} \tilde{\otimes} \ldots \tilde{\otimes} \textbf{x}^{(L)} = 
				\left[ \begin{array}{c}
					x_1^{(1)} \cdot \ldots \cdot x_1^{(L-1)} x_1^{(L)}\\\\
					x_1^{(1)} \cdot \ldots \cdot x_1^{(L-1)} x_2^{(L)}\\
					\vdots\\
					x_1^{(1)} \cdot \ldots \cdot x_1^{(L-1)} x_{I_L}^{(L)}\\\\
					x_1^{(1)} \cdot \ldots \cdot x_2^{(L-1)} x_1^{(L)}\\\\
					\vdots\\\\
					x_{I_1}^{(1)} \cdot \ldots \cdot x_{I_L-1}^{(L-1)} x_1^{(L)}\\
					\vdots\\
					x_{I_1}^{(1)} \cdot \ldots \cdot x_{I_L-1}^{(L-1)} x_{I_L}^{(L)}
				\end{array} \right].$$
			\end{lemma}
			
			The important thing to notice in this lemma is the ordering of the multi-index from top to bottom. It goes to $1 \ldots 11$ to $I_L \ldots I_{L-1} I_L$ following the numbering increasing order. In particular this lemma implies that $\textbf{x}^{(1)} \tilde{\otimes} \ldots \tilde{\otimes} \textbf{x}^{(L)} = vec(\textbf{x}^{(1)} \otimes \ldots \otimes \textbf{x}^{(L)})$, assuming we are using the same ordering of the multi-indexes. This explains why the notation $\otimes$ is commonly used for the Kronecker product.
			
			\begin{definition}
				Let two matrices $\textbf{A} \in \mathbb{K}^{k \times n}, \textbf{B} \in \mathbb{K}^{m \times n}$. The \emph{Khatri-Rao product} between $\textbf{A}$ and $\textbf{B}$ is defined by
				
				$$\textbf{A} \odot \textbf{B} = \left[ \textbf{A}_{:1} \tilde{\otimes} \textbf{B}_{:1}, \ \textbf{A}_{:2} \tilde{\otimes} \textbf{B}_{:2}, \ \ldots, \ \textbf{A}_{:n} \tilde{\otimes} \textbf{B}_{:n} \right].$$
			\end{definition}
	
			In the above definition, each product $\textbf{A}_{:i} \tilde{\otimes} \textbf{B}_{:i}$ is a Kronecker product of two column vectors. Thus, each $\textbf{A}_{:i} \tilde{\otimes} \textbf{B}_{:i}$ is a $ km \times 1 $ vector, so that $\textbf{A} \odot \textbf{B}$ is $ km \times n$ matrix. In particular, if $\textbf{W}^{(1)}, \ldots, \textbf{W}^{(L)}$ are the factor matrices of a CPD, then $\textbf{W}^{(1)} \odot \ldots \odot \textbf{W}^{(L)}$ is a matrix of shape $\displaystyle\prod_{\ell=1}^L I_\ell \times R$, where its $r$-th column is the coordinate representation of the $r$-th rank one term of the CPD.
			
			\begin{definition} \label{hadamard}
				Let two matrices $\textbf{A}, \textbf{B} \in \mathbb{K}^{m \times n}$. The \emph{Hadamard product} between $\textbf{A}$ and $\textbf{B}$ is defined by
				$$\textbf{A} \ast \textbf{B} = \left[ 
				\begin{array}{cccc}
					a_{11} b_{11} & a_{12} b_{12} & \ldots & a_{1n} b_{1n}\\
					a_{21} b_{21} & a_{22} b_{22} & \ldots & a_{2n} b_{2n}\\
					\vdots & \vdots & \ddots & \vdots\\
					a_{m1} b_{m1} & a_{m2} b_{m2} & \ldots & a_{mn} b_{mn}\\
				\end{array}
				\right].$$
			\end{definition}
			
			Note that the Hadamard product is nothing more than the coordinate-wise product between matrices. The next theorem summarizes some of the key properties of these three products.
			
			\begin{theorem} \label{special-products}
				Let $\textbf{A},\textbf{B},\textbf{C},\textbf{D}$ be matrices with the sizes necessary to have well defined operations below. The the following holds.
				\begin{enumerate}
					\item $\textbf{A} \tilde{\otimes} (\textbf{B}+\textbf{C}) = \textbf{A} \tilde{\otimes} \textbf{B} + \textbf{A} \tilde{\otimes} \textbf{C}$.
					\item $(\textbf{A}+\textbf{B}) \tilde{\otimes} \textbf{C} = \textbf{A} \tilde{\otimes} \textbf{C} + \textbf{B} \tilde{\otimes} \textbf{C}$.
					\item For all $\alpha \in \mathbb{K}, \ (\alpha \textbf{A}) \tilde{\otimes} \textbf{B} = \textbf{A} \tilde{\otimes} (\alpha \textbf{B}) = \alpha (\textbf{A} \tilde{\otimes} \textbf{B})$.
					\item $(\textbf{A} \tilde{\otimes} \textbf{B}) \tilde{\otimes} \textbf{C} = \textbf{A} \tilde{\otimes} (\textbf{B} \tilde{\otimes} \textbf{C})$.
					\item $(\textbf{A} \tilde{\otimes} \textbf{B}) (\textbf{C} \tilde{\otimes} \textbf{D}) = (\textbf{AC}) \tilde{\otimes} (\textbf{BD})$.
					\item $\textbf{A}^{-1} \tilde{\otimes} \textbf{B}^{-1} = (\textbf{A} \tilde{\otimes} \textbf{B})^{-1}$.
					\item $\textbf{A}^\dagger \tilde{\otimes} \textbf{B}^\dagger = (\textbf{A} \tilde{\otimes} \textbf{B})^\dagger$.
					\item $\textbf{A}^T \tilde{\otimes} \textbf{B}^T = (\textbf{A} \tilde{\otimes} \textbf{B})^T$.
					\item $(\textbf{A} \tilde{\otimes} \textbf{B})^\dagger = \textbf{A}^\dagger \tilde{\otimes} \textbf{B}^\dagger$.
					\item $(\textbf{A} \odot \textbf{B}) \odot \textbf{C} = \textbf{A} \odot (\textbf{B} \odot \textbf{C})$.
					\item $(\textbf{A} \odot \textbf{B})^T (\textbf{A} \odot \textbf{B}) = (\textbf{A}^T\textbf{A}) \ast (\textbf{B}^T\textbf{B})$.
					\item $(\textbf{A} \odot \textbf{B})^\dagger = ((\textbf{A}^T\textbf{A}) \ast (\textbf{B}^T\textbf{B}))^\dagger (\textbf{A} \odot \textbf{B})^T$.
					\item If $\textbf{a}, \textbf{b}$ are vectors, then $\textbf{a} \tilde{\otimes} \textbf{b} = \textbf{a} \odot \textbf{b}$.
					\item $(\textbf{A} \tilde{\otimes} \textbf{B}) \ast (\textbf{C} \tilde{\otimes} \textbf{D}) = (\textbf{A} \ast \textbf{C}) \tilde{\otimes} (\textbf{B} \ast \textbf{D})$.
				\end{enumerate}
			\end{theorem}	

	\section{Symmetric tensors}
		Let $S_L$ be the group of permutations of $L$ elements. Each element of $S_L$ will be denoted by $\sigma$. We can see $\sigma$ simply as a $L$-tuple of integers between 1 and $L$, without repetitions. Thus, we denote each entry of $\sigma$ by $\sigma(\ell)$. For example, suppose $L = 3$. A possible permutation $\sigma \in S_3$ is $\sigma = (3,2,1)$. In this case, $\sigma$ permutes the $3$-tuple $(1,2,3)$ to $(3,2,1)$. Thus we have that $\sigma(1) = 3$, $\sigma(2) = 2$, $\sigma(3) = 1$. Now, define the map $\pi_S:\mathbb{V}^{\otimes L} \to \mathbb{V}^{\otimes L}$ by
	
		$$\pi_S(\textbf{v}^{(1)} \otimes \ldots \otimes \textbf{v}^{(L)}) = \frac{1}{L!} \sum_{\sigma \in S_L} \textbf{v}^{(\sigma(1))} \otimes \ldots \otimes \textbf{v}^{(\sigma(L))}.$$
	
		We may interpret $\pi_S (\textbf{v}^{(1)} \otimes \ldots \otimes \textbf{v}^ {(L)}) $ as the average of the tensor products between the vectors $\textbf{v}^{(\ell)}$ considering all possible permutations of the indexes. Usually one denote $\textbf{v}^{(1)} \cdot \ldots \cdot \textbf{v}^{(L)} = \pi_S (\textbf{v}^{(1)} \otimes \ldots \otimes \textbf{v}^{(L)})$ and call this the \emph{symmetric tensor product} between the vectors $\textbf{v}^{(\ell)}$. In the case all vectors are equal, say equal to $\textbf{v}$, we denote $\textbf{v}^L = \pi_S (\underbrace{\textbf{v} \otimes \ldots \otimes \textbf{v}}_{L \text{ times}})$.	The image of $\pi_S$ is called the space of \emph{symmetric tensors}\footnote{If one want to specify the tensor order it is possible to refer to this space as the space of \emph{symmetric $L$-tensors} or \emph{symmetric $L$-order tensors}} of $\mathbb{V}$ and it is denoted by $S^L(\mathbb{V})$. Another way to characterize symmetric tensors is given by the following theorem.
	
		\begin{theorem}
			A tensor $\mathcal{T} \in \mathbb{V}^{\otimes L}$ is symmetric if, and only if, 
			
			$$\mathcal{T}(\textbf{v}^{(1)}, \ldots, \textbf{v}^{(L)}) = \mathcal{T}(\textbf{v}^{(\sigma(1))}, \ldots, \textbf{v}^{(\sigma(L)}))$$ 
for all $\sigma \in S_L$.
		\end{theorem}
		
		We finish our discussion about symmetric tensors giving a few more relevant results.
	
		\begin{theorem}
			Let $\dim(\mathbb{V}) = n$ and $\{ \textbf{e}^{(1)}, \ldots, \textbf{e}^{(n)} \}$ be a basis for $\mathbb{V}$. Then
			\begin{enumerate}
				\item $\{ e_{i_1}\cdot \ldots \cdot e_{i_L}: \ 1 \leq i_1 \leq \ldots \leq i_L \leq n \}$ is a basis of $S^L(\mathbb{V})$.
				\item $\displaystyle\dim S^L(\mathbb{V}) = \binom{L - 1 + n}{L}$.
			\end{enumerate}
		\end{theorem}
	
		\begin{theorem}
			Let $\mathcal{T} \in \mathbb{V}^{\otimes L}$ and consider a choice of basis for $\mathbb{V}$. If $\mathcal{T}$ is symmetric, then $\mathcal{T}_{i_1 \ldots i_L} = \mathcal{T}_{i_{\sigma(1)} \ldots i_{\sigma(L)}}$ for all $\sigma \in S_L$.
		\end{theorem}
		
		\begin{example}[Derivatives]
			Let $f:\mathbb{C}^m \to \mathbb{C}^n$ be a function $L$ times differentiable at $\textbf{a} \in \mathbb{C}^m$. The $L$-th derivative of $f$ at $\textbf{a}$ is a $L$-linear map $Df^{(L)}(\textbf{a}): \underbrace{\mathbb{C}^m \times \ldots \times \mathbb{C}^m}_{L \text{ times}} \to \mathbb{C}^n$. Because of the isomorphism
			
		$$\mathcal{L}_L(\mathbb{C}^m; \mathbb{C}^n) \cong (\mathbb{C}^{m})^\ast \otimes \ldots \otimes (\mathbb{C}^{m})^\ast \otimes \mathbb{C}^n,$$ 
we can write $Df^{(L)}(\textbf{a}) \in (\mathbb{C}^{m})^\ast \otimes \ldots \otimes (\mathbb{C}^{m})^\ast \otimes \mathbb{C}^n$ and consider this map as a mixed tensor of type $(L,1)$.  

		Let $\{ \textbf{e}_1, \ldots, \textbf{e}_m\}$ be the canonical basis of $\mathbb{C}^m$. Given any vectors $\dot{\textbf{x}}^{(1)}, \ldots, \dot{\textbf{x}}^{(L)} \in \mathbb{C}^m$, we have that	
		
		$$Df^{(L)}(\textbf{a})(\dot{\textbf{x}}^{(1)}, \ldots, \dot{\textbf{x}}^{(L)}) = \sum_{i_1, \ldots, i_L = 1}^m \dot{\textbf{x}}^{(1)}_{i_1} \ldots \dot{\textbf{x}}^{(L)}_{i_L} \cdot Df^{(L)}(\textbf{a})(\textbf{e}_{i_1}, \ldots, \textbf{e}_{i_k}) = $$
		$$ = \sum_{i_1, \ldots, i_L = 1}^m \dot{\textbf{x}}^{(1)}_{i_1} \ldots \dot{\textbf{x}}^{(L)}_{i_L} \cdot  
		\left[ \begin{array}{c}
			\displaystyle\frac{\partial^L f_1}{\partial x_{i_1} \ldots \partial x_{i_L}}(\textbf{a})\\
			\displaystyle\vdots\\
			\displaystyle\frac{\partial^L f_n}{\partial x_{i_1} \ldots \partial x_{i_L}}(\textbf{a})
		\end{array} \right].$$	
	
		Using the notation
		
		$$\frac{\partial^L f}{\partial x_{i_1} \ldots \partial x_{i_L}}(\textbf{a}) = \left[ \begin{array}{c}
			\displaystyle\frac{\partial^L f_1}{\partial x_{i_1} \ldots \partial x_{i_L}}(\textbf{a})\\
			\displaystyle\vdots\\
			\displaystyle\frac{\partial^L f_n}{\partial x_{i_1} \ldots \partial x_{i_L}}(\textbf{a})
		\end{array} \right]$$
and considering the dual basis $dx_i:\mathbb{C}^m \to \mathbb{C}$ such that $dx_i(\textbf{x}) = dx_i(x_1, \ldots, x_m) = x_i$, we have that

		$$Df^{(L)}(\textbf{a})(\dot{\textbf{x}}^{(1)}, \ldots, \dot{\textbf{x}}^{(L)}) = \sum_{i_1, \ldots, i_L = 1}^m dx_{i_1} \otimes \ldots \otimes dx_{i_L}(\dot{\textbf{x}}^{(1)}, \ldots, \dot{\textbf{x}}^{(L)}) \cdot \frac{\partial^L f}{\partial x_{i_1} \ldots \partial x_{i_L}}(\textbf{a}).$$ 	
	
		Notice that $Df^{(L)}(\textbf{a})$ a tensor is such that each coordinate $(i_1, \ldots, i_L, j)$ is given by $\displaystyle\frac{\partial^L f_j}{\partial x_{i_1} \ldots \partial x_{i_L}}(\textbf{a})$, for $1 \leq i_1, \ldots, i_L \leq m$ and $1 \leq j \leq n$. In a more simplified way, we can write 	
		
		$$Df^{(L)}(\textbf{a}) =  \sum_{i_1, \ldots, i_L = 1}^m dx_{i_1} \otimes \ldots \otimes dx_{i_L} \otimes \frac{\partial^L f}{\partial x_{i_1} \ldots \partial x_{i_L}}(\textbf{a}),$$
with the understanding that we are using the interpretation $Df^{(L)}(\textbf{a}) \in \mathcal{L}_L(\mathbb{C}^m; \mathbb{C}^n)$ and not $Df^{(L)}(\textbf{a}) \in \mathcal{L}(\mathbb{C}^m, \ldots, \mathbb{C}^m, \mathbb{C}^{n\ast}; \mathbb{C})$. This is a nice characterization of the derivative. Another elegant one is to expand

		$$\frac{\partial^L f}{\partial x_{i_1} \ldots \partial x_{i_L}}(\textbf{a}) = \left[ \begin{array}{c}
			\displaystyle\frac{\partial^L f_1}{\partial x_{i_1} \ldots \partial x_{i_L}}(\textbf{a})\\
			\displaystyle\vdots\\
			\displaystyle\frac{\partial^L f_n}{\partial x_{i_1} \ldots \partial x_{i_L}}(\textbf{a})
		\end{array} \right] = 
		\sum_{j=1}^n \displaystyle\frac{\partial^L f_j}{\partial x_{i_1} \ldots \partial x_{i_L}}(\textbf{a}) \cdot \textbf{e}_j$$
so that the previous formula becomes

		$$Df^{(L)}(\textbf{a}) =  \sum_{i_1, \ldots, i_L = 1}^m dx_{i_1} \otimes \ldots \otimes dx_{i_L} \otimes \left( \sum_{j=1}^n \displaystyle\frac{\partial^L f_j}{\partial x_{i_1} \ldots \partial x_{i_L}}(\textbf{a}) \cdot \textbf{e}_j \right) = $$
		
		$$ = \sum_{i_1, \ldots, i_L = 1}^m \sum_{j=1}^n \displaystyle\frac{\partial^L f_j}{\partial x_{i_1} \ldots \partial x_{i_L}}(\textbf{a}) \cdot dx_{i_1} \otimes \ldots \otimes dx_{i_L} \otimes \textbf{e}_j.$$
In particular, with this expression we can see that the coordinates of $Df^{(L)}(\textbf{a})$, as a $(L+1)$-th order tensor, are given by
	
	$$Df^{(L)}(\textbf{a})_{i_1 \ldots i_L j} = \displaystyle\frac{\partial^L f_j}{\partial x_{i_1} \ldots \partial x_{i_L}}(\textbf{a}).$$

It is also possible to consider the identification we are using so far and consider $Df^{(L)}(\textbf{a})$ as a tensor in $\mathbb{C}^{m} \otimes \ldots \otimes \mathbb{C}^{m} \otimes \mathbb{C}^n$. Then, instead of $dx_i$ we use the canonical basis vector $\textbf{e}_i$. 
		\end{example}
	
	\section{Antisymmetric tensors}
  		Given $\sigma \in S_L$, the \emph{sign} of $\sigma$ is denoted by $\text{sgn}(\sigma)$. We define $\text{sgn}(\sigma) = 1$ if $\sigma$ makes an even quantity of permutation and $\text{sgn}(\sigma) = -1$ if $\sigma$ makes an odd quantity of permutations. For example, let $L = 3$ and $\sigma = (3,2,1)$. Note that $\sigma(2) = 2$, then this number is not permuted, but the numbers 1 and 3 are permuted. From this it follows that $\sigma$ makes 2 permutations. Since 2 is even, we conclude that $\text{sgn}(\sigma) = 1$. Now let $\tilde{\sigma} = (3,1,2)$. In this case all three numbers are permuted. Since 3 is odd, we conclude that $\text{sgn}(\tilde{\sigma}) = -1$.
	
		Now define the map  $\pi_\Lambda:\mathbb{V}^{\otimes L} \to \mathbb{V}^{\otimes L}$ by
	
		$$\pi_\Lambda(\textbf{v}^{(1)} \otimes \ldots \otimes \textbf{v}^{(L)}) = \frac{1}{L!} \sum_{\sigma \in S_L} \text{sgn}(\sigma) \textbf{v}^{(\sigma(1))} \otimes \ldots \otimes \textbf{v}^{(\sigma(L))}.$$
	
		We may interpret $\pi_\Lambda(\textbf{v}^{(1)} \otimes \ldots \otimes \textbf{v}^{(L)})$ as the signed average of the tensor products between the vectors $\textbf{v}^{(\ell)} $ considering all possible permutations of the indexes. Usually one denote $\textbf{v}^{(1)} \wedge \ldots \wedge \textbf{v}^{(L)} = \pi_\Lambda(\textbf{v}^{(1)} \otimes \ldots \otimes \textbf{v}^{(L)})$ and can this the \emph{exterior product} or \emph{antisymmetric product} between the vectors $\textbf{v}^{(\ell)}$. The image of $\pi_\Lambda$ is called the space of \emph{antisymmetric tensors} or \emph{alternating tensors} of $V$ and it is denoted by $\Lambda^L(\mathbb{V})$. 
	
		\begin{theorem}
			$\textbf{v}^{(1)} \wedge \ldots \wedge \textbf{v}^{(L)} = 0 \iff$ the vectors $\textbf{v}^{(1)}, \ldots, \textbf{v}^{(L)}$ are linearly dependent.
		\end{theorem}
	
		\begin{theorem}
			$\textbf{v}^{(1)} \wedge \ldots \wedge \textbf{v}^{(\ell)} \wedge \ldots \wedge \textbf{v}^{(j)} \wedge \ldots \wedge \textbf{v}^{(L)} = - \textbf{v}^{(1)} \wedge \ldots \wedge \textbf{v}^{(j)} \wedge \ldots \wedge \textbf{v}^{(\ell)} \wedge \ldots \wedge \textbf{v}^{(L)}$ for all $i \neq j$.
		\end{theorem}
	
		\begin{theorem}
			Let $\dim(\mathbb{V}) = n$ and $\{ \textbf{e}_1, \ldots, \textbf{e}_n \}$ be a basis for $\mathbb{V}$. Then
			\begin{enumerate}
				\item $\{ \textbf{e}_{i_1}\wedge \ldots \wedge \textbf{e}_{i_L}: \ 1 \leq i_1 < \ldots < i_L \leq n \}$ is a basis of $\Lambda^L(V)$.
				\item $\displaystyle\dim \Lambda^L(\mathbb{V}) = \binom{n}{L}$.\bigskip
			\end{enumerate}
		\end{theorem}
	
		\begin{theorem}
			Let $\mathcal{T} \in \mathbb{V}^{\otimes L}$ and consider a choice of basis for $\mathbb{V}$. If $\mathcal{T}$ is antisymmetric, then $\mathcal{T}_{i_1 \ldots i_p \ldots i_q \ldots i_L} = -\mathcal{T}_{i_1 \ldots i_q \ldots i_p \ldots i_L}$ for all $p \neq q$.
		\end{theorem}